\documentclass[10pt]{amsart}
\usepackage{amsfonts, amsmath, amssymb, fullpage, amsthm, semantex, tikz, graphicx, hyperref, placeins, bbm, euscript, bbold}
\usetikzlibrary{cd}

\title{HOMFLY homology for oriented tangles}
\author{William Ballinger and Luke Conners}

\newcommand{\Hom}{{\rm Hom}}

\renewcommand{\to}{\rightarrow}

\newcommand{\id}{{\rm id}}

\let\hat=\widehat
\let\tilde=\widetilde

\def\Z{{\mathbb Z}}

\def\deg{{\mathrm{deg}}}

\renewcommand{\L}{\mathcal{L}}

\newcommand{\zip}{\operatorname{zip}}
\newcommand{\unzip}{\operatorname{unzip}}

\newcommand{\strand}[1]{\mathbb{1}_{#1}}

\newcommand{\Bim}{\mathrm{Bim}}

\newcommand{\SBim}{\mathrm{SBim}}

\renewcommand{\aa}{\mathbf{a}}
\newcommand{\bb}{\mathbf{b}}
\newcommand{\cc}{\mathbf{c}}
\newcommand{\dd}{\mathbf{d}}

\newcommand{\extsbim}{\mathrm{ExtSBim}}
\newcommand{\Ext}{\mathrm{Ext}}
\newcommand{\HHH}{\mathrm{HHH}}
\newcommand{\otan}{\mathcal{OTAN}}
\newcommand{\otanbal}{\mathcal{OTAN}_{bal}}
\newcommand{\funcat}{\mathrm{Fun}(K^b(\extsbim^{or}))}
\newcommand{\Br}{\mathrm{Br}}

\NewSymbolClass\MyBinaryOperator[
    define keys={
        {Lder}{command=\overset{L}},
        {Rder}{upper=R},
    },
]
\NewObject\MyBinaryOperator\tensor{\otimes}[
    define keys={
        {der}{Lder},
    },
]
\NewObject\MyBinaryOperator\fibre{\times}[
    define keys={
        {der}{Rder},
    },
]

\usepackage{chngcntr}
\counterwithin{figure}{section}

\tikzset{anchorbase/.style={baseline={([yshift=-0.5ex]current bounding box.center)}},
tinynodes/.style={font=\tiny,text height=0.75ex,text depth=0.15ex},
smallnodes/.style={font=\scriptsize,text height=0.75ex,text depth=0.15ex},
>={Latex[length=1mm, width=1.5mm]},
overcross/.style={line width=4pt,color=white},
1label/.style={very thick, black},
2label/.style={very thick, ourblue},
m1label/.style={very thick, double, black},
m2label/.style={very thick, double, ourblue},
glabel/.style={very thick, gray},
mglabel/.style={very thick, double, gray},
Klabel/.style={line width=5pt, mypurple},
webs/.style={line width=.9,color=black},
}

\tikzset{
    webs/.style={
        line width=1pt,
        draw=black
    },
    over/.style={
        preaction={
            draw=white,
            line width=5pt
        }
    }
}

\usetikzlibrary{decorations.markings,arrows.meta}

\tikzset{
    ->-/.style={
        postaction={decorate},
        decoration={
            markings,
            mark=at position #1 with {
                \arrow[sloped, xshift=1.1mm]{Stealth[length=2mm,width=2mm]}
            }
        }
    },
    -<-/.style={
        postaction={decorate},
        decoration={
            markings,
            mark=at position #1 with {
                \arrowreversed[sloped, xshift=-1.1mm]{Stealth[length=2mm,width=2mm]}
            }
        }
    },
    biarr/.style={
        postaction={decorate},
        decoration={
            markings,
            mark=at position .01 with {
                \arrowreversed{Stealth[length=2mm,width=2mm]}
            },
            mark=at position .9 with {
                \arrowreversed{Stealth[length=2mm,width=2mm]}
            }
        }
    },
    register/.style={webs,color=black},
    instruction/.style={webs,color=blue},
    shuffle/.style={webs,color=red},
    new/.style={webs,color=new_new}
}

\definecolor{new_new}{HTML}{489928}

\usepackage{adjustbox}

\newcommand{\MoveArrow}{\ensuremath{\;\Longleftrightarrow\;}}
\newcommand{\MoveCell}[2]{%
  \begin{minipage}[t]{0.47\textwidth}
    \vspace{0pt}
    \centering
    \textbf{\large #1}\par
    \vspace{0.4em}
    \begin{adjustbox}{max width=\linewidth,center}
      #2
    \end{adjustbox}
  \end{minipage}
}

\numberwithin{equation}{section}

\newtheorem{theorem}{Theorem}[section]

\newtheorem{proposition}[theorem]{Proposition}
\newtheorem{prop}[theorem]{Proposition}

\newtheorem{cor}[theorem]{Corollary}

\newtheorem{lemma}[theorem]{Lemma}

\theoremstyle{definition}

\newtheorem{remark}[theorem]{Remark}
\newtheorem{rem}[theorem]{Remark}

\newtheorem{notation}[theorem]{Notation}

\newtheorem{definition}[theorem]{Definition}
\newtheorem{defn}[theorem]{Definition}

\newtheorem{example}[theorem]{Example}

\newtheorem{warning}[theorem]{Warning}

\begin{document}

\begin{abstract}
    We introduce an extension of Khovanov and Rozansky's categorification of the HOMFLY polynomial to arbitrary oriented tangles. Unlike previous such extensions, ours is invariant under the full suite of oriented Reidemeister moves, including the ``misoriented" Reidemeister IIb. For oriented links $\mathcal{L}$, our construction recovers Khovanov and Rozansky's construction while removing the preliminary step of presenting $\mathcal{L}$ as a braid closure.
\end{abstract}

\maketitle

\tableofcontents

\section{Introduction} \label{sec: intro}

HOMFLY homology (also known as \textit{triply-graded} homology) is an invariant of oriented links, first constructed by Khovanov and Rozansky in 2005 \cite{KhR08}. The investigation of its structure and properties has driven a great deal of research over the past 20 years, elucidating deep connections across a wide variety of fields including algebraic geometry \cite{GNR21, OR18}, mathematical physics \cite{GS12, GGS18}, and higher algebra \cite{LMRSW24, DW25}. HOMFLY homology also specializes to a number of other interesting link invariants, coming equipped with spectral sequences to Khovanov homology, $\mathfrak{sl}_N$ Khovanov--Rozansky homology, and knot Floer homology \cite{DGR06, Ras15, BPRW22} and more direct specializations to the entire suite of Type A link polynomials, including the Jones, Alexander, and HOMFLY polynomials, giving it something of a universal status among all Type A quantum link invariants.

Unfortunately, this added complexity comes at a cost. Direct computations are exceedingly difficult, limited at present only to positive torus links \cite{HM19}, knots of small crossing number \cite{NS25}, and a few other very special cases. HOMFLY homology is also not known to be functorial under general link cobordisms (although this is known for \textit{braid} cobordisms by work of Elias and Krasner \cite{EK10a}), preventing its use for the spectacular applications to low-dimensional topology enjoyed by other link homology theories, such as Khovanov homology \cite{RW24, Ras10} and knot Floer homology \cite{LLP23}.

Both of these issues stem from the essentially nonlocal nature of existing constructions of HOMFLY homology. All of the other link invariants described above admit local descriptions in terms of skein relations (in the polynomial case) or cubes of resolutions (in the homological case) -- the invariant is specified for the simplest tangles, and a procedure is given to glue smaller tangles into larger ones. HOMFLY homology as constructed in \cite{KhR08} admits no such description; an arbitrary oriented link diagram must first be isotoped to a braid closure (potentially getting much larger in the process) before performing any algebraic computations.

The primary reason for this nonlocal behavior is that the invariant constructed in \cite{KhR08} does not respect the ``wrong-way" Reidemeister IIb move. This obstructs an extension of HOMFLY homology as previously defined to an invariant of oriented tangles; as we explain below, the existence of such a tangle invariant is a key feature of sufficiently local theories. Beyond theoretical considerations, a local assignment to tangles is also very useful computationally, enabling dramatically faster computations \cite{BN07} and a relatively straightforward strategy for establishing functoriality under link cobordisms \cite{BN05, ETW18}.

We address these weaknesses in this work, providing a new, essentially local construction of HOMFLY homology which is invariant under the full suite of oriented Reidemeister moves, including Reidemeister IIb, and which extends to a well-defined invariant of oriented tangles. In the special case of oriented links, our invariant recovers the original Khovanov--Rozansky construction.

\subsection{Tangle Invariants} \label{subsec: tangle_invariants}

Before describing our construction in detail, we briefly recall why one might be interested in tangle invariants and the usual procedure for constructing them.

Links are embedded in three-dimensional space; nevertheless almost all link invariants are constructed using two-dimensional projections called \textit{link diagrams}. The topological content of a link diagram $\mathcal{L}$ can be recovered from two sets of data: a collection of \textit{crossings} in which one segment of $\mathcal{L}$ crosses over another, and a collection of arcs connecting the endpoints of these crossings. Up to isotopy in the plane of projection, crossings are completely determined by their sign (i.e. which strand crosses above and which below), while the arcs connecting their endpoints are completely determined by the corresponding endpoints and their orientations.

With this description in mind, most link invariants are constructed according to the following recipe\footnote{This framework is equivalent to constructing a homomorphism from the \textit{planar algebra} of tangles in the sense of Vaughan Jones \cite{Jon21}.}:

\begin{enumerate}
    \item Assign \textit{local} data to each crossing in a way depending only on its sign;
    \item Assign \textit{global} gluing data to each arc connecting two crossings in a way depending only on its endpoints.
\end{enumerate}

Given any such assignment $f$ and any link diagram $\mathcal{L}$, there is a perfectly well-defined assignment $f(\mathcal{L})$ given by evaluating the crossings of $\mathcal{L}$ according to the local data of $f$ and gluing these evaluations along the arcs of $\mathcal{L}$ according to the global data of $f$. Essentially by definition, $f(\mathcal{L})$ is an invariant of the link diagram $\mathcal{L}$ up to isotopy in the plane of projection. To check that $f(\mathcal{L})$ is an invariant of the actual link represented by $\mathcal{L}$, one has to check that any two diagrams $\mathcal{L}$ and $\mathcal{L}'$ which represent isotopic links satisfy $f(\mathcal{L}) \cong f(\mathcal{L}')$ (for some suitable notion of equivalence $\cong$). Any two such diagrams are necessarily related by a finite sequence of Reidemeister moves (see Figure \ref{fig: reidemeisters_intro}), so it suffices to check that $f$ is invariant with respect to these moves.

These moves are \textit{local}, in the sense that they are performed only in some small region of any given diagram, leaving the rest of the diagram unchanged, and they ought to be checked locally. It's therefore profitable to think of $f$ as assigning invariants to the smaller \textit{tangle} diagrams appearing on either side of these moves, again by assigning the same local data to each crossing and gluing those assignments within these smaller diagrams using the same global data. If the invariants assigned by $f$ to the tangles appearing on either side of each Reidemeister move agree, then any two larger tangle diagrams $\tau$ and $\tau'$ related by such a move in a small region will also satisfy $f(\tau) \cong f(\tau')$ essentially by construction, and $f$ is therefore a tangle invariant.

\subsection{Previous Approaches}

In \cite{Web07}, Webster gave an interpretation of Khovanov and Rozansky's HOMFLY homology in terms of local and global data. Each crossing is assigned to a certain chain complex of Soergel bimodules\footnote{Both Khovanov--Rozansky's original construction in \cite{KhR08} and Webster's \cite{Web07} actually give a slightly different construction using matrix factorizations. This is equivalent to the construction described here, which uses the algebraic language of Soergel bimodules from Khovanov's construction \cite{Kh07}, as can be seen by always picking an appropriate resolution of each Soergel bimodule in the derived category $D^b(\SBim)$.} called its \textit{Rouquier complex}, and these crossings are glued by identifying (in a derived fashion) the actions of variables at the corresponding endpoints. (See Section \ref{sec: background} for a more precise explanation; we will not rely on the details of this construction in this introduction.) This gives a perfectly well-defined invariant $HHH_{alg}$ of oriented tangle diagrams, and one can check whether $HHH_{alg}$ is a \textit{bona fide} invariant of oriented tangles by checking its invariance under a complete list of oriented Reidemeister moves.

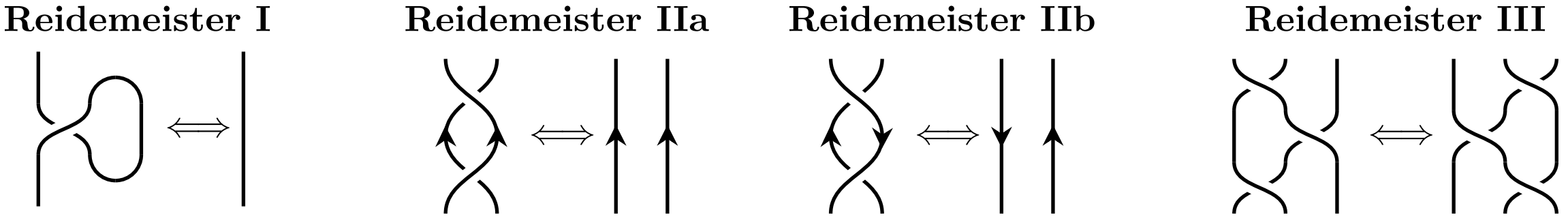
\begin{figure}[ht]
    \centering

    \makebox[\textwidth][c]{%
    \begin{minipage}[t]{0.24\textwidth}
        \centering
        \textbf{Reidemeister I}

        \medskip

        \begin{tikzpicture}[anchorbase]
            \draw[webs] (0,0) to (0,.5);
            \draw[webs] (0,1) to (0,1.5);
            \draw[webs] (.5,.5) to[out=90,in=270] (0,1);
            \draw[webs,over] (0,.5) to[out=90,in=270] (.5,1);
            \draw[webs] (.5,1) to[out=90,in=180] (.75,1.25);
            \draw[webs] (.75,1.25) to[out=0,in=90] (1,1);
            \draw[webs] (.5,.5) to[out=270,in=180] (.75,.25);
            \draw[webs] (.75,.25) to[out=0,in=270] (1,.5);
            \draw[webs] (1,.5) to (1,1);
        \end{tikzpicture}
        \MoveArrow
        \begin{tikzpicture}[anchorbase]
            \draw[webs] (0,0) to (0,1.5);
        \end{tikzpicture}
    \end{minipage}%
    \hfill
    \begin{minipage}[t]{0.20\textwidth}
        \centering
        \textbf{Reidemeister IIa}

        \medskip

        \begin{tikzpicture}[anchorbase]
            \draw[webs,->-=.99] (.5,0) to[out=90,in=270] (0,.75);
            \draw[webs] (0,.75) to[out=90,in=270] (.5,1.5);
            \draw[webs,over] (.5,.75) to[out=90,in=270] (0,1.5);
            \draw[webs,over,->-=.99] (0,0) to[out=90,in=270] (.5,.75);
        \end{tikzpicture}
        \MoveArrow
        \begin{tikzpicture}[anchorbase]
            \draw[webs,->-=.5] (0,0) to (0,1.5);
            \draw[webs,->-=.5] (.5,0) to (.5,1.5);
        \end{tikzpicture}
    \end{minipage}%
    \hfill
    \begin{minipage}[t]{0.20\textwidth}
        \centering
        \textbf{Reidemeister IIb}

        \medskip

        \begin{tikzpicture}[anchorbase]
            \draw[webs,->-=.99] (.5,0) to[out=90,in=270] (0,.75);
            \draw[webs] (0,.75) to[out=90,in=270] (.5,1.5);
            \draw[webs,over] (.5,.75) to[out=90,in=270] (0,1.5);
            \draw[webs,over,-<-=.99] (0,0) to[out=90,in=270] (.5,.75);
        \end{tikzpicture}
        \MoveArrow
        \begin{tikzpicture}[anchorbase]
            \draw[webs,-<-=.5] (0,0) to (0,1.5);
            \draw[webs,->-=.5] (.5,0) to (.5,1.5);
        \end{tikzpicture}
    \end{minipage}%
    \hfill
    \begin{minipage}[t]{0.28\textwidth}
        \centering
        \textbf{Reidemeister III}

        \medskip

        \begin{tikzpicture}[anchorbase]
            \draw[webs] (0,0) to[out=90,in=270] (.5,.5);
            \draw[webs,over] (.5,0) to[out=90,in=270] (0,.5);
            \draw[webs] (.5,.5) to[out=90,in=270] (1,1);
            \draw[webs,over] (1,.5) to[out=90,in=270] (.5,1);
            \draw[webs] (0,1) to[out=90,in=270] (.5,1.5);
            \draw[webs,over] (.5,1) to[out=90,in=270] (0,1.5);
            \draw[webs] (1,0) to (1,.5);
            \draw[webs] (0,.5) to (0,1);
            \draw[webs] (1,1) to (1,1.5);
        \end{tikzpicture}
        \MoveArrow
        \begin{tikzpicture}[anchorbase]
            \draw[webs] (.5,0) to[out=90,in=270] (1,.5);
            \draw[webs,over] (1,0) to[out=90,in=270] (.5,.5);
            \draw[webs] (0,.5) to[out=90,in=270] (.5,1);
            \draw[webs,over] (.5,.5) to[out=90,in=270] (0,1);
            \draw[webs] (.5,1) to[out=90,in=270] (1,1.5);
            \draw[webs,over] (1,1) to[out=90,in=270] (.5,1.5);
            \draw[webs] (0,0) to (0,.5);
            \draw[webs] (1,.5) to (1,1);
            \draw[webs] (0,1) to (0,1.5);
        \end{tikzpicture}
    \end{minipage}%
    }

    \caption{Reidemeister moves for oriented tangles. Each move should be considered with all consistent choices of orientation on each strand, as should reflections of these moves.}
    \label{fig: reidemeisters_intro}
\end{figure}

Khovanov and Rozansky verified that $HHH_{alg}$ respects Reidemeister I, IIa, and III\footnote{At least when all strands appearing in RIII are oriented in the same direction.}, but possibly \textit{not} RIIb (see Figure \ref{fig: reidemeisters_intro}). Much later, Abel proved in \cite{Abe17} that $HHH_{alg}$ in fact does \textit{not} respect Reidemeister $IIb$, giving a fatal obstruction to $HHH_{alg}$ being an oriented tangle invariant. One might hope that examples where this failure is relevant are few and far between, but in fact the problem is catastrophic; a recent result of Borodzik--Malashchuk \cite{BM26} shows that for \textit{any} two oriented knots $K_1$ and $K_2$, there exist diagrams $D_1$ of $K_1$ and $D_2$ of $K_2$ satisfying $HHH_{alg}(D_1) \cong HHH_{alg}(D_2)$.

Khovanov and Rozansky addressed this failure by restricting their allowable diagrams to \textit{braid closures}. A classic result of Alexander \cite{Al23} shows that any oriented link admits such a diagram, and a classic result of Markov \cite{Mar35} shows that any two braid closures representing the same oriented link can be related via isotopy in the plane and the Reidemeister moves $I$, $IIa$, and $III^+$. Because $HHH_{alg}$ is invariant under this more restrictive set of moves, given any oriented link $\mathcal{L}$ presented as a braid closure $\mathcal{L} \sim \hat{\beta}$, the assignment $\mathcal{L} \mapsto HHH_{alg}(\hat{\beta})$ is a genuine oriented link invariant.

This procedure has two obvious defects. First, since an arbitrary oriented tangle with nonempty boundary is not isotopic to an oriented braid closure, this correction does not give an oriented tangle invariant $HHH_{alg}(\tau)$. Second (and relatedly), even for an oriented link $\mathcal{L}$, one cannot compute $HHH_{alg}(\mathcal{L})$ directly from an arbitrary diagram and must first find a presentation of $\mathcal{L}$ as a braid closure. Addressing these defects has been a long-standing open problem, appearing as Problem 1.34 in the most recent Kirby list \cite{BKR26}.

\subsection{Tangle Homology}

\subsubsection{Version I: Coherent Tangle Diagrams}

One version of our solution to this problem essentially mimics that of Khovanov and Rozansky: we identify a class of oriented tangle diagrams on which $HHH_{alg}$ is well-behaved and which is sufficiently general to capture all oriented tangles up to isotopy. To that end, note that braid closures can be built iteratively from the empty tangle diagram via the following operations:

\begin{itemize}
    \item Add a single upwards-oriented strand to the left of an oriented tangle diagram;
    \item Attach a braid above or below a collection of parallel upwards-oriented strands;
    \item Close the leftmost top and bottom endpoints of an oriented tangle diagram in which both such endpoints are positively oriented.
\end{itemize}

\begin{figure}[ht]
    \centering
    \begin{tikzpicture}[scale=.75]
    
    \draw[webs] (1,0) to[out=90,in=270] (0,1);
    \draw[webs,over,->-=.2] (0,0) to[out=90,in=270] (1,1);
    \draw[webs,->-=.5] (2,-.5) to (2,1);

    \draw[webs] (2,1) to[out=90,in=270] (1,2);
    \draw[webs,over,->-=.01] (1,1) to[out=90,in=270] (2,2);
    \draw[webs,->-=.5] (0,1) to (0,2);

    \draw[webs,->-=.01] (1,2) to[out=90,in=270] (0,3);
    \draw[webs,over] (0,2) to[out=90,in=270] (1,3);
    \draw[webs,->-=.5] (2,2) to (2,5);

    \draw[webs] (0,3) to[out=90,in=0] (-.25,3.25);
    \draw[webs] (-.25,3.25) to[out=180,in=90] (-.5,3);
    \draw[webs,->-=.5] (-.5,3) to (-.5,0);
    \draw[webs] (-.5,0) to[out=270,in=180] (-.25,-.25);
    \draw[webs] (-.25,-.25) to[out=0,in=270] (0,0);

    \draw[webs,->-=.5] (1,-.5) to (1,0);
    \draw[webs,->-=.5] (1,3) to (1,3.5);

    \draw[webs,->-=.5] (-1,-.5) to (-1,3.5);
    \draw[webs,->-=.9] (-1,3.5) to[out=90,in=270] (1,5);
    \draw[webs,->-=.9,over] (1,3.5) to[out=90,in=270] (-1,5);
    \end{tikzpicture}
    \caption{An example of a coherent tangle diagram which is not a braid closure.}
    \label{fig: coh_tang}
\end{figure}
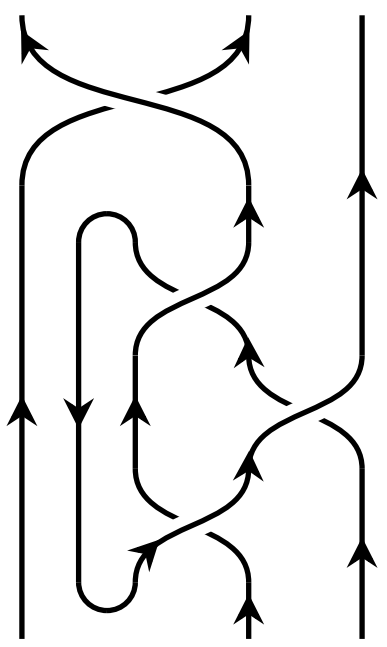

Rather than restricting to braid closures, we instead consider the collection of \textit{all} oriented tangle diagrams which can be obtained from the empty diagram by iteratively applying these operations. We refer to these diagrams as \textit{coherent}; see Figure \ref{fig: coh_tang} for an example of one such diagram which is not a braid closure.

\FloatBarrier

Our first novel result shows that $HHH_{alg}$ is well-behaved on this class of diagrams:

\begin{theorem} \label{thm: intro_invariance_coh}
    Let $\tau, \tau'$ be two coherent tangle diagrams representing isotopic oriented tangles. Then $HHH_{alg}(\tau) \cong HHH_{alg}(\tau')$.
\end{theorem}

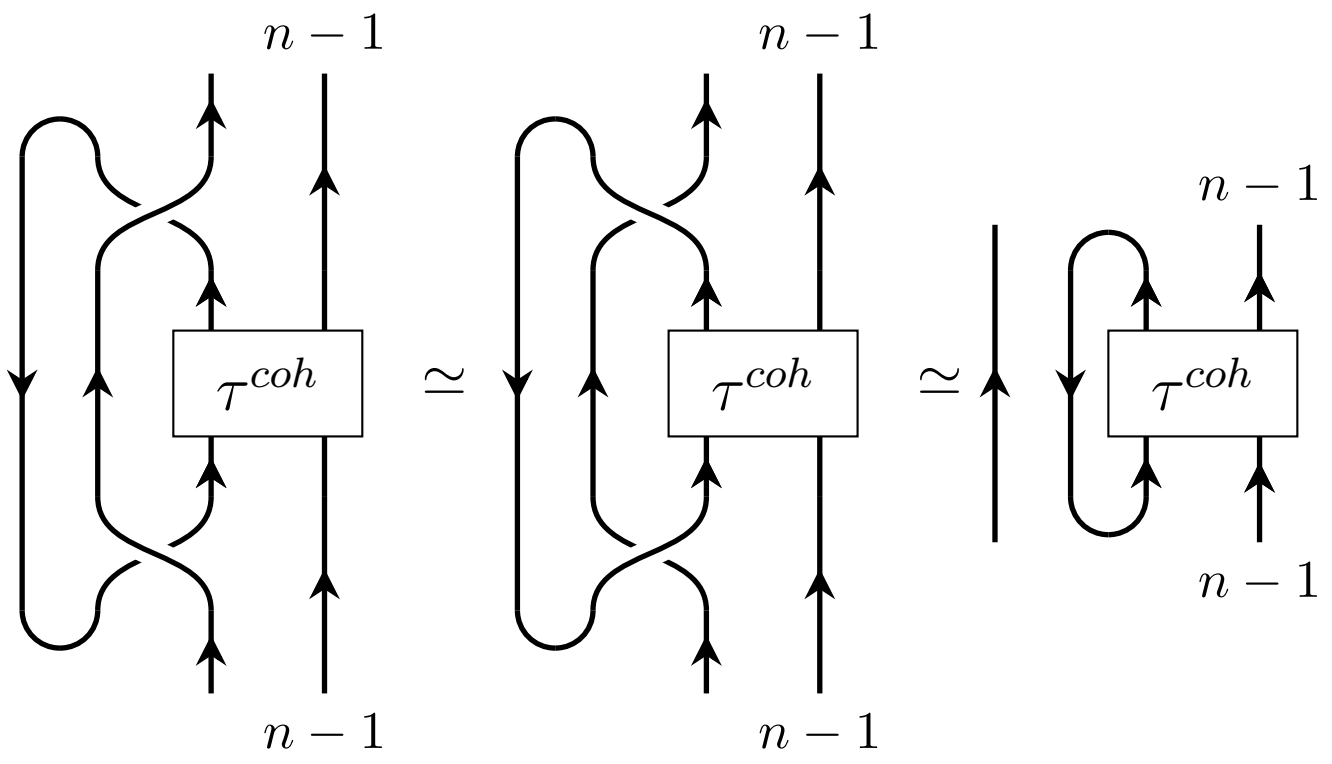
\begin{figure}[ht]
    \centering

    \begin{tikzpicture}[anchorbase]
        \draw[webs,->-=.1,->-=.9] (.25,0) to (.25,1.5);
        \draw[webs] (1,0) to (1,1.5);

        \filldraw[fill=white, draw=black]
            (0,.4) rectangle (1.25,1.1);

        \node at (.625,.75) {\Large{$\tau^{coh}$}};

        \draw[webs] (-.5,-.75) to[out=90,in=270] (.25,0);
        \draw[webs,over] (.25,-.75) to[out=90,in=270] (-.5,0);

        \draw[webs] (.25,1.5) to[out=90,in=270] (-.5,2.25);
        \draw[webs,over] (-.5,1.5) to[out=90,in=270] (.25,2.25);

        \draw[webs,->-=.5] (-.5,0) to (-.5,1.5);

        \draw[webs] (-.5,2.25) to[out=90,in=0] (-.75,2.5);
        \draw[webs] (-.75,2.5) to[out=180,in=90] (-1,2.25);
        \draw[webs,->-=.5] (-1,2.25) to (-1,-.75);
        \draw[webs] (-1,-.75) to[out=270,in=180] (-.75,-1);
        \draw[webs] (-.75,-1) to[out=0,in=270] (-.5,-.75);

        \draw[webs,-<-=.5] (.25,-.75) to (.25,-1.3);
        \draw[webs,->-=.5] (.25,2.25) to (.25,2.8);

        \draw[webs,-<-=.45] (1,0) to (1,-1.3) node[below]{$n-1$};
        \draw[webs,->-=.45] (1,1.5) to (1,2.8) node[above]{$n-1$};
    \end{tikzpicture}
    {\large{$\simeq$}} \ 
    \begin{tikzpicture}[anchorbase]
        \draw[webs,->-=.1,->-=.9] (.25,0) to (.25,1.5);
        \draw[webs] (1,0) to (1,1.5);

        \filldraw[fill=white, draw=black]
            (0,.4) rectangle (1.25,1.1);

        \node at (.625,.75) {\Large{$\tau^{coh}$}};

        \draw[webs] (.25,-.75) to[out=90,in=270] (-.5,0);
        \draw[webs,over] (-.5,-.75) to[out=90,in=270] (.25,0);
        
        \draw[webs] (-.5,1.5) to[out=90,in=270] (.25,2.25);
        \draw[webs,over] (.25,1.5) to[out=90,in=270] (-.5,2.25);
        
        \draw[webs,->-=.5] (-.5,0) to (-.5,1.5);

        \draw[webs] (-.5,2.25) to[out=90,in=0] (-.75,2.5);
        \draw[webs] (-.75,2.5) to[out=180,in=90] (-1,2.25);
        \draw[webs,->-=.5] (-1,2.25) to (-1,-.75);
        \draw[webs] (-1,-.75) to[out=270,in=180] (-.75,-1);
        \draw[webs] (-.75,-1) to[out=0,in=270] (-.5,-.75);

        \draw[webs,-<-=.5] (.25,-.75) to (.25,-1.3);
        \draw[webs,->-=.5] (.25,2.25) to (.25,2.8);

        \draw[webs,-<-=.45] (1,0) to (1,-1.3) node[below]{$n-1$};
        \draw[webs,->-=.45] (1,1.5) to (1,2.8) node[above]{$n-1$};
    \end{tikzpicture}
    {\large{$\simeq$}}
    \begin{tikzpicture}[anchorbase]
        \draw[webs,->-=.1,->-=.9] (.25,0) to (.25,1.5);
        \draw[webs,->-=.2,->-=.8] (1,-.3) node[below]{$n-1$} to (1,1.8) node[above]{$n-1$};

        \filldraw[fill=white, draw=black]
            (0,.4) rectangle (1.25,1.1);

        \node at (.625,.75) {\Large{$\tau^{coh}$}};

        \draw[webs] (.25,1.5) to[out=90,in=0] (0,1.75);
        \draw[webs] (0,1.75) to[out=180,in=90] (-.25,1.5);
        \draw[webs,->-=.5] (-.25,1.5) to (-.25,0);
        \draw[webs] (-.25,0) to[out=270,in=180] (0,-.25);
        \draw[webs] (0,-.25) to[out=0,in=270] (.25,0);

        \draw[webs,->-=.5] (-.75,-.3) to (-.75,1.8);
    \end{tikzpicture}

    \caption{The exchange move applied to an arbitrary coherent tangle diagram $\tau^{coh}$.}
    \label{fig: exchange_move_intro}
    \end{figure}

Our proof of Theorem \ref{thm: intro_invariance_coh} goes by two steps. The first is a Markov-type theorem giving a short list of moves which suffice to relate any two isotopic coherent tangle diagrams; these moves are listed in Theorem \ref{thm: markov_tangles_coherent}, and the interested reader should feel free to go read them now. Of these, the most interesting is the \textit{exchange move} depicted in Figure \ref{fig: exchange_move_intro}. The second step is to show that $HHH_{alg}$ respects each of these coherent Markov moves. This consists mostly of collecting existing results, with the glaring exception of the exchange move; the fact that $HHH_{alg}$ respects this move appears as Theorem \ref{thm: exchange} in the body of the paper, and its proof is more or less the only new algebraic input we require to run our construction.

\begin{figure}[ht]
\centering
\scalebox{.75}{%
    \begin{tikzpicture}[anchorbase]
        \filldraw[fill=white, draw=black]
            (0,-.5) rectangle (2,1.5);

        \node at (1,.5) {\Huge{$\tau$}};

        \draw[webs,-<-=.5] (0,.2) to (-.5,.2);
        \draw[webs,-<-=.5] (0,.8) to (-.5,.8);

        \draw[webs] (-.5,.2) to[out=180,in=90] (-1,-.5);
        \draw[webs] (-.5,.8) to[out=180,in=90] (-1.5,-.5);

        \draw[webs,-<-=.5] (-1,-.5) to (-1,-1.5);
        \draw[webs,-<-=.5] (-1.5,-.5) to (-1.5,-1.5);

        \draw[webs,->-=.5] (2,-.1) to (2.5,-.1);
        \draw[webs,->-=.5] (2,.3) to (2.5,.3);
        \draw[webs,-<-=.5] (2,.7) to (2.5,.7);
        \draw[webs,->-=.5] (2,1.1) to (2.5,1.1);

        \draw[webs] (2.5,-.1) to[out=0,in=90] (3,-.5);
        \draw[webs] (2.5,.3)  to[out=0,in=90] (3.5,-.5);
        \draw[webs] (2.5,.7)  to[out=0,in=90] (4,-.5);
        \draw[webs] (2.5,1.1)  to[out=0,in=90] (4.5,-.5);

        \draw[webs,->-=.5] (3,-.5) to (3,-1.5);
        \draw[webs,->-=.5] (3.5,-.5) to (3.5,-1.5);
        \draw[webs,-<-=.5] (4,-.5) to (4,-1.5);
        \draw[webs,-<-=.5] (4.5,-.5) to (4.5,-1.5);
    \end{tikzpicture}
    {\Huge{$\sim$}}
    \begin{tikzpicture}[anchorbase]

        \draw[shuffle,->-=.5] (0,3) to (0,4);
        \draw[shuffle,->-=.5] (1,3) to (1,4);
        \draw[shuffle,-<-=.5] (2,3) to (2,4);
        \draw[shuffle,-<-=.5] (3,3) to (3,4);
        \draw[shuffle,->-=.5] (5,3) to (5,4);
        \draw[shuffle,-<-=.5] (6,3) to (6,4);
        

        \draw[shuffle,-<-=.8]
            (2,4) to[out=90,in=270] (-2,7.5);

        \draw[shuffle,-<-=.8]
            (3,4) to[out=90,in=270] (-1,7.5);

        \draw[shuffle,-<-=.85]
            (6,4) to[out=90,in=270] (0,7.5);

        \draw[shuffle,->-=.85,over]
            (0,4) to[out=90,in=270] (6,7.5);

        \draw[shuffle,->-=.85,over]
            (1,4) to[out=90,in=270] (7,7.5);

        \draw[shuffle,->-=.8,over]
            (5,4) to[out=90,in=270] (8,7.5);


        \draw[register] (0,7.5) to[out=90,in=180] (.5,8);
        \draw[register] (-1,7.5) to[out=90,in=180] (.5,8.5);
        \draw[register] (-2,7.5) to[out=90,in=180] (.5,9);

        \draw[register] (6,7.5) to[out=90,in=0] (5.5,8);
        \draw[register] (7,7.5) to[out=90,in=0] (5.5,8.5);
        \draw[register] (8,7.5) to[out=90,in=0] (5.5,9);

        \draw[register,biarr] (.5,8) to (5.5,8);
        \draw[register,biarr] (.5,8.5) to (5.5,8.5);
        \draw[register,biarr] (.5,9) to (5.5,9);

        \filldraw[fill=white, draw=black]
            (1,7.5) rectangle (5,9.5);

        \node at (3,8.5) {\Huge{$\tau^{coh}$}};

    \end{tikzpicture}%
    }

    \caption{
        An arbitrary oriented tangle $\tau$ presented as a coherent tangle diagram $\tau^{coh}$ with an attached shuffle braid. Note that our convention for the orientation of coherent diagrams is rotated for ease of drawing.
    }
    \label{fig: intro_coh_shuffle}
\end{figure}
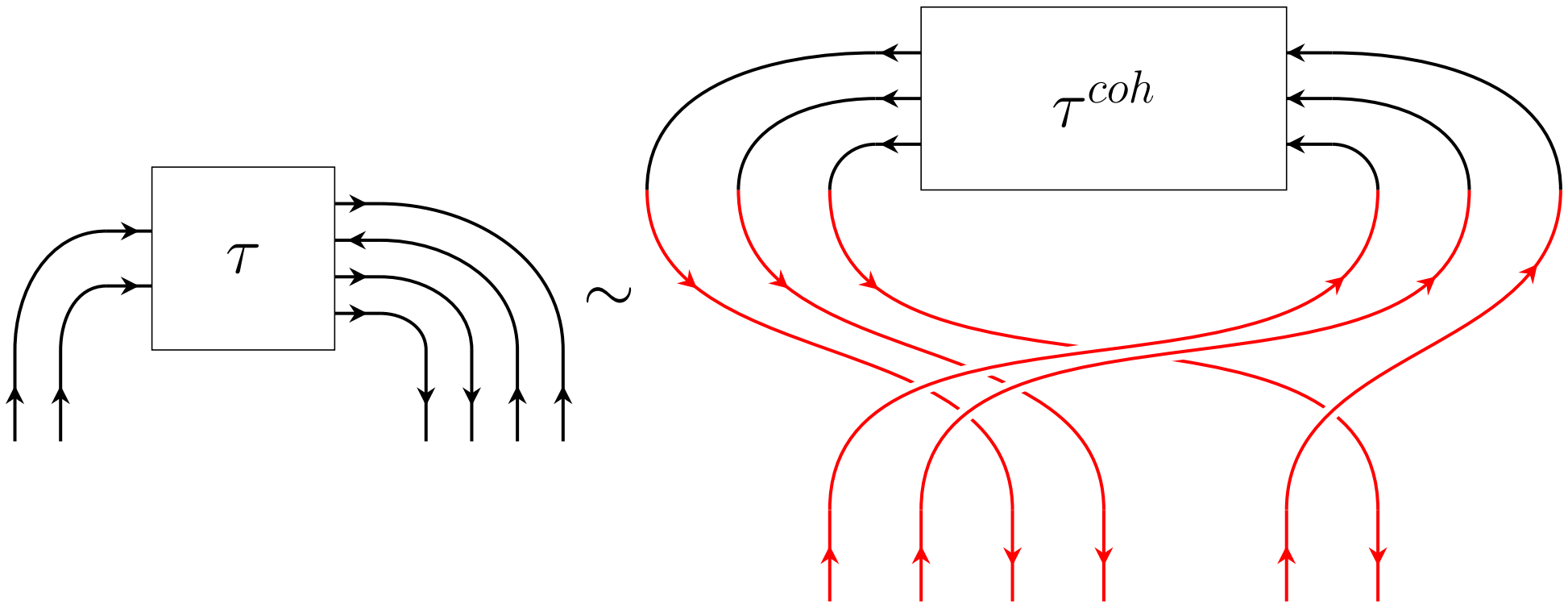

It is immediately clear that not all oriented tangle diagrams are coherent ones. The most obvious obstruction is the orientation of the boundary: all of the endpoints in any coherent tangle diagram must be upwards-oriented, while an arbitrary oriented tangle diagram can have both upwards- and downwards-oriented endpoints. Up to isotopy, this turns out to be the \textit{only} such obstruction; in fact any oriented tangle admits a diagram consisting of a coherent tangle diagram together with a shuffle braid depending only on the orientation of its boundary. See Figure \ref{fig: intro_coh_shuffle} for an example.

Because the shuffle braid depends only on the orientation of $\tau$ at the boundary, given any two oriented tangle diagrams $\tau$ and $\tau'$, we have $\tau \sim \tau'$ if and only if $\tau^{coh} \sim \tau'^{coh}$. As a consequence, we immediately obtain our desired tangle invariant (appearing as Theorem \ref{thm: main_invariance_HHH} in the body of the paper):

\begin{theorem} \label{thm: intro_invariance}
    The assignment $HHH \colon \tau \mapsto HHH_{alg}(\tau^{coh})$ is an invariant of oriented tangles in the thickened disk $\mathbb{D}^2 \times (-\epsilon, \epsilon)$ together with a marked point $p \in \partial \mathbb{D}^2$.
\end{theorem}

\begin{remark}
    The appearance of a marked point $p$ on the boundary may appear a bit strange at first glance, but it is unavoidable from our construction. To justify it, observe that Figure \ref{fig: intro_coh_shuffle} implicitly contains a choice of a ``leftmost" and ``rightmost" boundary point of $\tau$. After embedding $\tau$ in the thickened disk, such a choice is equivalent to a marked point $p$, which one can imagine as sitting above the tangle diagrams there.
\end{remark}

\begin{remark}
    While this paper was in development, Higgins proposed a construction of HOMFLY homology for the class of \textit{rational} tangles by considering a stable limit of the chain complexes assigned to them by $\mathfrak{sl}_N$ homology \cite{Hig26}. It would be very interesting to investigate whether the complexes assigned there agree with our invariant constructed here.
\end{remark}

\FloatBarrier

\subsubsection{Version II: Tangle Functors}

\textit{A priori}, the attentive reader might protest that our construction of $HHH$ suffers the same defect as Khovanov and Rozansky's $HHH_{alg}$. After all, we need to first isotope the entirety of any given tangle diagram $\tau$ to a coherent diagram $\tau^{coh}$ before applying $HHH_{alg}$, and at first glance, this is far from a local construction on crossings together with global gluing. Nevertheless, we can retroactively describe our invariant in terms of local crossing data and global gluing data, as we now outline.

The local data is easy: to any crossing $\tau$ embedded in a disk with marked boundary point, we can simply assign the complex $HHH(\tau)$ described above. The global data is subtler. Given an arc $\mathcal{A}$ connecting two disks containing tangle diagrams $\tau$ and $\tau'$, one would like to assign some functor $\Phi(\mathcal{A})$ which can be applied to the tensor product $HHH(\tau) \otimes HHH(\tau')$ to recover $HHH(\tau \cup \tau' \cup \mathcal{A})$. Further, $\Phi(\mathcal{A})$ should be independent of both the isotopy class of $\mathcal{A}$ within $\mathbb{D}^2$ and the tangles $\tau$ and $\tau'$, depending only on their boundaries. More broadly, given \textit{any} tangle diagram $\nu$ and \textit{any} way to attach $\nu$ to two disks containing diagrams $\tau$ and $\tau'$, we could hope to assign a functor $\Phi(\nu)$ satisfying the same criteria. This functor should again be independent of the isotopy class of $\nu$; two diagrams $\nu$ and $\nu'$ representing isotopic tangles should be assigned to naturally equivalent functors.

Put another way, we should hope to see a map $\Phi$ from the collection of oriented tangle diagrams to the category of endofunctors on the target of $HHH$ encoding all of the global gluing operations we might wish to perform, and $\Phi$ should take diagrams representing isotopic tangles to naturally isomorphic functors. Our final main result gives exactly such a construction (see Theorem \ref{thm: main_invariance} in the body of the paper for a more precise statement):

\begin{theorem} \label{thm: tangle_functor_intro}
    There is a well-defined functor $\Phi \colon \otan_{bal} \to \funcat$ which implements gluing operations on $HHH$.
\end{theorem}

Here $\otanbal$ denotes the category of \textit{balanced} (i.e. with equal numbers of upwards- and downwards-oriented endpoints on each boundary line) oriented tangles embedded in the thickened unit square, and $\funcat$ denotes a certain category of endofunctors on the universal target of our invariant; see Sections \ref{subsec: otan} and \ref{subsec: general_framework}, respectively, for further details.

In fact, our construction of $HHH(\tau)$ for a given diagram $\tau$ in the body of the paper proceeds directly by applying the functor $\Phi(\tilde{\tau})$ to the ground ring $\mathbb{Z}$, thought of as gluing the bent-around tangle diagram $\tilde{\tau}$ to the empty diagram. That this procedure recovers the above description in terms of coherent diagrams is then a nontrivial result, stated as Theorem \ref{thm: coherence} below.

We construct $\Phi$ essentially by a scanning argument. Given any (balanced) oriented tangle diagram $\tau$, we begin by expressing $\tau$ as a composition of elementary tangle diagrams (i.e. caps, cups, and crossings). We then imagine connecting $\tau$ via some shuffle braid to a pre-existing coherent diagram, thought of as encoding some algebraic gadget $C$ in our target category $K^b(\extsbim)$. For each such elementary diagram $\tau'$, we then specify some fixed isotopy through the shuffle braid. This has the effect of applying some topological operations on the coherent diagram, but in all cases, the resulting diagram will still be coherent. We can then translate these topological operations using $HHH_{alg}$ into some corresponding functor, and we assign exactly this functor as $\Phi(\tau')$; we assign as $\Phi(\tau)$ simply the composition of these functors for each elementary tangle constituting $\tau$.

\begin{remark} \label{rem: nice_construction_intro}
    This version of our construction has a few nice consequences. First, our Markov theorem for coherent tangle diagrams falls out immediately: each such diagram can be built by the same sort of scanning argument applied to a tangle diagram which is already coherent, and we can just check how the coherent diagrams produced in this way change under the application of each of the Reidemeister moves before scanning. Second, we automatically recover a description of our invariant in terms of local crossing data and global gluing operations. This endows the target category $K^b(\extsbim)$ with a \textit{canopolis} structure in the sense of Bar-Natan \cite{BN05}; we defer a detailed discussion of that structure for future work.
    
    Throughout this scanning procedure, a few steps introduce additional crossings not present in the original diagram, but only one introduces additional cups and caps: a sideways pointing crossing. As a consequence, if a tangle $\tau$ has Morse width $b$, the coherent tangle $\tau^{coh}$ that we produce has Morse width at most $b+1$. Going all the way to a braid can inflate the width arbitrarily, so this is a potentially significant aid to computation as the complexity of computing similar invariants tends to scale badly with diagram width. In particular, we can compute the HOMFLY homology of a knot with Morse width $b$ while never mentioning Soergel bimodules on more than $b+1$ strands.
\end{remark}

\subsection{Structure}

Section \ref{sec: background} consists of background regarding the topological setting of oriented tangles, the algebraic setting of Soergel bimodules, and Khovanov--Rozansky's $HHH_{alg}$. With the exception of Theorem \ref{thm: exchange} verifying the invariance of $HHH_{alg}$ under the exchange move, this section can safely be skipped by experts and referred back to for notation. Section \ref{sec: tangle_functors} is dedicated to the construction of the functor $\Phi$ of Theorem \ref{thm: tangle_functor_intro}. In Section \ref{sec: HHH}, we extract our oriented tangle invariant $HHH$ from this functor and describe its properties. We give a topological spin on these results in Appendix \ref{app: markov}, where our Markov-style theorem for coherent tangle diagrams is stated as Theorem \ref{thm: markov_tangles_coherent}.

\subsection{Acknowledgements}

We would like to thank Anna Beliakova, Eugene Gorsky, Jonathan Higgins, Mee Seong Im, Mikhail Khovanov, and Lev Rozansky for helpful conversations during the preparation of this work. This collaboration began as a result of discussions during the 2025 AMS Spring Eastern Sectional. Both authors were partially supported by the Simons Collaboration grant on New Structures in Low-Dimensional Topology, and LC was partially supported by a University of Zurich Postdoc Grant.

\subsubsection{LLM Use}

The authors used large language models for figure formatting and copyediting, and the argument contained in Figure~\ref{fig: ExchangeViaStabilization} was extracted from an example produced by Claude Sonnet 4.6.

\section{Background} \label{sec: background}

\subsection{Oriented Tangles} \label{subsec: otan}

\begin{definition} \label{def: balanced_tangles}
    Let $\aa \in \{+1,-1\}^m$ and $\bb \in \{+1, -1\}^n$ be two sequences. Then an \textit{oriented $(\aa, \bb)$-tangle} $\tau$ is a proper smooth embedding of a finite collection of disjoint arcs and circles into $\mathbb{R} \times [0,1] \times (-\epsilon, \epsilon)$ satisfying:

    \begin{itemize}
        \item The boundary points of the arcs map bijectively onto the sets $\{1, 2, \dots, m\} \times \{1\} \times \{0\}$ and $\{1, 2, \dots, n\} \times \{0\} \times \{0\}$, called respectively the \textit{top} and \textit{bottom} boundaries of $\tau$;

        \item Arranging $\aa$ (resp. $\bb$) along the top (resp. bottom) boundary of $\tau$, the positive and negative entries correspond with the positively and negatively oriented points.
    \end{itemize}

    We will also call $\tau$ an \textit{oriented $(m,n)$-tangle}.
\end{definition}

\begin{definition} \label{def: otan}
    The \textit{oriented tangle category} $\otan$ consists of the following data:

    \begin{itemize}
        \item \textbf{Objects:} Orientation sequences.

        \item \textbf{Morphisms from $\aa$ to $\bb$:} Isotopy classes rel. boundary of oriented $(\bb, \aa)$-tangles.

        \item \textbf{Composition:} Given an oriented $(\aa, \bb)$-tangle $\tau$ (i.e. a morphism from $\bb$ to $\aa$) and an oriented $(\bb, \cc)$-tangle $\tau'$ (i.e. a morphism from $\cc$ to $\bb$), their composition $\tau \star \tau'$ is the oriented $(\aa, \cc)$-tangle (i.e. the morphism from $\cc$ to $\aa$) given by vertically concatenating $\tau$ above $\tau'$ along their common boundary $\bb$ and rescaling the vertical coordinate appropriately\footnote{Note that this composition respects isotopy classes.}.
    \end{itemize}

    We will refer to all morphisms in $\otan$ as \textit{oriented tangles}. We equip $\otan$ with a monoidal structure, denoted $- \sqcup -$, given on objects and morphisms by horizontal concatenation and rescaling the horizontal coordinate appropriately\footnote{Note that this monoidal structure respects isotopy classes.}.
\end{definition}

By arranging a general oriented tangle in Morse position, we see that $\otan$ is generated up to monoidal product and composition by oriented crossings, caps, and cups. In fact, we can cut down on this generating set by restricting to \textit{vertical} crossings, in which both strands are oriented the same direction, by choosing once and for all a decomposition of the other, \textit{horizontal} crossings into cups, caps, and vertical crossings. See Figures \ref{fig: handles}, \ref{fig: vertical_crossings}, and \ref{fig: horizontal_crossing_decompositions} below.

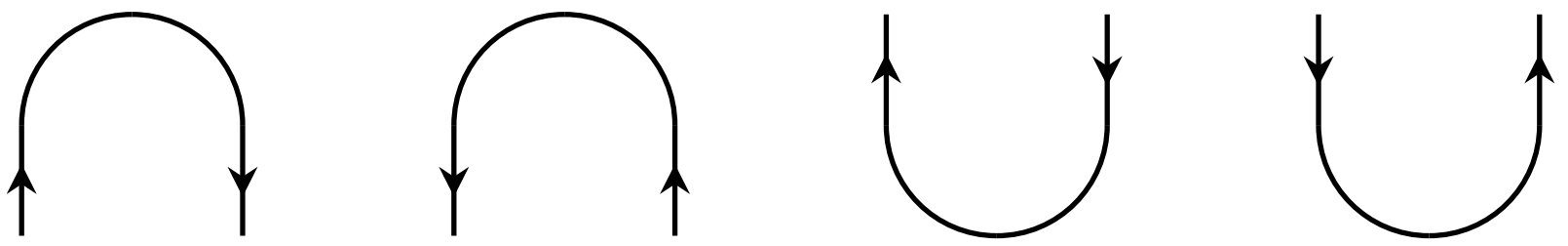
\begin{figure}[ht]
    \centering
    \begin{tikzpicture}
        \draw[webs,->-=.5] (0,0) to (0,.75);
        \draw[webs] (0,.75) to[out=90,in=180] (.75,1.5);
        \draw[webs] (.75,1.5) to[out=0,in=90] (1.5,.75);
        \draw[webs,->-=.5] (1.5,.75) to (1.5,0);
    \end{tikzpicture}
    \hspace{1cm}
    \begin{tikzpicture}
        \draw[webs,-<-=.5] (0,0) to (0,.75);
        \draw[webs] (0,.75) to[out=90,in=180] (.75,1.5);
        \draw[webs] (.75,1.5) to[out=0,in=90] (1.5,.75);
        \draw[webs,-<-=.5] (1.5,.75) to (1.5,0);
    \end{tikzpicture}
    \hspace{1cm}
    \begin{tikzpicture}
        \begin{scope}[yscale=-1]
        \draw[webs,-<-=.5] (0,0) to (0,.75);
        \draw[webs] (0,.75) to[out=90,in=180] (.75,1.5);
        \draw[webs] (.75,1.5) to[out=0,in=90] (1.5,.75);
        \draw[webs,-<-=.5] (1.5,.75) to (1.5,0);
        \end{scope}
    \end{tikzpicture}
    \hspace{1cm}
    \begin{tikzpicture}
        \begin{scope}[yscale=-1]
        \draw[webs,->-=.5] (0,0) to (0,.75);
        \draw[webs] (0,.75) to[out=90,in=180] (.75,1.5);
        \draw[webs] (.75,1.5) to[out=0,in=90] (1.5,.75);
        \draw[webs,->-=.5] (1.5,.75) to (1.5,0);
        \end{scope}
    \end{tikzpicture}

    \caption{A clockwise cup, counterclockwise cup, clockwise cap, and counterclockwise cap. Note that cups open downwards in our convention, as we read tangles from top to bottom.}
    \label{fig: handles}
\end{figure}

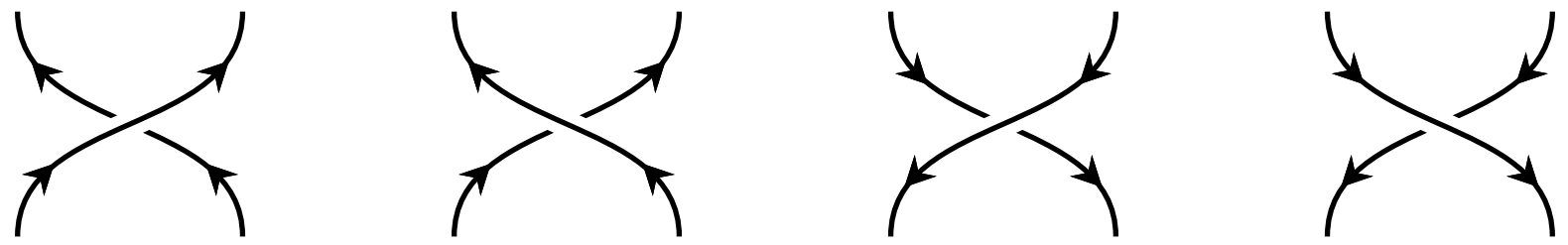
\begin{figure}[ht]
    \centering

    \begin{tikzpicture}
        \draw[webs,->-=.2,->-=.8] (1.5,0) to[out=90,in=270] (0,1.5);
        \draw[webs,->-=.2,->-=.8, over] (0,0) to[out=90,in=270] (1.5,1.5);
    \end{tikzpicture}
    \hspace{1cm}
    \begin{tikzpicture}
        \draw[webs,->-=.2,->-=.8] (0,0) to[out=90,in=270] (1.5,1.5);
        \draw[webs,->-=.2,->-=.8, over] (1.5,0) to[out=90,in=270] (0,1.5);
    \end{tikzpicture}
    \hspace{1cm}
    \begin{tikzpicture}
        \draw[webs,->-=.2,->-=.8] (0,1.5) to[in=90,out=270] (1.5,0);
        \draw[webs,->-=.2,->-=.8, over] (1.5,1.5) to[in=90,out=270] (0,0);
    \end{tikzpicture}
    \hspace{1cm}
    \begin{tikzpicture}
        \draw[webs,->-=.2,->-=.8] (1.5,1.5) to[in=90,out=270] (0,0);
        \draw[webs,->-=.2,->-=.8, over] (0,1.5) to[in=90,out=270] (1.5,0);
    \end{tikzpicture}

    \caption{Vertical crossings}
    \label{fig: vertical_crossings}
\end{figure}

\begin{figure}[ht]
    \centering

    \begin{tikzpicture}[anchorbase]
        \draw[webs, ->-=.2,->-=.8] (0,1.5) to[out=270, in=90] (1.5,0);
        \draw[webs, ->-=.2,->-=.8, over] (0,0) to[out=90,in=270] (1.5,1.5);
    \end{tikzpicture}
    {\Large{$:=$}}
    \begin{tikzpicture}[anchorbase]
    \draw[webs,->-=.5] (1,0) to (1,.5);
    \draw[webs,->-=.5] (0,1.5) to (0,.5);
    \draw[webs,->-=.5] (.5,1) to (.5,1.5);
    \draw[webs,->-=.5] (1.5,1) to (1.5,0);
    
    \draw[webs,->-=.99] (0,.5) to[out=270,in=180] (.25,.25);
    \draw[webs] (.25,.25) to[out=0,in=270] (.5,.5);

    \draw[webs] (.5,.5) to[out=90,in=270] (1,1);
    \draw[webs,over] (1,.5) to[out=90,in=270] (.5,1);
    
    \draw[webs,->-=.99] (1,1) to[out=90,in=180] (1.25,1.25);
    \draw[webs] (1.25,1.25) to[out=0,in=90] (1.5,1);

    \end{tikzpicture}
    \quad \quad
    \begin{tikzpicture}[anchorbase]
    \begin{scope}[xscale=-1,yscale=1]
        \draw[webs, ->-=.2,->-=.8] (0,1.5) to[out=270, in=90] (1.5,0);
        \draw[webs, ->-=.2,->-=.8, over] (0,0) to[out=90,in=270] (1.5,1.5);
    \end{scope}
    \end{tikzpicture}
    {\Large{$:=$}}
    \begin{tikzpicture}[anchorbase]
    \begin{scope}[xscale=-1,yscale=1]
    \draw[webs,->-=.5] (1,0) to (1,.5);
    \draw[webs,->-=.5] (0,1.5) to (0,.5);
    \draw[webs,->-=.5] (.5,1) to (.5,1.5);
    \draw[webs,->-=.5] (1.5,1) to (1.5,0);
    
    \draw[webs,->-=.99] (0,.5) to[out=270,in=180] (.25,.25);
    \draw[webs] (.25,.25) to[out=0,in=270] (.5,.5);

    \draw[webs] (.5,.5) to[out=90,in=270] (1,1);
    \draw[webs,over] (1,.5) to[out=90,in=270] (.5,1);
    
    \draw[webs,->-=.99] (1,1) to[out=90,in=180] (1.25,1.25);
    \draw[webs] (1.25,1.25) to[out=0,in=90] (1.5,1);

    \end{scope}
    \end{tikzpicture}
    \\
    \vspace{1cm}
    \begin{tikzpicture}[anchorbase]
        \draw[webs, -<-=.2,-<-=.8] (0,1.5) to[out=270, in=90] (1.5,0);
        \draw[webs, -<-=.2,-<-=.8, over] (0,0) to[out=90,in=270] (1.5,1.5);
    \end{tikzpicture}
    {\Large{$:=$}}
    \begin{tikzpicture}[anchorbase]
    \draw[webs,-<-=.5] (1,0) to (1,.5);
    \draw[webs,-<-=.5] (0,1.5) to (0,.5);
    \draw[webs,-<-=.5] (.5,1) to (.5,1.5);
    \draw[webs,-<-=.5] (1.5,1) to (1.5,0);
    
    \draw[webs,-<-=.99] (0,.5) to[out=270,in=180] (.25,.25);
    \draw[webs] (.25,.25) to[out=0,in=270] (.5,.5);

    \draw[webs] (.5,.5) to[out=90,in=270] (1,1);
    \draw[webs,over] (1,.5) to[out=90,in=270] (.5,1);
    
    \draw[webs,-<-=.99] (1,1) to[out=90,in=180] (1.25,1.25);
    \draw[webs] (1.25,1.25) to[out=0,in=90] (1.5,1);

    \end{tikzpicture}
    \quad \quad
    \begin{tikzpicture}[anchorbase]
    \begin{scope}[xscale=-1,yscale=1]
        \draw[webs, -<-=.2,-<-=.8] (0,1.5) to[out=270, in=90] (1.5,0);
        \draw[webs, -<-=.2,-<-=.8, over] (0,0) to[out=90,in=270] (1.5,1.5);
    \end{scope}
    \end{tikzpicture}
    {\Large{$:=$}}
    \begin{tikzpicture}[anchorbase]
    \begin{scope}[xscale=-1,yscale=1]
    \draw[webs,-<-=.5] (1,0) to (1,.5);
    \draw[webs,-<-=.5] (0,1.5) to (0,.5);
    \draw[webs,-<-=.5] (.5,1) to (.5,1.5);
    \draw[webs,-<-=.5] (1.5,1) to (1.5,0);
    
    \draw[webs,-<-=.99] (0,.5) to[out=270,in=180] (.25,.25);
    \draw[webs] (.25,.25) to[out=0,in=270] (.5,.5);

    \draw[webs] (.5,.5) to[out=90,in=270] (1,1);
    \draw[webs,over] (1,.5) to[out=90,in=270] (.5,1);
    
    \draw[webs,-<-=.99] (1,1) to[out=90,in=180] (1.25,1.25);
    \draw[webs] (1.25,1.25) to[out=0,in=90] (1.5,1);

    \end{scope}
    \end{tikzpicture}

    \caption{Building horizontal crossings out of vertical crossings and handles}
    \label{fig: horizontal_crossing_decompositions}
\end{figure}
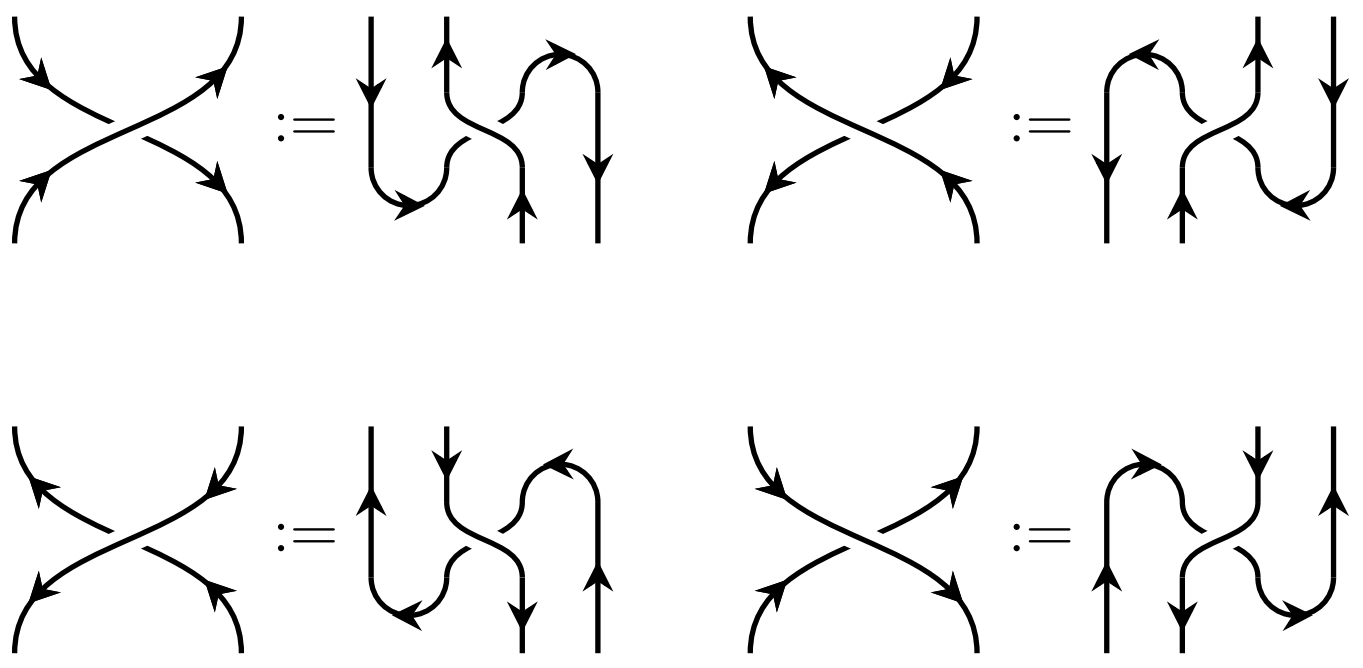

\FloatBarrier

We call a presentation of an oriented tangle as a composition of tensor products of (vertical or horizontal) crossings and handles an \textit{oriented tangle diagram}. Just as for link diagrams, there is a finite list of Reidemeister moves relating oriented tangle diagrams which represent isotopic oriented tangles. There are many such lists in the literature; here, we will use the one from Sawin's expository article \cite{Saw95}.

\begin{theorem} \label{thm: reidemeister_tangles}
    Two oriented tangle diagrams represent the same morphism in $\otan$ (i.e. isotopic oriented tangles) if and only if they are related by a finite sequence of the following local moves:
    
    \begin{figure}[ht]
    \centering

    \renewcommand{\arraystretch}{1.2}
    \begin{tabular}{@{}c@{\hspace{1.5em}}c@{}}

    \MoveCell{Distant interchange}{
        \begin{tikzpicture}[anchorbase]
            \draw[webs] (0,0) to (0,2);
            \draw[webs] (1.5,0) to (1.5,2);

            \filldraw[fill=white, draw=black] (-.5,1.1) rectangle (.5,1.9);
            \node at (0,1.5) {\Large{$\tau$}};

            \filldraw[fill=white, draw=black] (1,.1) rectangle (2,.9);
            \node at (1.5,.5) {\Large{$\tau'$}};
        \end{tikzpicture}
        \MoveArrow
        \begin{tikzpicture}[anchorbase]
            \draw[webs] (0,0) to (0,2);
            \draw[webs] (1.5,0) to (1.5,2);

            \filldraw[fill=white, draw=black] (-.5,.1) rectangle (.5,.9);
            \node at (0,.5) {\Large{$\tau$}};

            \filldraw[fill=white, draw=black] (1,1.1) rectangle (2,1.9);
            \node at (1.5,1.5) {\Large{$\tau'$}};
        \end{tikzpicture}
    }
    &
    \MoveCell{Reidemeister I}{
        \begin{tikzpicture}[anchorbase]
            \draw[webs] (0,0) to (0,.5);
            \draw[webs] (0,1) to (0,1.5);
            \draw[webs] (.5,.5) to[out=90,in=270] (0,1);
            \draw[webs,over] (0,.5) to[out=90,in=270] (.5,1);
            \draw[webs] (.5,1) to[out=90,in=180] (.75,1.25);
            \draw[webs] (.75,1.25) to[out=0,in=90] (1,1);
            \draw[webs] (.5,.5) to[out=270,in=180] (.75,.25);
            \draw[webs] (.75,.25) to[out=0,in=270] (1,.5);
            \draw[webs] (1,.5) to (1,1);
        \end{tikzpicture}
        \MoveArrow
        \begin{tikzpicture}[anchorbase]
            \draw[webs] (0,0) to (0,1.5);
        \end{tikzpicture}
        \MoveArrow
        \begin{tikzpicture}[anchorbase]
            \begin{scope}[xscale=-1]
            \draw[webs] (0,0) to (0,.5);
            \draw[webs] (0,1) to (0,1.5);
            \draw[webs] (.5,.5) to[out=90,in=270] (0,1);
            \draw[webs,over] (0,.5) to[out=90,in=270] (.5,1);
            \draw[webs] (.5,1) to[out=90,in=180] (.75,1.25);
            \draw[webs] (.75,1.25) to[out=0,in=90] (1,1);
            \draw[webs] (.5,.5) to[out=270,in=180] (.75,.25);
            \draw[webs] (.75,.25) to[out=0,in=270] (1,.5);
            \draw[webs] (1,.5) to (1,1);
            \end{scope}
        \end{tikzpicture}
    }

    \\[1.5em]

    \MoveCell{Zigzag}{
        \begin{tikzpicture}[anchorbase]
            \draw[webs] (0,0) to (0,.75);
            \draw[webs] (0,.75) to[out=90,in=180] (.25,1);
            \draw[webs] (.25,1) to[out=0,in=90] (.5,.75);
            \draw[webs] (.5,.75) to[out=270,in=180] (.75,.5);
            \draw[webs] (.75,.5) to[out=0,in=270] (1,.75);
            \draw[webs] (1,.75) to (1,1.5);
        \end{tikzpicture}
        \MoveArrow
        \begin{tikzpicture}[anchorbase]
            \draw[webs] (0,0) to (0,1.5);
        \end{tikzpicture}
        \MoveArrow
        \begin{tikzpicture}[anchorbase]
            \draw[webs] (0,0) to (0,.75);
            \draw[webs] (0,.75) to[out=90,in=0] (-.25,1);
            \draw[webs] (-.25,1) to[out=180,in=90] (-.5,.75);
            \draw[webs] (-.5,.75) to[out=270,in=0] (-.75,.5);
            \draw[webs] (-.75,.5) to[out=180,in=270] (-1,.75);
            \draw[webs] (-1,.75) to (-1,1.5);
        \end{tikzpicture}
    }
    &
    \MoveCell{Reidemeister II}{
        \begin{tikzpicture}[anchorbase]
            \draw[webs] (.5,0) to[out=90,in=270] (0,.75);
            \draw[webs,over] (0,0) to[out=90,in=270] (.5,.75);
            \draw[webs] (0,.75) to[out=90,in=270] (.5,1.5);
            \draw[webs,over] (.5,.75) to[out=90,in=270] (0,1.5);
        \end{tikzpicture}
        \MoveArrow
        \begin{tikzpicture}[anchorbase]
            \draw[webs] (0,0) to (0,1.5);
            \draw[webs] (.5,0) to (.5,1.5);
        \end{tikzpicture}
    }

    \\[1.5em]

    \MoveCell{Dial moves}{
        \begin{tabular}{@{}c@{}}
        \begin{tikzpicture}[anchorbase]
            \draw[webs] (0,0) to[out=90,in=270] (.75,.75);
            \draw[webs,over] (.75,0) to[out=90,in=270] (0,.75);
            \draw[webs] (0,.75) to (0,1.5);
            \draw[webs] (.75,.75) to[out=90,in=180] (1.125,1.125);
            \draw[webs] (1.125,1.125) to[out=0,in=90] (1.5,.75);
            \draw[webs] (1.5,.75) to (1.5,0);
        \end{tikzpicture}
        \MoveArrow
        \begin{tikzpicture}[anchorbase]
            \draw[webs] (0,0) to (0,.75);
            \draw[webs] (.75,.75) to[out=270,in=90] (1.5,0);
            \draw[webs,over] (.75,0) to[out=90,in=270] (1.5,.75);
            \draw[webs] (0,.75) to[out=90,in=180] (.375,1.125);
            \draw[webs] (.375,1.125) to[out=0,in=90] (.75,.75);
            \draw[webs] (1.5,.75) to (1.5,1.5);
        \end{tikzpicture}
        \\[1em]
        \begin{tikzpicture}[anchorbase]
            \draw[webs] (.75,0) to[out=90,in=270] (0,.75);
            \draw[webs,over] (0,0) to[out=90,in=270] (.75,.75);
            \draw[webs] (0,.75) to (0,1.5);
            \draw[webs] (.75,.75) to[out=90,in=180] (1.125,1.125);
            \draw[webs] (1.125,1.125) to[out=0,in=90] (1.5,.75);
            \draw[webs] (1.5,.75) to (1.5,0);
        \end{tikzpicture}
        \MoveArrow
        \begin{tikzpicture}[anchorbase]
            \draw[webs] (0,0) to (0,.75);
            \draw[webs] (.75,0) to[out=90,in=270] (1.5,.75);
            \draw[webs,over] (.75,.75) to[out=270,in=90] (1.5,0);
            \draw[webs] (0,.75) to[out=90,in=180] (.375,1.125);
            \draw[webs] (.375,1.125) to[out=0,in=90] (.75,.75);
            \draw[webs] (1.5,.75) to (1.5,1.5);
        \end{tikzpicture}
        \end{tabular}
    }
    &
    \MoveCell{Reidemeister III}{
        \begin{tikzpicture}[anchorbase]
            \draw[webs] (0,0) to[out=90,in=270] (.5,.5);
            \draw[webs,over] (.5,0) to[out=90,in=270] (0,.5);
            \draw[webs] (.5,.5) to[out=90,in=270] (1,1);
            \draw[webs,over] (1,.5) to[out=90,in=270] (.5,1);
            \draw[webs] (0,1) to[out=90,in=270] (.5,1.5);
            \draw[webs,over] (.5,1) to[out=90,in=270] (0,1.5);
            \draw[webs] (1,0) to (1,.5);
            \draw[webs] (0,.5) to (0,1);
            \draw[webs] (1,1) to (1,1.5);
        \end{tikzpicture}
        \MoveArrow
        \begin{tikzpicture}[anchorbase]
            \draw[webs] (.5,0) to[out=90,in=270] (1,.5);
            \draw[webs,over] (1,0) to[out=90,in=270] (.5,.5);
            \draw[webs] (0,.5) to[out=90,in=270] (.5,1);
            \draw[webs,over] (.5,.5) to[out=90,in=270] (0,1);
            \draw[webs] (.5,1) to[out=90,in=270] (1,1.5);
            \draw[webs,over] (1,1) to[out=90,in=270] (.5,1.5);
            \draw[webs] (0,0) to (0,.5);
            \draw[webs] (1,.5) to (1,1);
            \draw[webs] (0,1) to (0,1.5);
        \end{tikzpicture}
    }

    \end{tabular}

    \caption{A list of Reidemeister moves for oriented tangles. Each move is considered with all consistent choices of orientation on each strand.}
    \label{fig: reidemeister_moves}
\end{figure}
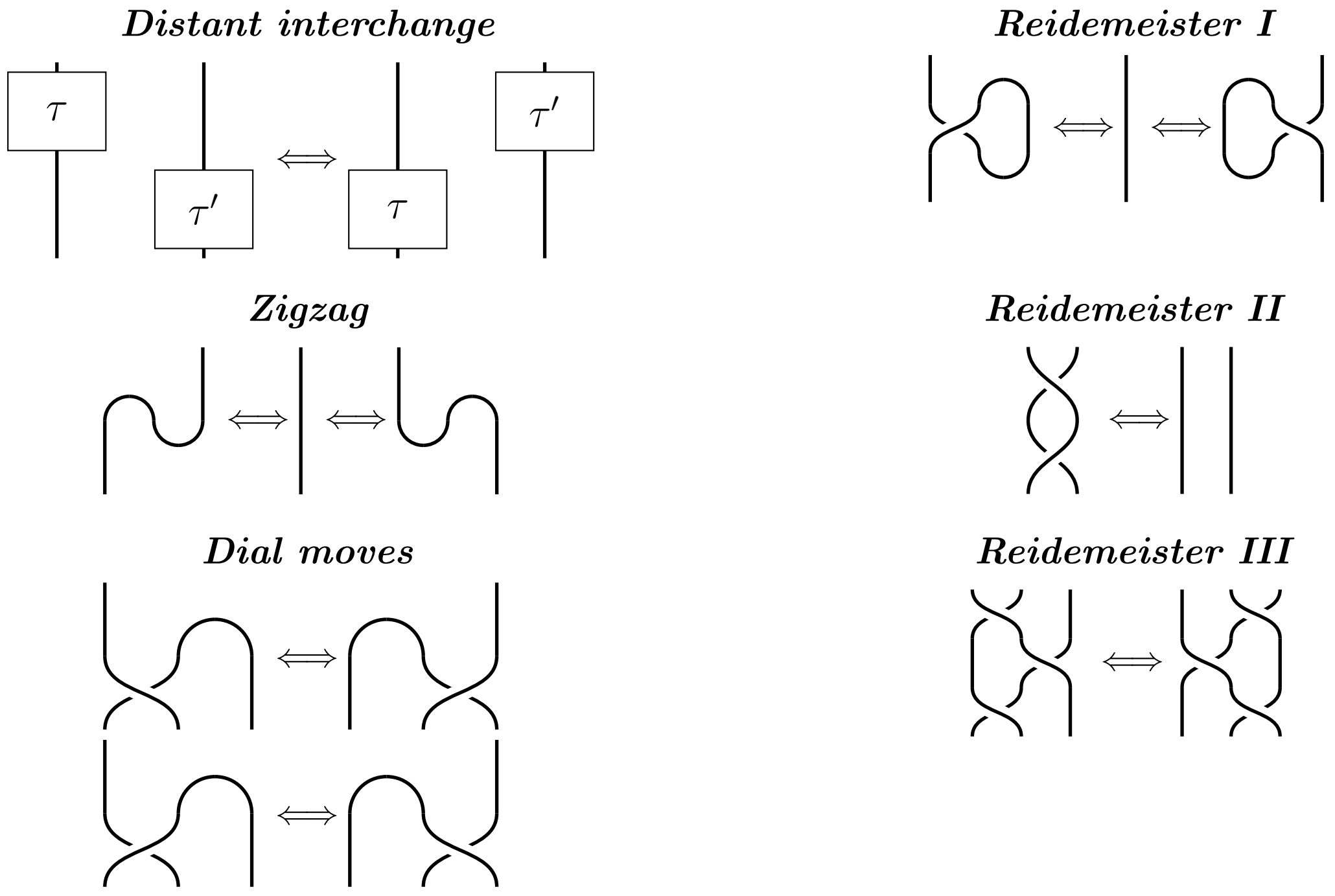
    
\end{theorem}

\FloatBarrier

We will occasionally make use of the following involution:

\begin{defn} \label{def: rotated_tangles}
    Given an oriented tangle $\tau$, we denote by $r(\tau)$ the oriented tangle obtained by rotating $\tau$ by $\pi$. Given a sequence $\aa = (a_1, \dots, a_m)$, we will also write $r(\aa) = (-a_m, \dots, -a_1)$, so that if $\tau$ is an oriented $(\aa, \bb)$-tangle, then $r(\tau)$ is an oriented $(r(\bb), r(\aa))$-tangle.
\end{defn}

We will be particularly interested in a few special classes of oriented tangle diagrams. The first class, consisting of \textit{braid diagrams}, will be useful in recalling the usual construction of HOMFLY homology.

\begin{definition} \label{def: braid}
    An \textit{$n$-strand braid diagram} is an oriented $(n,n)$-tangle diagram built entirely from vertical crossings between upwards-oriented strands. We denote by $\mathrm{Br}_n$ the \textit{$n$-strand braid group}, whose elements are $n$-strand braid diagrams considered up to isotopy.
\end{definition}

Note that any $n$-strand braid diagram has domain and codomain objects $(+1, \dots, +1)$. When considering oriented tangle diagrams with negatively oriented endpoints, we will typically restrict our attention to the class of \textit{balanced} oriented tangles, which we now define.

\begin{definition} \label{def: balanced_sequence}
    We call a sequence $\aa = (a_1, \dots, a_n) \in \{+1,-1\}^n$ \textit{balanced} if $\sum_{i = 1}^n a_i = 0$.
\end{definition}

\begin{definition} \label{def: otan_bal}
    The \textit{balanced oriented tangle category} $\otanbal$ is the full subcategory of $\otan$ on balanced sequences. We refer to all morphisms in $\otanbal$ as \textit{balanced oriented tangles}.
\end{definition}

\subsection{Soergel Bimodules} \label{subsec: soergel_bimodules}

In this section we recall the category of \textit{Soergel bimodules} introduced in \cite{Soe07}. For each $n \geq 1$, let $R_n = \mathbb{Z}[x_1, \dots, x_n]$ be a polynomial ring in $n$ variables (we set $R_0 = \mathbb{Z}$ for notational convenience). We endow $R_n$ with a $\mathbb{Z}$-grading by placing each variable $x_i$ in degree $2$; we refer to this grading throughout as the \textit{quantum} or \textit{internal} grading. To distinguish from other grading groups introduced later, we will denote by $q$ the generator of the quantum grading group and denote this grading group as $\mathbb{Z}_q$. We write $\mathbb{Z}_q$-gradings multiplicatively as powers of $q$, so that e.g. we have $\mathrm{deg}(x_i) = q^2$.

\begin{definition} \label{def: bimodule_1_cat}
    For each $n \geq 0$, we denote by $\Bim_n$ the category whose objects are $\mathbb{Z}_q$-graded $R_n$-bimodules and morphisms are $\mathbb{Z}_q$-graded modules of all homogeneous graded bimodule homomorphisms. We equip $\Bim_n$ with the standard monoidal structure given by $- \otimes_{R_n} -$ and denote the monoidal unit with respect to this structure (i.e. $R_n$ considered as a bimodule over itself by multiplication) by $\mathbb{1}_n$.
\end{definition}

\begin{notation} \label{not: grading_shift}
    The categories $\Bim_n$ are equipped with endofunctors which act on objects by formally shifting the quantum grading. We denote these functors also by powers of $q$, where $q^m$ denotes the functor which shifts the degree of each element \textit{up} by $q^m$. Explicitly, given a graded $R_n$-bimodule $M$ and an element $s \in M$ of degree $q^d$, the corresponding element $s \in q^m M$ has degree $q^{d+m}$.
\end{notation}

For each $1 \leq i \leq n - 1$, we denote by $B_i$ the graded $R_n$-bimodule

\[
B_i = q^{-1} R_n \otimes_{R_n^i} R_n,
\]
where $R_n^i$ is the subring of $R_n$ consisting of polynomials which are symmetric in the variables $x_i, x_{i+1}$. The category $\SBim_n$ of \textit{Soergel bimodules on $n$ strands} is the full subcategory of $\Bim_n$ generated by $R_n$ and $B_i$ for $1 \leq i \leq n - 1$ under (possibly infinite) direct sums, tensor product, grading shifts, and taking direct summands. Note that $\SBim_n$ inherits a monoidal structure from $\Bim_n$ with monoidal unit $\mathbb{1}_n$.


\begin{rem} \label{rem: explicit_soergel_bimodules}
It is convenient to consider $R_n$-bimodules as (left) modules over the enveloping algebra 
\[
R_n^e := R_n \otimes_\mathbb{Z} R_n^{op} \cong \mathbb{Z}[x_1, \dots, x_n, x_1', \dots, x_n'],
\]
where $\mathrm{deg}(x_i') = \mathrm{deg}(x_i) = q^2$, acting on the left via the variables $x_i$ and on the right via $x_i'$. Within this description, we have an explicit description of $B_i$ as the quotient

\[
B_i = q^{-1} \cfrac{\mathbb{Z}[x_1, \dots, x_n, x_1', \dots, x_n']}{\langle x_i + x_{i+1} - (x_i' + x_{i+1}'), \ x_ix_{i+1} - x_i'x_{i+1}' , x_j - x_j' \text{ (for }j \neq i , i+1)\rangle}.
\]

Note in particular that the shifted identity Soergel bimodule $q^{-1} R_n \in \SBim_n$ is a further quotient of $B_i \in \SBim_n$ by the collection $\langle x_j - x_j' \rangle_{1 \leq j \leq n}$. We denote the corresponding quotient map by $\unzip \colon B_i \to q^{-1} R_n$. As this terminology suggests, there is a corresponding map $\zip \colon q R_n \to B_i$ defined via $\zip(1) := x_i - x_{i+1}'$ and extended bilinearly. It is a useful exercise to check directly that, accounting for explicit grading shift, $\zip$ and $\unzip$ are grading-preserving bimodule homomorphisms. 
\end{rem}

The collection of monoidal ($1$-)categories $\SBim_n$ for all $n \geq 0$ naturally assembles into a monoidal $2$-category, as we now describe.

\begin{defn} \label{def: SBim_2_cat}
    We denote by $\SBim := \bigsqcup_{n \geq 0} \SBim_n$ the $2$-category consisting of the following data:

    \begin{itemize}
        \item Objects: Non-negative integers $n \geq 0$.

        \item Morphism categories from $m$ to $n$: $\SBim_m$ if $m = n$, empty otherwise. We denote (vertical) composition of $2$-morphisms within these $1$-morphism categories by $- \circ -$ (this is just the usual composition of graded bimodule homomorphisms).

        \item Horizontal composition of $1$- and $2$-morphisms, denoted $- \star-$, given by tensor product over the common coefficient ring $R_n$.
    \end{itemize}

    We will generally refer to all $1$-morphisms in $\SBim$ as \textit{Soergel bimodules}. We equip $\SBim$ with a monoidal structure $-\boxtimes -$ given on objects by $m \boxtimes n := m + n$ and on $1$- and $2$-morphisms by $- \otimes_\mathbb{Z} -$.
\end{defn}

There is a convenient diagrammatic description\footnote{Strictly speaking this description applies only to \textit{Bott-Samelson} bimodules, i.e. direct sums of tensor products of shifts of the bimodules $B_i$, and not more general Soergel bimodules which might only be direct summands of these. The distinction will not be relevant for this work.} of $\SBim$ paralleling the Type A webs of Cautis--Kamnitzer--Morrison \cite{CKM14}. Within this framework, Soergel bimodules on $n$ strands are depicted as certain upwards-oriented trivalent graphs embedded in $[0,1] \times [0,1]$ with $n$ boundary points lying on each of the intervals $[0,1] \times \{0\}$ and $[0,1] \times \{1\}$. We depict the bimodules $B_i \in \SBim_n$ as follows:

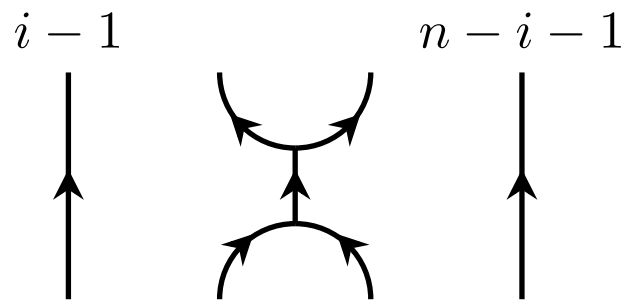
\begin{figure}[ht]
    \centering

    {\Large{$B_i := $}}
    \begin{tikzpicture}[anchorbase]
        \draw[webs,->-=.5] (0,0) to (0,1.5) node[above]{$i - 1$};
        \draw[webs,->-=.5] (1,0) to[out=90,in=180] (1.5,.5);
        \draw[webs,->-=.5] (2,0) to[out=90,in=0] (1.5,.5);
        \draw[webs,->-=.5] (1.5,.5) to (1.5,1);
        \draw[webs,->-=.5] (1.5,1) to[out=180,in=270] (1,1.5);
        \draw[webs,->-=.5] (1.5,1) to[out=0,in=270] (2,1.5);
        \draw[webs,->-=.5] (3,0) to (3,1.5) node[above]{$n - i - 1$};
    \end{tikzpicture}
\end{figure}

\FloatBarrier

Here and throughout, a single strand labeled with a non-negative integer $c \geq 0$ denotes a cable of $c$-parallel strands. We depict horizontal composition $\star$ within the graphical calculus as vertical concatenation along the common $n$ boundary points, written so that $M \star N$ has $M$ above $N$. Note that this gives a natural graphical description of $\mathbb{1}_n$ as $n$ parallel upwards-oriented vertical strands. We depict the monoidal product $\boxtimes$ on $\SBim$ as horizontal concatenation, written so that $M \boxtimes N$ has $M$ to the left of $N$. 

To extract topological information from $\SBim$, we will need to consider various categories of chain complexes over it. These homological enhancements will generally preserve objects, occurring at the level of $1$- and $2$-morphisms.

\begin{definition} \label{def: homotopy_cat}
We denote by $\mathrm{Ch}^b(\SBim)$ the $2$-category of (bounded) chain complexes over $\SBim$, whose objects are the same as those of $\SBim$ and morphism categories from $n$ to $n$ are given by $\mathrm{Ch}^b(\SBim_n)$, the category of (bounded) chain complexes over the additive category $\SBim_n$. We denote by $K^b(\SBim)$ the bounded homotopy $2$-category of Soergel bimodules, obtained from $\mathrm{Ch}^b(\SBim)$ by declaring homotopic $2$-morphisms (i.e. chain maps between complexes of Soergel bimodules) to be equivalent. We equip both $\mathrm{Ch}^b(\SBim)$ and $K^b(\SBim)$ with the monoidal structure induced from $\SBim$.

Note that $\mathrm{Ch}^b(\SBim)$ and $K^b(\SBim)$ carry an additional $\mathbb{Z}$-grading indexing the chain groups, which we call the \textit{homological} grading. We denote by $t$ the generator of this grading group, again writing $\mathbb{Z}_t$ for the group and denoting all gradings and formal grading shifts multiplicatively as powers of $t$.
\end{definition}

\begin{remark} \label{rem: trivial_homological_inclusion}
    There is a natural (fully faithful) inclusion $\SBim \hookrightarrow K^b(\SBim)$ given on objects as the identity, on $1$-morphisms $M$ (i.e. Soergel bimodules) by considering $M$ as a chain complex concentrated in homological degree $0$, and on $2$-morphisms by considering the corresponding chain maps. We will make implicit use of this inclusion throughout by considering chain groups of complexes in $K^b(\SBim)$ as explicit homological shifts of bimodules in $\SBim$. We will also extend the  diagrammatics for $\SBim$ along this inclusion, depicting more general chain complexes $C \in \mathrm{Ch}^b(\SBim_n)$ by depicting $C$ as a coupon spanning $n$ parallel upwards-oriented strands.
\end{remark}

\subsection{The Ext-Soergel Category}

Next, we recall a derived version of the categories $\SBim_n$ and $K^b(\SBim_n)$. Consider the bounded derived category $D^b(\Bim_n)$; here objects are chain complexes of $\Z_q$-graded $R_n$-bimodules, and morphisms are generated by chain maps over $\Bim_n$ and formal inverses of quasi-isomorphisms. There is once again a copy of $\Z$ which tracks the homological index of chain groups in $D^b(\Bim_n)$; we call this the \textit{Hochschild} or \textit{derived} grading. We denote by $a$ the generator of this grading group, again writing $\mathbb{Z}_a$ for the group and denoting gradings multiplicatively as powers of $a$. As in Notation \ref{not: grading_shift}, we also use powers of $a$ to denote formal grading shifts of $\mathbb{Z}_a$-graded objects.

As stated, all morphisms in $D^b(\Bim_n)$ preserve Hochschild degrees, and so morphism spaces are only $\Z_q$-graded. We prefer to work with a modest enhancement of these categories in which morphism spaces are $\Z_a \times \Z_q$-graded:

\begin{definition} \label{def: graded_derived_1cat}
    For each $n \geq 0$, let $D^b(\Bim_n)_{gr}$ denote the category whose objects are bounded chain complexes of $\Z_q$-graded $R_n$-bimodules, and whose morphism spaces are defined as follows:

    \[
    \Hom_{D^b(\Bim_n)_{gr}}(M, N) := \bigoplus_{i \in \Z} a^i \Hom_{D^b(\Bim_n)}(M, a^{-i}N).
    \]
\end{definition}

Explicitly, $D^b(\Bim_n)_{gr}$ has the same objects as $D^b(\Bim_n)$, but all higher (hyper-)Exts are included as morphisms with $a$-degree equal to their index. We now introduce the relevant subcategory of $D^b(\Bim_n)_{gr}$ in which we work:

\begin{defn} \label{def: ext_soergel_1_cat}
    For each $n \geq 0$, we denote by $\extsbim_n$ the full subcategory of $D^b(\Bim_n)_{gr}$ generated by finite direct sums of formal $\mathbb{Z}_a$-grading shifts of objects of $\SBim_n$ with trivial differential. We call $\extsbim_n$ the category of \textit{$n$-strand Ext-Soergel bimodules}.
\end{defn}

More explicitly, an $n$-strand Ext-Soergel bimodule is a finite direct sum of the form

\[
M = \bigoplus_j a^j M_j
\]
where each $M_j$ is an $n$-strand Soergel bimodule. All $n$-strand Ext-Soergel bimodules have trivial differential, so their morphism spaces are straightforward to describe using Ext groups: given $M = \bigoplus_j a^j M_j$ and $N = \bigoplus_k a^k N_k$, we have

\[
\Hom_{\extsbim_n}(M, N) = \bigoplus_{i \in \Z} a^i \Ext^i_{R^e_n}(M, N) := \bigoplus_{i \in \Z} \bigoplus_j \bigoplus_k  a^{i-j+k} \Ext^i_{R^e_n}(M_j, N_k).
\]
Here the overall $a$ shift is determined by observing that $\Ext^i_{R^e_n}(a^jM_j, a^kN_k) = \Ext^{i+j-k}_{R^e_n}(M_j, N_k)$ for all $M_j, N_k \in \SBim_n$.

\begin{warning} \label{warn: non_triangulated}
    It is common in the existing link homology literature to work instead with the larger category $D^b(\SBim_n)$, which is the full \textit{triangulated} subcategory of $D^b(\Bim_n)$ generated by $n$-strand Soergel bimodules. Objects there and in $\extsbim_n$ have the same underlying chain groups, but complexes in $D^b(\SBim_n)$ are allowed to have nontrivial differential. One consequence is that the category $\extsbim_n$ is not triangulated, but merely additive, leading to some (surmountable) technical difficulties (see e.g. Lemma \ref{lem: trace_ext_sbim}). We restrict to the smaller category $\extsbim_n$ because the exchange move (Theorem \ref{thm: exchange}) does not hold for arbitrary complexes in $D^b(\SBim_n)$ (see e.g. Example 5.6 from \cite{Abe17} for a counterexample involving a ``virtual crossing").
\end{warning}

As was the case for the underived categories $\SBim_n$, the $n$-strand Ext-Soergel categories assemble into a monoidal $2$-category as follows:

\begin{defn} \label{def: ext_soergel_2_cat}
    We denote by $\extsbim$ the $2$-category consisting of the following data:

    \begin{itemize}
        \item Objects: Non-negative integers $n \geq 0$.

        \item $1$-morphisms from $m$ to $n$: $n$-strand Ext-Soergel bimodules (i.e. finite direct sums of formal $\mathbb{Z}_a$-grading shifts of objects of $\SBim_n$) if $m = n$, empty otherwise.

        \item $2$-morphisms from $M$ to $N$ given by the $\mathbb{Z}_a \times \mathbb{Z}_q$-graded module

        \[
        \Hom_{\extsbim_n}(M,N): = \bigoplus_j a^j \Ext^j_{R_n^e}(M,N)
        \]
        with (vertical) composition $- \circ -$ inherited from $D^b(\Bim_n)$.

        \item Horizontal composition of $1$- and $2$-morphisms, denoted $- \star -$, given by tensor product over the common coefficient ring $R_n$.
    \end{itemize}

    We will generally refer to all $1$-morphisms in $\extsbim$ as \textit{Ext-Soergel bimodules}. We denote by $\extsbim_n$ the endomorphism category of the object $n \in \extsbim$, called the category of \textit{$n$-strand Ext-Soergel bimodules}. We equip $\extsbim$ with a monoidal structure $- \boxtimes -$ given on objects by $m \boxtimes n := m + n$ and on $1$- and $2$-morphisms by $- \otimes_\mathbb{Z} -$.
\end{defn}

\begin{remark} \label{rem: left_right_free}
    One could modify Definition \ref{def: ext_soergel_2_cat} by instead declaring horizontal composition to be the \textit{derived} tensor product over $R_n$; however, since all $n$-strand Soergel bimodules are free as left and right $R_n$-modules (cf. \cite{EK10}), the derived tensor product and ordinary tensor product agree in this setting. (There is a similar monoidal $2$-category $D^b(\Bim) = \bigsqcup_{n \geq 0} D^b(\Bim_n)$ in which horizontal composition is really the derived tensor product and \textit{not} the ordinary tensor product).
\end{remark}

For several of our computations, we will need a much more concrete description of $2$-morphisms in this category, and in particular of their horizontal compositions. This will take the form of a particular chain complex computing $\Hom_{\extsbim_n}(M,N)$. 

\begin{definition} \label{def: koszul_complex}
    For each $n \geq 0$, we denote by $K_n \in D^b(\Bim_n)$ the complex of free $R_n$-bimodules given as the dg algebra
    $$K_n = \mathbb{Z}[x_1,\dots,x_n, x_1',\dots,x_n', \eta_1,\dots,\eta_n]; \quad d(x_i) = d(x_i') = 0, \ \ d(\eta_i) = x_i - x_i'.$$

    Here $x_i$ and $x_i'$ commute with everything and the $\eta_i$ anticommute; the $R_n$-bimodule structure is by letting $R_n$ act on the left via multiplication by $x_i$ and on the right via multiplication by $x_i'$. $K_n$ is $\mathbb{Z}_a \times \mathbb{Z}_q$-graded by putting $x_i$ and $x_i'$ in degree $q^2$ as usual and $\eta_i$ in degree $a^{-1}q^2$, so the differential has degree $a$ as expected.
\end{definition}

Note that $K_n$ is simply the Koszul resolution of the identity $R_n$-bimodule $\mathbb{1}_n$. More precisely, there is a quasi-isomorphism of complexes of $R_n$-bimodules $\epsilon \colon K_n \to \mathbb{1}_n$ given in $a$-degree $0$ by the quotient by the ideal $\langle x_i - x_i'\rangle_{1 \leq i \leq n}$ and in all other $a$-degrees by $0$. We can therefore use $K_n$ to compute the Ext groups of Definition \ref{def: ext_soergel_2_cat} as follows:

\begin{prop} \label{prop: Hom_Complex}
    For any $n$-strand Ext-Soergel bimodules $M$ and $N$, we have an isomorphism of $\mathbb{Z}_a \times \mathbb{Z}_q$-graded modules
    
    \[
    \Hom_{\extsbim_n}(M,N) \cong H^\bullet(\Hom_{R_n^e}(K_n \star M, N)).
    \]

    (Here $\Hom_{R^e_n}(K_n \star M, N)$ denotes the chain complex of $R_n$-bimodules obtained by applying the functor $\Hom_{R^e_n}(-, N)$ to each chain group in $K_n \star M$.)
\end{prop}

\begin{proof}
    Recall from Remark \ref{rem: left_right_free} that any $n$-strand Soergel bimodule $M \in \SBim_n$ (and therefore any $n$-strand Ext-Soergel bimodule as well) is free as both a left and a right $R_n$-module. It follows that $K_n \star M$ is a free resolution of $\mathbb{1}_n \star M \cong M$ as an $R_n$-bimodule, and the desired result is then the standard construction of Ext groups using projective resolutions in the domain.
\end{proof}

Formally, the isomorphism of Proposition \ref{prop: Hom_Complex} realizes a morphism $f \colon M \to N$ as (the homology class of) the composite $\tilde{f} \circ (\epsilon \star \id_M)^{-1}$, where the latter is a formal inverse to the quasi-isomorphism $\epsilon \star \id_M$ realizing $K_n \star M$ as a projective resolution of $M$. We will typically identify each such morphism $f$ with its image $\tilde{f}$ under this isomorphism, calling $\tilde{f}$ a \textit{chain model} of $f$. Note that chain models are exactly characterized by the equation $\tilde{f} = f \circ (\epsilon \star \id_M)$; in particular, $\epsilon \star \id_M \colon K_n \star M \to M$ is a chain model for the identity $\id_M \in \Hom_{\extsbim_n}(M, M)$.

\begin{remark}
    The resolution $K_n \star M$ is typically not a minimal resolution of $M$, and certain calculations will be easier in a minimal resolution -- we will freely do this when it helps. However, working with the canonical resolution $K_n \star M$ becomes absolutely essential when we start taking horizontal compositions (tensor products) of morphisms: all calculations live in $K_n \star M_1 \star M_2$, so we need only to consider the ordinary tensor product $M_1 \star M_2$, not to simplify the (very infinite) tensor product of their minimal resolutions.
\end{remark}

Our next task is to identify the various structural operations of $2$-morphisms in $\extsbim$ (i.e. horizontal composition, vertical composition, and tensor products) with operations on chain models. Here we are aided by the existence of a coalgebra structure on $K_n$:

\begin{proposition} \label{prop: dg_coalg}
    There exist chain maps of $R_n$-bimodules $\epsilon: K_n \to \mathbb{1}_n$ and $\Delta: K_n \to K_n \star K_n$ making $K_n$ into a (strictly) counital, coassociative coalgebra in the chain category $\mathrm{Ch}^b(\Bim_n)$.
\end{proposition}

\begin{proof}
    The counit $\epsilon$ is the chain map involved in the projective resolution of $\mathbb{1}_n$; it sends $1$ to $1$ and any monomial in the $\eta_i$ to $0$, and is extended to the rest of $K_n$ by the fact that it must be a bimodule map (so e.g. $x_i$ and $x_i'$ are both sent to $x_i$).
    
    To describe the comultiplication $\Delta$, first observe that $K_n \star K_n$ itself carries an associative multiplication (with signs defined via the usual Koszul sign rule) given on pure tensors by the rule $(x_1 \otimes y_1) \cdot (x_2 \otimes y_2) := (-1)^{\deg_a(y_1) \deg_a(x_2)} (x_1x_2) \otimes (y_1 y_2)$. Using this multiplication, the action of $x_i$ on $K_n \star K_n$ is given by multiplication by $x_i \otimes 1$, while the action of $x_i'$ is given by multiplication by $1 \otimes x_i'$. It is a straightforward exercise to check that this multiplication satisfies the (graded) Leibniz rule $d(xy) = d(x)y + (-1)^{\deg_a(x)}x d(y)$.
    
    Now $\Delta$ sends $\eta_i$ to $\eta_i \otimes 1 + 1 \otimes \eta_i$, is extended to monomials in the $\eta_i$ by declaring $\Delta(\eta_i \eta_j) = \Delta(\eta_i) \Delta(\eta_j)$, and again is extended to the rest of $K_n$ using the fact that it is a bimodule map. For example, 
    \begin{equation*}
        \Delta(x_1x_2' \eta_1\eta_2) = x_1 \eta_1\eta_2 \otimes x_2' + x_1\eta_1 \otimes \eta_2 x_2' - x_1\eta_2 \otimes \eta_1 x_2' + x_1 \otimes \eta_1\eta_2 x_2'.
    \end{equation*}
    It is straightforward to check that $\Delta$ indeed defines a chain map, and it remains to check that $\Delta$ is coassociative and counital. To cut down on the amount of work required, it is helpful to first observe that since $\epsilon$ and $\Delta$ are both multiplication-preserving bimodule maps, these properties need only be checked for the generators $\eta_i$ and $1$. We leave this as an easy exercise.
\end{proof}

\begin{rem}
    We mention for the homologically devout that the coalgebra structure on $K_n$ mostly follows from general principles. Indeed, let $\epsilon \colon P \to \mathbb{1}$ be any bounded projective resolution of the monoidal unit in a monoidal dg category $\mathcal{C}$ in which tensor products of projective objects are projective (for example, the category of dg modules over a commutative ring). In particular, $\epsilon$ is a quasi-isomorphism, as is the tensor product $\mu \colon P \otimes P \xrightarrow{\epsilon \otimes \id} \mathbb{1} \otimes P \xrightarrow{\cong} P$. Now both $P \otimes P$ and $P$ are bounded complexes of projective objects, so by standard homological algebra (see e.g. \cite{Wei94}), $\mu$ is moreover a homotopy equivalence. It therefore admits a homotopy inverse $\Delta \colon P \to P \otimes P$, and one can show that $(P, \epsilon, \Delta)$ is a coalgebra object in the homotopy category $H^0(\mathcal{C})$. The advantage of working with the explicit structure maps of Proposition \ref{prop: dg_coalg} is that they turn $K_n$ into a \textit{strictly} counital coassociative coalgebra, in which the required coherence axioms hold on the nose rather than up to homotopy.
\end{rem}

We are now equipped to compute a chain model for vertical composition:

\begin{proposition} \label{prop: ext_horizontal_composites}
    Let $f \colon M_0 \to M_1$ and $g \colon M_1 \to M_2$ be morphisms in $\extsbim_n$ with chain models $\tilde{f} \colon K_n \star M_0 \to M_1$ and $\tilde{g} \colon K_n \star M_1 \to M_2$. Then the following is a chain model for the vertical composite $g \circ f:$

    \begin{equation*}
        K_n \star M_0 \xrightarrow{\Delta \star \id} K_n \star K_n \star M_0 \xrightarrow{\id \star \tilde{f}} K_n \star M_1 \xrightarrow{\tilde{g}} M_2.
    \end{equation*}
\end{proposition}

\begin{proof}
    Denote this composition by $\tilde{g \circ f}$; we need to verify that $\tilde{g \circ f} = (g \circ f) \circ (\epsilon \star \id_{M_0})$. Consider the following diagram:

    \begin{center}
    \begin{tikzcd} [row sep = large]
        K_n \star M_0 \arrow[r, "\Delta \star \id"] \arrow[d, "\epsilon \star \id"'] \arrow[dr, "\tilde{f}"] & K_n \star K_n \star M_0 \arrow[r, "\id \star \tilde{f}"] & K_n \star M_1 \arrow[d, "\tilde{g}"] \arrow[dl, "\epsilon \star \id"'] \\
        M_0 \arrow[r, "f"] & M_1 \arrow[r, "g"] & M_2
    \end{tikzcd}
    \end{center}

    The bottom-left and bottom-right triangles commute by the definition of chain models. Meanwhile, a direct computation gives

    \[
    (\epsilon \star \id_{M_1}) \circ (\id_{K_n} \star \tilde{f}) \circ (\Delta \star \id_{M_0}) = \tilde{f} \circ (((\epsilon \star \id_{K_n}) \circ \Delta) \star \id_{M_0}) = \tilde{f},
    \]
    where the second equality uses counitality of the coalgebra $K_n$. Hence the whole diagram commutes.
\end{proof}

To construct a chain model for horizontal composition $g \star f$ requires a brief technical detour:

\begin{lemma} \label{lem: koszul_exchange}
    For each $n$-strand Ext-Soergel bimodule $M$, there is a unique-up-to-homotopy chain homotopy equivalence $s_M \colon K_n \star M \to M \star K_n$ making the following diagram commute:

    \vspace{1em}
    \begin{center}
    \begin{tikzcd}
        K_n \star M \arrow[r, "s_M"] \arrow[dd, "\epsilon \star \id_M"'] & M \star K_n \arrow[dd, "\id_M \star \epsilon"] \\
        & \\
        M \arrow[r, "\id_M"] & M
    \end{tikzcd}
    \end{center}

    Further, the maps $s_M$ are natural up to homotopy in $M$, in the sense that for any morphism $f \colon M \to N$ of $n$-strand Ext-Soergel bimodules, the following diagram commutes up to homotopy:

    \vspace{1em}
    \centering
    \begin{tikzcd}
        K_n \star M \arrow[r, "s_M"] \arrow[dd, "\id_{K_n} \star f"'] & M \star K_n \arrow[dd, "f \star \id_{K_n}"] \\
        & \\
        K_n \star N \arrow[r, "s_N"] & N \star K_n
    \end{tikzcd}
\end{lemma}

\begin{proof}
    Recall that $K_n$ is a complex of free $R_n$-bimodules and $M$ is separately free as a left or as a right $R_n$ module, so both $K_n \star M$ and $M \star K_n$ are bounded complexes of free $R_n$-bimodules. By the usual comparison theorem for projective resolutions (see e.g. Chapter 2 of \cite{Wei94}), there exists a unique-up-to-homotopy chain map $s_M \colon K_n \star M \to M \star K_n$ making the first diagram (strictly) commute.
    
    Now, observe that $s_M$ is a chain model for the map $\tilde{s}_M = (\id_M \star \epsilon)^{-1} \colon M \to M \star K_n$ in the derived category of bimodules, where $(\id_M \star \epsilon)^{-1}$ is again a formal inverse to the quasi-isomorphism $\id_M \star \epsilon$. Naturality of $\tilde{s}_M$ in the derived category can be easily verified using the following diagram, all but the square cell of which automatically commute:

    
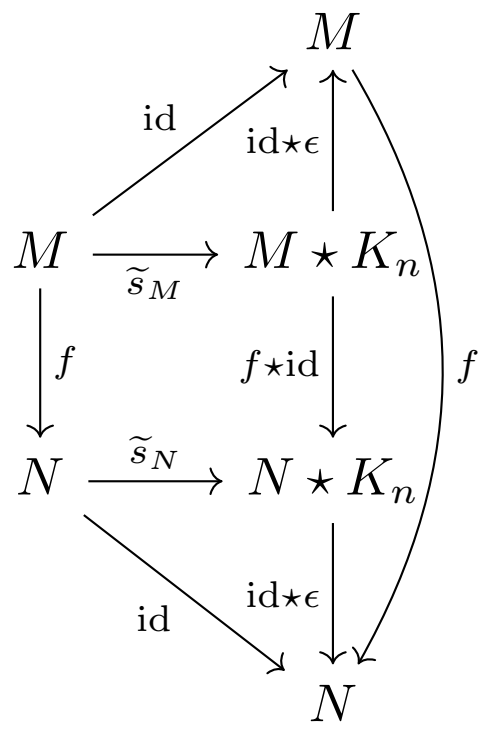
\begin{center}
    \begin{tikzcd} [row sep=large]
        & M \arrow[ddd, bend left, "f"] \\
        M \arrow[ur, "\id"] \arrow[r, "\tilde{s}_M"'] \arrow[d, "f"] & M \star K_n \arrow[u, "\id \star \epsilon"] \arrow[d, "f \star \id"'] \\
        N \arrow[dr, "\id"'] \arrow[r, "\tilde{s}_N"] & N \star K_n \arrow[d, "\id \star \epsilon"'] \\
        & N
    \end{tikzcd}
    \end{center}

    We conclude by observing that both $(f \star \id_{K_n}) \circ s_M$ and $s_N \circ (\id_{K_n} \star f)$ are chain models for the morphism $(f \star \id) \circ \tilde{s}_M = \tilde{s}_N \circ f \colon M \to N \star K_n$ in the derived category, and such chain models are unique up to homotopy.
\end{proof}

The map $s_M$ of Lemma \ref{lem: koszul_exchange} is easy to characterize, but is nontrivial to calculate for particular $M$. We will make essential use of this map for $M = B = B_1$ on two strands in later calculations, so we record it here.

\begin{lemma}\label{lem: Bswap}
    The map $s_B: K_2 \star B \to B \star K_2$ acts as follows:
    \begin{align*}
        1 \star 1 &\mapsto 1 \star 1 && 1 \star (x_1 - x_2') \mapsto (x_1 - x_2') \star 1\\
        \eta_1 \star 1 &\mapsto 1 \star \eta_2 && \eta_1 \star (x_1 - x_2') \mapsto (x_1-x_2') \star \eta_1 \\
        \eta_2 \star 1 &\mapsto 1\star \eta_1 && \eta_2 \star (x_1 - x_2') \mapsto (x_1-x_2')\star \eta_2\\
        \eta_1\eta_2 \star 1 &\mapsto 1\star (-\eta_1\eta_2) && \eta_1\eta_2 \star (x_1 - x_2') \mapsto (x_1-x_2') \star \eta_1\eta_2
    \end{align*}
    (These eight classes generate $K_2 \star B$ freely as a bimodule.)
\end{lemma}
\begin{proof}
    The first line guarantees that this $s_B$ commutes with the counit on $K_2$, so all that needs to be checked is that it really is a chain map. For the entries in the second column, this follows from the fact that $x_1(x_1-x_2') = (x_1-x_2')x_1'$. In the first, we have for example that
    \begin{align*}
        d(\eta_1 \star 1) = (x_1-x_1')\star 1 &= x_1 \star 1 - 1 \star x_1 \\
        &= x_1(1\star 1) - 1 \star (x_1-x_2') - (1 \star 1)x_2'\\
        &\mapsto x_1(1 \star 1) - (x_1-x_2')\star 1 - (1 \star 1)x_2'\\
        &=x_2' \star 1 - 1 \star x_2' = 1 \star (x_2 - x_2') = d(1 \star \eta_2)
    \end{align*}
\end{proof}

We are now equipped to give a chain model for horizontal composition:

\begin{prop} \label{prop: ext_vertical_composites}
    Let $f \colon M_0 \to N_0$ and $g \colon M_1 \to N_1$ be morphisms in $\extsbim_n$ with chain models $\tilde{f} \colon K_n \star M_0 \to N_0$ and $\tilde{g} \colon K_n \star M_1 \to N_1$. Then the following is a chain model for the horizontal composite $g \star f$:

    \[
    K_n \star (M_1 \star M_0) \xrightarrow{\Delta \star \id \star \id} K_n \star K_n \star M_1 \star M_0 \xrightarrow{\id \star s_{M_1} \star \id} (K_n \star M_1) \star (K_n \star M_0) \xrightarrow{\tilde{g} \star \tilde{f}} N_1 \star N_0.
    \]
\end{prop}

\begin{proof}
    Denote this composition by $\tilde{g \star f}$; we need to verify that $\tilde{g \star f} = (g \star f) \circ (\epsilon \star \id_{M_1} \star \id_{M_0})$. Recall that $g \star f = (\id \star f) \circ (g \star \id)$, and consider the following diagram:

    
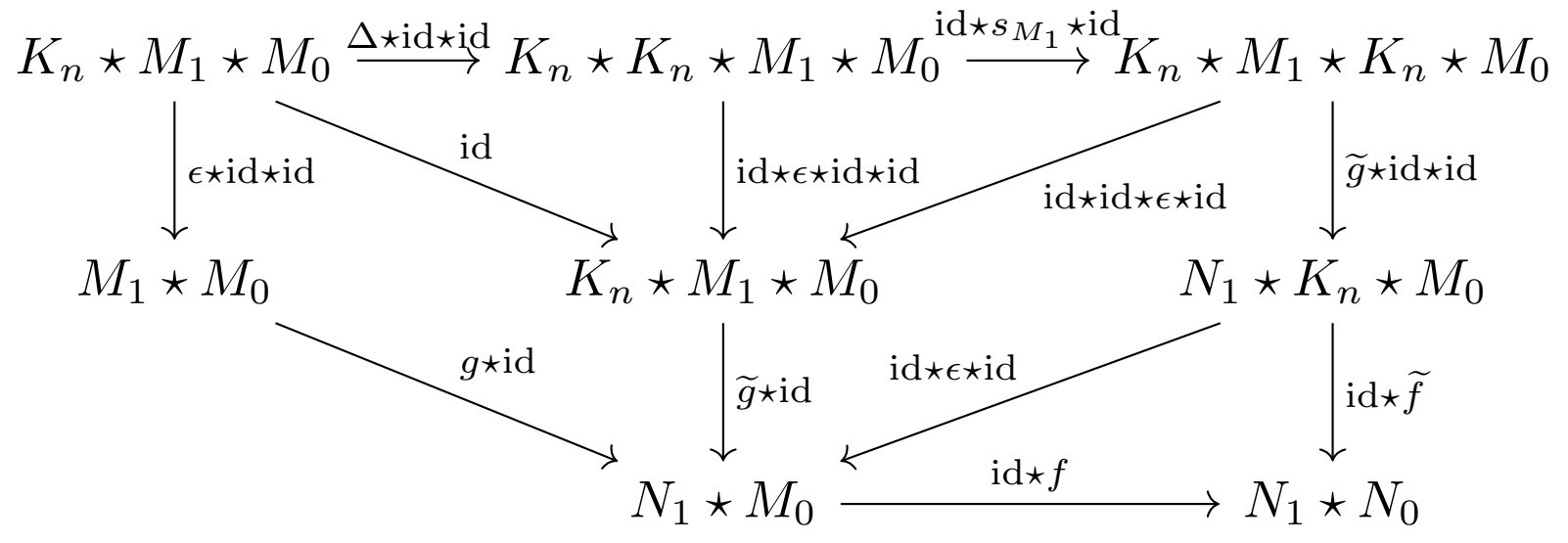
\begin{center}
        \begin{tikzcd}[row sep = large]
        K_n \star M_1 \star M_0 \arrow[r, "\Delta \star \id \star \id"] \arrow[d, "\epsilon \star \id \star \id"] \arrow[dr, "\id"] & K_n \star K_n \star M_1 \star M_0 \arrow[r, "\id \star s_{M_1} \star \id"] \arrow[d, "\id \star \epsilon \star \id \star \id"] & K_n \star M_1 \star K_n \star M_0 \arrow[d, "\tilde{g} \star \id \star \id"] \arrow[dl, "\id \star \id \star \epsilon \star \id"] \\
        M_1 \star M_0 \arrow[dr, "g \star \id"] & K_n \star M_1 \star M_0 \arrow[d, "\tilde{g} \star \id"] & N_1 \star K_n \star M_0 \arrow[d, "\id \star \tilde{f}"] \arrow[dl, "\id \star \epsilon \star \id"'] \\
        & N_1 \star M_0 \arrow[r, "\id \star f"] & N_1 \star N_0
        \end{tikzcd}
    \end{center}
    The entire diagram commutes by the definition of chain models, the counitality of $K_n$, and Lemma \ref{lem: koszul_exchange}.
\end{proof}

Finally, we turn to external tensor product $\boxtimes$. First, observe that for each $m, n \geq 0$, we have $K_m \boxtimes K_n = K_{m+n}$ by construction.

\begin{prop} \label{prop: ext_tensors}
    Let $f \colon M_0 \to N_0$ be a morphism in $\extsbim_m$ and $g \colon M_1 \to N_1$ be a morphism in $\extsbim_n$ with chain models $\tilde{f} \colon K_m \star M_0 \to N_0$ and $\tilde{g} \colon K_n \star M_1 \to N_1$. Then the following is a chain model for the tensor product $f \boxtimes g$:
    \begin{equation*}
        K_{m+n} \star (M_0 \boxtimes M_1) \xrightarrow{=} (K_m \boxtimes K_n) \star (M_0 \boxtimes M_1) \xrightarrow{\cong} (K_m \star M_0) \boxtimes (K_n \star M_1) \xrightarrow{\tilde{f} \boxtimes \tilde{g}} N_0 \boxtimes N_1.
    \end{equation*}
\end{prop}

\begin{proof}
    Analogous to the proofs of Propositions \ref{prop: ext_horizontal_composites} and \ref{prop: ext_vertical_composites}.
\end{proof}

As in Definition \ref{def: homotopy_cat}, we will again consider the bounded chain $2$-category $\mathrm{Ch}^b(\extsbim)$ and bounded homotopy $2$-category $K^b(\extsbim)$, whose objects are the same as those of $\extsbim$ and morphism categories are given by passing to the bounded chain categories and bounded homotopy categories of the additive categories $\extsbim_n$. We also equip both $\mathrm{Ch}^b(\extsbim)$ and $K^b(\extsbim)$ with the monoidal structures induced from $\extsbim$. Note that $\mathrm{Ch}^b(\extsbim)$ and $K^b(\extsbim)$ carry three distinct gradings, indexed by $\mathbb{Z}_a \times \mathbb{Z}_q \times \mathbb{Z}_t$.

\begin{remark} \label{rem: ext_embedding}
    There is a natural (non-full) inclusion $\SBim \hookrightarrow \extsbim$ given on objects as the identity, on $1$-morphisms $M$ (i.e. Soergel bimodules) by concentrating $M$ in Hochschild degree $0$, and on $2$-morphisms (i.e. bimodule homomorphisms) by considering the corresponding chain maps\footnote{This functor respects horizontal composition of $1$-morphisms precisely because of Remark \ref{rem: left_right_free}. In particular, the same prescription does \textit{not} define a $2$-functor $\Bim \hookrightarrow D^b(\Bim)$.}. This is an additive functor on morphism categories, so it immediately extends to an inclusion $K^b(\SBim) \hookrightarrow K^b(\extsbim)$ on the corresponding homotopy $2$-categories. As in Remark \ref{rem: trivial_homological_inclusion}, we will also extend the diagrammatic notation along this inclusion by again depicting an arbitrary complex of Ext-Soergel bimodules $C \in \mathrm{Ch}^b(\extsbim_n)$ as a coupon spanning $n$ parallel upwards-oriented strands.
\end{remark}

\begin{example}\label{ex: exthomnames}
    In $\extsbim_n$, the endomorphism ring of $\strand{n}$ is the (super)polynomial ring $\Z[x_1,\dots,x_n,\theta_1,\dots,\theta_n]$, with even variables $x$ in degree $q^2$ and odd variables $\theta_i$ in degree $aq^{-2}$. The map $\theta_i$ has a chain model $K_n \to \strand{n}$ given by the odd derivation $\epsilon \circ \partial / \partial \eta_i$. 

    Many (it would be reasonable to conjecture ``all" here) other morphisms of $\extsbim_n$ may be built out of tensor products and compositions of the $\theta_i$ and morphisms of $\SBim_n$. Most importantly for our purposes, we have an endomorphism $\theta_{i,i+1}$ of $B_i$ of degree $aq^{-4}$, given as the composite
    \begin{equation*}
        B_i \rightarrow B_i \star B_i \rightarrow B_i \star \strand{n} \star B_i \xrightarrow{1 \star \theta_i \star 1} B_i \star \strand{n} \star B_i \rightarrow B_i \star B_i \rightarrow B_i
    \end{equation*}
    of the coproduct on $B_i$, a unitor, $\theta_i$, and then another unitor and the product on $B_i$. Similarly, for any $M$, we have an endomorphism $\theta_i$ given by identifying $M$ with $\strand{n} \star M$ and acting by $\theta_i \star 1$, and on the other side an endomorphism $\theta_i'$ given by $1 \star \theta_i$ acting on $M \star \strand{n}$.
\end{example}

\begin{lemma}\label{lem:thetaidentities}
    The following identities hold:
    \begin{enumerate}
        \item $\unzip \circ \theta_{i,i+1} \circ \zip = \theta_i - \theta_{i+1}$, as endomorphisms of $\strand{n}$
        \item $(x_i-x_{i+1}) \theta_{i,i+1} = \theta_i - \theta_{i+1}$, as endomorphisms of $B_i$
        \item $(x_i-x_i')\theta_{i,i+1} = \theta_i - \theta_i'$, as endomorphisms of $B_i$
    \end{enumerate}
\end{lemma}
\begin{proof}
    To see the first identity, it is easier to work with an alternate chain model for $\theta_{i,i+1}$: there is a homotopy equivalence between $K_n \star B_i$ and the Koszul complex on $x_ix_{i+1} - x_i' x_{i+1}'$, $x_i+x_{i+1} - x_i' - x_{i+1}'$, and all other $x_j-x_j'$, and under this equivalence $\theta_{i,i+1}$ is taken to the odd derivation dual to the nullhomotopy of $x_ix_{i+1} - x_i'x_{i+1}'$. One now needs only to work out models for $\zip$ and $\unzip$ as maps between $K_n$ and this Koszul complex to see the result. The second two identities follow from this one and standard identities in $\SBim$, most importantly the fact that $\zip \circ \unzip = x_i - x_{i+1}'$.
\end{proof}

\subsection{Rouquier Complexes} \label{subsec: rouquier}

\begin{figure}[ht]
\centering

\Large{$\Br_4 \ni \sigma_1 := $}
\begin{tikzpicture}[anchorbase]
    \draw[webs,->-=.2,->-=.8] (1.5,0) to[out=90,in=270] (.5,1.5);
    \draw[webs,->-=.2,->-=.8, over] (.5,0) to[out=90,in=270] (1.5,1.5);
    \draw[webs,->-=.5] (2,0) to (2,1.5);
    \draw[webs,->-=.5] (2.5,0) to (2.5,1.5);
\end{tikzpicture}
\Large{; \quad $\Br_4 \ni \sigma_1^{-1} := $}
\begin{tikzpicture}[anchorbase]
    \draw[webs,->-=.2,->-=.8] (.5,0) to[out=90,in=270] (1.5,1.5);
    \draw[webs,->-=.2,->-=.8,over] (1.5,0) to[out=90,in=270] (.5,1.5);
    \draw[webs,->-=.5] (2,0) to (2,1.5);
    \draw[webs,->-=.5] (2.5,0) to (2.5,1.5);
\end{tikzpicture}

\caption{The positive and negative crossings $\sigma_1^{\pm} \in \Br_4$. Note that strand indices are read from left to right along the bottom boundary.}
\label{fig: crossings_reference}
\end{figure}
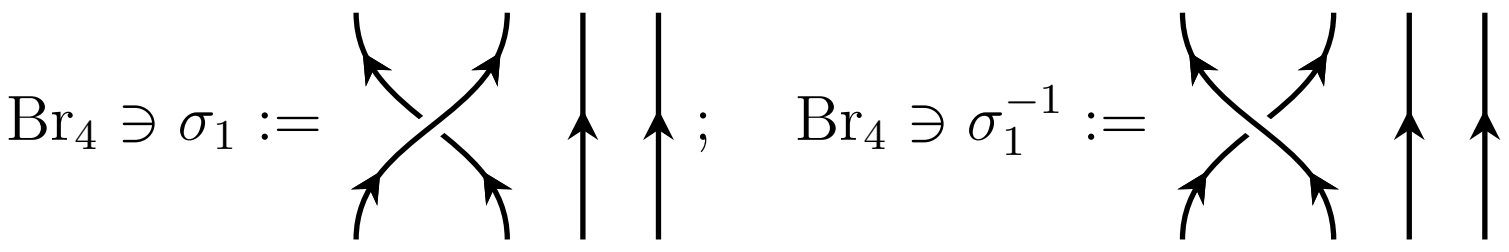

Recall the $n$-strand braid group $\Br_n$ of Definition \ref{def: braid}. We denote by $\sigma_i$ the $i^{th}$ \textit{positive} crossing, in which the $i^{th}$ strand passes over the $(i+1)^{st}$ strand, and by $\sigma_i^{-1}$ the corresponding \textit{negative crossing}; see Figure \ref{fig: crossings_reference} for an example. Given any braid $\beta \in \Br_n$, we can write $\beta$ as a product of upwards-oriented vertical crossings; we refer to the resulting expression $\beta = \sigma_{i_1}^{s_1} \star \dots \star\sigma_{i_l}^{s_l}$ as a \textit{braid word}.

Following Rouquier \cite{Rou04}, we now recall how to assign to each $n$-strand braid word $\beta$ a chain complex $F(\beta) \in \mathrm{Ch}^b(\SBim_n)$ of Soergel bimodules.

\begin{definition} \label{def: rouquier_complex}
The \textit{Rouquier complex} assigned to a single crossing $\beta = \sigma_i^{\pm} \in \Br_n$ is the complex of Soergel bimodules $F(\beta) \in \mathrm{Ch}^b(\SBim_n)$ depicted below:

\[
F(\sigma_i) := B_i \xrightarrow{\unzip} tq^{-1} R_n; \quad F(\sigma_i^{-1}) := t^{-1}q R_n \xrightarrow{\zip} B_i
\]

Given an arbitrary $n$-strand braid word $\beta = \sigma_{i_1}^{s_1} \star \dots \star \sigma_{i_l}^{s_l}$, we further set $F(\beta) := F(\sigma_{i_1}^{s_1}) \star \dots \star F(\sigma_{i_l}^{s_l}) \in \mathrm{Ch}^b(\SBim_n)$.
\end{definition}

\begin{figure}[ht]
    \centering
    
    \Large{$F(\sigma_3 \star \sigma_2 \star \sigma_3) \star C = $}
    \begin{tikzpicture}[anchorbase]
        \draw[webs,->-=.1,->-=.9] (0,0) to (0,1.5);
        \draw[webs,->-=.1,->-=.9] (.5,0) to (.5,1.5);
        \draw[webs,->-=.1,->-=.9] (1,0) to (1,1.5);
        \draw[webs,->-=.1,->-=.9] (1.5,0) to (1.5,1.5);
        \draw[webs,->-=.1,->-=.9] (2,0) to (2,1.5);

        \draw[webs,->-=.67] (0,1.5) to (0,3.75);
        \draw[webs,->-=.67] (2,1.5) to (2,3.75);

        \draw[webs] (1.5,3) to[out=90,in=270] (1,3.75);
        \draw[webs,over] (1,3) to[out=90,in=270] (1.5,3.75);
        \draw[webs] (1.5,1.5) to[out=90,in=270] (1,2.25);
        \draw[webs,over] (1,1.5) to[out=90,in=270] (1.5,2.25);
        \draw[webs,->-=.99] (1.5,2.25) to (1.5,3);
        \draw[webs,->-=.99] (1,2.25) to[out=90,in=270] (.5,3);
        \draw[webs,over,->-=.99] (.5,2.25) to[out=90,in=270] (1,3);
        \draw[webs] (.5,1.5) to (.5,2.25);
        \draw[webs] (.5,3) to (.5,3.75);

        \filldraw[fill=white, draw=black]
            (-.25,.4) rectangle (2.25,1.1);

        \node at (1,.75) {\Large{$C$}};
        
    \end{tikzpicture}
    \Large{;} \quad \Large{$F(\sigma_2 \star \sigma_3 \star \sigma_2) \star C = $}
    \begin{tikzpicture}[anchorbase]
        \draw[webs,->-=.1,->-=.9] (0,0) to (0,1.5);
        \draw[webs,->-=.1,->-=.9] (.5,0) to (.5,1.5);
        \draw[webs,->-=.1,->-=.9] (1,0) to (1,1.5);
        \draw[webs,->-=.1,->-=.9] (1.5,0) to (1.5,1.5);
        \draw[webs,->-=.1,->-=.9] (2,0) to (2,1.5);

        \draw[webs,->-=.67] (0,1.5) to (0,3.75);
        \draw[webs,->-=.67] (2,1.5) to (2,3.75);

        \draw[webs] (1,3) to[out=90,in=270] (.5,3.75);
        \draw[webs,over] (.5,3) to[out=90,in=270] (1,3.75);
        \draw[webs] (1,1.5) to[out=90,in=270] (.5,2.25);
        \draw[webs,over] (.5,1.5) to[out=90,in=270] (1,2.25);
        \draw[webs,->-=.99] (.5,2.25) to (.5,3);
        \draw[webs,->-=.99] (1.5,2.25) to[out=90,in=270] (1,3);
        \draw[webs,over,->-=.99] (1,2.25) to[out=90,in=270] (1.5,3);
        \draw[webs] (1.5,1.5) to (1.5,2.25);
        \draw[webs] (1.5,3) to (1.5,3.75);

        \filldraw[fill=white, draw=black]
            (-.25,.4) rectangle (2.25,1.1);

        \node at (1,.75) {\Large{$C$}};
        
    \end{tikzpicture}
 
    \caption{A diagrammatic representation of the compositions $F(\sigma_3 \star \sigma_2 \star \sigma_3) \star C$ and $F(\sigma_2 \star \sigma_3 \star \sigma_2) \star C$, where $C \in \mathrm{Ch}^b(\extsbim_5)$ is an arbitrary complex of Ext-Soergel bimodules. These two complexes become isomorphic in the homotopy category $K^b(\extsbim_5)$.}
    \label{fig: rouquier_diagram}
\end{figure}
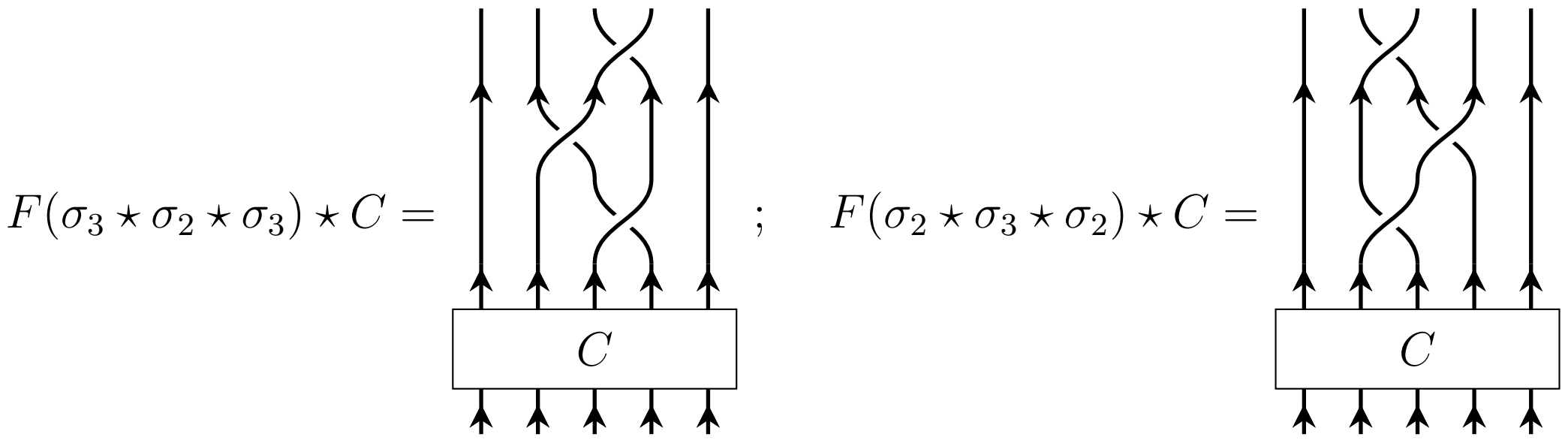

We can use Rouquier complexes to further extend the diagrammatics for Soergel bimodules introduced in Section \ref{subsec: soergel_bimodules}, depicting a Rouquier complex $F(\beta)$ by $\beta$ itself for some braid word $\beta$. See Figure \ref{fig: rouquier_diagram} for an example. There are obvious isotopies suggested by these diagrammatics corresponding to Reidemeister II and III moves in the braid words\footnote{To really complete the picture, one should also wish to slide coupons labeled by arbitrary complexes $C \in \mathrm{Ch}^b(\extsbim)$ over or under strands. This is an important step on the way to showing that Rouquier complexes equip $K^b(\extsbim)$ with the structure of a braided monoidal $2$-category. We neither prove nor rely on this statement, but see \cite{LMRSW24, SW24, EH24} for results in this direction.}; in Figure \ref{fig: rouquier_diagram}, these correspond to (the isotopy implied by) the braid relation $\sigma_3 \star \sigma_2 \star \sigma_3 \sim \sigma_2 \star \sigma_3 \star \sigma_2$. Ideally we would have the freedom to apply such isotopies while preserving the isomorphism class of the underlying chain complex. This is not true on the nose in the chain $2$-category $\mathrm{Ch}^b(\extsbim)$, but it becomes true upon passing to the homotopy $2$-category $K^b(\extsbim)$.

\begin{prop} \label{prop: rouquier_canonicity}
    Let $\beta, \beta'$ be two $n$-strand braid words representing the same $n$-strand braid. Then there are natural equivalences of functors

    \[
    F(\beta) \star - \cong F(\beta') \star -\colon K^b(\extsbim_n) \to K^b(\extsbim_n)
    \]
    and
    \[
    - \star F(\beta) \cong - \star F(\beta') \colon K^b(\extsbim_n) \to K^b(\extsbim_n). 
    \]
\end{prop}

\begin{proof}
    This was shown by Rouquier in \cite{Rou04} when restricted to the homotopy category of Soergel bimodules; since the homotopy equivalences involved are local, the result immediately extends to Ext-Soergel bimodules.
\end{proof}

In light of Proposition \ref{prop: rouquier_canonicity}, we will from now on work directly with the homotopy $2$-category $K^b(\extsbim)$ rather than its underlying chain $2$-category $\mathrm{Ch}^b(\extsbim)$. Given an $n$-strand braid $\beta \in \Br_n$ represented by a particular braid word $\underline{\beta}$, we will typically refer to the complex $F(\underline{\beta}) \in K^b(\SBim_n)$ as the \textit{Rouquier complex associated to $\beta$} and denote this complex by $F(\beta)$. It is easy to see directly from Definition \ref{def: rouquier_complex} that Rouquier complexes are compatible with the monoidal structure on the $2$-category $K^b(\extsbim)$, in the sense that $F(\beta) \boxtimes F(\beta') = F(\beta \sqcup \beta')$ for any braids $\beta$ and $\beta'$.

\subsection{Inclusion and Partial Trace} \label{subsec: inclusion_trace}

We now introduce the final ingredients required to extract topological information from Soergel bimodules: the \textit{inclusion} and \textit{partial trace} functors. Almost all of the results in this section are drawn directly from existing literature, especially Section 3 of \cite{Hog18}. The one exception is that our statement of Proposition \ref{prop: trace-like} seems to be more general than previous results, but this form is well-known to experts (and can be distilled from the proofs of similar existing results; see e.g. Lemma 4.11 of \cite{RT21}).

\begin{definition} \label{def: adding_strand}
    We denote by $I \colon K^b(\extsbim) \to K^b(\extsbim)$ the \textit{inclusion} $2$-functor $1 \boxtimes -$. Explicitly, $I$ acts on objects by sending $n$ to $1 + n$ for each $n \geq 0$, on $1$-morphisms by sending $M$ to $\strand{1} \otimes_\mathbb{Z} M$, and on $2$-morphisms by sending $f \colon M \to N$ to $\id_{\strand{1}} \otimes f \colon \strand{1} \otimes_\Z M \to \strand{1} \otimes_\Z N$.
\end{definition}

\begin{figure}[ht]
    \centering

    \Large{$I(C) := $}
    \begin{tikzpicture}[anchorbase]
        \draw[webs,->-=.5] (0,0) to (0,1.5);
        \draw[webs,->-=.1,->-=.9] (1,0) node[below]{$n$} to (1,1.5) node[above]{$n$};
        \filldraw[fill=white, draw=black]
            (.5,.4) rectangle (1.5,1.1);

        \node at (1,.75) {\Large{$C$}};
    \end{tikzpicture}

    \caption{A diagrammatic depiction of the complex $I(C) \in K^b(\extsbim_{n+1})$ for an arbitrary complex $C \in K^b(\extsbim_n)$.}
    \label{fig: inclusion}
\end{figure}
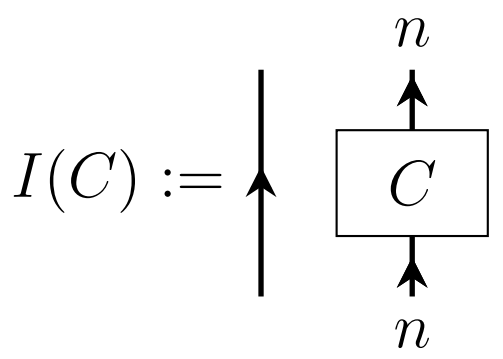

The diagrammatic description of $\SBim$, together with its extension to $K^b(\extsbim)$ via coupons, already encodes the monoidal structure, and it therefore gives a diagrammatic interpretation of the inclusion functor at the level of $1$-morphisms: we obtain $I(C)$ from $C$ by adding a single upwards-oriented strand on the left. See Figure \ref{fig: inclusion}.

By passing to the ordinary categories $K^b(\extsbim_n)$, we can think of the inclusion $2$-functor as a collection of ordinary functors $I_n \colon K^b(\extsbim_n) \to K^b(\extsbim_{n+1})$ for each $n \geq 0$. We now describe adjoints to these functors due to Hogancamp in \cite{Hog18}.

\begin{definition} \label{def: partial_trace}
    For each $n \geq 0$, we denote by $Tr_n \colon D^b(\Bim_{n+1}) \to D^b(\Bim_n)$ the \textit{partial trace functor}

    \[
    Tr_n(C) := C \xrightarrow{x_1 - x_1'} aq^{-2} C.
    \]

    Here the $R_n^e$ action on the right-hand side is given by the restriction of the $R_{n+1}^e$-action along the inclusion $\iota \colon R_n \to R_{n+1}$ sending $x_i$ to $x_{i+1}$ for each $1 \leq i \leq n$. Note that this shifts the index of each variable by one position.
\end{definition}

\begin{remark} \label{rem: partial_trace_well_defined}
    The attentive reader might have some concerns about Definition \ref{def: partial_trace} being a well-defined functor on $D^b(\Bim_{n+1})$, but in fact formation of mapping cones along central, closed morphisms preserves quasi-isomorphisms; see e.g. \cite{Con25}, Remark 2.3.
\end{remark}

The following key fact is necessary for all of our subsequent material:

\begin{lemma} \label{lem: trace_ext_sbim}
    Partial trace restricts to a functor $Tr_n \colon \extsbim_{n+1} \to \extsbim_n$.
\end{lemma}

Before giving the proof, we recall a key structural fact about Soergel bimodules, categorifying a recursive description of the Hecke algebra. This fact is exploited for similar purposes to ours in \cite{Ras15} (see Lemma~4.1 and Proposition~4.8), which itself cites Wu \cite{Wu06} for an elementary proof.
\begin{lemma} \label{lem: singleB1}
    Each $1$-morphism $M$ of $\extsbim$ is a direct sum of summands of $1$-morphisms of the form either $I(N)$ or $I(N) \star B_1 \star I(N')$, for $N,N'$ $1$-morphisms on one fewer strand than $M$.
\end{lemma}

\begin{proof}
    Every 1-morphism $M$ is a direct sum of summands of words in the $B_i$, and using the decompositions $B_i \star B_i \cong qB_i \oplus q^{-1} B_i$ and $B_i \star B_{i+1} \star B_i \oplus B_{i+1} \cong B_{i+1} \star B_i \star B_{i+1} \oplus B_i$ such a word is itself a direct sum of summands of words indexed by reduced words for elements of $S_n$, which may be chosen to use $B_1$ at most once.
\end{proof}

\begin{proof}[Proof of Lemma~\ref{lem: trace_ext_sbim}]
    By Lemma~\ref{lem: singleB1}, we may restrict to showing that traces of $I(N)$ or $I(N) \star B_1 \star I(N')$ land in $\extsbim$, and by Proposition~\ref{prop: trace_local} it suffices to show that $Tr_n(\strand{n+1})$ and $Tr_n(B_1)$ both lie in $\extsbim_n$. This is now a standard computation, dating back to at least \cite{KhR08}: both of these are (infinite) direct sums of shifts of $\strand{n}$, which may be identified with $\strand{n}[x,\theta]$ and $q\strand{n}[x,\theta^{(-1)}],$ respectively, where the old $x_1$ acts by $x$ in either case, $\theta$ records the action of the old $\theta_1$, and $\theta^{(-1)}$ records the action of the old $\theta_{1,2}$.
\end{proof}

\begin{proposition} \label{prop: trace_inclusion_adjoint}
    There is an adjunction $I_n \dashv Tr_n$ for each $n \geq 0$.
\end{proposition}

\begin{proof}
    See \cite{Hog18}, Proposition 3.6.
\end{proof}

\begin{example}
    We can now fully work out the action of $Tr_1$. On one-morphisms, as stated before $Tr_1(\strand{2}) = \strand{1}[x,\theta]$ and $Tr_1(B_1) = q\strand{1}[x,\theta^{(-1)}]$ are direct sums of copies of $\strand{1}$. To understand the morphisms, note that the adjunction tells us that summands of $Tr_1(B_1)$ are really indexed by summands of $\Hom(\strand{2},B_1),$ where $\zip$ is identified with $1$ and $\theta_{1,2} \circ \zip$ is identified with $\theta^{(-1)}$. On both of these objects, $Tr_1(x_1)$, $Tr_1(x_2)$, and $Tr_1(\theta_2)$ act as $x$, $x_1$, and $\theta_1$, respectively. For the other morphisms, we have:
    \begin{align*}
        Tr_1(\theta_1: \strand{2} \to \strand{2}) &= \theta \\
        Tr_1(\theta_1: B_1 \to B_1) &= (x-x_1)\theta^{(-1)} + \theta_1 \text{ (by Lemma \ref{lem:thetaidentities} (2))}\\
        Tr_1(\theta_{1,2}: B_1 \to B_1) &= \theta^{(-1)} \\
        Tr_1(\zip: \strand{2} \to B_1) &= \left(1 \mapsto 1, \theta \mapsto (x - x_1) \theta^{(-1)} + \theta_1 \right) \\
        Tr_1(\unzip: B_1 \to \strand{2}) &= \left(1 \mapsto x-x_1, \theta^{(-1)} \mapsto \theta - \theta_1 \right) \text{ (by Lemma \ref{lem:thetaidentities} (1))}
    \end{align*}
    where all rows should be extended to intertwine $x_1$ with $x$,  $x_2$ with $x_1$, and $\theta_2$ with $\theta_1$, and all right actions agree with left actions after taking the trace (this last feature is specific to $Tr_1$, since $\extsbim_1$ is symmetric monoidal unlike all larger $\extsbim_n$).
\end{example}

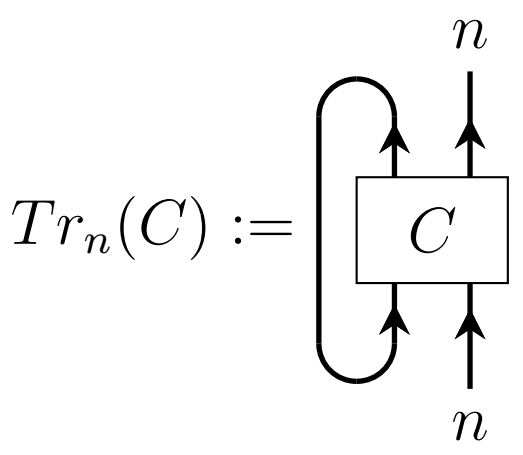
\begin{figure}[ht]
    \centering

    \Large{$Tr_n(C) := $}
    \begin{tikzpicture}[anchorbase]
        \draw[webs,->-=.1,->-=.9] (.5,0) to (.5,1.5);
        \draw[webs,->-=.2,->-=.8] (1,-.3) node[below]{$n$} to (1,1.8) node[above]{$n$};
        \draw[webs] (.5,1.5) to[out=90,in=0] (.25,1.75);
        \draw[webs] (.25,1.75) to[out=180,in=90] (0,1.5);
        \draw[webs] (0,1.5) to (0,0);
        \draw[webs] (0,0) to[out=270,in=180] (.25,-.25);
        \draw[webs] (.25,-.25) to[out=0,in=270] (.5,0);
        \filldraw[fill=white, draw=black]
            (.25,.4) rectangle (1.25,1.1);

        \node at (.75,.75) {\Large{$C$}};
    \end{tikzpicture}

    \caption{A diagrammatic depiction of the complex $Tr_n(C) \in K^b(\extsbim_n)$ for an arbitrary complex $C \in K^b(\extsbim_{n+1})$.}
    \label{fig: partial_trace}
\end{figure}

We can again incorporate the functor $Tr$ into the diagrammatic calculus\footnote{The primary motivation for this diagrammatic notation, as well as the name ``partial trace", is that $Tr$ categorifies the Markov trace used in Jones's construction of the HOMFLY polynomial (\cite{Jon87}).} as the closure of the two leftmost boundary points; see Figure \ref{fig: partial_trace}.

\FloatBarrier

There are again obvious isotopies suggested by these diagrammatics, and it is once again reasonable to insist that these isotopies should respect the isomorphism class of the underlying complex. We conclude this section by recalling a few such well-known results, mostly following Section 3 of \cite{Hog18}. We begin with an algebraic version of the distant interchange move:

\begin{prop} \label{prop: trace_local}
    Let $M, N \in K^b(\extsbim_n)$ and $C \in K^b(\extsbim_{n+1})$ be given. Then there is a natural (in $M$, $N$, and $C$) equivalence in $K^b(\extsbim_n)$

    \[
    Tr_n \left( (\strand{1} \boxtimes M) \star C \star (\strand{1} \boxtimes N)\right) \simeq M \star Tr_n(C) \star N.
    \]

    Diagrammatically:

    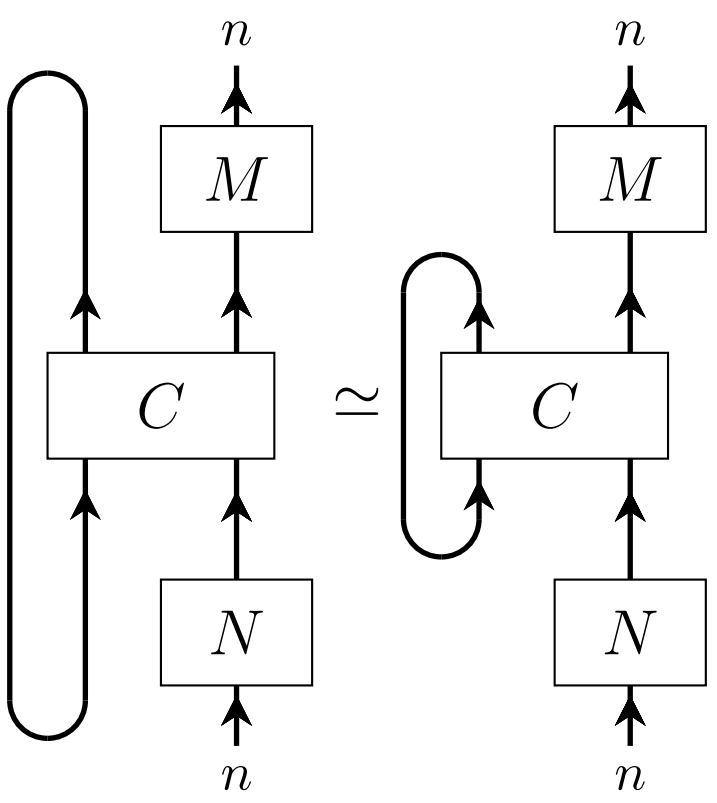
\begin{figure}[ht]
    \centering

    \begin{tikzpicture}[anchorbase]
        \draw[webs,->-=.33,->-=.67] (0,-1.2) to (0,2.7);
        \draw[webs,->-=.05,->-=.35,->-=.65,->-=.95] (1,-1.5) node[below]{$n$} to (1,3) node[above]{$n$};
        
        \draw[webs] (0,2.7) to[out=90,in=0] (-.25,2.95);
        \draw[webs] (-.25,2.95) to[out=180,in=90] (-.5,2.7);
        \draw[webs] (-.5,2.7) to (-.5,-1.2);
        \draw[webs] (-.5,-1.2) to[out=270,in=180] (-.25,-1.45);
        \draw[webs] (-.25,-1.45) to[out=0,in=270] (0,-1.2);
        
        \filldraw[fill=white, draw=black]
            (-.25,.4) rectangle (1.25,1.1);

        \node at (.5,.75) {\Large{$C$}};
        
        \filldraw[fill=white, draw=black]
            (.5,1.9) rectangle (1.5,2.6);

        \node at (1,2.25) {\Large{$M$}};
        
        \filldraw[fill=white, draw=black]
            (.5,-.4) rectangle (1.5,-1.1);

        \node at (1,-.75) {\Large{$N$}};
    \end{tikzpicture}
    {\Large{$\simeq$}}
    \begin{tikzpicture}[anchorbase]
        \draw[webs,->-=.1,->-=.9] (0,0) to (0,1.5);
        \draw[webs,->-=.05,->-=.35,->-=.65,->-=.95] (1,-1.5) node[below]{$n$} to (1,3) node[above]{$n$};
        
        \draw[webs] (0,1.5) to[out=90,in=0] (-.25,1.75);
        \draw[webs] (-.25,1.75) to[out=180,in=90] (-.5,1.5);
        \draw[webs] (-.5,1.5) to (-.5,0);
        \draw[webs] (-.5,0) to[out=270,in=180] (-.25,-.25);
        \draw[webs] (-.25,-.25) to[out=0,in=270] (0,0);
        
        \filldraw[fill=white, draw=black]
            (-.25,.4) rectangle (1.25,1.1);

        \node at (.5,.75) {\Large{$C$}};
        
        \filldraw[fill=white, draw=black]
            (.5,1.9) rectangle (1.5,2.6);

        \node at (1,2.25) {\Large{$M$}};
        
        \filldraw[fill=white, draw=black]
            (.5,-.4) rectangle (1.5,-1.1);

        \node at (1,-.75) {\Large{$N$}};
    \end{tikzpicture}

    \caption{A diagrammatic version of the distant interchange isotopy.}
    \label{fig: trace_local}
    \end{figure}

    \FloatBarrier
\end{prop}

\begin{proof}
    Clear from the definition of $Tr$.
\end{proof}

Next, we have an algebraic version of the Reidemeister I move:

\begin{prop} \label{prop: markov}
    For each $n \geq 0$, there are natural equivalences of functors

    \[
    Tr_n \circ (I_n(-) \star F(\sigma_1)) \cong \ tq^{-1} \mathrm{Id} \colon K^b(\extsbim_n) \to K^b(\extsbim_n)
    \]
    and
    \[
    Tr_n \circ (I_n(-) \star F(\sigma_1^{-1})) \cong aq^{-3} \ \mathrm{Id} \colon K^b(\extsbim_n) \to K^b(\extsbim_n).
    \]

    Diagrammatically, for each $C \in K^b(\extsbim_n)$:

    \begin{figure}[ht]
    \centering

    \begin{tikzpicture}[anchorbase]
        \draw[webs,->-=.8] (.25,0) to (.25,1.8);
        \draw[webs,->-=.15,->-=.885] (1,-1.3) node[below]{$n-1$} to (1,1.8) node[above]{$n-1$};

        \filldraw[fill=white, draw=black]
            (0,.4) rectangle (1.25,1.1);

        \node at (.625,.75) {\Large{$C$}};

        \draw[webs,->-=.2,->-=.8] (.25,-1) to[out=90,in=270] (-.5,0);
        \draw[webs,->-=.2,->-=.8,over] (-.5,-1) to[out=90,in=270] (.25,0);

        \draw[webs,->-=.5] (-.5,0) to (-.5,1.5);
        \draw[webs] (-.5,1.5) to[out=90,in=0] (-.75,1.75);
        \draw[webs] (-.75,1.75) to[out=180,in=90] (-1,1.5);
        \draw[webs] (-1,1.5) to (-1,-1);
        \draw[webs] (-1,-1) to[out=270,in=180] (-.75,-1.25);
        \draw[webs] (-.75,-1.25) to[out=0,in=270] (-.5,-1);

        \draw[webs] (.25,-1.3) to (.25,-1);
    \end{tikzpicture}
    {\Large{$\simeq tq^{-1}$}} \ 
    \begin{tikzpicture}[anchorbase]
        \draw[webs,->-=.2,->-=.8] (1,0) node[below]{$n$} to (1,2.5) node[above]{$n$};
        \filldraw[fill=white, draw=black]
            (.5,.9) rectangle (1.5,1.6);

        \node at (1,1.25) {\Large{$C$}};
    \end{tikzpicture}
    {\Large{;}} \quad 
    \begin{tikzpicture}[anchorbase]
        \draw[webs,->-=.8] (.25,0) to (.25,1.8);
        \draw[webs,->-=.15,->-=.885] (1,-1.3) node[below]{$n-1$} to (1,1.8) node[above]{$n-1$};

        \filldraw[fill=white, draw=black]
            (0,.4) rectangle (1.25,1.1);

        \node at (.625,.75) {\Large{$C$}};

        \draw[webs,->-=.2,->-=.8] (-.5,-1) to[out=90,in=270] (.25,0);
        \draw[webs,->-=.2,->-=.8,over] (.25,-1) to[out=90,in=270] (-.5,0);

        \draw[webs,->-=.5] (-.5,0) to (-.5,1.5);
        \draw[webs] (-.5,1.5) to[out=90,in=0] (-.75,1.75);
        \draw[webs] (-.75,1.75) to[out=180,in=90] (-1,1.5);
        \draw[webs] (-1,1.5) to (-1,-1);
        \draw[webs] (-1,-1) to[out=270,in=180] (-.75,-1.25);
        \draw[webs] (-.75,-1.25) to[out=0,in=270] (-.5,-1);

        \draw[webs] (.25,-1.3) to (.25,-1);
    \end{tikzpicture}
    {\Large{$\simeq aq^{-3}$}} \ 
    \begin{tikzpicture}[anchorbase]
        \draw[webs,->-=.2,->-=.8] (1,0) node[below]{$n$} to (1,2.5) node[above]{$n$};
        \filldraw[fill=white, draw=black]
            (.5,.9) rectangle (1.5,1.6);

        \node at (1,1.25) {\Large{$C$}};
    \end{tikzpicture}
    
    \caption{A diagrammatic version of the Reidemeister I move.}
    \end{figure}
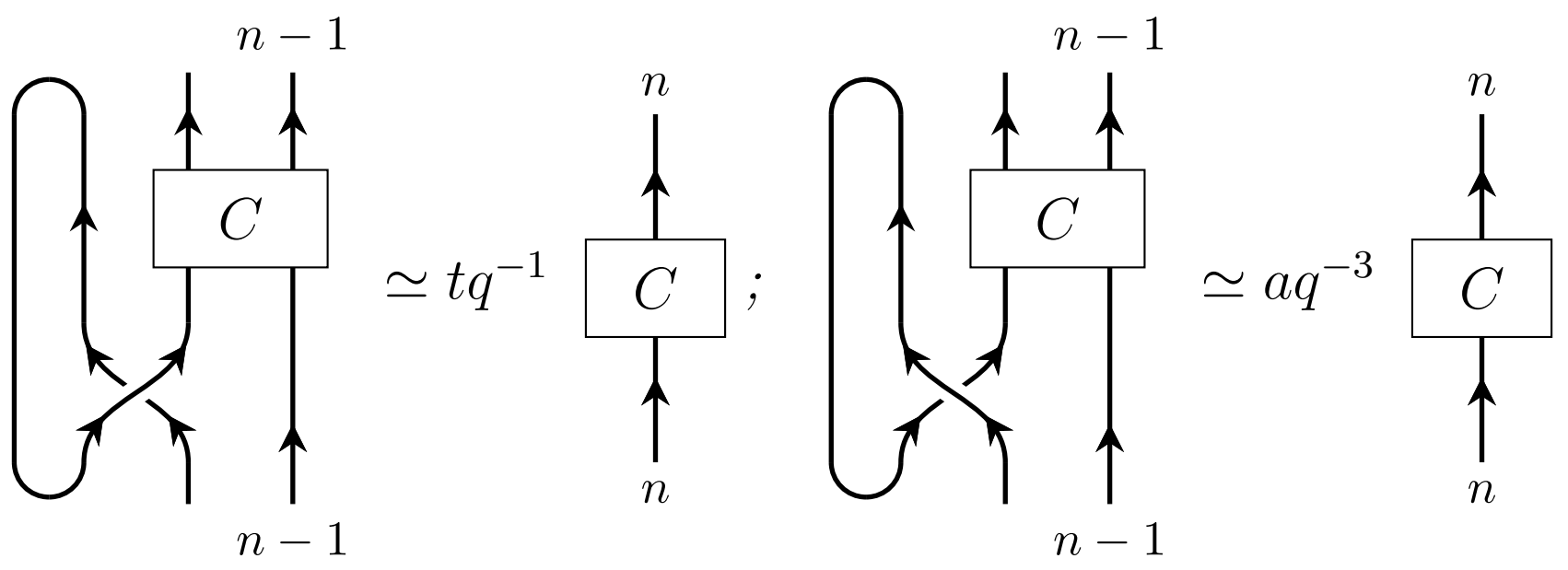

    \FloatBarrier
\end{prop}

\begin{proof}
    By Proposition \ref{prop: trace_local}, it suffices to compute $Tr_1(F(\sigma_1^{\pm 1})) \in K^b(\extsbim_1)$. This is done e.g. in Proposition 3.10 of \cite{Hog18}; renormalizing the shifts there to account for our slightly different presentation of the positive and negative crossings gives exactly the desired result.
\end{proof}

Finally, consider two complexes $M \in K^b(\extsbim_m)$, $C \in K^b(\extsbim_n)$ for some $m \leq n$. There is an obvious isotopy between the diagrammatic representations of $Tr^m(M \star C)$ and $Tr^m(C \star M)$ given by sliding the coupon labeled $M$ through the $m$ parallel closure strands, and this isotopy once again has an algebraic version:

\begin{prop} \label{prop: trace-like}
    Fix $0 \leq m \leq n$, and let $M \in K^b(\extsbim_m)$, $C \in K^b(\extsbim_n)$ be given. Then there is a natural (in $M$ and $C$) isomorphism in $K^b(\extsbim_{n - m})$

    \[
    Tr^m((M \boxtimes \strand{n-m}) \star C) \simeq Tr^m(C \star (M \boxtimes \strand{n-m})).
    \]

    Diagrammatically:

    \begin{figure}[ht]
        \centering

        \begin{tikzpicture}[anchorbase]
            \draw[webs,->-=.05,->-=.5,->-=.95] (0,0) to (0,2.5);
        \draw[webs,->-=.1,->-=.7] (1,-.8) node[below]{$n - m$} to (1,3.3) node[above]{$n - m$};
        
        \draw[webs] (0,2.5) to[out=90,in=0] (-.5,3);
        \draw[webs] (-.5,3) to[out=180,in=90] (-1,2.5);
        \draw[webs] (-1,2.5) to node[left]{$m$} (-1,0);
        \draw[webs] (-1,0) to[out=270,in=180] (-.5,-.5);
        \draw[webs] (-.5,-.5) to[out=0,in=270] (0,0);
        
        \filldraw[fill=white, draw=black]
            (-.25,.4) rectangle (1.25,1.1);

        \node at (.5,.75) {\Large{$C$}};
        
        \filldraw[fill=white, draw=black]
            (-.5,1.5) rectangle (.5,2.2);

        \node at (0,1.85) {\Large{$M$}};
        \end{tikzpicture}
        {\Large{$\simeq$}}
        \begin{tikzpicture}[anchorbase]
            \draw[webs,->-=.05,->-=.5,->-=.95] (0,0) to (0,2.5);
        \draw[webs,->-=.3,->-=.9] (1,-.8) node[below]{$n - m$} to (1,3.3) node[above]{$n - m$};
        
        \draw[webs] (0,2.5) to[out=90,in=0] (-.5,3);
        \draw[webs] (-.5,3) to[out=180,in=90] (-1,2.5);
        \draw[webs] (-1,2.5) to node[left]{$m$} (-1,0);
        \draw[webs] (-1,0) to[out=270,in=180] (-.5,-.5);
        \draw[webs] (-.5,-.5) to[out=0,in=270] (0,0);

        \filldraw[fill=white, draw=black]
            (-.25,1.5) rectangle (1.25,2.2);

        \node at (.5,1.85) {\Large{$C$}};
        
        \filldraw[fill=white, draw=black]
            (-.5,.4) rectangle (.5,1.1);

        \node at (0,.75) {\Large{$M$}};
        \end{tikzpicture}
    \end{figure}

    \FloatBarrier
\end{prop}

\begin{proof}
    This is essentially a cyclicity statement for a relative version of Hochschild homology, and we proceed with this motivation in mind. We begin by building an alternate model of $Tr^m((M \boxtimes \strand{n-m}) \star C)$ using the Koszul complex $K_m$. Recall from the proof of Proposition \ref{prop: Hom_Complex} that $((K_m \star M) \boxtimes \strand{n-m}) \star C$ is a free $R_m$-bimodule. In particular, we may consider the naive (i.e. underived) quotient by the ideal $I_m := \langle x_j - x_j' \ | \ 1 \leq j \leq m \rangle$ generated by differences of the variables involved in the left and right actions. By inspecting the differentials on each side, we immediately obtain an isomorphism

    \[
    Tr^m((M \boxtimes \strand{n-m}) \star C) \cong a^mq^{-2m} \left(((K_m \star M) \boxtimes \strand{n-m}) \star C \right) / I_m.
    \]

    We now apply the natural homotopy equivalence $s_M \colon K_m \star M \simeq M \star K_m$ of Lemma \ref{lem: koszul_exchange}\footnote{Note that $s_M$ is strictly natural as a morphism in the derived category of bimodules, so it automatically extends to chain complexes over the derived category.} to the right-hand side to obtain an isomorphism
    
    \[
    Tr^m((M \boxtimes \strand{n-m}) \star C) \cong a^mq^{-2m} \left((M \boxtimes \strand{n-m}) \star ( K_m \boxtimes \strand{n-m}) \star C \right) / I_m.
    \]

    Now the right-hand side is exactly a shift of the $0^{th}$ Hochschild homology of $(M \boxtimes \strand{n-m}) \star (K_m \boxtimes \strand{n-m}) \star C$, considered as a complex of $R_m$-bimodules. This functor is well-known to be invariant under cyclic permutations up to natural isomorphism, so we immediately obtain
    
    \[
    a^mq^{-2m} \left((M \boxtimes \strand{n-m}) \star ( K_m \boxtimes \strand{n-m}) \star C \right) / I_m \cong a^mq^{-2m} \left(( K_m \boxtimes \strand{n-m}) \star C \star (M \boxtimes \strand{n-m})\right) / I_m.
    \]

    Finally, again applying the first step of this argument, we obtain

    \[
    a^mq^{-2m} ((K_m \boxtimes \strand{n-m}) \star C \star (M \boxtimes \strand{n-m}))/I_m \cong Tr^m(C \star (M \boxtimes \strand{n-m})).
    \]
\end{proof}

\subsection{Coherent Tangle Diagrams} \label{subsec: coherent_diagrams}

We can obtain various oriented tangle diagrams by starting with the empty diagram and iteratively including new strands, concatenating with braid diagrams, and closing up endpoints. By iteratively applying the corresponding functors $I$, $F(\beta) \star -$, $- \star F(\beta)$, and $Tr$ to the complex $\Z = \strand{0} \in K^b(\extsbim_0)$ according to the diagrammatics described above, we can interpret any oriented tangle diagram obtained in this way as encoding a complex in $K^b(\extsbim)$. We summarize this idea in the following definition:

\begin{defn} \label{def: coherent_tangle_diagram}
    We call an oriented tangle diagram $\tau$ \textit{coherent} if it can be obtained from the empty tangle diagram via a finite sequence of the following operations:
    \begin{enumerate}
        \item Adding a single upwards-oriented vertical strand to the left of an oriented tangle diagram $\tau'$;

        \item Closing the leftmost top and bottom endpoints of an oriented tangle diagram $\tau'$ in which both such endpoints are positively oriented;

        \item Vertically composing an oriented tangle diagram $\tau'$ above or below with a single crossing between upwards-oriented strands.
    \end{enumerate}

    Given a coherent tangle diagram $\tau$ with $2n$ endpoints, we say it \textit{encodes} the complex $C \in K^b(\extsbim_n)$ which is represented by $\tau$ diagrammatically.
\end{defn}

There are a few special classes of coherent tangle diagrams worth singling out. One class consists of braid diagrams; any such $n$-strand diagram $\beta$ can be obtained from the empty diagram by first including $n$ parallel upwards strands, then concatenating with $\beta$. Of course $\beta$ encodes the complex $F(\beta) \in K^b(\extsbim_n)$. Another class consists of \textit{braid closures}; these are those diagrams obtained from $n$-strand braid diagrams $\beta$ by closing all $n$ pairs of endpoints from left to right. Having no endpoints, each such braid closure encodes a complex $C = Tr^n(F(\beta)) \in K^b(\extsbim_0)$.

The following theorem is due to Khovanov in \cite{Kh07}, where it is proven by comparison to an earlier formulation by Khovanov and Rozansky in \cite{KhR08}.

\begin{theorem} \label{thm: hhh_alg_invariance}
    Let $\L$ and $\L'$ be two braid closures depicting isotopic oriented links. Then $\L$ and $\L'$ encode homotopy equivalent complexes in $K^b(\extsbim_0)$ up to an overall regrading.
\end{theorem}

\begin{proof}
    Suppose $\L$ is the closure of an $m$-strand braid word $\beta \in \Br_m$ and $\L'$ is the closure of an $n$-strand braid word $\beta' \in \Br_n$. By a classical theorem of Markov \cite{Mar35}, $\beta$ and $\beta'$ are related by a finite sequence of braid relations, conjugations, and stabilizations. By Propositions \ref{prop: rouquier_canonicity}, \ref{prop: trace-like}, and \ref{prop: markov}, the complexes encoded by the diagrams in any step in this sequence are homotopy equivalent up to an overall regrading. In particular, the complexes encoded by $\L$ and $\L'$ are homotopy equivalent up to an overall regrading.
\end{proof}

\begin{remark} \label{rem: alg_hhh_regrading}
    The regrading involved in Theorem \ref{thm: hhh_alg_invariance} can be made perfectly explicit. One way to do this is via a global renormalization after computing $HHH^{alg}$. For this, we first formally extend to a half-integer grading\footnote{One needs to be somewhat careful about managing signs: the functors $a^{\pm 1}$ and $t^{\pm 1}$ both negate the signs on all differentials, so we cannot work with grading shifts by $a^{1/2}$ or $t^{1/2}$ alone but only in the combination $a^{1/2}t^{1/2}$, which we may arbitrarily take to be odd.} in $\Z_a$ and $\Z_t$ and set $\epsilon := a^{1/2}q^{-1}t^{-1/2}$, $\rho := a^{-1/2}q^2t^{-1/2}$. Given a braid word $\beta$, denote by $e(\beta)$ its \textit{braid exponent}, given by the number of positive crossings minus the number of negative crossings. Then after choosing a presentation $\L = \hat{\beta}$ of an oriented link $\L$ as the closure of an $n$-strand braid $\beta$, tracking all regradings explicitly through the above argument gives that $HHH^{norm}(\L) := \epsilon^{e(\beta)} \rho^n HHH^{alg}(\L)$ is a genuine link invariant.

    Alternatively, we can renormalize the Rouquier complexes $F(\beta)$ and the partial trace directly to eliminate the grading shifts in Proposition \ref{prop: markov} by redefining $F^{norm}(\beta) := \epsilon^{e(\beta)} F(\beta)$ for each braid word $\beta$ and $Tr^{norm}(C) := \rho Tr(C)$ for each complex $C \in K^b(\extsbim_{\geq 1})$. It is easy to check that $e(\beta)$ is invariant under the braid relations, so $F^{norm}(\beta)$ is once again a braid invariant. An easy computation reveals $\epsilon \rho = t^{-1}q$ and $\epsilon^{-1} \rho = a^{-1}q^3$. We can now check that the Markov moves hold without corrections; for example, the positive Markov move becomes

    \[
    Tr_n^{norm} (I_n(C) \star F^{norm}(\sigma_1)) = \epsilon \rho \ Tr_n (I_n(C) \star F(\sigma_1)) \simeq (t^{-1}q) (tq^{-1}) C = C.
    \]

    This \textit{local} renormalization has a huge advantage over the global renormalization: it can easily be applied to any coherent tangle diagram to encode a renormalized complex, not just to a closed braid. We prefer to work with the original, unnormalized Rouquier complex and trace mainly to avoid managing two distinct conventions or introducing half-integer gradings from the start.
\end{remark}

By a classical theorem of Alexander \cite{Al23}, \textit{any} oriented link $\L$ can be presented as the closure of a braid $\beta$ and therefore interpreted as encoding a complex in $K^b(\extsbim_0)$. Theorem \ref{thm: hhh_alg_invariance} guarantees that this complex does not depend on the choice of $\beta$. We are therefore led to the following definition, originally due to Khovanov in \cite{Kh07} following Khovanov and Rozansky in \cite{KhR08}:

\begin{definition} \label{def: hhh_algebraic}
    Let $\L$ be an oriented link. Then the \textit{algebraic HOMFLY homology} of $\L$ is the complex $HHH_{alg}(\L) \in K^b(\extsbim_0)$ encoded by any presentation of $\L$ as a braid closure.
\end{definition}

\begin{remark} \label{rem: complex_vs_homology}
    It is common in the literature to instead refer to the HOMFLY homology of $\mathcal{L}$ as the $\Z_a \times \Z_q \times \Z_t$-graded module extracted from what we denote $HHH_{alg}(\L)$ by taking homology in the $t$-direction; hence the name ``triply-graded homology". We prefer to work directly with the complex $HHH_{alg}(\L)$ for closer agreement with the complexes assigned to oriented tangles with non-empty boundary. Over a coefficient field $HHH_{alg}(\L)$ is already homotopy equivalent to its homology (considered as a complex with $0$ differential), so this distinction is largely irrelevant.
\end{remark}

If we drop the restriction to braid closures and allow ourselves to pass through \textit{arbitrary} coherent tangle diagrams, there is another isotopy that suggests itself, displayed in Figure \ref{fig: exchange_move}. We refer to this isotopy as the \textit{exchange move}; it is an open-ended analogue of the exchange move for braid closures first introduced by Birman and Menasco in \cite{BM90}. At the diagrammatic level, this isotopy is a composition of two dial moves and a ``wrong-way" Reidemeister II move, in which a free upwards-oriented strand slides over or under the downwards-oriented component of closure. To the best of our knowledge, the following algebraic version of this move is new:

\begin{theorem} \label{thm: exchange}
    For each $n \geq 1$, there are natural equivalences of functors from $K^b(\extsbim_n)$ to $K^b(\extsbim_n)$ of the form

    \[
    Tr_n \circ (F(\sigma_1) \star I_n(-) \star F(\sigma_1^{-1})) \cong Tr_n \circ (F(\sigma_1^{-1}) \star I_n(-) \star F(\sigma_1)) \cong I_{n-1} \circ Tr_{n-1}.
    \]

    Diagrammatically, for each $C \in K^b(\extsbim_n)$:

    \begin{figure}[ht]
    \centering

    \begin{tikzpicture}[anchorbase]
        \draw[webs,->-=.1,->-=.9] (.25,0) to (.25,1.5);
        \draw[webs] (1,0) to (1,1.5);

        \filldraw[fill=white, draw=black]
            (0,.4) rectangle (1.25,1.1);

        \node at (.625,.75) {\Large{$C$}};

        \draw[webs] (-.5,-.75) to[out=90,in=270] (.25,0);
        \draw[webs,over] (.25,-.75) to[out=90,in=270] (-.5,0);

        \draw[webs] (.25,1.5) to[out=90,in=270] (-.5,2.25);
        \draw[webs,over] (-.5,1.5) to[out=90,in=270] (.25,2.25);

        \draw[webs,->-=.5] (-.5,0) to (-.5,1.5);

        \draw[webs] (-.5,2.25) to[out=90,in=0] (-.75,2.5);
        \draw[webs] (-.75,2.5) to[out=180,in=90] (-1,2.25);
        \draw[webs] (-1,2.25) to (-1,-.75);
        \draw[webs] (-1,-.75) to[out=270,in=180] (-.75,-1);
        \draw[webs] (-.75,-1) to[out=0,in=270] (-.5,-.75);

        \draw[webs,-<-=.5] (.25,-.75) to (.25,-1.3);
        \draw[webs,->-=.5] (.25,2.25) to (.25,2.8);

        \draw[webs,-<-=.45] (1,0) to (1,-1.3) node[below]{$n-1$};
        \draw[webs,->-=.45] (1,1.5) to (1,2.8) node[above]{$n-1$};
    \end{tikzpicture}
    {\large{$\simeq$}} \ 
    \begin{tikzpicture}[anchorbase]
        \draw[webs,->-=.1,->-=.9] (.25,0) to (.25,1.5);
        \draw[webs] (1,0) to (1,1.5);

        \filldraw[fill=white, draw=black]
            (0,.4) rectangle (1.25,1.1);

        \node at (.625,.75) {\Large{$C$}};

        \draw[webs] (.25,-.75) to[out=90,in=270] (-.5,0);
        \draw[webs,over] (-.5,-.75) to[out=90,in=270] (.25,0);
        
        \draw[webs] (-.5,1.5) to[out=90,in=270] (.25,2.25);
        \draw[webs,over] (.25,1.5) to[out=90,in=270] (-.5,2.25);
        
        \draw[webs,->-=.5] (-.5,0) to (-.5,1.5);

        \draw[webs] (-.5,2.25) to[out=90,in=0] (-.75,2.5);
        \draw[webs] (-.75,2.5) to[out=180,in=90] (-1,2.25);
        \draw[webs] (-1,2.25) to (-1,-.75);
        \draw[webs] (-1,-.75) to[out=270,in=180] (-.75,-1);
        \draw[webs] (-.75,-1) to[out=0,in=270] (-.5,-.75);

        \draw[webs,-<-=.5] (.25,-.75) to (.25,-1.3);
        \draw[webs,->-=.5] (.25,2.25) to (.25,2.8);

        \draw[webs,-<-=.45] (1,0) to (1,-1.3) node[below]{$n-1$};
        \draw[webs,->-=.45] (1,1.5) to (1,2.8) node[above]{$n-1$};
    \end{tikzpicture}
    {\large{$\simeq$}}
    \begin{tikzpicture}[anchorbase]
        \draw[webs,->-=.1,->-=.9] (.25,0) to (.25,1.5);
        \draw[webs,->-=.2,->-=.8] (1,-.3) node[below]{$n-1$} to (1,1.8) node[above]{$n-1$};

        \filldraw[fill=white, draw=black]
            (0,.4) rectangle (1.25,1.1);

        \node at (.625,.75) {\Large{$C$}};

        \draw[webs] (.25,1.5) to[out=90,in=0] (0,1.75);
        \draw[webs] (0,1.75) to[out=180,in=90] (-.25,1.5);
        \draw[webs] (-.25,1.5) to (-.25,0);
        \draw[webs] (-.25,0) to[out=270,in=180] (0,-.25);
        \draw[webs] (0,-.25) to[out=0,in=270] (.25,0);

        \draw[webs,->-=.5] (-.75,-.3) to (-.75,1.8);
    \end{tikzpicture}

    \caption{The exchange move.}
    \label{fig: exchange_move}
    \end{figure}
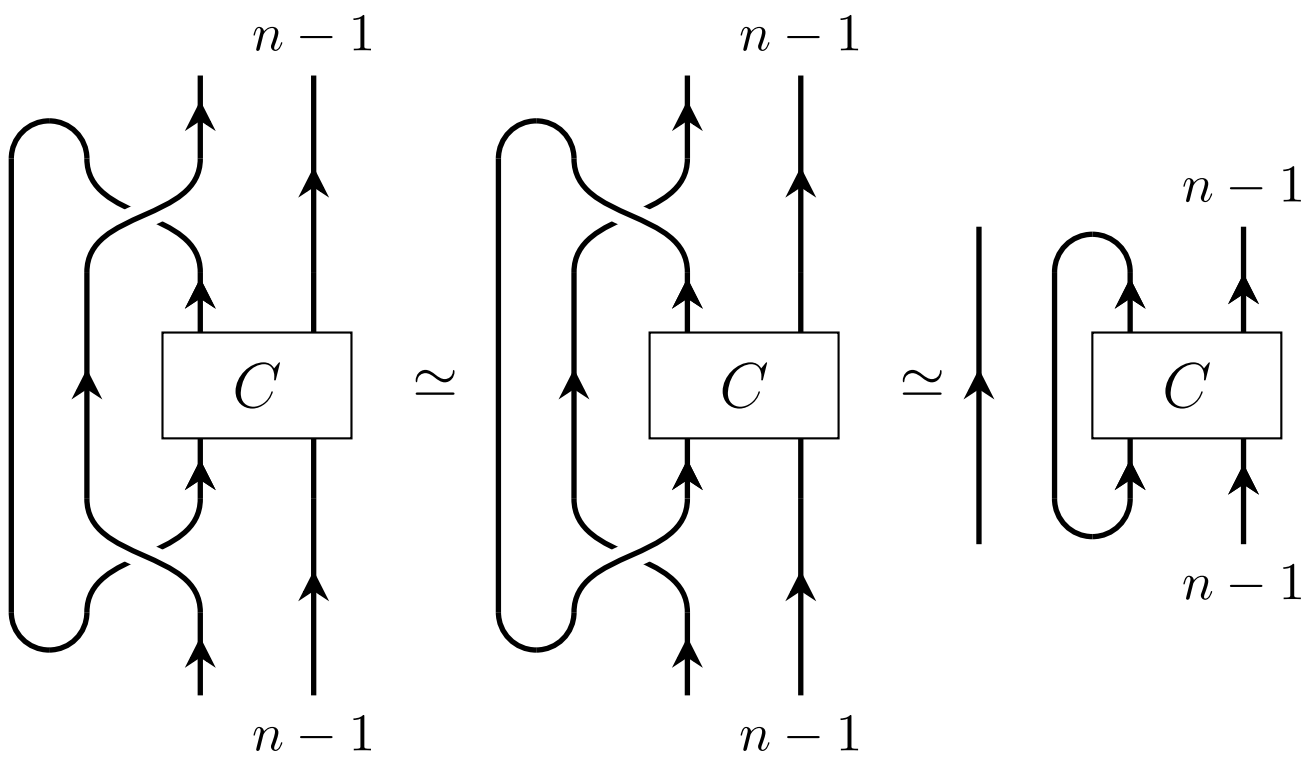

    \FloatBarrier
\end{theorem}

\begin{proof}
    We proceed in two steps: first, we construct a natural transformation $$I_{n-1} \circ Tr_{n-1} \to Tr_n \circ (F(\sigma_1^\epsilon) \star I_n(-) \star F(\sigma_1^{-\epsilon})),$$ for $\epsilon = \pm 1$, which turns out to exist on the whole derived category of bimodules, not just $\extsbim$; second, we show that when restricted to $\extsbim$, this is pointwise a homotopy equivalence, so invertible as a natural transformation on $K^b(\extsbim)$.

    To build the map, we chain together a series of (co)unit maps from the adjunction of Proposition \ref{prop: trace_inclusion_adjoint}, braid isotopies, and the cyclicity map from Proposition~\ref{prop: trace-like}, following Figure~\ref{fig: exchange_proof}:
    \begin{align*}
        I_{n-1} \circ Tr_{n-1}(-) &\rightarrow I_{n-1} \circ Tr_{n-1} \circ Tr_n \circ I_n(-) &\text{(unit of the adjunction)} \\
        &\rightarrow I_{n-1} \circ Tr_{n-1} \circ Tr_n \circ (I_n(-) \star F(\sigma_1^{-\epsilon}\sigma_1^\epsilon)) &\text{(Reidemeister II)} \\
        &\rightarrow I_{n-1} \circ Tr_{n-1} \circ Tr_n \circ (F(\sigma_1^\epsilon)\star I_n(-) \star F(\sigma_1^{-\epsilon})) &\text{(Proposition~\ref{prop: trace-like})} \\
        &\rightarrow Tr_n \circ (F(\sigma_1^\epsilon)\star I_n(-) \star F(\sigma_1^{-\epsilon})) &\text{(counit of the adjunction)}
    \end{align*}

    \begin{figure} 
    \centering
    \scalebox{.9}{%
    \begin{tikzpicture}[anchorbase]
    \draw[webs,->-=.1,->-=.9] (.25,0) to (.25,1.5);
        \draw[webs,->-=.2,->-=.8] (1,-.3) node[below]{$n-1$} to (1,1.8) node[above]{$n-1$};

        \filldraw[fill=white, draw=black]
            (0,.4) rectangle (1.25,1.1);

        \node at (.625,.75) {\Large{$C$}};

        \draw[webs] (.25,1.5) to[out=90,in=0] (0,1.75);
        \draw[webs] (0,1.75) to[out=180,in=90] (-.25,1.5);
        \draw[webs] (-.25,1.5) to (-.25,0);
        \draw[webs] (-.25,0) to[out=270,in=180] (0,-.25);
        \draw[webs] (0,-.25) to[out=0,in=270] (.25,0);

        \draw[webs,->-=.5] (-.75,-.3) to (-.75,1.8);
    \end{tikzpicture}
    {\Large{$\rightarrow$}} \ 
    \begin{tikzpicture}[anchorbase]
        \draw[webs,->-=.05,->-=.95] (.25,0) to (.25,1.5);
        \draw[webs,->-=.28,->-=.72] (1,-.8) node[below]{$n-1$} to (1,2.3) node[above]{$n-1$};
    
        \filldraw[fill=white, draw=black]
            (0,.4) rectangle (1.25,1.1);
        \node at (.625,.75) {\Large{$C$}};
    
        \draw[webs] (.25,1.5) to[out=90,in=0] (-.5,2.25);
        \draw[webs] (-.5,2.25) to[out=180,in=90] (-1.25,1.5);
        \draw[webs] (-1.25,1.5) to (-1.25,0);
        \draw[webs] (-1.25,0) to[out=270,in=180] (-.5,-.75);
        \draw[webs] (-.5,-.75) to[out=0,in=270] (.25,0);
    
        \draw[webs,->-=.5] (-.25,0) to (-.25,1.5);
        \draw[webs] (-.25,1.5) to[out=90,in=0] (-.5,1.75);
        \draw[webs] (-.5,1.75) to[out=180,in=90] (-.75,1.5);
        \draw[webs] (-.75,1.5) to (-.75,0);
        \draw[webs] (-.75,0) to[out=270,in=180] (-.5,-.25);
        \draw[webs] (-.5,-.25) to[out=0,in=270] (-.25,0);
        
        \draw[webs,->-=.5] (-1.75,-.8) to (-1.75,2.3);
    \end{tikzpicture}
    {\Large{$\simeq$}} \ 
    \begin{tikzpicture}[anchorbase]
        \draw[webs,->-=.1,->-=.9] (.25,0) to (.25,1.5);
        \draw[webs,->-=.2,->-=.8] (1,-2.3) node[below]{$n-1$} to (1,2.3) node[above]{$n-1$};
    
        \filldraw[fill=white, draw=black]
            (0,.4) rectangle (1.25,1.1);
        \node at (.625,.75) {\Large{$C$}};
    
        \draw[webs] (.25,1.5) to[out=90,in=0] (-.625,2.25);
        \draw[webs] (-.625,2.25) to[out=180,in=90] (-1.5,1.5);
        \draw[webs] (-1.5,1.5) to (-1.5,-1.5);
        \draw[webs] (-1.5,-1.5) to[out=270,in=180] (-.625,-2.25);
        \draw[webs] (-.625,-2.25) to[out=0,in=270] (.25,-1.5);
    
        \draw[webs,->-=.01] (-.5,-1.5) to[out=90,in=270] (.25,-.75);
        \draw[webs,over,->-=.01] (.25,-1.5) to[out=90,in=270] (-.5,-.75);
    
        \draw[webs,->-=.01] (.25,-.75) to[out=90,in=270] (-.5,0);
        \draw[webs,over,->-=.01] (-.5,-.75) to[out=90,in=270] (.25,0);
    
        \draw[webs,->-=.5] (-.5,0) to (-.5,1.5);
        \draw[webs] (-.5,1.5) to[out=90,in=0] (-.75,1.75);
        \draw[webs] (-.75,1.75) to[out=180,in=90] (-1,1.5);
        \draw[webs] (-1,1.5) to (-1,-1.5);
        \draw[webs] (-1,-1.5) to[out=270,in=180] (-.75,-1.75);
        \draw[webs] (-.75,-1.75) to[out=0,in=270] (-.5,-1.5);
    
        \draw[webs,->-=.5] (-2,-2.3) to (-2,2.3);
    \end{tikzpicture}
    {\Large{$\simeq$}} \ 
    \begin{tikzpicture}[anchorbase]
        \draw[webs,->-=.1,->-=.9] (.25,0) to (.25,1.5);
        \draw[webs,->-=.2,->-=.8] (1,-1.55) node[below]{$n-1$} to (1,3.05) node[above]{$n-1$};
    
        \filldraw[fill=white, draw=black]
            (0,.4) rectangle (1.25,1.1);
        \node at (.625,.75) {\Large{$C$}};
    
        \draw[webs,->-=.99] (-.5,1.5) to[out=90,in=270] (.25,2.25);
        \draw[webs,over,->-=.99] (.25,1.5) to[out=90,in=270] (-.5,2.25);
    
        \draw[webs,->-=.01] (.25,-.75) to[out=90,in=270] (-.5,0);
        \draw[webs,over,->-=.01] (-.5,-.75) to[out=90,in=270] (.25,0);
    
        \draw[webs] (-.5,2.25) to[out=90,in=0] (-.75,2.5);
        \draw[webs] (-.75,2.5) to[out=180,in=90] (-1,2.25);
        \draw[webs] (-1,2.25) to (-1,-.75);
        \draw[webs] (-1,-.75) to[out=270,in=180] (-.75,-1);
        \draw[webs] (-.75,-1) to[out=0,in=270] (-.5,-.75);
    
        \draw[webs,->-=.5] (-.5,0) to (-.5,1.5);
        \draw[webs] (.25,2.25) to[out=90,in=0] (-.625,3);
        \draw[webs] (-.625,3) to[out=180,in=90] (-1.5,2.25);
        \draw[webs] (-1.5,2.25) to (-1.5,-.75);
        \draw[webs] (-1.5,-.75) to[out=270,in=180] (-.625,-1.5);
        \draw[webs] (-.625,-1.5) to[out=0,in=270] (.25,-.75);
    
        \draw[webs,->-=.5] (-2,-1.55) to (-2,3.05);
    \end{tikzpicture}
    {\Large{$\rightarrow$}} \ 
    \begin{tikzpicture}[anchorbase]
        \draw[webs,->-=.1,->-=.9] (.25,0) to (.25,1.5);
        \draw[webs] (1,0) to (1,1.5);

        \filldraw[fill=white, draw=black]
            (0,.4) rectangle (1.25,1.1);

        \node at (.625,.75) {\Large{$C$}};

        \draw[webs] (.25,-.75) to[out=90,in=270] (-.5,0);
        \draw[webs,over] (-.5,-.75) to[out=90,in=270] (.25,0);
        
        \draw[webs] (-.5,1.5) to[out=90,in=270] (.25,2.25);
        \draw[webs,over] (.25,1.5) to[out=90,in=270] (-.5,2.25);
        
        \draw[webs,->-=.5] (-.5,0) to (-.5,1.5);

        \draw[webs] (-.5,2.25) to[out=90,in=0] (-.75,2.5);
        \draw[webs] (-.75,2.5) to[out=180,in=90] (-1,2.25);
        \draw[webs] (-1,2.25) to (-1,-.75);
        \draw[webs] (-1,-.75) to[out=270,in=180] (-.75,-1);
        \draw[webs] (-.75,-1) to[out=0,in=270] (-.5,-.75);

        \draw[webs,-<-=.5] (.25,-.75) to (.25,-1.3);
        \draw[webs,->-=.5] (.25,2.25) to (.25,2.8);

        \draw[webs,-<-=.45] (1,0) to (1,-1.3) node[below]{$n-1$};
        \draw[webs,->-=.45] (1,1.5) to (1,2.8) node[above]{$n-1$};
    \end{tikzpicture}
    }
    \caption{Intermediate stages of the exchange map. The first and fourth maps are the unit and counit, respectively, of the adjunction of Proposition \ref{prop: trace_inclusion_adjoint}, the second is a braid relation, and the third is the equivalence from Proposition \ref{prop: trace-like}.}
    \label{fig: exchange_proof}
\end{figure}
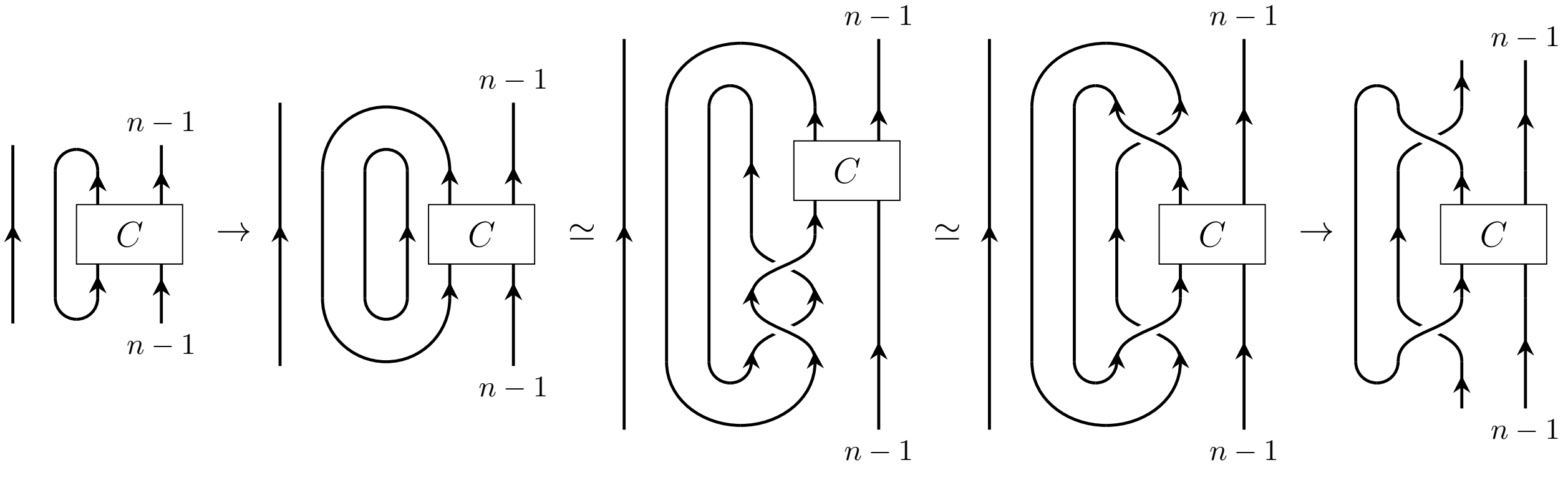

    Write $\psi_C$ for the component of this natural transformation on a one-morphism $C$ of $\extsbim_n$. It now suffices to show that each $\psi_C$ is a homotopy equivalence, since a natural transformation that is invertible pointwise is invertible overall. In fact, it suffices to consider only $C$ an individual bimodule, rather than a complex: $\psi$ is a natural transformation between two functors $K^b(\extsbim_n) \to K^b(\extsbim_n)$, but each of $\psi$ and the (co)domain functors actually originates as a strict natural transformation / functor $\extsbim_n \to \operatorname{Ch}^b(\extsbim_n)$, from which one formally extracts first an endofunctor of $\operatorname{Ch}^b(\extsbim_n)$ and then of $K^b(\extsbim_n)$. This means that, if $C$ is a complex with component bimodules $C_i$, then $\psi_C$ acquires the structure of a map between two double complexes, with columns given by the $\psi_{C_i}$. One now needs only to use that a map between bounded double complexes, restricting to a homotopy equivalence on each column, is a homotopy equivalence on the total complex.
    
    Now, for single bimodules $M$, whether $\psi_M$ is a homotopy equivalence is unchanged by passing to direct sums or summands, so it suffices to consider only $M$ having the two forms of Lemma~\ref{lem: singleB1}. Since also each step in the construction of $\psi_M$ visibly commutes with tensor products on the rightmost $n-1$ strands (i.e. by $I(\extsbim_{n-1})$, above or below), we need only show that $\psi_{\strand{n}}$ and $\psi_{B_1}$ are homotopy equivalences.

    To understand $\psi_{\strand{n}}$, we may as well take $n=1$ because the other strands do not participate in the calculation. Then a single Reidemeister II move identifies the last frame of Figure~\ref{fig: exchange_proof} with $Tr_1 \circ I_1 (\strand{1})$, so overall we have a map $I_0 \circ Tr_0(\strand{1}) \to Tr_1 \circ I_1(\strand{1})$, both sides of which are isomorphic to $\strand{1}[x,\theta]$. As a final simplification, the counit map in the last step of Figure~\ref{fig: exchange_proof} commutes with the added Reidemeister II move, so overall we need to check that a map
    \begin{equation*}
        I_0 \circ Tr_0(\strand{1}) \to I_0 \circ Tr_0 \circ Tr_1 (\strand{2}) \to I_0 \circ Tr_0 \circ Tr_1 (\strand{2}) \to Tr_1(\strand{2}) \tag{*}
    \end{equation*}
    is an isomorphism, where the first map is a unit, the middle one is built from a Reidemeister II move, a cyclicity map, and a second Reidemeister II move, and the final map is a counit. Visually, the middle map swaps the two unknots in the picture; our main task will be to see that this also occurs algebraically \footnote{This could be seen purely formally, using the fact that $x$ and $\theta$ slide past crossings, but the direct computation below is not so bad.}. 

    Unpacking the Reidemeister II isomorphisms, the middle map is a sum of two terms. The first is $-1$ (factoring through the $R$ portion of the Reidemeister II complex), and the other (coming from $B_1 \star B_1$) factors as ($I_0$ applied to)
    \begin{equation*}
        Tr_0\circ Tr_1(\strand{2}) \to Tr_0 \circ Tr_1( B_1 \star B_1) \to Tr_0\circ Tr_1(B_1\star B_1) \to Tr_0\circ Tr_1(\strand{2})
    \end{equation*}
    where the first and last maps are the cup and cap of the self-duality of $B_1$, and the middle one is the cyclicity map. To compute this explicitly, note that every term is the full trace of some Soergel bimodule. The full trace involves taking the ($a$-direction) cone on $x_1-x_1'$ and $x_2-x_2'$; we may factor this as the two-step process of first tensoring by the resolution bimodule $K_2$ and then taking the underived trace, simply identifying $x_1$ with $x_1'$ and $x_2$ with $x_2'$. We will represent this final underived trace as a $\star$ at the end of each term, imagining the formulae drawn in an annulus which we may freely rotate, with the last $\star$ serving as a tensor product between the initial and final bimodules involved. Under this interpretation, what we must calculate is a map
    \begin{equation*}
        K_2 \star \to K_2 \star B_1 \star B_1 \star \to B_1  \star K_2 \star B_1 \star \to B_1 \star K_2 \star B_1 \star \to K_2 \star B_1 \star B_1 \star \to K_2 \star,
    \end{equation*}
    where in sequence we first birth a pair of $B_1$ using the duality, commute $K_2$ past $B_1$ using the map of Lemma~\ref{lem: Bswap}, cycle the newly-initial $B_1$ around to the end, and then finally kill the pair of $B_1$'s. Starting with an arbitrary monomial $ x_1^ax_2^b\eta_1^\alpha \eta_2^\beta \star $ in $K_2 \star,$ we get
    \begin{align*}
        x_1^ax_2^b\eta_1^\alpha \eta_2^\beta \star &\mapsto x_1^ax_2^b\eta_1^\alpha \eta_2^\beta \star (x_1 \star 1 - 1 \star x_2') \star \\
        &= x_1^ax_2^b\eta_1^\alpha \eta_2^\beta \star ((x_1-x_2')\star 1 + 1 \star (x_2-x_2')) \star\\
        &\mapsto \left(x_1^a x_2^b (x_1-x_2') \star \eta_1^\alpha \eta_2^\beta \star 1 \star \right)+ \left( x_1^ax_2^b \star \eta_2^\alpha \eta_1^\beta \star (x_2-x_2') \star\right) \\
        &\mapsto \left(\eta_1^\alpha \eta_2^\beta \star 1 \star x_1^a x_2^b (x_1-x_2') \star \right) + \left(\eta_2^\alpha \eta_1^\beta \star (x_2-x_2') \star x_1^ax_2^b \star \right) \\
        &= \left(\eta_1^\alpha \eta_2^\beta \star 1 \star(x_1-x_2') x_1'^a x_2'^b\star \right) + \left(\eta_2^\alpha \eta_1^\beta \star x_1^bx_2^a (x_2-x_2') \star 1 \star \right) \\
        &\mapsto \eta_1^\alpha \eta_2^\beta x_1^a x_2^b + \eta_2^\alpha \eta_1^\beta x_1^b x_2^a,
    \end{align*}
    where throughout we've used the fact that $\star$ is a tensor product over the polynomial ring, so $x_i' \star = \star x_i$, and that on $B_1$, we have $(x_1-x_1')(x_1-x_2') = 0$. Subtracting $1$, then, the middle map in (*), above, is $I_0$ applied to the endomorphism of $K_2 \star$ sending $x_1^ax_2^b\eta_1^\alpha \eta_2^\beta \star$ to $ x_1^b x_2^a\eta_1^\beta \eta_2^\alpha \star,$ visibly swapping the two components. Overall, this shows that (*) is the identity, using the zigzag identity for the adjunction between $Tr$ and $I$: we may interpret (*) as a map
    $$ \strand{1} \otimes \HHH(O_1) \to \strand{1} \otimes \HHH(O_1) \otimes \HHH(O_2) \to \strand{1} \otimes \HHH(O_2) \otimes \HHH(O_1) \to \strand{1} \otimes \HHH(O_1),$$ where $O_1$ and $O_2$ are two unknots and, as we have just computed, the middle map is the tensor swap. The factor $\HHH(O_1)$ does not participate in the first or last steps, so the whole map is the tensor product of $\HHH(O_1)$ and the composite of a unit and counit in the adjunction. 

    We now need only to show that $\psi_{B_1}$ is a homotopy equivalence, where as before we might as well take $n=2$. Let us first show that the codomain of $\psi_{B_1}$ is in fact homotopy equivalent to the domain; the fact that $\psi_{B_1}$ is a homotopy equivalence will then follow almost formally by considering $\psi_{\sigma_1}$.

    So, consider the complex $Tr_2(\sigma_1^{\epsilon} \star B_2 \star \sigma_1^{-\epsilon})$ (let's take $\epsilon = 1$ below; the case $\epsilon = -1$ follows by dualizing the complex). This will include a term $Tr_2(B_1 \star B_2 \star B_1)$ involving two copies of $B_1$, so to make taking the trace easier we will expand $B_1 \star B_2 \star B_1 \cong B_1 \oplus B_{121}$, giving overall the functor $Tr_2$ applied to
    
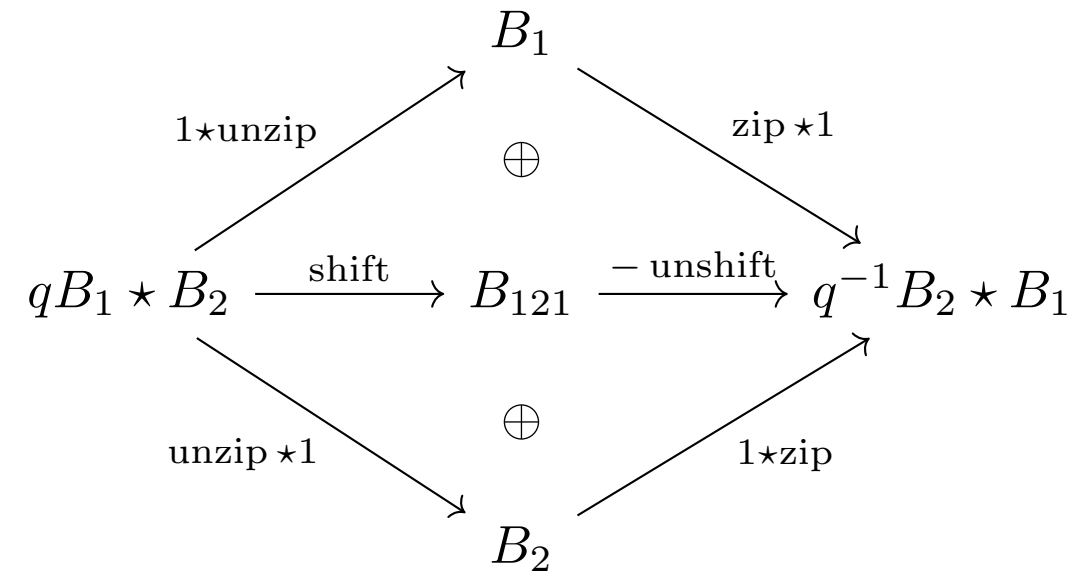
\begin{equation*}
        \begin{tikzcd}[row sep=small, column sep=large]
         & B_1 \arrow[ddr, "\zip \star 1"] & \\
         & \oplus & \\
        qB_1 \star B_2
          \arrow[uur, "1 \star \unzip"]
          \arrow[r, "\operatorname{shift}"]
          \arrow[ddr, "\unzip \star 1"']
         & B_{121} \arrow[r, "-\operatorname{unshift}"]
         & q^{-1}B_2 \star B_1 \\
         & \oplus & \\
         & B_2 \arrow[uur, "1 \star \zip"'] &
        \end{tikzcd}
    \end{equation*}
    where $\operatorname{shift}$ and $\operatorname{unshift}$ are the unique-up-to-sign maps generating the appropriate space of degree $q$ maps.
    \begin{lemma}
        $$Tr_2(B_{121}) \cong q^2 B_1[x,\theta^{(-2)}],$$ with $\theta^{(-2)}$ in degree $aq^{-6}$. 
        The maps $\operatorname{shift}$ and $\operatorname{unshift}$ induce
        \begin{align*}
            Tr_2(\operatorname{shift}: B_1 \star B_2 \to B_{121}) &= \left( 1 \mapsto 1, \theta^{(-1)} \mapsto (x-x_2)\theta^{(-2)} + \theta_{1,2}\right) \\
            Tr_2(\operatorname{unshift}: B_{121} \to B_1 \star B_2) &= \left(1 \mapsto x-x_2, \theta^{(-2)} \mapsto \theta^{(-1)} - \theta_{1,2} \right)
        \end{align*}
        (extended to be $x$-linear, as before), and the $\operatorname{shift}$ and $\operatorname{unshift}$ maps involving $B_2 \star B_1$ are the same, but with $x_2$ replaced by $x_2'$.
    \end{lemma}
    \begin{proof}
        The bimodule $B_{121}$ is the quotient of $R_3^e$ by differences of the elementary symmetric functions $e_1,e_2,$ and $e_3$ in the variables $x$ and $x'$, which we may re-express as $$B_{121} = q^{-3} R_3^e / \left(e_1(x) - e_1(x'), e_2(x) - e_2(x'), (x_1-x_1')(x_2-x_1')(x_3-x_1')\right).$$ Similarly, $$B_1 \star B_2 = q^{-2} R_3^e / \left(e_1(x) - e_1(x'), e_2(x) - e_2(x'), (x_1-x_1')(x_2-x_1')\right).$$ The maps $\operatorname{shift}$ and $\operatorname{unshift}$ are the quotient and multiplication by $x_1'-x_3$, respectively. We may compute using free resolutions of $B_{121}$ and $B_1 \star B_2$ given as Koszul complexes for these same relations; since the relations are nearly identical, $\operatorname{shift}$ and $\operatorname{unshift}$ lift to diagonal maps between the complexes with matrix entries all $1$ or $x_1' - x_3$. The result follows, after taking care to note that the endomorphism of the resolution for $B_1 \star B_2$ moving along the arrows of the differential carrying $(x_1-x_1')(x_2-x_1')$ is, in our earlier naming, $\theta_{1,2} \star 1 - 1 \star \theta_{2,3},$ which becomes $\theta - \theta_{1,2}$ after taking the trace.
    \end{proof}
    This tells us that the lower portion of the diagram for $Tr_2(\sigma_1 \star B_2 \star \sigma_1^{-1})$ may be contracted away, canceling the $q^2B_1[x]$ pieces of $Tr_2(qB_1 \star B_2)$ and $Tr_2(B_{121})$, the $q^2\theta^{(-1)}B_1[x]$ piece of $Tr_2(qB_1 \star B_2)$ with the $\theta B_1[x]$ piece of $Tr_2(B_2)$, the $q^2\theta^{(-2)} B_1[x]$ piece of $Tr_2(B_{121})$ with the $\theta^{(-1)}B_1[x]$ piece of $Tr_2(q^{-1} B_2 \star B_1)$, and finally the $B_1[x]$ pieces of $Tr_2(B_2)$ and $Tr_2(q^{-1}B_2 \star B_1)$. This leaves only $Tr_2(B_1) = q\strand{2}[x,\theta^{(-1)}]$, which is isomorphic to the domain of $\psi_{B_1}$ as desired.

    Finally, using a Reidemeister III move and a stabilization, we see that the domain and codomain of $\psi_{\sigma_1}$ are homotopy equivalent to $\strand{2}$, and we have just calculated that the domain and codomain of $\psi_{\strand{2}}$ are homotopy equivalent to $\strand{2}[x,\theta]$, that the domain and codomain of $\psi_{B_1}$ are homotopy equivalent to $q\strand{2}[x,\theta^{(-1)}],$ and that (under these equivalences) $\psi_{\strand{2}}$ is the identity. However, $\sigma_1$ is a mapping cone $B_1 \xrightarrow{\unzip} tq^{-1} \strand{2}$, so both the domain and codomain of $\psi$ must take $\unzip$ to the inclusion of a direct summand if we are to have only a single copy of $\strand{2}$ survive after applying the functors. Therefore, the naturality square for $\psi$ tells us that, after these homotopy equivalences, $\psi_{B_1}$ must be an isomorphism.
\end{proof}

\begin{remark}
    The fact that the map of Figure~\ref{fig: exchange_proof} should be an isomorphism can be made believable by noticing that, if we interpret the adjunction maps topologically as the birth of a circle and a saddle, this traces out a cobordism obtained from the ``wrong-way" RII cobordism by introducing a canceling pair of handles and applying a movie move.
\end{remark}

\begin{remark}
    Just as for $B_1$, the class $\theta^{(-2)}$ in $Tr_2(B_{121})$ in fact comes from an endomorphism $\theta_{1,2,3}$ of $B_{121}$ defined similarly to $\theta_{1,2}$: First open a bubble consisting of a 1- and 2- colored strand in the central 3-colored strand of $B_{121}$, then act by $\theta$ on the 1-colored strand inside the bubble and then collapse the bubble. This of course generalizes further, and most computations in $\extsbim$ are nicely expressible in terms of these divided-difference-of-$\theta$ endomorphisms.
\end{remark}

\begin{remark} \label{rem: exchange_renormalization}
    It is easy to check that the exchange move still preserves gradings if we use the normalized Rouquier complexes and partial trace of Remark \ref{rem: alg_hhh_regrading}.
\end{remark}

\FloatBarrier

\section{Tangle Functors} \label{sec: tangle_functors}

\subsection{General Framework} \label{subsec: general_framework}

In this section we take the first steps towards defining our oriented tangle invariant by assigning to each balanced oriented tangle diagram a functor between homotopy categories of Ext-Soergel bimodules. Before we begin, we will find it useful to artificially extend the objects of $\extsbim$ to encode orientation data while preserving the underlying morphism categories.

\begin{definition} \label{def: oriented_functor_cat}
    We denote by $\funcat$ the ($1$-)category consisting of the following data:

    \begin{itemize}
        \item Objects: Balanced orientation sequences $\aa = (a_1, \dots, a_{2n}) \in \{\pm 1\}^{2n}$ for $n \geq 0$.

        \item Morphisms from $\aa = (a_1, \dots, a_{2m})$ to $\bb = (b_1, \dots, b_{2n})$: Exact functors from $K^b(\extsbim_m)$ to $K^b(\extsbim_n)$, considered up to natural equivalence.
    \end{itemize}
\end{definition}

\begin{remark} \label{rem: oriented_identities}
    For any two balanced orientation sequences $\aa, \bb \in \{\pm 1\}^{2n}$, there is a distinguished isomorphism $\Psi_{\aa, \bb} \colon \aa \to \bb$ in $\funcat$ given by the identity functor on $K^b(\extsbim_n)$.
\end{remark}

Our goal for the remainder of this section is to construct a contravariant functor $\Phi \colon \otanbal \to \funcat$. On objects, $\Phi$ will send a balanced orientation sequence $\aa \in \{+1,-1\}^{2m}$ to the same sequence $\aa$; on morphisms, it will take a balanced oriented $(2m, 2n)$-tangle $\tau$ to an exact functor $\Phi(\tau) \colon K^b(\extsbim_m) \to K^b(\extsbim_n)$. We construct the data of $\Phi$ by assigning a functor to each vertical crossing and oriented handle in \S \ref{subsec: actions}. We check that this data gives a well-defined functor in \S \ref{subsec: invariance}.

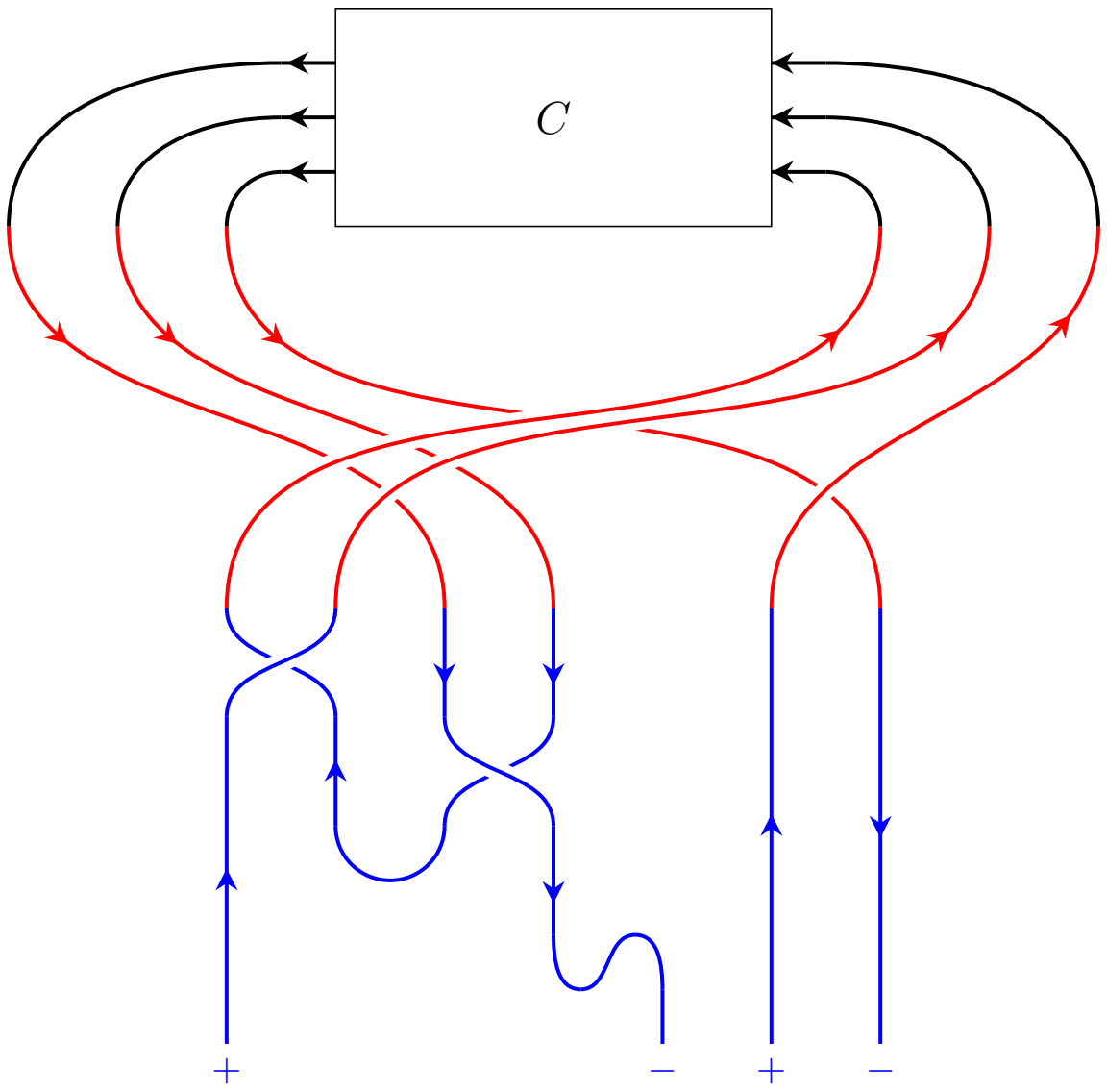
\begin{figure}[ht]
\centering
\scalebox{.75}{%

    \begin{tikzpicture}[anchorbase]


        \draw[instruction,->-=.5] (0,0) node[below]{$+$} to (0,3);
        \draw[instruction,->-=.5] (6,4) to (6,0) node[below]{$-$};
        \draw[instruction,->-=.5] (5,0) node[below]{$+$} to (5,4);

        \draw[instruction] (4,0) node[below]{$-$} to (4,.5);
        \draw[instruction] (4,.5) to[out=90,in=0] (3.75,1);
        \draw[instruction] (3.75,1) to[out=180,in=0] (3.25,.5);
        \draw[instruction] (3.25,.5) to[out=180,in=270] (3,1);

        \draw[instruction,-<-=.4] (3,1) to (3,2);

        \draw[instruction] (2,2) to[out=90,in=270] (3,3);
        \draw[instruction,over] (3,2) to[out=90,in=270] (2,3);

        \draw[instruction] (2,2) to[out=270,in=0] (1.5,1.5);
        \draw[instruction] (1.5,1.5) to[out=180,in=270] (1,2);

        \draw[instruction,->-=.5] (1,2) to (1,3);

        \draw[instruction] (1,3) to[out=90,in=270] (0,4);
        \draw[instruction,over] (0,3) to[out=90,in=270] (1,4);

        \draw[instruction,-<-=.4] (2,3) to (2,4);
        \draw[instruction,-<-=.4] (3,3) to (3,4);


        \draw[shuffle,-<-=.8]
            (2,4) to[out=90,in=270] (-2,7.5);

        \draw[shuffle,-<-=.8]
            (3,4) to[out=90,in=270] (-1,7.5);

        \draw[shuffle,-<-=.85]
            (6,4) to[out=90,in=270] (0,7.5);

        \draw[shuffle,->-=.85,over]
            (0,4) to[out=90,in=270] (6,7.5);

        \draw[shuffle,->-=.85,over]
            (1,4) to[out=90,in=270] (7,7.5);

        \draw[shuffle,->-=.8,over]
            (5,4) to[out=90,in=270] (8,7.5);


        \draw[register] (0,7.5) to[out=90,in=180] (.5,8);
        \draw[register] (-1,7.5) to[out=90,in=180] (.5,8.5);
        \draw[register] (-2,7.5) to[out=90,in=180] (.5,9);

        \draw[register] (6,7.5) to[out=90,in=0] (5.5,8);
        \draw[register] (7,7.5) to[out=90,in=0] (5.5,8.5);
        \draw[register] (8,7.5) to[out=90,in=0] (5.5,9);

        \draw[register,biarr] (.5,8) to (5.5,8);
        \draw[register,biarr] (.5,8.5) to (5.5,8.5);
        \draw[register,biarr] (.5,9) to (5.5,9);

        \filldraw[fill=white, draw=black]
            (1,7.5) rectangle (5,9.5);

        \node at (3,8.5) {\Large{$C$}};

    \end{tikzpicture}%
    }

    \caption{
        A visual representation of a balanced oriented $(6,4)$-tangle $\tau$, drawn in \textcolor{blue}{blue},
        acting on a complex $C \in K^b(\extsbim_3)$. These are separated by a shuffle braid, drawn in \textcolor{red}{red}. We assign to $\tau$ an exact functor $\Phi(\tau) \colon K^b(\extsbim_3) \to K^b(\extsbim_2)$.
    }
    \label{fig: general_framework}
\end{figure}

We will define $\Phi$ for each elementary tangle algebraically, but it will be useful to keep Figure \ref{fig: general_framework} in mind as topological motivation. This figure represents the action of a balanced oriented tangle $\tau$, drawn in \textcolor{blue}{blue}, on a complex $C \in K^b(\extsbim)$. We will sometimes refer to $\tau$ as the \textit{instruction tangle} and $C$ as the \textit{register complex}. The register complex and the instruction tangle are separated by a \textit{shuffle braid}, consisting of horizontal crossings and drawn in \textcolor{red}{red}. This braid is purely cosmetic and serves only to connect the positively and negatively oriented points at the boundary of $C$ with the positively and negatively oriented input points of $\tau$.

Our assignment of functors to elementary tangles is motivated by the idea that, topologically, one can slide an instruction tangle $\tau$ through the shuffle braid in such a way that, after a specified isotopy, the result parses as a diagrammatic description of a functor $\Phi(\tau)$ applied to $C$. After specifying some conventions regarding the action of cups and caps, the functor we assign to any tangle can be recovered from this picture.

\begin{remark} \label{rem: rotated_diagrams}
    The register complex $C$ in Figure \ref{fig: general_framework} is rotated relative to the diagrammatics introduced in Section \ref{sec: background}. We retain this rotated scheme for the remainder of this work. This is useful for several reasons. First, it dramatically simplifies the role of the shuffle braid in Figure \ref{fig: general_framework}; second, it affords a clear marker for what is a genuine oriented tangle diagram and what is a diagram encoding some complex of Ext-Soergel bimodules when both are present.
\end{remark}

\begin{remark} \label{rem: shuffle_convention}
    We take the convention that upwards-oriented strands in the shuffle braid pass over downwards-oriented strands. The opposite convention is also possible, and this choice would lead to the assignment of a different functor to oriented caps and cups. This difference is best viewed in light of Lemma \ref{lem: trivial_horizontals} below; under the other set of conventions, the second two horizontal crossings of Figure \ref{fig: horizontal_crossing_decompositions} would act trivially rather than the first two.
\end{remark}

\begin{warning} \label{warning: non_local}
    We wish to stress that $\Phi$ is \textit{not} a monoidal functor (indeed, we have not even equipped its target category $\funcat$ with a monoidal structure). In particular, there is no sense in which $\Phi$ is local with respect to horizontal concatenation $(-) \sqcup (-)$ on $\otanbal$. We are therefore in principle required to define $\Phi$ for each elementary tangle in the presence of any possible configuration of parallel upwards- and downwards-oriented vertical strands on either side.

    In practice, we can cut down on the amount of data needed to specify $\Phi$ by making use of the distinguished isomorphisms $\Psi_{\aa, \bb}$ from Remark \ref{rem: oriented_identities}. We define directly an assignment of $\Phi$ only to elementary tangles with all downwards-oriented vertical strands appearing to the left of all upwards-oriented vertical strands in any given region in which several parallel vertical strands are present, and we implicitly extend this assignment to arbitrary arrangements of parallel vertical strands by (pre- or post-)composition with the isomorphisms $\Psi_{\aa, \bb}$. Since $\Psi_{\aa, \bb}$ are identity functors after unwinding definitions, this has no effect on any explicit computations in $\funcat$. Nevertheless, the reader should keep in mind that the figures in Section \ref{subsec: actions} are really shorthand for assignments in full generality. See Lemma \ref{lem: trivial_horizontals} for another perspective on this convention.
\end{warning}

\FloatBarrier

\subsection{Action of Crossings and Handles} \label{subsec: actions}

\subsubsection{Crossings} \label{subsubsec: crossings}

\begin{defn} \label{def: upwards_crossing}
    Let $\tau$ consist of a single crossing $\sigma_r^{\pm}$ between the $r$ and $(r+1)^{st}$ upwards-oriented strands of the output. Then we set $\Phi(\tau) := - \star F(\sigma_r^{\pm})$.

    We depict this assignment diagrammatically as follows:

    
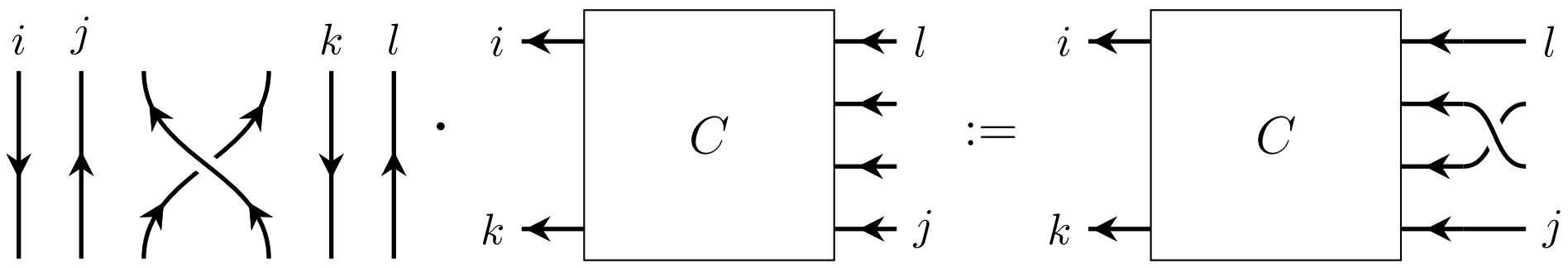
\begin{figure}[ht]
        \centering

        \begin{tikzpicture}[anchorbase]
            \draw[webs,-<-=.5] (0,0) to (0,1.5) node[above]{$i$};
            \draw[webs,->-=.5] (.5,0) to (.5,1.5) node[above]{$j$};

            \draw[webs,->-=.2,->-=.8] (1,0) to[out=90,in=270] (2,1.5);
            \draw[webs,over,->-=.2,->-=.8] (2,0) to[out=90,in=270] (1,1.5);

            \draw[webs,-<-=.5] (2.5,0) to (2.5,1.5) node[above]{$k$};
            \draw[webs,->-=.5] (3,0) to (3,1.5) node[above]{$l$};
        \end{tikzpicture}
        {\huge{$\cdot$}}
        \begin{tikzpicture}[anchorbase]
            \draw[webs,biarr] (0,0) node[left]{$k$} to (3,0) node[right]{$j$};
            \draw[webs,biarr] (0,1.5) node[left]{$i$} to (3,1.5) node[right]{$l$};
            \draw[webs,-<-=.9] (1,.5) to (3,.5);
            \draw[webs,-<-=.9] (1,1) to (3,1);

            \filldraw[fill=white, draw=black] (.5,-.25) rectangle (2.5,1.75);
            \node at (1.5,.75) {\Large{$C$}};
        \end{tikzpicture}
        {\Large{$:=$}}
        \begin{tikzpicture}[anchorbase]
            \draw[webs,biarr] (0,0) node[left]{$k$} to (3,0);
            \draw[webs,biarr] (0,1.5) node[left]{$i$} to (3,1.5);
            \draw[webs,-<-=.9] (1,.5) to (3,.5);
            \draw[webs,-<-=.9] (1,1) to (3,1);

            \filldraw[fill=white, draw=black] (.5,-.25) rectangle (2.5,1.75);
            \node at (1.5,.75) {\Large{$C$}};

            \draw[webs] (3,.5) to[out=0,in=180] (3.5,1);
            \draw[webs,over] (3,1) to[out=0,in=180] (3.5,.5);
            \draw[webs] (3,0) to (3.5,0) node[right]{$j$};
            \draw[webs] (3,1.5) to (3.5,1.5) node[right]{$l$};
        \end{tikzpicture}  
    \end{figure}
\end{defn}

\FloatBarrier

\begin{defn} \label{def: downwards_crossing}
    Let $\tau$ consist of a single crossing $\sigma_i^{\pm}$ between the $i$ and $(i+1)^{st}$ downwards-oriented strands of the output. Then we set $\Phi(\tau) := F(\sigma_i^{\pm}) \star -$.

    We depict this assignment diagrammatically as follows:

    \begin{figure}[ht]
        \centering

        \begin{tikzpicture}[anchorbase]
            \draw[webs,-<-=.5] (0,0) to (0,1.5) node[above]{$i$};
            \draw[webs,->-=.5] (.5,0) to (.5,1.5) node[above]{$j$};

            \draw[webs,-<-=.2,-<-=.8] (1,0) to[out=90,in=270] (2,1.5);
            \draw[webs,over,-<-=.2,-<-=.8] (2,0) to[out=90,in=270] (1,1.5);

            \draw[webs,-<-=.5] (2.5,0) to (2.5,1.5) node[above]{$k$};
            \draw[webs,->-=.5] (3,0) to (3,1.5) node[above]{$l$};
        \end{tikzpicture}
        {\huge{$\cdot$}}
        \begin{tikzpicture}[anchorbase]
            \draw[webs,biarr] (0,0) node[left]{$k$} to (3,0) node[right]{$j$};
            \draw[webs,biarr] (0,1.5) node[left]{$i$} to (3,1.5) node[right]{$l$};
            \draw[webs,-<-=.1] (0,.5) to (2,.5);
            \draw[webs,-<-=.1] (0,1) to (2,1);

            \filldraw[fill=white, draw=black] (.5,-.25) rectangle (2.5,1.75);
            \node at (1.5,.75) {\Large{$C$}};
        \end{tikzpicture}
        {\Large{$:=$}}
        \begin{tikzpicture}[anchorbase]
            \draw[webs,biarr] (0,0) to (3,0) node[right]{$j$};
            \draw[webs,biarr] (0,1.5) to (3,1.5) node[right]{$l$};
            \draw[webs,-<-=.1] (0,.5) to (2,.5);
            \draw[webs,-<-=.1] (0,1) to (2,1);

            \filldraw[fill=white, draw=black] (.5,-.25) rectangle (2.5,1.75);
            \node at (1.5,.75) {\Large{$C$}};

            \draw[webs] (0,.5) to[out=180,in=0] (-.5,1);
            \draw[webs,over] (0,1) to[out=180,in=18] (-.5,.5);
            \draw[webs] (0,0) to (-.5,0) node[left]{$k$};
            \draw[webs] (0,1.5) to (-.5,1.5) node[left]{$i$};
        \end{tikzpicture}  
    \end{figure}
\end{defn}

\FloatBarrier

\begin{rem} \label{rem: vertical_crossings_motivation}
    In the context of Figure \ref{fig: general_framework}, we view vertical crossings in the instruction tangle as sliding directly through the shuffle braid to meet the register complex, either on the right for upwards-oriented crossings or on the left for downwards-oriented crossings. See Figure \ref{fig: upwards_crossing_figure} for an example of this perspective applied to an upwards-oriented crossing.

    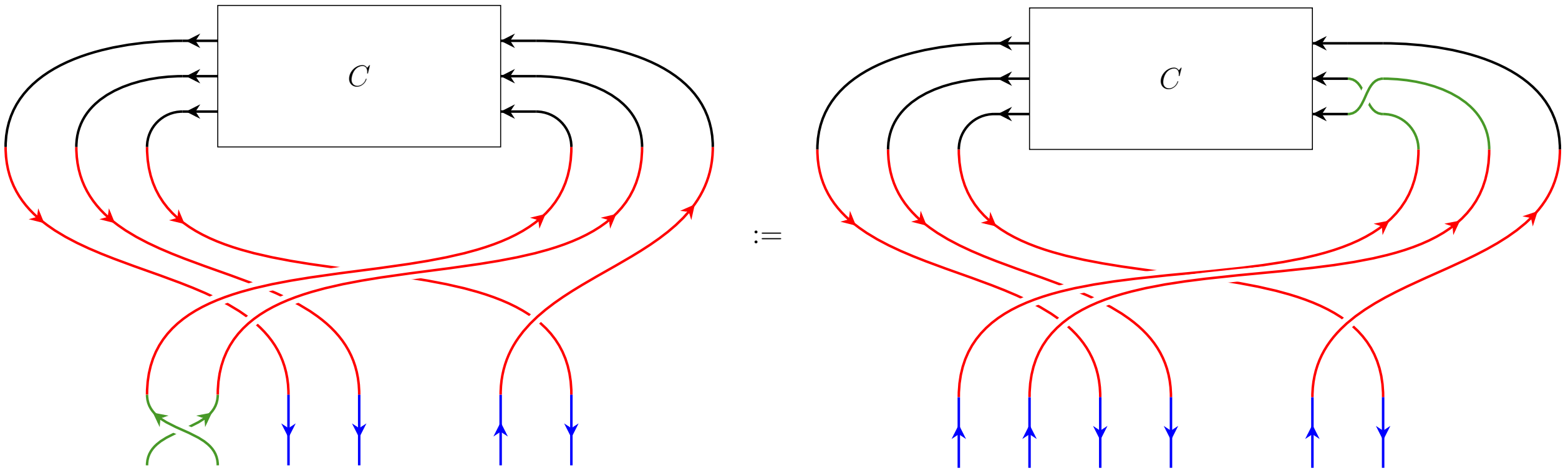
\begin{figure}[ht]
    \centering

    \scalebox{0.75}{%
    \begin{tikzpicture}[anchorbase]

    
        \draw[instruction,->-=.5] (2,4) to (2,3);
        \draw[instruction,->-=.5] (3,4) to (3,3);
        \draw[instruction,->-=.5] (5,3) to (5,4);
        \draw[instruction,->-=.5] (6,4) to (6,3);
        \draw[new,->-=.75] (0,3) to[out=90,in=270] (1,4);
        \draw[new, over, ->-=.75] (1,3) to[out=90,in=270] (0,4);


        \draw[shuffle,-<-=.8]
            (2,4) to[out=90,in=270] (-2,7.5);

        \draw[shuffle,-<-=.8]
            (3,4) to[out=90,in=270] (-1,7.5);

        \draw[shuffle,-<-=.85]
            (6,4) to[out=90,in=270] (0,7.5);

        \draw[shuffle,->-=.85,over]
            (0,4) to[out=90,in=270] (6,7.5);

        \draw[shuffle,->-=.85,over]
            (1,4) to[out=90,in=270] (7,7.5);

        \draw[shuffle,->-=.8,over]
            (5,4) to[out=90,in=270] (8,7.5);


        \draw[register] (0,7.5) to[out=90,in=180] (.5,8);
        \draw[register] (-1,7.5) to[out=90,in=180] (.5,8.5);
        \draw[register] (-2,7.5) to[out=90,in=180] (.5,9);

        \draw[register] (6,7.5) to[out=90,in=0] (5.5,8);
        \draw[register] (7,7.5) to[out=90,in=0] (5.5,8.5);
        \draw[register] (8,7.5) to[out=90,in=0] (5.5,9);

        \draw[register,biarr] (.5,8) to (5.5,8);
        \draw[register,biarr] (.5,8.5) to (5.5,8.5);
        \draw[register,biarr] (.5,9) to (5.5,9);

        \filldraw[fill=white, draw=black]
            (1,7.5) rectangle (5,9.5);

        \node at (3,8.5) {\Large{$C$}};

    \end{tikzpicture}
    \quad {\Large{$:=$}} \quad
    \begin{tikzpicture}[anchorbase]

    
        \draw[instruction,->-=.5] (2,4) to (2,3);
        \draw[instruction,->-=.5] (3,4) to (3,3);
        \draw[instruction,->-=.5] (5,3) to (5,4);
        \draw[instruction,->-=.5] (6,4) to (6,3);
        \draw[instruction,->-=.5] (0,3) to (0,4);
        \draw[instruction, ->-=.5] (1,3) to (1,4);


        \draw[shuffle,-<-=.8]
            (2,4) to[out=90,in=270] (-2,7.5);

        \draw[shuffle,-<-=.8]
            (3,4) to[out=90,in=270] (-1,7.5);

        \draw[shuffle,-<-=.85]
            (6,4) to[out=90,in=270] (0,7.5);

        \draw[shuffle,->-=.85,over]
            (0,4) to[out=90,in=270] (6.5,7.5);

        \draw[shuffle,->-=.85,over]
            (1,4) to[out=90,in=270] (7.5,7.5);

        \draw[shuffle,->-=.8,over]
            (5,4) to[out=90,in=270] (8.5,7.5);


        \draw[register] (0,7.5) to[out=90,in=180] (.5,8);
        \draw[register] (-1,7.5) to[out=90,in=180] (.5,8.5);
        \draw[register] (-2,7.5) to[out=90,in=180] (.5,9);

        \draw[new] (6.5,7.5) to[out=90,in=0] (6,8);
        \draw[new] (7.5,7.5) to[out=90,in=0] (6,8.5);
        \draw[register] (8.5,7.5) to[out=90,in=0] (6,9);

        \draw[new] (5.5,8.5) to[out=0,in=180] (6,8);
        \draw[new, over] (5.5,8) to[out=0,in=180] (6,8.5);
        \draw[register] (5.5,9) to (6,9);

        \draw[register,biarr] (.5,8) to (5.5,8);
        \draw[register,biarr] (.5,8.5) to (5.5,8.5);
        \draw[register,biarr] (.5,9) to (5.5,9);

        \filldraw[fill=white, draw=black]
            (1,7.5) rectangle (5,9.5);

        \node at (3,8.5) {\Large{$C$}};

    \end{tikzpicture}%
    }

    \caption{
        An upwards-oriented crossing sliding from the instruction tangle to the register complex through the shuffle braid.
    }
    \label{fig: upwards_crossing_figure}
\end{figure}
\end{rem}

\subsubsection{Handles} \label{subsubsec: handles}

\begin{defn} \label{def: cups}
    Let $\tau$ consist of a single clockwise- or counterclockwise-oriented cup to the left of $k$ downwards-oriented strands and to the right of $j$ upwards-oriented strands. Then we set 
    \[
    \Phi(\tau) := F(\sigma_k\sigma_{k-1} \dots \sigma_1) \star (I(-)) \star F(\sigma_1\sigma_2 \dots \sigma_j).
    \]

    We depict this assignment diagrammatically\footnote{Here we have chosen a clockwise orientation, but the functor assigned to a counterclockwise cup is identical.} as follows:

    \begin{figure}[ht]
        \centering

        \begin{tikzpicture}[anchorbase]
            \draw[webs,-<-=.5] (0,0) to (0,1.5) node[above]{$i$};
            \draw[webs,->-=.5] (.5,0) to (.5,1.5) node[above]{$j$};

            \draw[webs,->-=.5] (1,0) to (1,.5);
            \draw[webs] (1,.5) to[out=90,in=180] (1.25,.75);
            \draw[webs] (1.25,.75) to[out=0,in=90] (1.5,.5);
            \draw[webs,->-=.5] (1.5,.5) to (1.5,0);

            \draw[webs,-<-=.5] (2,0) to (2,1.5) node[above]{$k$};
            \draw[webs,->-=.5] (2.5,0) to (2.5,1.5) node[above]{$l$};
        \end{tikzpicture}
        {\huge{$\cdot$}}
        \begin{tikzpicture}[anchorbase]
            \draw[webs,biarr] (0,0) node[left]{$k$} to (3,0) node[right]{$j$};
            \draw[webs,biarr] (0,1) node[left]{$i$} to (3,1) node[right]{$l$};

            \filldraw[fill=white, draw=black] (.5,-.25) rectangle (2.5,1.25);
            \node at (1.5,.5) {\Large{$C$}};
        \end{tikzpicture}
        {\Large{$:=$}}
        \begin{tikzpicture}[anchorbase]
            \draw[webs,biarr] (0,0) to (3,0);
            \draw[webs,biarr] (0,1) to (3,1);

            \filldraw[fill=white, draw=black] (.5,-.25) rectangle (2.5,1.25);
            \node at (1.5,.5) {\Large{$C$}};

            \draw[webs,-<-=.5] (0,-.5) to (3,-.5);

            \draw[webs] (3,-.5) to[out=0,in=180] (4,.5);
            \draw[webs,over] (3,0) to (4,0) node[right]{$j$};
            \draw[webs] (3,1) to (4,1) node[right]{$l$};

            \draw[webs] (0,0) to (-1,0) node[left]{$k$};
            \draw[webs] (0,1) to (-1,1) node[left]{$i$};
            \draw[webs,over] (0,-.5) to[out=180,in=0] (-1,.5);
        \end{tikzpicture}  
    \end{figure}
    
\end{defn}

\FloatBarrier

\begin{defn} \label{def: caps}
    Let $\tau$ consist of a single clockwise- or counterclockwise-oriented cap to the left of $k$ downwards-oriented strands and to the right of $j$ upwards-oriented strands. Then we set 
    \[
    \Phi(\tau) := Tr \left( F(\sigma_1^{-1} \sigma_2^{-1} \dots \sigma_k^{-1}) \star - \star F(\sigma_j^{-1} \sigma_{j-1}^{-1} \dots \sigma_1^{-1})\right).
    \]

    We depict this assignment diagrammatically\footnote{Here we have chosen a counterclockwise orientation, but as was the case for cups, the functor assigned to a clockwise cap is identical.} as follows:

    \begin{figure}[ht]
        \centering

        \begin{tikzpicture}[anchorbase]
            \draw[webs,-<-=.5] (0,0) to (0,1.5) node[above]{$i$};
            \draw[webs,->-=.5] (.5,0) to (.5,1.5) node[above]{$j$};

            \draw[webs,->-=.5] (1,1.5) to (1,1);
            \draw[webs] (1,1) to[out=270,in=180] (1.25,.75);
            \draw[webs] (1.25,.75) to[out=0,in=270] (1.5,1);
            \draw[webs,->-=.5] (1.5,1) to (1.5,1.5);

            \draw[webs,-<-=.5] (2,0) to (2,1.5) node[above]{$k$};
            \draw[webs,->-=.5] (2.5,0) to (2.5,1.5) node[above]{$l$};
        \end{tikzpicture}
        {\huge{$\cdot$}}
        \begin{tikzpicture}[anchorbase]
            \draw[webs,biarr] (0,0) node[left]{$k$} to (3,0) node[right]{$j$};
            \draw[webs,biarr] (0,1) node[left]{$i$} to (3,1) node[right]{$l$};
            \draw[webs,biarr] (0,.5) to (3,.5);

            \filldraw[fill=white, draw=black] (.5,-.25) rectangle (2.5,1.25);
            \node at (1.5,.5) {\Large{$C$}};
        \end{tikzpicture}
        {\Large{$:=$}}
        \begin{tikzpicture}[anchorbase]
            \draw[webs,biarr] (0,0) to (3,0);
            \draw[webs,biarr] (0,1) to (3,1);
            \draw[webs,biarr] (0,.5) to (3,.5);

            \filldraw[fill=white, draw=black] (.5,-.25) rectangle (2.5,1.25);
            \node at (1.5,.5) {\Large{$C$}};

            \draw[webs] (0,0) to (-1,0) node[left]{$k$};
            \draw[webs] (0,1) to (-1,1) node[left]{$i$};
            \draw[webs,over] (0,.5) to[out=180,in=90] (-.5,-.5);
            
            \draw[webs] (3,.5) to[out=0,in=90] (3.5,-.5);
            \draw[webs,over] (3,0) to (4,0) node[right]{$j$};
            \draw[webs] (3,1) to (4,1) node[right]{$l$};
            
            \draw[webs] (3.5,-.5) to[out=270,in=0] (3.25,-.75);
            \draw[webs] (-.5,-.5) to[out=270,in=180] (-.25,-.75);
            \draw[webs] (-.25,-.75) to (3.25,-.75);
        \end{tikzpicture}  
    \end{figure}
    
\end{defn}

\FloatBarrier

\begin{rem} \label{rem: handles_motivation}
    As in Remark \ref{rem: vertical_crossings_motivation}, in the context of Figure \ref{fig: general_framework}, we view cups and caps in the instruction tangle as sliding through the shuffle braid to meet the register complex. This requires either the addition (resp. subtraction) of a single strand to the register complex in the case of a cup (resp. cap). We adopt the convention that this strand is added on the bottom. The sign of the crossings appended to the register complex is fixed by our previous convention for shuffle braids set in Remark \ref{rem: shuffle_convention}. See Figure \ref{fig: handles_figure} for one example each of this perspective applied to a cap and a cup.

    \begin{figure}[ht]
    \centering

    \scalebox{0.75}{%
    \begin{tikzpicture}[anchorbase]

    
        \draw[instruction,->-=.5] (2,4) to (2,3);
        \draw[instruction,->-=.5] (3,4) to (3,3);
        \draw[instruction,->-=.5] (5,3) to (5,4);
        \draw[instruction,->-=.5] (6,4) to (6,3);
        \draw[instruction,->-=.5] (0,3) to (0,4);
        \draw[instruction, ->-=.5] (1,3) to (1,4);

        \draw[new, -<-=.99] (4.5,3) to[out=90,in=0] (4,3.5);
        \draw[new] (4,3.5) to[out=180,in=90] (3.5,3);


        \draw[shuffle,-<-=.8]
            (2,4) to[out=90,in=270] (-2,7.5);

        \draw[shuffle,-<-=.8]
            (3,4) to[out=90,in=270] (-1,7.5);

        \draw[shuffle,-<-=.85]
            (6,4) to[out=90,in=270] (0,7.5);

        \draw[shuffle,->-=.85,over]
            (0,4) to[out=90,in=270] (6,7.5);

        \draw[shuffle,->-=.85,over]
            (1,4) to[out=90,in=270] (7,7.5);

        \draw[shuffle,->-=.8,over]
            (5,4) to[out=90,in=270] (8,7.5);


        \draw[register] (0,7.5) to[out=90,in=180] (.5,8);
        \draw[register] (-1,7.5) to[out=90,in=180] (.5,8.5);
        \draw[register] (-2,7.5) to[out=90,in=180] (.5,9);

        \draw[register] (6,7.5) to[out=90,in=0] (5.5,8);
        \draw[register] (7,7.5) to[out=90,in=0] (5.5,8.5);
        \draw[register] (8,7.5) to[out=90,in=0] (5.5,9);

        \draw[register,biarr] (.5,8) to (5.5,8);
        \draw[register,biarr] (.5,8.5) to (5.5,8.5);
        \draw[register,biarr] (.5,9) to (5.5,9);

        \filldraw[fill=white, draw=black]
            (1,7.5) rectangle (5,9.5);

        \node at (3,8.5) {\Large{$C$}};

    \end{tikzpicture}%
    \quad {\Large{$:=$}} \quad
    \begin{tikzpicture}[anchorbase]


        \draw[instruction,->-=.5] (2,4) to (2,3);
        \draw[instruction,->-=.5] (3,4) to (3,3);
        \draw[instruction,->-=.5] (5,3) to (5,4);
        \draw[instruction,->-=.5] (6,4) to (6,3);
        \draw[instruction,->-=.5] (0,3) to (0,4);
        \draw[instruction, ->-=.5] (1,3) to (1,4);

        \draw[new, -<-=.5] (4.5,3) to (4.5,4);
        \draw[new, ->-=.5] (3.5,3) to (3.5,4);


        \draw[shuffle,-<-=.8]
            (2,4) to[out=90,in=270] (-3,8);

        \draw[shuffle,-<-=.8]
            (3,4) to[out=90,in=270] (-2,8);

        \draw[shuffle,-<-=.85]
            (6,4) to[out=90,in=270] (-1,8);

        \draw[new,-<-=.8]
            (4.5,4) to[out=90,in=270] (-1.5,8);

        \draw[shuffle,->-=.85,over]
            (0,4) to[out=90,in=270] (7,8);

        \draw[shuffle,->-=.85,over]
            (1,4) to[out=90,in=270] (8,8);

        \draw[shuffle,->-=.8,over]
            (5,4) to[out=90,in=270] (9,8);

        \draw[new,->-=.8,over]
            (3.5,4) to[out=90,in=270] (8.5,8);


        \draw[register] (-1,8) to[out=90,in=180] (-.5,8.5);
        \draw[new] (-1.5,8) to[out=90,in=180] (-.5,8.75);
        \draw[register] (-2,8) to[out=90,in=180] (-.5,9);
        \draw[register] (-3,8) to[out=90,in=180] (-.5,9.5);

        \draw[register] (7,8) to[out=90,in=0] (6.5,8.5);
        \draw[register] (8,8) to[out=90,in=0] (6.5,9);
        \draw[new] (8.5,8) to[out=90,in=0] (6.5,9.25);
        \draw[register] (9,8) to[out=90,in=0] (6.5,9.5);

        \draw[new] (6.5,9.25) to[out=180,in=0] (5.5,7.75);
        \draw[register, over] (6.5,9) to (5.5,9);
        \draw[register, over] (6.5,8.5) to (5.5,8.5);
        \draw[register] (6.5,9.5) to (5.5,9.5);

        \draw[register] (.5,8.5) to (-.5,8.5);
        \draw[register] (.5,9) to (-.5,9);
        \draw[register] (.5,9.5) to (-.5,9.5);
        \draw[new, over] (.5,7.75) to[out=180,in=0] (-.5,8.75);

        \draw[register,biarr] (.5,8.5) to (5.5,8.5);
        \draw[register,biarr] (.5,9) to (5.5,9);
        \draw[register,biarr] (.5,9.5) to (5.5,9.5);
        \draw[new,biarr] (.5,7.75) to (5.5,7.75);

        \filldraw[fill=white, draw=black]
            (1,8) rectangle (5,10);

        \node at (3,9) {\Large{$C$}};

    \end{tikzpicture}%
    }
    \\
    \vspace{1em}
    \scalebox{0.75}{%
    \begin{tikzpicture}[anchorbase]

    
        \draw[instruction,->-=.5] (2,4) to (2,3);
        \draw[instruction,->-=.5] (3,4) to (3,3);
        \draw[instruction,->-=.5] (5,3) to (5,4);
        \draw[instruction,->-=.5] (6,4) to (6,3);
        \draw[instruction,->-=.5] (0,3) to (0,4);
        \draw[instruction, ->-=.5] (1,3) to (1,4);

        \draw[new, ->-=.99] (3,3) to[out=270,in=180] (4,2.5);
        \draw[new] (4,2.5) to[out=0,in=270] (5,3);


        \draw[shuffle,-<-=.8]
            (2,4) to[out=90,in=270] (-2,7.5);

        \draw[shuffle,-<-=.8]
            (3,4) to[out=90,in=270] (-1,7.5);

        \draw[shuffle,-<-=.85]
            (6,4) to[out=90,in=270] (0,7.5);

        \draw[shuffle,->-=.85,over]
            (0,4) to[out=90,in=270] (6,7.5);

        \draw[shuffle,->-=.85,over]
            (1,4) to[out=90,in=270] (7,7.5);

        \draw[shuffle,->-=.8,over]
            (5,4) to[out=90,in=270] (8,7.5);


        \draw[register] (0,7.5) to[out=90,in=180] (.5,8);
        \draw[register] (-1,7.5) to[out=90,in=180] (.5,8.5);
        \draw[register] (-2,7.5) to[out=90,in=180] (.5,9);

        \draw[register] (6,7.5) to[out=90,in=0] (5.5,8);
        \draw[register] (7,7.5) to[out=90,in=0] (5.5,8.5);
        \draw[register] (8,7.5) to[out=90,in=0] (5.5,9);

        \draw[register,biarr] (.5,8) to (5.5,8);
        \draw[register,biarr] (.5,8.5) to (5.5,8.5);
        \draw[register,biarr] (.5,9) to (5.5,9);

        \filldraw[fill=white, draw=black]
            (1,7.5) rectangle (5,9.5);

        \node at (3,8.5) {\Large{$C$}};

    \end{tikzpicture}%
    \quad {\Large{$:=$}} \quad
    \begin{tikzpicture}[anchorbase]


        \draw[instruction,->-=.5] (2,4) to (2,3);
        \draw[instruction,->-=.5] (6,4) to (6,3);
        \draw[instruction,->-=.5] (0,3) to (0,4);
        \draw[instruction, ->-=.5] (1,3) to (1,4);


        \draw[shuffle,-<-=.8]
            (2,4) to[out=90,in=270] (-3,8);

        \draw[shuffle,-<-=.85]
            (6,4) to[out=90,in=270] (-1,8);

        \draw[shuffle,->-=.85,over]
            (0,4) to[out=90,in=270] (7,8);

        \draw[shuffle,->-=.85,over]
            (1,4) to[out=90,in=270] (8,8);


        \draw[register] (-1,8) to[out=90,in=180] (-.5,8.5);
        \draw[register] (-3,8) to[out=90,in=180] (-.5,9.5);

        \draw[register] (7,8) to[out=90,in=0] (6.5,8.5);
        \draw[register] (8,8) to[out=90,in=0] (6.5,9);

        \draw[new] (5.5,9.5) to[out=0,in=90] (6.5,8);
        \draw[register, over] (6.5,8.5) to (5.5,8.5);
        \draw[register, over] (6.5,9) to (5.5,9);

        \draw[new] (6.5,8) to[out=270,in=0] (6, 7.5);
        \draw[new, -<-=.5] (6,7.5) to (.5,7.5);
        \draw[new] (.5,7.5) to[out=180,in=270] (0,8);

        \draw[register] (.5,8.5) to (-.5,8.5);
        \draw[new, over] (.5,9) to[out=180,in=90] (0,8);
        \draw[register] (.5,9.5) to (-.5,9.5);

        \draw[register,biarr] (.5,8.5) to (5.5,8.5);
        \draw[register,biarr] (.5,9) to (5.5,9);
        \draw[register,biarr] (.5,9.5) to (5.5,9.5);

        \filldraw[fill=white, draw=black]
            (1,8) rectangle (5,10);

        \node at (3,9) {\Large{$C$}};

    \end{tikzpicture}%
    }

    \caption{
        A cup and a cap sliding from the instruction tangle to the register complex through a shuffle braid. Note that the number of strands in the register complex increases or decreases by $1$.
    }
    \label{fig: handles_figure}
\end{figure}
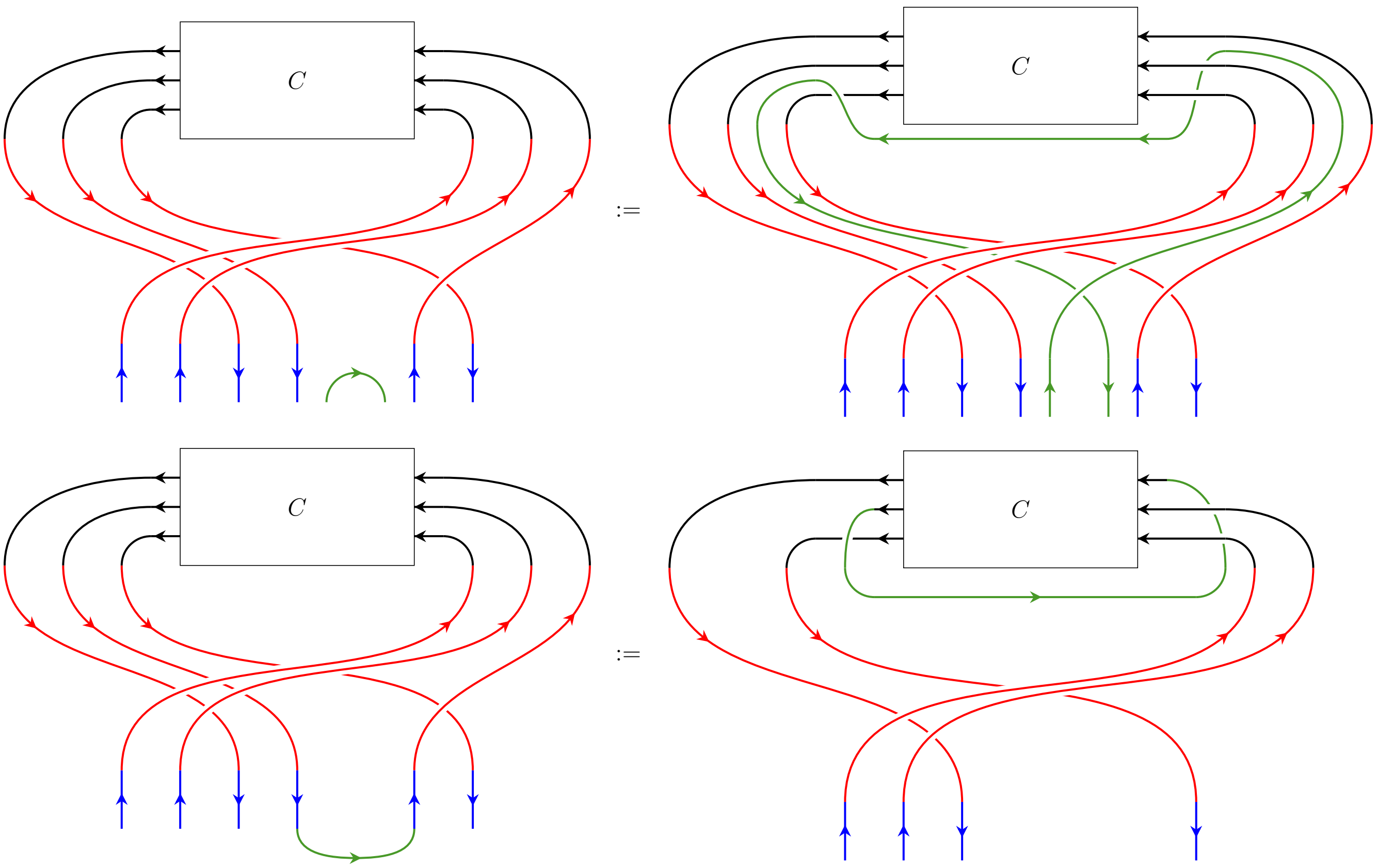
\end{rem}

\begin{remark} \label{remark: tangle_regrading}
    From this point onwards, we will exclusively use the renormalized crossing and partial trace functors described in Remark \ref{rem: alg_hhh_regrading} so that the Markov moves of Proposition \ref{prop: markov} hold with no regradings.
\end{remark}

\FloatBarrier

\subsubsection{Trivial Horizontal Crossings}

Before checking that $\Phi$ is a well-defined functor, we pause to point out an extremely useful result which allows us to disregard many horizontal crossings when computing in $\funcat$.

\begin{lemma} \label{lem: trivial_horizontals}
    All horizontal crossings in which an upwards-oriented strand crosses over a downwards-oriented strand are sent by $\Phi$ to the distinguished isomorphisms $\Psi_{\aa, \bb}$ in $\funcat$.
\end{lemma}

\begin{proof}
    The relevant crossings are the first two from Figure \ref{fig: horizontal_crossing_decompositions}. We illustrate only the action of the first crossing; the second is similar. Suppose the horizontal crossing in question occurs to the left of $i$ downwards and $j$ upwards strands and to the right of $k$ downwards and $l$ upwards strands. We compute its action on a complex $C$ directly; see Figure \ref{fig: horizontal_crossing_trivial} below.

    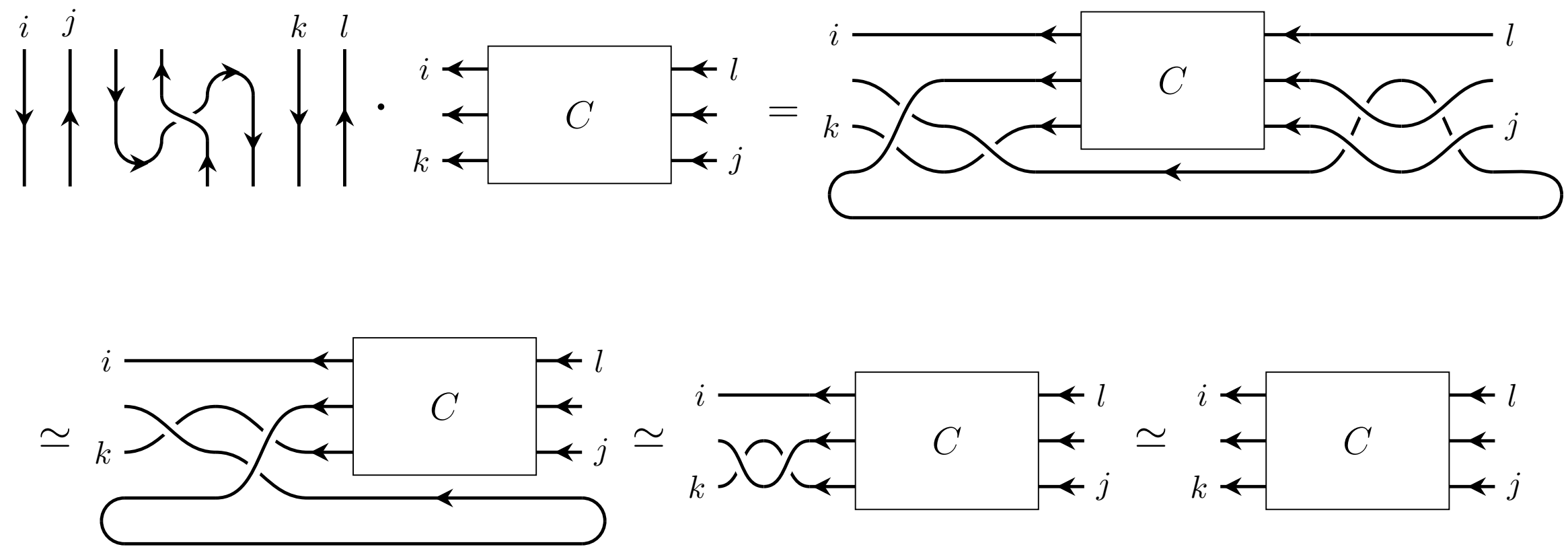
\begin{figure}[ht]
    \centering

    \scalebox{.9}{%
    \begin{tikzpicture}[anchorbase]
    \draw[webs,-<-=.5] (0,0) to (0,1.5) node[above]{$i$};
    \draw[webs,->-=.5] (0.5,0) to (0.5,1.5) node[above]{$j$};
    \draw[webs,-<-=.5] (3,0) node[below]{\textcolor{white}{$k$}} to (3,1.5) node[above]{$k$};
    \draw[webs,->-=.5] (3.5,0) to (3.5,1.5) node[above]{$l$};

    \begin{scope}[xshift=1cm]
    \draw[webs,->-=.5] (0,1.5) to (0,.5);
    \draw[webs,->-=.5] (.5,1) to (.5,1.5);
        
    \draw[webs] (.5,.5) to[out=90,in=270] (1,1);
    \draw[webs, over] (1,.5) to[out=90,in=270] (.5,1);
        
    \draw[webs,->-=.99] (0,.5) to[out=270,in=180] (.25,.25);
    \draw[webs] (.25,.25) to[out=0,in=270] (.5,.5);
    \draw[webs,->-=.99] (1,1) to[out=90,in=180] (1.25,1.25);
    \draw[webs] (1.25,1.25) to[out=0,in=90] (1.5,1);
    \draw[webs,->-=.5] (1,0) to (1,.5);
    \draw[webs,-<-=.5] (1.5,0) to (1.5,1);
    \end{scope}
    \end{tikzpicture}%
    }
    {\huge{$\cdot$}}
    \scalebox{.9}{%
    \begin{tikzpicture}[anchorbase]
    \draw[webs,biarr] (0,0) node[left]{$k$} to (3,0) node[right]{$j$};
    \draw[webs,biarr] (0,.5) to (3,.5);
    \draw[webs,biarr] (0,1) node[left]{$i$} to (3,1) node[right]{$l$};
    \filldraw[fill=white, draw=black]
        (.5,-.25) rectangle (2.5,1.25);
    \node at (1.5,.5) {\Large{$C$}};
    \end{tikzpicture}%
    }
    {\Large{$=$}}
    \scalebox{.9}{%
    \begin{tikzpicture}[anchorbase]
    \draw[webs,biarr] (0,0) to (3,0);
    \draw[webs,biarr] (0,.5) to (3,.5);
    \draw[webs,biarr] (0,1) to (3,1);
    \filldraw[fill=white, draw=black]
        (.5,-.25) rectangle (2.5,1.25);
    \node at (1.5,.5) {\Large{$C$}};

    \draw[webs] (0,0) to[out=180,in=0] (-1,-.5);
    \draw[webs, over] (0,-.5) to[out=180,in=0] (-1,0);
    \draw[webs] (0,.5) to (-1,.5);
    \draw[webs] (-1,0) to[out=180,in=0] (-2,.5);
    \draw[webs] (-1,-.5) to[out=180,in=0] (-2,0) node[left]{$k$};
    \draw[webs, over] (-1,.5) to[out=180,in=0] (-2,-.5);
    \draw[webs] (0,1) to (-2,1) node[left]{$i$};

    \draw[webs,-<-=.5] (0,-.5) to (3,-.5);

    \draw[webs] (3,-.5) to[out=0,in=180] (4,.5);
    \draw[webs, over] (3,0) to[out=0,in=180] (4,-.5);
    \draw[webs, over] (3,.5) to[out=0,in=180] (4,0);
    \draw[webs] (4,.5) to[out=0,in=180] (5,-.5);
    \draw[webs, over] (4,0) to[out=0,in=180] (5,.5);
    \draw[webs, over] (4,-.5) to[out=0,in=180] (5,0) node[right]{$j$};
    \draw[webs] (3,1) to (5,1) node[right]{$l$};

    \draw[webs] (5,-.5) to[out=0,in=90] (5.75,-.75);
    \draw[webs] (5.75,-.75) to[out=270,in=0] (5.5,-1);
    \draw[webs] (-2,-.5) to[out=180,in=90] (-2.25,-0.75);
    \draw[webs] (-2.25,-0.75) to[out=270,in=180] (-2,-1);

    \draw[webs] (-2,-1) to (5.5,-1);
    \end{tikzpicture}%
    }
    \\
    \vspace{1cm}
    {\Large{$\simeq$}}
    \scalebox{.9}{%
    \begin{tikzpicture}[anchorbase]
    \draw[webs,biarr] (0,0) to (3,0) node[right]{$j$};
    \draw[webs,biarr] (0,.5) to (3,.5);
    \draw[webs,biarr] (0,1) to (3,1) node[right]{$l$};
    \filldraw[fill=white, draw=black]
        (.5,-.25) rectangle (2.5,1.25);
    \node at (1.5,.5) {\Large{$C$}};

    \draw[webs] (0,-.5) to[out=180,in=0] (-1,0);
    \draw[webs] (0,0) to[out=180,in=0] (-1,.5);
    \draw[webs, over] (0,.5) to[out=180,in=0] (-1,-.5);
    \draw[webs] (-1,.5) to[out=180,in=0] (-2,0) node[left]{$k$};
    \draw[webs, over] (-1,0) to[out=180,in=0] (-2,.5);
    \draw[webs] (0,1) to (-2,1) node[left]{$i$};
    \draw[webs] (-1,-.5) to (-2,-.5);

    \draw[webs,-<-=.5] (0,-.5) to (3,-.5);

    \draw[webs] (3,-.5) to[out=0,in=90] (3.25,-.75);
    \draw[webs] (3.25,-.75) to[out=270,in=0] (3,-1);
    \draw[webs] (-2,-.5) to[out=180,in=90] (-2.25,-0.75);
    \draw[webs] (-2.25,-0.75) to[out=270,in=180] (-2,-1);

    \draw[webs] (-2,-1) to (3,-1);
    \end{tikzpicture}%
    }
    {\Large{$\simeq$}}
    \scalebox{.9}{%
    \begin{tikzpicture}[anchorbase]
    \draw[webs,biarr] (0,0) to (3,0) node[right]{$j$};
    \draw[webs,biarr] (0,.5) to (3,.5);
    \draw[webs,biarr] (0,1) to (3,1) node[right]{$l$};
    \filldraw[fill=white, draw=black]
        (.5,-.25) rectangle (2.5,1.25);
    \node at (1.5,.5) {\Large{$C$}};

    \draw[webs] (0,0) to[out=180,in=0] (-.5,.5);
    \draw[webs,over] (0,.5) to[out=180,in=0] (-.5,0);
    \draw[webs] (-.5,.5) to[out=180,in=0] (-1,0) node[left]{$k$};
    \draw[webs, over] (-.5,0) to[out=180,in=0] (-1,.5);
    \draw[webs] (0,1) to (-1,1) node[left]{$i$};
    \end{tikzpicture}%
    }
    {\Large{$\simeq$}}
    \scalebox{.9}{%
    \begin{tikzpicture}[anchorbase]
    \draw[webs,biarr] (0,0) node[left]{$k$} to (3,0) node[right]{$j$};
    \draw[webs,biarr] (0,.5) to (3,.5);
    \draw[webs,biarr] (0,1) node[left]{$i$} to (3,1) node[right]{$l$};
    \filldraw[fill=white, draw=black]
        (.5,-.25) rectangle (2.5,1.25);
    \node at (1.5,.5) {\Large{$C$}};
    \end{tikzpicture}%
    }
    
    \caption{The action of a trivial horizontal crossing. Here the first and third isotopies are braid relations, while the second is an application of Propositions \ref{prop: trace_local} and \ref{prop: markov}.}
    \label{fig: horizontal_crossing_trivial}
\end{figure}
\end{proof}

\begin{remark} \label{rem: trivial_horizontals}
    Lemma \ref{lem: trivial_horizontals} has a natural interpretation from the perspective of Figure \ref{fig: general_framework}; the horizontal crossings which act trivially are exactly those which can be absorbed into the shuffle braid with no impact on the functor assigned via $\Phi$.
\end{remark}

\FloatBarrier

\begin{remark}
    Suppose $\aa, \bb, \cc, \dd$ are four balanced sequences with $\mathrm{len}(\aa) = 2l$, $\mathrm{len}(\bb) = \mathrm{len}(\cc) = 2m$, and $\mathrm{len}(\dd) = 2n$. Given a balanced oriented $(\aa,\bb)$-tangle $\tau$ and a balanced oriented $(\cc,\dd)$-tangle $\tau'$, unless $\bb = \cc$, there is no canonical vertical concatenation $\tau \star \tau'$ in $\otanbal$. Nevertheless, for any complex $C \in K^b(\extsbim_l)$, we can still consider the composition $\Phi(\tau') \circ \Phi(\tau) (C)$. Lemma \ref{lem: trivial_horizontals} gives a natural tangle-theoretic interpretation of this composition by interposing an oriented $(\bb, \cc)$-tangle constructed entirely from horizontal crossings of upwards-oriented strands over downwards-oriented strands.
\end{remark}

\subsection{Invariance} \label{subsec: invariance}

This section is devoted to the proof of our main invariance theorem:

\begin{theorem} \label{thm: main_invariance}
    The assignment $\Phi$ constructed in \S \ref{subsec: actions} gives a well-defined functor from $\otanbal$ to $\funcat$.
\end{theorem}

\begin{proof}
    We need to verify that $\Phi$ takes isotopic balanced oriented tangle diagrams to naturally equivalent functors. By Theorem \ref{thm: reidemeister_tangles}, it suffices to check that $\Phi$ preserves the local moves in Figure \ref{fig: reidemeister_moves}. We check the distant interchange moves in \S \ref{subsec: distant_interchange}, zigzag moves in \S \ref{subsec: zigzag}, dial moves in \S \ref{subsec: dial_moves}, Reidemeister I moves in \S \ref{subsec: reidemeister_I}, Reidemeister II moves in \S \ref{subsec: reidemeister_II}, and Reidemeister III moves in \S \ref{subsec: reidemeister_III}.
\end{proof}

\begin{notation} \label{not: equivalence_symbols}
    In our proof of Theorem \ref{thm: main_invariance}, we will wish to distinguish among a variety of different notions of equivalence. We will use the symbol $\simeq$ to denote two complexes in $\extsbim$ which are (naturally in some complex $C$ from which they are built) homotopy equivalent. In practice, this equivalence will always arise through a combination of braid relations and applications of Propositions \ref{prop: trace_local}, \ref{prop: trace-like} and \ref{prop: markov} and Theorem \ref{thm: exchange}.
    
    For two balanced oriented tangles $\tau$ and $\tau'$, we will write $\tau \sim \tau'$ when $\Phi(\tau)$ is naturally equivalent to $\Phi(\tau')$. In light of Theorem \ref{thm: main_invariance}, any two isotopic balanced oriented tangles are related in this way, and we will typically use $\sim$ to denote isotopy after we have concluded the proof of this theorem. Nevertheless, we wish to stress that for the duration of the proof, the symbol $\sim$ indicates some algebraic statement which must be justified.
\end{notation}

\FloatBarrier

\subsubsection{Distant Interchange} \label{subsec: distant_interchange}

It suffices to check the cases in which both $\tau$ and $\tau'$ are a single vertical crossing, cup, or cap, possibly separated by some vertical strands. If both $\tau$ and $\tau'$ are vertical crossings of the same orientation, the desired equivalence follows immediately from the distant interchange of Rouquier complexes. Similarly, if $\tau$ and $\tau'$ are vertical crossings of opposite orientations, the desired equivalence follows from associativity of the tensor product $\star$.

If one of $\tau, \tau'$ is a vertical crossing and the other is a cap or a cup, the desired equivalence follows immediately from Proposition \ref{prop: trace_local}. This leaves only the cases in which $\tau$ and $\tau'$ consist of caps and cups (again, potentially separated by some identity strands). Recall that the action of clockwise and counterclockwise cups/caps on a given complex are identical, so we need only consider a single orientation. If $\tau$ and $\tau'$ are both cups, the desired equivalence follows from a sequence of braid relations; see Figure \ref{fig: cup_interchange} below. On the other hand, if $\tau$ and $\tau'$ are both caps, the desired equivalence follows from a sequence of braid relations together with Proposition \ref{prop: trace-like}; see Figure \ref{fig: cap_interchange} below.

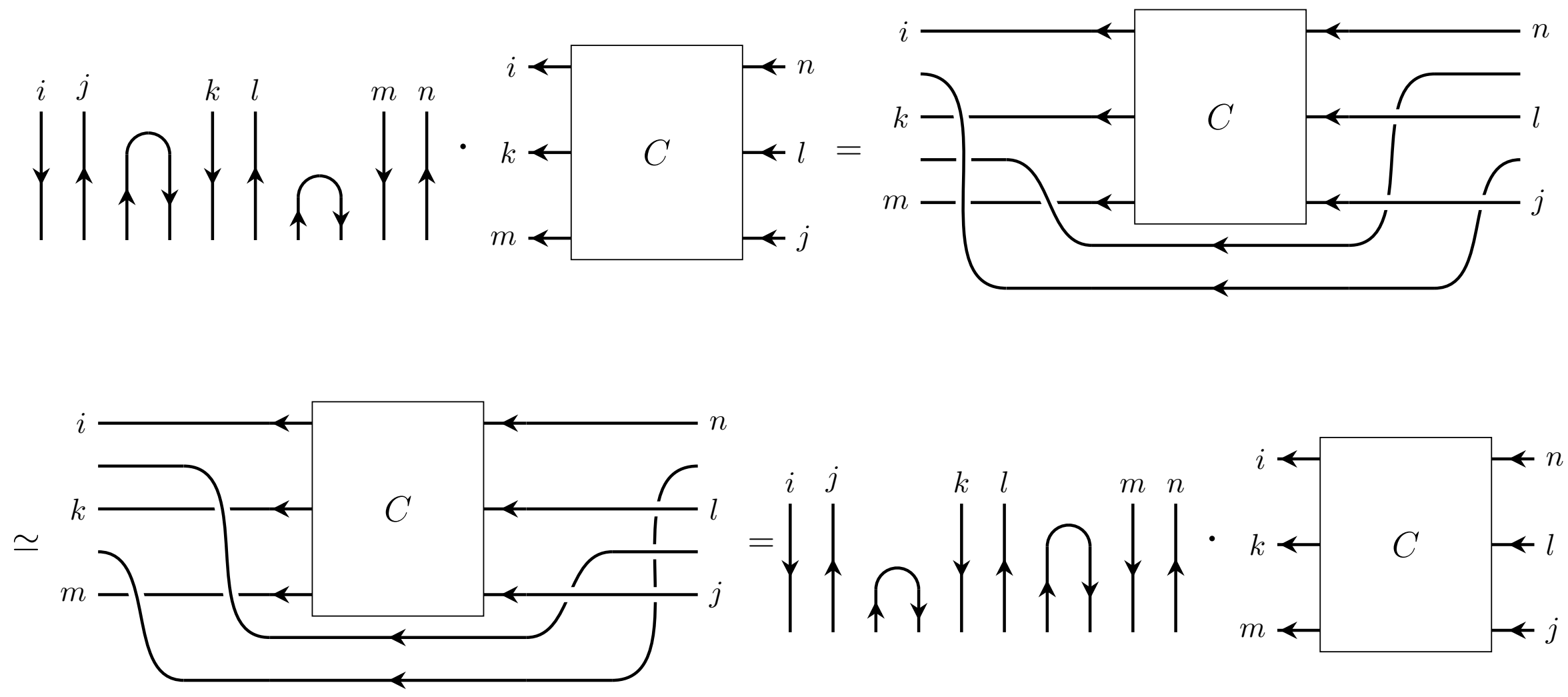
\begin{figure}[ht]
    \centering

\scalebox{.85}{%
\begin{tikzpicture}[anchorbase]
    \draw[webs,-<-=.5] (0,0) to (0,1.5) node[above]{$i$};
    \draw[webs,->-=.5] (.5,0) to (.5,1.5) node[above]{$j$};

    \draw[webs,->-=.5] (1,0) to (1,1);
    \draw[webs] (1,1) to[out=90,in=180] (1.25,1.25);
    \draw[webs] (1.25,1.25) to[out=0,in=90] (1.5,1);
    \draw[webs,->-=.5] (1.5,1) to (1.5,0);

    \draw[webs,->-=.5] (3,0) to (3,.5);
    \draw[webs] (3,.5) to[out=90,in=180] (3.25,.75);
    \draw[webs] (3.25,.75) to[out=0,in=90] (3.5,.5);
    \draw[webs,->-=.5] (3.5,.5) to (3.5,0);

    \draw[webs,-<-=.5] (2,0) to (2,1.5) node[above]{$k$};
    \draw[webs,->-=.5] (2.5,0) to (2.5,1.5) node[above]{$l$};

    \draw[webs,-<-=.5] (4,0) to (4,1.5) node[above]{$m$};
    \draw[webs,->-=.5] (4.5,0) to (4.5,1.5) node[above]{$n$};
\end{tikzpicture}
{\huge{$\cdot$}}
\begin{tikzpicture}[anchorbase]
    \draw[webs,biarr] (0,0) node[left]{$m$} to (3,0) node[right]{$j$};
    \draw[webs,biarr] (0,1) node[left]{$k$} to (3,1) node[right]{$l$};
    \draw[webs,biarr] (0,2) node[left]{$i$} to (3,2) node[right]{$n$};
    \filldraw[fill=white, draw=black]
        (.5,-.25) rectangle (2.5,2.25);
    \node at (1.5,1) {\Large{$C$}};
\end{tikzpicture}
{\Large{$=$}}
\begin{tikzpicture}[anchorbase]
    \draw[webs,biarr] (0,0) to (3,0);
    \draw[webs,biarr] (0,1) to (3,1);
    \draw[webs,biarr] (0,2) to (3,2);
    \filldraw[fill=white, draw=black]
        (.5,-.25) rectangle (2.5,2.25);
    \node at (1.5,1) {\Large{$C$}};

    \draw[webs,-<-=.5] (0,-.5) to (3,-.5);
    \draw[webs,-<-=.5] (0,-1) to (3,-1);

    \draw[webs] (-1,.5) to (-2,.5);
    \draw[webs] (0,-1) to (-1,-1);
    \draw[webs] (0,0) to (-2,0) node[left]{$m$};
    \draw[webs] (0,1) to (-2,1) node[left]{$k$};
    \draw[webs] (0,2) to (-2,2) node[left]{$i$};
    \draw[webs,over] (0,-.5) to[out=180,in=0] (-1,.5);
    \draw[webs,over] (-1,-1) to[out=180,in=0] (-2,1.5);

    \draw[webs] (3,-.5) to[out=0,in=180] (4,1.5);
    \draw[webs] (4,1.5) to (5,1.5);
    \draw[webs] (3,-1) to (4,-1);
    \draw[webs] (4,-1) to[out=0,in=180] (5,.5);
    \draw[webs,over] (3,0) to (5,0) node[right]{$j$};
    \draw[webs,over] (3,1) to (5,1) node[right]{$l$};
    \draw[webs,over] (3,2) to (5,2) node[right]{$n$};
\end{tikzpicture}
}
\\
\vspace{1cm}
\scalebox{.85}{%
{\Large{$\simeq$}}
\begin{tikzpicture}[anchorbase]
    \draw[webs,biarr] (0,0) to (3,0);
    \draw[webs,biarr] (0,1) to (3,1);
    \draw[webs,biarr] (0,2) to (3,2);
    \filldraw[fill=white, draw=black]
        (.5,-.25) rectangle (2.5,2.25);
    \node at (1.5,1) {\Large{$C$}};

    \draw[webs,-<-=.5] (0,-.5) to (3,-.5);
    \draw[webs,-<-=.5] (0,-1) to (3,-1);

    \draw[webs] (-1,1.5) to (-2,1.5);
    \draw[webs] (0,-1) to (-1,-1);
    \draw[webs] (0,0) to (-2,0) node[left]{$m$};
    \draw[webs] (0,1) to (-2,1) node[left]{$k$};
    \draw[webs] (0,2) to (-2,2) node[left]{$i$};
    \draw[webs,over] (0,-.5) to[out=180,in=0] (-1,1.5);
    \draw[webs,over] (-1,-1) to[out=180,in=0] (-2,.5);

    \draw[webs] (3,-.5) to[out=0,in=180] (4,.5);
    \draw[webs] (3,-1) to (4,-1);
    \draw[webs] (4,-1) to[out=0,in=180] (5,1.5);
    \draw[webs,over] (4,.5) to (5,.5);
    \draw[webs,over] (3,0) to (5,0) node[right]{$j$};
    \draw[webs,over] (3,1) to (5,1) node[right]{$l$};
    \draw[webs,over] (3,2) to (5,2) node[right]{$n$};
\end{tikzpicture}
{\Large{$=$}}\begin{tikzpicture}[anchorbase]
    \draw[webs,-<-=.5] (0,0) to (0,1.5) node[above]{$i$};
    \draw[webs,->-=.5] (.5,0) to (.5,1.5) node[above]{$j$};

    \draw[webs,->-=.5] (1,0) to (1,.5);
    \draw[webs] (1,.5) to[out=90,in=180] (1.25,.75);
    \draw[webs] (1.25,.75) to[out=0,in=90] (1.5,.5);
    \draw[webs,->-=.5] (1.5,.5) to (1.5,0);

    \draw[webs,->-=.5] (3,0) to (3,1);
    \draw[webs] (3,1) to[out=90,in=180] (3.25,1.25);
    \draw[webs] (3.25,1.25) to[out=0,in=90] (3.5,1);
    \draw[webs,->-=.5] (3.5,1) to (3.5,0);

    \draw[webs,-<-=.5] (2,0) to (2,1.5) node[above]{$k$};
    \draw[webs,->-=.5] (2.5,0) to (2.5,1.5) node[above]{$l$};

    \draw[webs,-<-=.5] (4,0) to (4,1.5) node[above]{$m$};
    \draw[webs,->-=.5] (4.5,0) to (4.5,1.5) node[above]{$n$};
\end{tikzpicture}
{\huge{$\cdot$}}
\begin{tikzpicture}[anchorbase]
    \draw[webs,biarr] (0,0) node[left]{$m$} to (3,0) node[right]{$j$};
    \draw[webs,biarr] (0,1) node[left]{$k$} to (3,1) node[right]{$l$};
    \draw[webs,biarr] (0,2) node[left]{$i$} to (3,2) node[right]{$n$};
    \filldraw[fill=white, draw=black]
        (.5,-.25) rectangle (2.5,2.25);
    \node at (1.5,1) {\Large{$C$}};
\end{tikzpicture}
}

\caption{Distant interchange of two cups. The equivalence is a sequence of braid relations.}
\label{fig: cup_interchange}
\end{figure}

\begin{figure}[ht]
    \centering

\scalebox{.9}{%
\begin{tikzpicture}[anchorbase]
    \draw[webs,-<-=.5] (0,0) to (0,1.5) node[above]{$i$};
    \draw[webs,->-=.5] (.5,0) to (.5,1.5) node[above]{$j$};

    \draw[webs,->-=.5] (1,1.5) to (1,.5);
    \draw[webs] (1,.5) to[out=270,in=180] (1.25,.25);
    \draw[webs] (1.25,.25) to[out=0,in=270] (1.5,.5);
    \draw[webs,->-=.5] (1.5,.5) to (1.5,1.5);

    \draw[webs,->-=.5] (3,1.5) to (3,1);
    \draw[webs] (3,1) to[out=270,in=180] (3.25,.75);
    \draw[webs] (3.25,.75) to[out=0,in=270] (3.5,1);
    \draw[webs,->-=.5] (3.5,1) to (3.5,1.5);

    \draw[webs,-<-=.5] (2,0) to (2,1.5) node[above]{$k$};
    \draw[webs,->-=.5] (2.5,0) to (2.5,1.5) node[above]{$l$};

    \draw[webs,-<-=.5] (4,0) to (4,1.5) node[above]{$m$};
    \draw[webs,->-=.5] (4.5,0) to (4.5,1.5) node[above]{$n$};
\end{tikzpicture}
{\huge{$\cdot$}}
\begin{tikzpicture}[anchorbase]
    \draw[webs,biarr] (0,0) node[left]{$m$} to (3,0) node[right]{$j$};
    \draw[webs,biarr] (0,.5) to (3,.5);
    \draw[webs,biarr] (0,1) node[left]{$k$} to (3,1) node[right]{$l$};
    \draw[webs,biarr] (0,1.5) to (3,1.5);
    \draw[webs,biarr] (0,2) node[left]{$i$} to (3,2) node[right]{$n$};
    \filldraw[fill=white, draw=black]
        (.5,-.25) rectangle (2.5,2.25);
    \node at (1.5,1) {\Large{$C$}};
\end{tikzpicture}
{\Large{$=$}}
\begin{tikzpicture}[anchorbase]
\draw[webs,biarr] (0,0) to (3,0);
    \draw[webs,biarr] (0,.5) to (3,.5);
    \draw[webs,biarr] (0,1) to (3,1);
    \draw[webs,biarr] (0,1.5) to (3,1.5);
    \draw[webs,biarr] (0,2) to (3,2);
    \filldraw[fill=white, draw=black]
        (.5,-.25) rectangle (2.5,2.25);
    \node at (1.5,1) {\Large{$C$}};
    
    \draw[webs] (0,0) to (-2,0) node[left]{$m$};
    \draw[webs] (0,1) to (-2,1) node[left]{$k$};
    \draw[webs] (0,2) to (-2,2) node[left]{$i$};
    \draw[webs,over] (0,.5) to[out=180,in=90] (-1,-.25);
    \draw[webs] (0,1.5) to (-1,1.5);
    \draw[webs,over] (-1,1.5) to[out=180,in=90] (-1.75,-.75);

    \draw[webs] (3,1.5) to[out=0,in=90] (4,-.25);
    \draw[webs,over] (3,.5) to (4,.5);
    \draw[webs] (4,.5) to[out=0,in=90] (4.5,-.75);
    \draw[webs,over] (3,2) to (5,2) node[right]{$n$};
    \draw[webs,over] (3,1) to (5,1) node[right]{$l$};
    \draw[webs,over] (3,0) to (5,0) node[right]{$j$};

    \draw[webs] (-1,-.25) to[out=270,in=180] (-.75,-.5);
    \draw[webs] (-1.75,-.75) to[out=270,in=180] (-1.5,-1);
    \draw[webs] (4,-.25) to[out=270,in=0] (3.75,-.5);
    \draw[webs] (4.5,-.75) to[out=270,in=0] (4.25,-1);
    \draw[webs] (-1.5,-1) to (4.25,-1);
    \draw[webs] (-.75,-.5) to (3.75,-.5);
    
\end{tikzpicture}%
}
\\
\vspace{1cm}
\scalebox{.9}{%
{\Large{$\simeq$}}
\begin{tikzpicture}[anchorbase]
\draw[webs,biarr] (0,0) to (3,0);
    \draw[webs,biarr] (0,.5) to (3,.5);
    \draw[webs,biarr] (0,1) to (3,1);
    \draw[webs,biarr] (0,1.5) to (3,1.5);
    \draw[webs,biarr] (0,2) to (3,2);
    \filldraw[fill=white, draw=black]
        (.5,-.25) rectangle (2.5,2.25);
    \node at (1.5,1) {\Large{$C$}};
    
    \draw[webs] (0,0) to (-2,0) node[left]{$m$};
    \draw[webs] (0,1) to (-2,1) node[left]{$k$};
    \draw[webs] (0,2) to (-2,2) node[left]{$i$};
    \draw[webs] (0,.5) to (-1,.5);
    \draw[webs,over] (0,1.5) to[out=180,in=0] (-1,-.5);
    \draw[webs,over] (-1,.5) to[out=180,in=0] (-2,-1);
    \draw[webs,over] (-1,-.5) to (-2,-.5);

    \draw[webs] (3,.5) to[out=0,in=180] (4,-.5);
    \draw[webs] (3,1.5) to (4,1.5);
    \draw[webs] (4,1.5) to[out=0,in=180] (5,-1);
    \draw[webs,over] (4,-.5) to (5,-.5);
    \draw[webs,over] (3,2) to (5,2) node[right]{$n$};
    \draw[webs,over] (3,1) to (5,1) node[right]{$l$};
    \draw[webs,over] (3,0) to (5,0) node[right]{$j$};

    \draw[webs] (-2,-1) to[out=180,in=90] (-2.25,-1.25);
    \draw[webs] (-2.25,-1.25) to[out=270,in=180] (-2,-1.5);
    \draw[webs] (5,-1) to[out=0,in=90] (5.25,-1.25);
    \draw[webs] (5.25,-1.25) to[out=270,in=0] (5,-1.5);
    \draw[webs] (-2,-1.5) to (5,-1.5);

    \draw[webs] (-2,-.5) to[out=180,in=90] (-2.75,-1.25);
    \draw[webs] (-2.75,-1.25) to[out=270,in=180] (-2,-2);
    \draw[webs] (5,-.5) to[out=0,in=90] (5.75,-1.25);
    \draw[webs] (5.75,-1.25) to[out=270,in=0] (5,-2);
    \draw[webs] (-2,-2) to (5,-2);
\end{tikzpicture}%
{\Large{$\simeq$}}
\begin{tikzpicture}[anchorbase]
\draw[webs,biarr] (0,0) to (3,0);
    \draw[webs,biarr] (0,.5) to (3,.5);
    \draw[webs,biarr] (0,1) to (3,1);
    \draw[webs,biarr] (0,1.5) to (3,1.5);
    \draw[webs,biarr] (0,2) to (3,2);
    \filldraw[fill=white, draw=black]
        (.5,-.25) rectangle (2.5,2.25);
    \node at (1.5,1) {\Large{$C$}};

    \draw[webs] (0,0) to (-2,0) node[left]{$m$};
    \draw[webs] (0,1) to (-2,1) node[left]{$k$};
    \draw[webs] (0,2) to (-2,2) node[left]{$i$};
    \draw[webs] (0,.5) to (-1,.5);
    \draw[webs,over] (0,1.5) to[out=180,in=90] (-1,-.25);
    \draw[webs,over] (-1,.5) to[out=180,in=90] (-2,-.75);

    \draw[webs] (3,.5) to[out=0,in=90] (4,-.25);
    \draw[webs] (3,1.5) to (4,1.5);
    \draw[webs] (4,1.5) to[out=0,in=90] (4.75,-.75);
    \draw[webs,over] (3,2) to (5,2) node[right]{$n$};
    \draw[webs,over] (3,1) to (5,1) node[right]{$l$};
    \draw[webs,over] (3,0) to (5,0) node[right]{$j$};

    \draw[webs] (-1,-.25) to[out=270,in=180] (-.75,-.5);
    \draw[webs] (4,-.25) to[out=270,in=0] (3.75,-.5);
    \draw[webs] (-.75,-.5) to (3.75,-.5);

    \draw[webs] (-2,-.75) to[out=270,in=180] (-1.75,-1);
    \draw[webs] (4.75,-.75) to[out=270,in=0] (4.5,-1);
    \draw[webs] (-1.75,-1) to (4.5,-1);
\end{tikzpicture}%
}
\\
\vspace{1cm}
\scalebox{.9}{%
{\Large{$=$}}
\begin{tikzpicture}[anchorbase]
    \draw[webs,-<-=.5] (0,0) to (0,1.5) node[above]{$i$};
    \draw[webs,->-=.5] (.5,0) to (.5,1.5) node[above]{$j$};

    \draw[webs,->-=.5] (1,1.5) to (1,1);
    \draw[webs] (1,1) to[out=270,in=180] (1.25,.75);
    \draw[webs] (1.25,.75) to[out=0,in=270] (1.5,1);
    \draw[webs,->-=.5] (1.5,1) to (1.5,1.5);

    \draw[webs,->-=.5] (3,1.5) to (3,.5);
    \draw[webs] (3,.5) to[out=270,in=180] (3.25,.25);
    \draw[webs] (3.25,.25) to[out=0,in=270] (3.5,.5);
    \draw[webs,->-=.5] (3.5,.5) to (3.5,1.5);

    \draw[webs,-<-=.5] (2,0) to (2,1.5) node[above]{$k$};
    \draw[webs,->-=.5] (2.5,0) to (2.5,1.5) node[above]{$l$};

    \draw[webs,-<-=.5] (4,0) to (4,1.5) node[above]{$m$};
    \draw[webs,->-=.5] (4.5,0) to (4.5,1.5) node[above]{$n$};
\end{tikzpicture}
{\huge{$\cdot$}}
\begin{tikzpicture}[anchorbase]
    \draw[webs,biarr] (0,0) node[left]{$m$} to (3,0) node[right]{$j$};
    \draw[webs,biarr] (0,.5) to (3,.5);
    \draw[webs,biarr] (0,1) node[left]{$k$} to (3,1) node[right]{$l$};
    \draw[webs,biarr] (0,1.5) to (3,1.5);
    \draw[webs,biarr] (0,2) node[left]{$i$} to (3,2) node[right]{$n$};
    \filldraw[fill=white, draw=black]
        (.5,-.25) rectangle (2.5,2.25);
    \node at (1.5,1) {\Large{$C$}};
\end{tikzpicture}%
}

\caption{Distant interchange of two caps. The first equivalence is a sequence of braid relations, and the second is an application of Proposition \ref{prop: trace-like}.}
\label{fig: cap_interchange}
\end{figure}
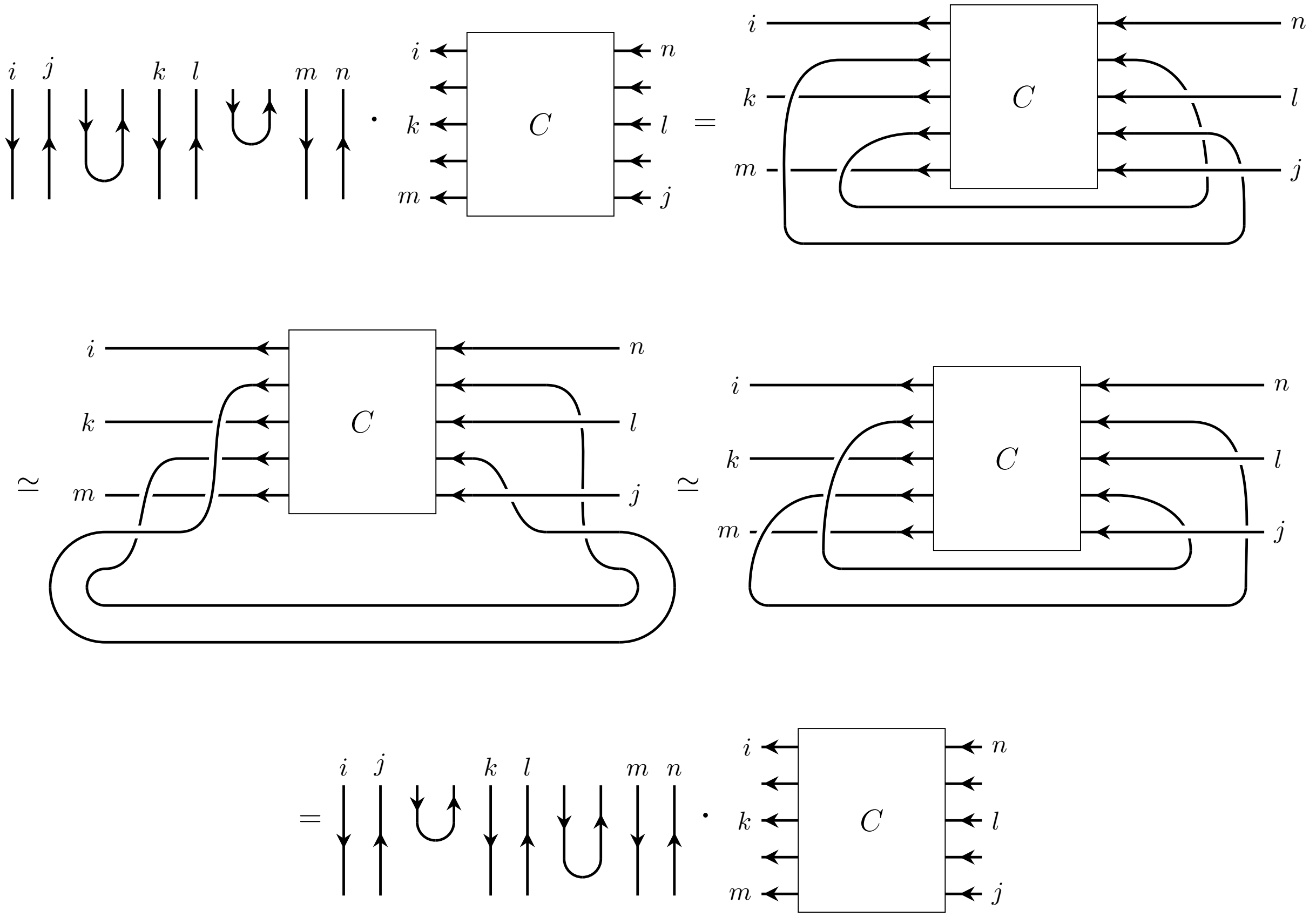

\FloatBarrier

It remains to consider distant interchange of caps past cups. Here we again compute directly; see Figure \ref{fig: caps_past_cups} below. We illustrate only the case in which $\tau$ is a cup and $\tau'$ is a cap; the other case is exactly analogous.

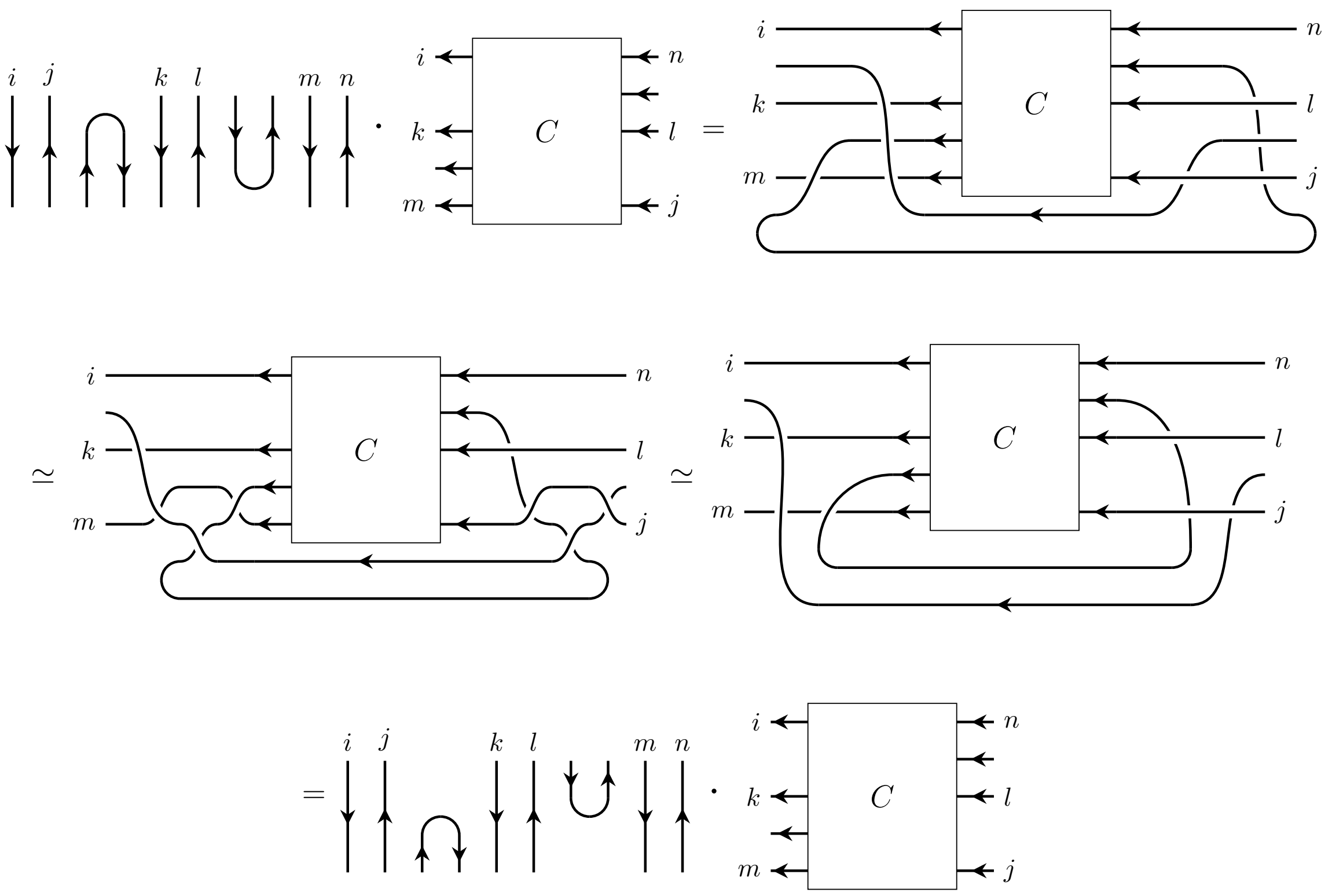
\begin{figure}[ht]
    \centering

    \scalebox{.85}{
    \begin{tikzpicture}[anchorbase]
    \draw[webs,-<-=.5] (0,0) to (0,1.5) node[above]{$i$};
    \draw[webs,->-=.5] (.5,0) to (.5,1.5) node[above]{$j$};

    \draw[webs,->-=.5] (1,0) to (1,1);
    \draw[webs] (1,1) to[out=90,in=180] (1.25,1.25);
    \draw[webs] (1.25,1.25) to[out=0,in=90] (1.5,1);
    \draw[webs,->-=.5] (1.5,1) to (1.5,0);

    \draw[webs,->-=.5] (3,1.5) to (3,.5);
    \draw[webs] (3,.5) to[out=270,in=180] (3.25,.25);
    \draw[webs] (3.25,.25) to[out=0,in=270] (3.5,.5);
    \draw[webs,->-=.5] (3.5,.5) to (3.5,1.5);

    \draw[webs,-<-=.5] (2,0) to (2,1.5) node[above]{$k$};
    \draw[webs,->-=.5] (2.5,0) to (2.5,1.5) node[above]{$l$};

    \draw[webs,-<-=.5] (4,0) to (4,1.5) node[above]{$m$};
    \draw[webs,->-=.5] (4.5,0) to (4.5,1.5) node[above]{$n$};
    \end{tikzpicture}
    {\huge{$\cdot$}}
    \begin{tikzpicture}[anchorbase]
        \draw[webs,biarr] (0,0) node[left]{$m$} to (3,0) node[right]{$j$};
        \draw[webs,-<-=.1] (0,.5) to (2,.5);
        \draw[webs,biarr] (0,1) node[left]{$k$} to (3,1) node[right]{$l$};
        \draw[webs,-<-=.9] (1,1.5) to (3,1.5);
        \draw[webs,biarr] (0,2) node[left]{$i$} to (3,2) node[right]{$n$};

        \filldraw[fill=white, draw=black]
        (.5,-.25) rectangle (2.5,2.25);
        \node at (1.5,1) {\Large{$C$}};
    \end{tikzpicture}
    {\Large{$=$}}
    \begin{tikzpicture}[anchorbase]
        \draw[webs,biarr] (0,0) to (3,0);
        \draw[webs,-<-=.1] (0,.5) to (2,.5);
        \draw[webs,biarr] (0,1) to (3,1);
        \draw[webs,-<-=.9] (1,1.5) to (3,1.5);
        \draw[webs,biarr] (0,2) to (3,2);

        \filldraw[fill=white, draw=black]
        (.5,-.25) rectangle (2.5,2.25);
        \node at (1.5,1) {\Large{$C$}};

        \draw[webs,-<-=.5] (0,-.5) to (3,-.5);

        \draw[webs] (0,0) to (-2,0) node[left]{$m$};
        \draw[webs] (0,1) to (-2,1) node[left]{$k$};
        \draw[webs] (0,2) to (-2,2) node[left]{$i$};
        \draw[webs] (0,.5) to (-1,.5);
        \draw[webs,over] (0,-.5) to[out=180,in=0] (-1,1.5);
        \draw[webs,over] (-1,.5) to[out=180,in=0] (-2,-.5);
        \draw[webs] (-1,1.5) to (-2,1.5);

        \draw[webs] (3,-.5) to[out=0,in=180] (4,.5);
        \draw[webs] (3,1.5) to (4,1.5);
        \draw[webs] (4,1.5) to[out=0,in=180] (5,-.5);
        \draw[webs,over] (4,.5) to (5,.5);
        \draw[webs,over] (3,0) to (5,0) node[right]{$j$};
        \draw[webs,over] (3,1) to (5,1) node[right]{$l$};
        \draw[webs,over] (3,2) to (5,2) node[right]{$n$};

        \draw[webs] (-2,-.5) to[out=180,in=90] (-2.25,-.75);
        \draw[webs] (-2.25,-.75) to[out=270,in=180] (-2,-1);
        \draw[webs] (5,-.5) to[out=0,in=90] (5.25,-.75);
        \draw[webs] (5.25,-.75) to[out=270,in=0] (5,-1);
        \draw[webs] (-2,-1) to (5,-1);
    \end{tikzpicture}
    }
    \\
    \vspace{1cm}
    \scalebox{.85}{
    {\Large{$\simeq$}}
    \begin{tikzpicture}[anchorbase]
        \draw[webs,biarr] (0,0) to (3,0);
        \draw[webs,-<-=.1] (0,.5) to (2,.5);
        \draw[webs,biarr] (0,1) to (3,1);
        \draw[webs,-<-=.9] (1,1.5) to (3,1.5);
        \draw[webs,biarr] (0,2) to (3,2);

        \filldraw[fill=white, draw=black]
        (.5,-.25) rectangle (2.5,2.25);
        \node at (1.5,1) {\Large{$C$}};

        \draw[webs,-<-=.5] (0,-.5) to (3,-.5);
        
        \draw[webs] (0,1) to (-2,1) node[left]{$k$};
        \draw[webs] (0,2) to (-2,2) node[left]{$i$};
        \draw[webs] (0,0) to[out=180,in=0] (-.5,.5);
        \draw[webs,over] (0,.5) to[out=180,in=0] (-.5,0);
        \draw[webs] (0,-.5) to (-.5,-.5);
        \draw[webs] (-.5,0) to[out=180,in=0] (-1,-.5);
        \draw[webs,over] (-.5,-.5) to[out=180,in=0] (-1,0);
        \draw[webs] (-.5,.5) to (-1,.5);
        \draw[webs] (-1,.5) to[out=180,in=0] (-1.5,0);
        \draw[webs,over] (-1,0) to[out=180,in=0] (-2,1.5);
        \draw[webs] (-1.5,0) to (-2,0) node[left]{$m$};

        \draw[webs] (3,1.5) to[out=0,in=180] (4,0);
        \draw[webs] (3,0) to (3.5,0);
        \draw[webs,over] (3.5,0) to[out=0,in=180] (4,.5);
        \draw[webs] (3,-.5) to (4,-.5);
        \draw[webs] (4,0) to[out=0,in=180] (4.5,-.5);
        \draw[webs,over] (4,-.5) to[out=0,in=180] (4.5,0);
        \draw[webs,over] (3,1) to (5,1) node[right]{$l$};
        \draw[webs,over] (3,2) to (5,2) node[right]{$n$};
        \draw[webs] (4,.5) to (4.5,.5);
        \draw[webs] (4.5,0) to[out=0,in=180] (5,.5);
        \draw[webs,over] (4.5,.5) to[out=0,in=180] (5,0) node[right]{$j$};

        \draw[webs] (-1,-.5) to[out=180,in=90] (-1.25,-.75);
        \draw[webs] (-1.25,-.75) to[out=270,in=180] (-1,-1);
        \draw[webs] (4.5,-.5) to[out=0,in=90] (4.75,-.75);
        \draw[webs] (4.75,-.75) to[out=270,in=0] (4.5,-1);
        \draw[webs] (-1,-1) to (4.5,-1);
    \end{tikzpicture}
    {\Large{$\simeq$}}
    \begin{tikzpicture}[anchorbase]
        \draw[webs,biarr] (0,0) to (3,0);
        \draw[webs,-<-=.1] (0,.5) to (2,.5);
        \draw[webs,biarr] (0,1) to (3,1);
        \draw[webs,-<-=.9] (1,1.5) to (3,1.5);
        \draw[webs,biarr] (0,2) to (3,2);

        \filldraw[fill=white, draw=black]
        (.5,-.25) rectangle (2.5,2.25);
        \node at (1.5,1) {\Large{$C$}};
        
        \draw[webs] (0,1) to (-2,1) node[left]{$k$};
        \draw[webs] (0,2) to (-2,2) node[left]{$i$};
        \draw[webs] (0,0) to (-2,0) node[left]{$m$};
        \draw[webs,over] (0,.5) to[out=180,in=90] (-1,-.5);

        \draw[webs] (3,1.5) to[out=0,in=90] (4,-.5);
        \draw[webs] (4,-1.25) to[out=0,in=180] (5,.5);
        \draw[webs,over] (3,1) to (5,1) node[right]{$l$};
        \draw[webs,over] (3,2) to (5,2) node[right]{$n$};
        \draw[webs,over] (3,0) to (5,0) node[right]{$j$};

        \draw[webs,-<-=.5] (-1,-1.25) to (4,-1.25);
        \draw[webs,over] (-1,-1.25) to[out=180,in=0] (-2,1.5);

        \draw[webs] (-1,-.5) to[out=270,in=180] (-.75,-.75);
        \draw[webs] (4,-.5) to[out=270,in=0] (3.75,-.75);
        \draw[webs] (-.75,-.75) to (3.75,-.75);
    \end{tikzpicture}
    }
    \\
    \vspace{1cm}
    \scalebox{.85}{
    {\Large{$=$}}
    \begin{tikzpicture}[anchorbase]
        \draw[webs,-<-=.5] (0,0) to (0,1.5) node[above]{$i$};
        \draw[webs,->-=.5] (.5,0) to (.5,1.5) node[above]{$j$};

        \draw[webs,->-=.5] (1,0) to (1,.5);
        \draw[webs] (1,.5) to[out=90,in=180] (1.25,.75);
        \draw[webs] (1.25,.75) to[out=0,in=90] (1.5,.5);
        \draw[webs,->-=.5] (1.5,.5) to (1.5,0);

        \draw[webs,->-=.5] (3,1.5) to (3,1);
        \draw[webs] (3,1) to[out=270,in=180] (3.25,.75);
        \draw[webs] (3.25,.75) to[out=0,in=270] (3.5,1);
        \draw[webs,->-=.5] (3.5,1) to (3.5,1.5);

        \draw[webs,-<-=.5] (2,0) to (2,1.5) node[above]{$k$};
        \draw[webs,->-=.5] (2.5,0) to (2.5,1.5) node[above]{$l$};

        \draw[webs,-<-=.5] (4,0) to (4,1.5) node[above]{$m$};
        \draw[webs,->-=.5] (4.5,0) to (4.5,1.5) node[above]{$n$};
    \end{tikzpicture}
    {\huge{$\cdot$}}
    \begin{tikzpicture}[anchorbase]
        \draw[webs,biarr] (0,0) node[left]{$m$} to (3,0) node[right]{$j$};
        \draw[webs,-<-=.1] (0,.5) to (2,.5);
        \draw[webs,biarr] (0,1) node[left]{$k$} to (3,1) node[right]{$l$};
        \draw[webs,-<-=.9] (1,1.5) to (3,1.5);
        \draw[webs,biarr] (0,2) node[left]{$i$} to (3,2) node[right]{$n$};
        
        \filldraw[fill=white, draw=black]
        (.5,-.25) rectangle (2.5,2.25);
        \node at (1.5,1) {\Large{$C$}};
    \end{tikzpicture}
    }

    \caption{Distant interchange of caps past cups. The first equivalence is a sequence of braid relations, and the second is an application of Theorem \ref{thm: exchange}.}
    \label{fig: caps_past_cups}
\end{figure}

\begin{remark} \label{rem: caps_cups_exchange}
    We pause here to point out that the distant interchange of caps past cups is the only explicit use of our new ingredient, Theorem \ref{thm: exchange}, in our invariance proof. There are implicit uses in further cases; Lemma \ref{lem: bending_around} requires an interchange of a cap past a cup, and this lemma is later applied to cut down on the number of required cases in checking Reidemeister II and III moves. This use is only cosmetic for Reidemeister III: while it saves some work, each case can be checked directly without ever appealing to Theorem \ref{thm: exchange}.
    
    On the other hand, checking the ``wrong-way" Reidemeister II move which is not covered trivially by Lemma \ref{lem: trivial_horizontals} requires the use of Theorem \ref{thm: exchange} in an essential way. This move has been the fundamental obstruction to previous attempts at extending HOMFLY homology to an oriented tangle invariant. Here we bypass this obstruction at the cost of having to check locality properties explicitly. From either perspective, it is exactly Theorem \ref{thm: exchange} which comes to the rescue, and it does so exactly once.
\end{remark}

\FloatBarrier

\subsubsection{Zigzag} \label{subsec: zigzag}

Here there are 4 distinct cases given by a choice of orientation and a kink on the left or right. We proceed by computing the action of these tangles on a complex directly; see Figure \ref{fig: zigzag_invariance} below. Only one case is illustrated here; the others are completely analogous.

\begin{figure}[ht]
    \centering

    \scalebox{.9}{%
    \begin{tikzpicture}[anchorbase]
    \draw[webs,-<-=.5] (0,0) to (0,1.5) node[above]{$i$};
    \draw[webs,->-=.5] (.5,0) to (.5,1.5) node[above]{$j$};
    \draw[webs,->-=.5] (1,0) to (1,.75);
    \draw[webs] (1,.75) to[out=90,in=180] (1.25,1);
    \draw[webs,->-=.99] (1.25,1) to[out=0,in=90] (1.5,.75);
    \draw[webs] (1.5,.75) to[out=270,in=180] (1.75,.5);
    \draw[webs] (1.75,.5) to[out=0,in=270] (2,.75);
    \draw[webs,->-=.5] (2,.75) to (2,1.5);
    \draw[webs,-<-=.5] (2.5,0) to (2.5,1.5) node[above]{$k$};
    \draw[webs,->-=.5] (3,0) to (3,1.5) node[above]{$l$};
    \end{tikzpicture}%
    }
    {\huge{$\cdot$}}
    \scalebox{.9}{%
    \begin{tikzpicture}[anchorbase]
    \draw[webs,biarr] (0,0) node[left]{$k$} to (3,0) node[right]{$j$};
    \draw[webs,->-=.1] (3,.5) to (1,.5);
    \draw[webs,biarr] (0,1) node[left]{$i$} to (3,1) node[right]{$l$};
    \filldraw[fill=white, draw=black]
        (.5,-.25) rectangle (2.5,1.25);
    \node at (1.5,.5) {\Large{$C$}};
    \end{tikzpicture}%
    }
    {\Large{$=$}}
    \scalebox{.9}{%
    \begin{tikzpicture}[anchorbase]
    
    \draw[webs,biarr] (0,0) to (3,0);
    \draw[webs,->-=.1] (3,.5) to (1,.5);
    \draw[webs,biarr] (0,1) to (3,1);
    \filldraw[fill=white, draw=black]
        (.5,-.25) rectangle (2.5,1.25);
    \node at (1.5,.5) {\Large{$C$}};

    \draw[webs,-<-=.5] (0,-.5) to (3,-.5);

    \draw[webs] (0,0) to[out=180,in=0] (-.5,-.5);
    \draw[webs,over] (0,-.5) to[out=180,in=0] (-.5,0);
    \draw[webs] (-.5,-.5) to[out=180,in=0] (-1,0) node[left]{$k$};
    \draw[webs,over] (-.5,0) to[out=180,in=0] (-1,-.5);
    \draw[webs] (0,1) to (-1,1) node[left]{$i$};

    \draw[webs] (3,-.5) to[out=0,in=180] (3.5,0);
    \draw[webs,over] (3,0) to[out=0,in=180] (3.5,-.5);
    \draw[webs] (3,.5) to (3.5,.5);
    \draw[webs] (3.5,.5) to[out=0,in=180] (4.5,-.5);
    \draw[webs,over] (3.5,0) to[out=0,in=180] (4,.5);
    \draw[webs] (4,.5) to (4.5,.5);
    \draw[webs] (3.5,-.5) to (4,-.5);
    \draw[webs,over] (4,-.5) to[out=0,in=180] (4.5,0) node[right]{$j$};
    \draw[webs] (3,1) to (4.5,1) node[right]{$l$};

    \draw[webs] (-1,-.5) to[out=180,in=90] (-1.25,-.75);
    \draw[webs] (-1.25,-.75) to[out=270,in=180] (-1,-1);
    \draw[webs] (4.5,-.5) to[out=0,in=90] (4.75,-.75);
    \draw[webs] (4.75,-.75) to[out=270,in=0] (4.5,-1);
    \draw[webs] (-1,-1) to (4.5,-1);
    \end{tikzpicture}%
    }
    \\
    \vspace{1cm}
    {\Large{$\simeq$}}
    \scalebox{.9}{%
    \begin{tikzpicture}[anchorbase]

    \draw[webs,biarr] (0,0) to (3,0);
    \draw[webs,->-=.1] (3,.5) to (1,.5);
    \draw[webs,biarr] (0,1) to (3,1);
    \filldraw[fill=white, draw=black]
        (.5,-.25) rectangle (2.5,1.25);
    \node at (1.5,.5) {\Large{$C$}};

    \draw[webs,-<-=.5] (0,-.5) to (3,-.5);
    \draw[webs] (0,-.5) to (-1,-.5);
    \draw[webs] (0,0) to (-1,0) node[left]{$k$};
    \draw[webs] (0,1) to (-1,1) node[left]{$i$};

    \draw[webs] (3,.5) to[out=0,in=180] (3.5,0);
    \draw[webs,over] (3,0) to[out=0,in=180] (3.5,.5);
    \draw[webs] (3.5,.5) to (4,.5);
    \draw[webs] (3.5,0) to[out=0,in=180] (4,-.5);
    \draw[webs] (3,-.5) to (3.5,-.5);
    \draw[webs,over] (3.5,-.5) to[out=0,in=180] (4,0);
    \draw[webs] (4,0) to[out=0,in=180] (4.5,.5);
    \draw[webs,over] (4,.5) to[out=0,in=180] (4.5,0) node[right]{$j$};
    \draw[webs] (4,-.5) to (4.5,-.5);
    \draw[webs] (3,1) to (4.5,1) node[right]{$l$};

    \draw[webs] (-1,-.5) to[out=180,in=90] (-1.25,-.75);
    \draw[webs] (-1.25,-.75) to[out=270,in=180] (-1,-1);
    \draw[webs] (4.5,-.5) to[out=0,in=90] (4.75,-.75);
    \draw[webs] (4.75,-.75) to[out=270,in=0] (4.5,-1);
    \draw[webs] (-1,-1) to (4.5,-1);
    
    \end{tikzpicture}%
    }
    {\Large{$\simeq$}}
    \scalebox{.9}{%
    \begin{tikzpicture}[anchorbase]
    
    \draw[webs,biarr] (0,0) node[left]{$k$} to (3,0);
    \draw[webs,->-=.1] (3,.5) to (1,.5);
    \draw[webs,biarr] (0,1) node[left]{$i$} to (3,1);
    \filldraw[fill=white, draw=black]
        (.5,-.25) rectangle (2.5,1.25);
    \node at (1.5,.5) {\Large{$C$}};

    \draw[webs] (3,.5) to[out=0,in=180] (3.5,0);
    \draw[webs,over] (3,0) to[out=0,in=180] (3.5,.5);
    \draw[webs] (3.5,0) to[out=0,in=180] (4,.5);
    \draw[webs,over] (3.5,.5) to[out=0,in=180] (4,0) node[right]{$j$};
    \draw[webs] (3,1) to (4,1) node[right]{$l$};
    
    \end{tikzpicture}%
    }
    {\Large{$\simeq$}}
    \scalebox{.9}{%
    \begin{tikzpicture}[anchorbase]
    \draw[webs,biarr] (0,0) node[left]{$k$} to (3,0) node[right]{$j$};
    \draw[webs,->-=.1] (3,.5) to (1,.5);
    \draw[webs,biarr] (0,1) node[left]{$i$} to (3,1) node[right]{$l$};
    \filldraw[fill=white, draw=black]
        (.5,-.25) rectangle (2.5,1.25);
    \node at (1.5,.5) {\Large{$C$}};
    \end{tikzpicture}%
    }
    \caption{Invariance under a zigzag move. The first and third equivalences are braid relations, while the second is an application of Propositions \ref{prop: trace_local} and \ref{prop: markov}.}
    \label{fig: zigzag_invariance}
\end{figure}
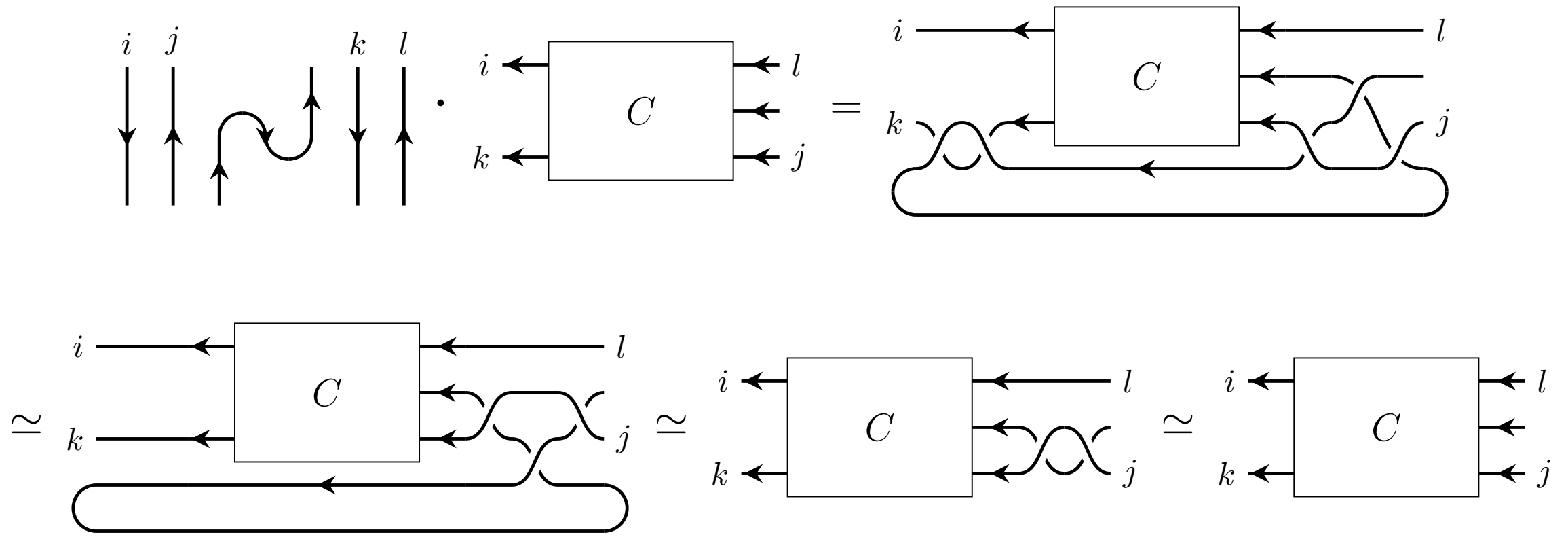

\FloatBarrier

\subsubsection{Dial Moves} \label{subsec: dial_moves}

There are eight\footnote{In \cite{Saw95}, only the four cases in which the vertical strand passes over the cup are listed. We were unable to derive the other four cases, in which the vertical strand passes under the cup, from these and the remaining relations. We instead opt to check these cases explicitly, as our preference for caution outweighs our fear of redundancy.} cases, given by a choice of orientation on each strand and the sign of the crossing. In each case, one side contains a vertical crossing while the other contains a horizontal crossing.

In each of the four cases in which the vertical strand passes over the cup, the desired relation follows immediately from resolving the horizontal crossing as in Figure \ref{fig: horizontal_crossing_decompositions} and applying zigzag and distant interchange moves (indeed, this motivated our particular choice of decomposition for each crossing). We illustrate with a single example in Figure \ref{fig: trivial_dial_move} below; the other three cases follow analogously.

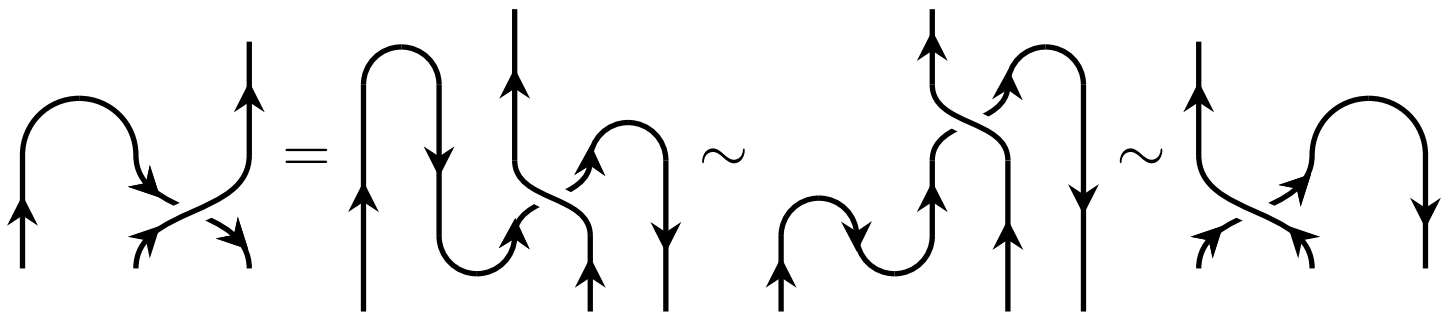
\begin{figure}[ht]
\centering
\begin{tikzpicture}[anchorbase]
\draw[webs,->-=.5] (0,0) to (0,.75);
\draw[webs,->-=.8,->-=.2] (.75,.75) to[out=270,in=90] (1.5,0);
\draw[webs,->-=.2,over] (.75,0) to[out=90,in=270] (1.5,.75);
\draw[webs] (0,.75) to[out=90,in=180] (.375,1.125);
\draw[webs] (.375,1.125) to[out=0,in=90] (.75,.75);
\draw[webs,->-=.5] (1.5,.75) to (1.5,1.5);
\end{tikzpicture}
{\Large{$=$}}
\begin{tikzpicture}[anchorbase]
\draw[webs,->-=.5] (0,0) to (0,1.5);
\draw[webs] (0,1.5) to[out=90,in=180] (.25,1.75);
\draw[webs] (.25,1.75) to[out=0,in=90] (.5,1.5);
\draw[webs,->-=.5] (.5,1.5) to (.5,.5);
\draw[webs] (.5,.5) to[out=270,in=180] (.75,.25);
\draw[webs] (.75,.25) to[out=0,in=270] (1,.5);
\draw[webs,->-=.01,->-=.99] (1,.5) to[out=90,in=270] (1.5,1);
\draw[webs] (1.5,1) to[out=90,in=180] (1.75,1.25);
\draw[webs] (1.75,1.25) to[out=0,in=90] (2,1);
\draw[webs,->-=.5] (2,1) to (2,0);
\draw[webs,->-=.5] (1.5,0) to (1.5,.5);
\draw[webs,over] (1.5,.5) to[out=90,in=270] (1,1);
\draw[webs,->-=.5] (1,1) to (1,2);
\end{tikzpicture}
{\Large{$\sim$}}
\begin{tikzpicture}[anchorbase]
\draw[webs,->-=.5] (0,0) to (0,.5);
\draw[webs] (0,.5) to[out=90,in=180] (.25,.75);
\draw[webs,->-=.99] (.25,.75) to[out=0,in=90] (.5,.5);
\draw[webs] (.5,.5) to[out=270,in=180] (.75,.25);
\draw[webs] (.75,.25) to[out=0,in=270] (1,.5);
\draw[webs,->-=.5] (1,.5) to (1,1);
\draw[webs,->-=.99] (1,1) to[out=90,in=270] (1.5,1.5);
\draw[webs] (1.5,1.5) to[out=90,in=180] (1.75,1.75);
\draw[webs] (1.75,1.75) to[out=0,in=90] (2,1.5);
\draw[webs,->-=.5] (2,1.5) to (2,0);
\draw[webs,->-=.5] (1.5,0) to (1.5,1);
\draw[webs,over] (1.5,1) to[out=90,in=270] (1,1.5);
\draw[webs,->-=.5] (1,1.5) to (1,2);
\end{tikzpicture}
{\Large{$\sim$}}
\begin{tikzpicture}[anchorbase]
\draw[webs,->-=.2,->-=.8] (0,0) to[out=90,in=270] (.75,.75);
\draw[webs,over,->-=.2] (.75,0) to[out=90,in=270] (0,.75);
\draw[webs,->-=.5] (0,.75) to (0,1.5);
\draw[webs] (.75,.75) to[out=90,in=180] (1.125,1.125);
\draw[webs] (1.125,1.125) to[out=0,in=90] (1.5,.75);
\draw[webs,->-=.5] (1.5,.75) to (1.5,0);
\end{tikzpicture}

\caption{A dial move in which the vertical strand passes over the cup. The first equality is from Figure \ref{fig: horizontal_crossing_decompositions}, while the two isotopies follow from a distant interchange and zigzag move, respectively.}
\label{fig: trivial_dial_move}
\end{figure}

The remaining cases must be checked directly. In two such cases, the horizontal crossing acts trivially by Lemma \ref{lem: trivial_horizontals}; we illustrate one such case in Figure \ref{fig: nontrivial_dial_move_trivial_crossing} below.

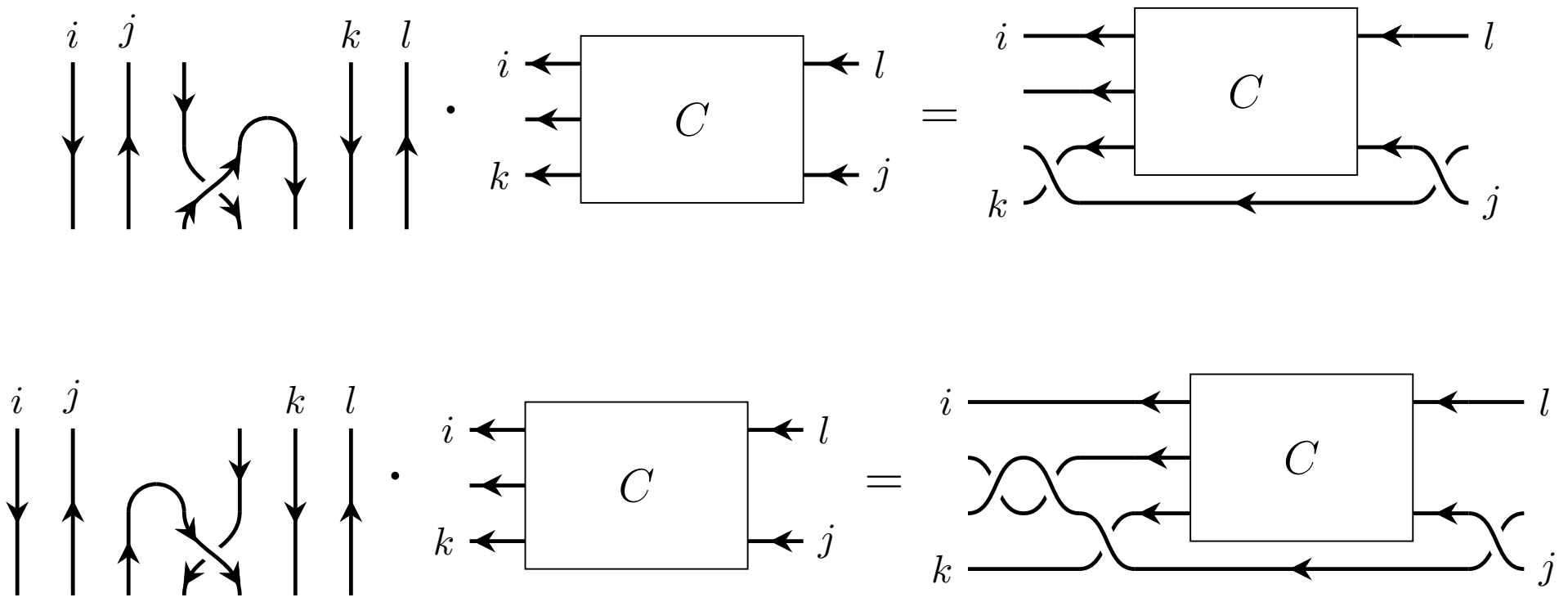
\begin{figure}[h!t]
\centering

\scalebox{.9}{%
\begin{tikzpicture}[anchorbase]
\draw[webs,-<-=.5] (0,0) to (0,1.5) node[above]{$i$};
\draw[webs,->-=.5] (.5,0) to (.5,1.5) node[above]{$j$};
\draw[webs,-<-=.2] (1.5,0) to[out=90,in=270] (1,.75);
\draw[webs,->-=.2,->-=.8,over] (1,0) to[out=90,in=270] (1.5,.75);
\draw[webs] (1.5,.75) to[out=90,in=180] (1.75,1);
\draw[webs] (1.75,1) to[out=0,in=90] (2,.75);
\draw[webs,->-=.5] (2,.75) to (2,0);
\draw[webs,-<-=.5] (1,.75) to (1,1.5);
\draw[webs,-<-=.5] (2.5,0) to (2.5,1.5) node[above]{$k$};
\draw[webs,->-=.5] (3,0) to (3,1.5) node[above]{$l$};
\end{tikzpicture}%
}
{\huge{$\cdot$}}
\scalebox{.9}{%
\begin{tikzpicture}[anchorbase]
\draw[webs,biarr] (0,0) node[left]{$k$} to (3,0) node[right]{$j$};
    \draw[webs,-<-=.1] (0,.5) to (2,.5);
    \draw[webs,biarr] (0,1) node[left]{$i$} to (3,1) node[right]{$l$};
    \filldraw[fill=white, draw=black]
        (.5,-.25) rectangle (2.5,1.25);
    \node at (1.5,.5) {\Large{$C$}};
\end{tikzpicture}%
}
{\Large{$=$}}
\scalebox{.9}{%
\begin{tikzpicture}[anchorbase]
\draw[webs,biarr] (0,0) to (3,0);
    \draw[webs,-<-=.1] (0,.5) to (2,.5);
    \draw[webs,biarr] (0,1) to (3,1);
    \filldraw[fill=white, draw=black]
        (.5,-.25) rectangle (2.5,1.25);
    \node at (1.5,.5) {\Large{$C$}};

\draw[webs,-<-=.5] (0,-.5) to (3,-.5);
\draw[webs] (0,0) to[out=180,in=0] (-.5,-.5) node[left]{$k$};
\draw[webs,over] (0,-.5) to[out=180,in=0] (-.5,0);
\draw[webs] (0,.5) to (-.5,.5);
\draw[webs] (0,1) to (-.5,1) node[left]{$i$};

\draw[webs] (3,-.5) to[out=0,in=180] (3.5,0);
\draw[webs,over] (3,0) to[out=0,in=180] (3.5,-.5) node[right]{$j$};
\draw[webs] (3,1) to (3.5,1) node[right]{$l$};
\end{tikzpicture}%
}
\\
\vspace{1cm}
\scalebox{.9}{%
\begin{tikzpicture}[anchorbase]
\draw[webs,-<-=.5] (0,0) to (0,1.5) node[above]{$i$};
\draw[webs,->-=.5] (.5,0) to (.5,1.5) node[above]{$j$};

\draw[webs,->-=.5] (1,0) to (1,.75);
\draw[webs] (1,.75) to[out=90,in=180] (1.25,1);
\draw[webs] (1.25,1) to[out=0,in=90] (1.5,.75);
\draw[webs,-<-=.2] (1.5,0) to[out=90,in=270] (2,.75);
\draw[webs, over,-<-=.2,-<-=.8] (2,0) to[out=90,in=270] (1.5,.75);
\draw[webs,-<-=.5] (2,.75) to (2,1.5);

\draw[webs,-<-=.5] (2.5,0) to (2.5,1.5) node[above]{$k$};
\draw[webs,->-=.5] (3,0) to (3,1.5) node[above]{$l$};
\end{tikzpicture}%
}
{\huge{$\cdot$}}
\scalebox{.9}{%
\begin{tikzpicture}[anchorbase]
\draw[webs,biarr] (0,0) node[left]{$k$} to (3,0) node[right]{$j$};
    \draw[webs,-<-=.1] (0,.5) to (2,.5);
    \draw[webs,biarr] (0,1) node[left]{$i$} to (3,1) node[right]{$l$};
    \filldraw[fill=white, draw=black]
        (.5,-.25) rectangle (2.5,1.25);
    \node at (1.5,.5) {\Large{$C$}};
\end{tikzpicture}%
}
{\Large{$=$}}
\scalebox{.9}{%
\begin{tikzpicture}[anchorbase]
\draw[webs,biarr] (0,0) to (3,0);
    \draw[webs,-<-=.1] (0,.5) to (2,.5);
    \draw[webs,biarr] (0,1) to (3,1);
    \filldraw[fill=white, draw=black]
        (.5,-.25) rectangle (2.5,1.25);
    \node at (1.5,.5) {\Large{$C$}};

\draw[webs,-<-=.5] (0,-.5) to (3,-.5);
\draw[webs] (0,0) to[out=180,in=0] (-.5,-.5);
\draw[webs,over] (0,-.5) to[out=180,in=0] (-.5,0);
\draw[webs] (0,.5) to (-.5,.5);
\draw[webs] (-.5,.5) to[out=180,in=0] (-1,0);
\draw[webs,over] (-.5,0) to[out=180,in=0] (-1,.5);
\draw[webs] (-1,0) to[out=180,in=0] (-1.5,.5);
\draw[webs,over] (-1,.5) to[out=180,in=0] (-1.5,0);
\draw[webs] (-.5,-.5) to (-1.5,-.5) node[left]{$k$};
\draw[webs] (0,1) to (-1.5,1) node[left]{$i$};

\draw[webs] (3,-.5) to[out=0,in=180] (3.5,0);
\draw[webs,over] (3,0) to[out=0,in=180] (3.5,-.5) node[right]{$j$};
\draw[webs] (3,1) to (3.5,1) node[right]{$l$};
\end{tikzpicture}%
}

\caption{A dial move in which the vertical strand passes under the cup and the only horizontal crossing acts trivially. Note that the resulting complexes are naturally homotopy equivalent via a sequence of braid relations.}
\label{fig: nontrivial_dial_move_trivial_crossing}
\end{figure} The other two cases require a bit more work, but we can again check directly using the decompositions of Figure \ref{fig: horizontal_crossing_decompositions}. We illustrate one such case in Figure \ref{fig: nontrivial_dial_move_nontrivial_crossing} below, and the other is analogous.


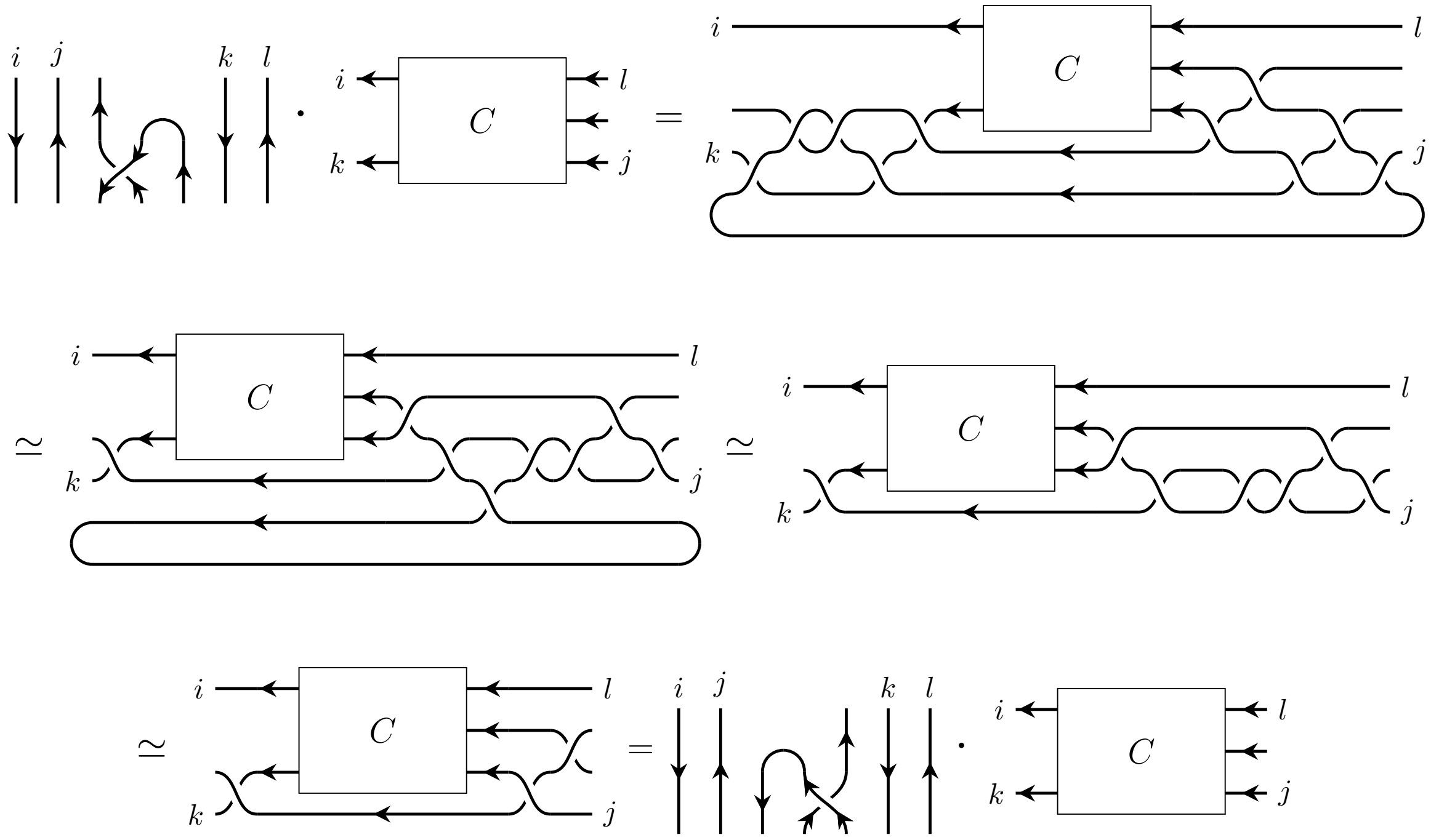
\begin{figure}[ht]
\centering

\scalebox{.9}{%
\begin{tikzpicture}[anchorbase]
\draw[webs,-<-=.5] (0,0) to (0,1.5) node[above]{$i$};
\draw[webs,->-=.5] (.5,0) to (.5,1.5) node[above]{$j$};

\draw[webs,->-=.2] (1.5,0) to[out=90,in=270] (1,.75);
\draw[webs,-<-=.2,-<-=.8,over] (1,0) to[out=90,in=270] (1.5,.75);
\draw[webs] (1.5,.75) to[out=90,in=180] (1.75,1);
\draw[webs] (1.75,1) to[out=0,in=90] (2,.75);
\draw[webs,-<-=.5] (2,.75) to (2,0);
\draw[webs,->-=.5] (1,.75) to (1,1.5);

\draw[webs,-<-=.5] (2.5,0) to (2.5,1.5) node[above]{$k$};
\draw[webs,->-=.5] (3,0) to (3,1.5) node[above]{$l$};
\end{tikzpicture}%
}
{\huge{$\cdot$}}
\scalebox{.9}{%
\begin{tikzpicture}[anchorbase]
\draw[webs,biarr] (0,0) node[left]{$k$} to (3,0) node[right]{$j$};
    \draw[webs,->-=.1] (3,.5) to (1,.5);
    \draw[webs,biarr] (0,1) node[left]{$i$} to (3,1) node[right]{$l$};
    \filldraw[fill=white, draw=black]
        (.5,-.25) rectangle (2.5,1.25);
    \node at (1.5,.5) {\Large{$C$}};
\end{tikzpicture}%
}
{\Large{$=$}}
\scalebox{.9}{%
\begin{tikzpicture}[anchorbase]
\draw[webs,biarr] (0,0) to (3,0);
    \draw[webs,->-=.1] (3,.5) to (1,.5);
    \draw[webs,biarr] (0,1) to (3,1);
    \filldraw[fill=white, draw=black]
        (.5,-.25) rectangle (2.5,1.25);
    \node at (1.5,.5) {\Large{$C$}};

\draw[webs,-<-=.5] (0,-.5) to (3,-.5);
\draw[webs,-<-=.5] (0,-1) to (3,-1);

\draw[webs] (0,0) to[out=180,in=0] (-.5,-.5);
\draw[webs,over] (0,-.5) to[out=180,in=0] (-.5,0);
\draw[webs] (-.5,0) to (-1,0);
\draw[webs] (0,-1) to (-.5,-1);
\draw[webs] (-.5,-.5) to[out=180,in=0] (-1,-1);
\draw[webs,over] (-.5,-1) to[out=180,in=0] (-1,-.5);
\draw[webs] (-1,-.5) to[out=180,in=0] (-1.5,0);
\draw[webs,over] (-1,0) to[out=180,in=0] (-1.5,-.5);
\draw[webs] (-1.5,-.5) to[out=180,in=0] (-2,0);
\draw[webs,over] (-1.5,0) to[out=180,in=0] (-2,-.5);
\draw[webs] (-1,-1) to (-2,-1);
\draw[webs] (-2,-1) to[out=180,in=0] (-2.5,-.5) node[left]{$k$};
\draw[webs,over] (-2,-.5) to[out=180,in=0] (-2.5,-1);
\draw[webs] (0,1) to (-2.5,1) node[left]{$i$};
\draw[webs] (-2,0) to (-2.5,0);

\draw[webs] (3,-.5) to[out=0,in=180] (3.5,0);
\draw[webs,over] (3,0) to[out=0,in=180] (3.5,-.5);
\draw[webs] (3,.5) to (3.5,.5);
\draw[webs] (3.5,0) to[out=0,in=180] (4,.5);
\draw[webs,over] (3.5,.5) to[out=0,in=180] (4,0);
\draw[webs] (3,-1) to (4,-1);
\draw[webs] (3.5,-.5) to (4,-.5);
\draw[webs] (4,-1) to[out=0,in=180] (4.5,-.5);
\draw[webs,over] (4,-.5) to[out=0,in=180] (4.5,-1);
\draw[webs] (4,0) to (4.5,0);
\draw[webs] (4.5,-.5) to[out=0,in=180] (5,0);
\draw[webs,over] (4.5,0) to[out=0,in=180] (5,-.5);
\draw[webs] (4.5,-1) to (5,-1);
\draw[webs] (5,-.5) to[out=0,in=180] (5.5,-1);
\draw[webs,over] (5,-1) to[out=0,in=180] (5.5,-.5) node[right]{$j$};
\draw[webs] (5,0) to (5.5,0);
\draw[webs] (4,.5) to (5.5,.5);
\draw[webs] (3,1) to (5.5,1) node[right]{$l$};

\draw[webs] (-2.5,-1) to[out=180,in=90] (-2.75,-1.25);
\draw[webs] (-2.75,-1.25) to[out=270,in=180] (-2.5,-1.5);
\draw[webs] (5.5,-1) to[out=0,in=90] (5.75,-1.25);
\draw[webs] (5.75,-1.25) to[out=270,in=0] (5.5,-1.5);
\draw[webs] (-2.5,-1.5) to (5.5,-1.5);
\end{tikzpicture}%
}
\\
\vspace{1cm}
{\Large{$\simeq$}}
\scalebox{.9}{%
\begin{tikzpicture}[anchorbase]
\draw[webs,biarr] (0,0) to (3,0);
    \draw[webs,->-=.1] (3,.5) to (1,.5);
    \draw[webs,biarr] (0,1) to (3,1);
    \filldraw[fill=white, draw=black]
        (.5,-.25) rectangle (2.5,1.25);
    \node at (1.5,.5) {\Large{$C$}};

\draw[webs,-<-=.5] (0,-.5) to (3,-.5);
\draw[webs,-<-=.5] (0,-1) to (3,-1);

\draw[webs] (0,0) to[out=180,in=0] (-.5,-.5) node[left]{$k$};
\draw[webs,over] (0,-.5) to[out=180,in=0] (-.5,0);
\draw[webs] (0,1) to (-.5,1) node[left]{$i$};
\draw[webs] (0,-1) to (-.5,-1);

\draw[webs] (3,.5) to[out=0,in=180] (3.5,0);
\draw[webs,over] (3,0) to[out=0,in=180] (3.5,.5);
\draw[webs] (3,-.5) to (3.5,-.5);
\draw[webs] (3.5,-.5) to[out=0,in=180] (4,0);
\draw[webs,over] (3.5,0) to[out=0,in=180] (4,-.5);
\draw[webs] (3,-1) to (4,-1);
\draw[webs] (4,-1) to[out=0,in=180] (4.5,-.5);
\draw[webs,over] (4,-.5) to[out=0,in=180] (4.5,-1);
\draw[webs] (4,0) to (4.5,0);
\draw[webs] (4.5,0) to[out=0,in=180] (5,-.5);
\draw[webs,over] (4.5,-.5) to[out=0,in=180] (5,0);
\draw[webs] (5,0) to[out=0,in=180] (5.5,-.5);
\draw[webs,over] (5,-.5) to[out=0,in=180] (5.5,0);
\draw[webs] (3.5,.5) to (5.5,.5);
\draw[webs] (5.5,0) to[out=0,in=180] (6,.5);
\draw[webs,over] (5.5,.5) to[out=0,in=180] (6,0);
\draw[webs] (5.5,-.5) to (6,-.5);
\draw[webs] (6,-.5) to[out=0,in=180] (6.5,0);
\draw[webs,over] (6,0) to[out=0,in=180] (6.5,-.5) node[right]{$j$};
\draw[webs] (3,1) to (6.5,1) node[right]{$l$};
\draw[webs] (6,.5) to (6.5,.5);
\draw[webs] (4.5,-1) to (6.5,-1);

\draw[webs] (-.5,-1) to[out=180,in=90] (-.75,-1.25);
\draw[webs] (-.75,-1.25) to[out=270,in=180] (-.5,-1.5);
\draw[webs] (6.5,-1) to[out=0,in=90] (6.75,-1.25);
\draw[webs] (6.75,-1.25) to[out=270,in=0] (6.5,-1.5);
\draw[webs] (-.5,-1.5) to (6.5,-1.5);
\end{tikzpicture}%
}
{\Large{$\simeq$}}
\scalebox{.9}{%
\begin{tikzpicture}[anchorbase]
\draw[webs,biarr] (0,0) to (3,0);
    \draw[webs,->-=.1] (3,.5) to (1,.5);
    \draw[webs,biarr] (0,1) to (3,1);
    \filldraw[fill=white, draw=black]
        (.5,-.25) rectangle (2.5,1.25);
    \node at (1.5,.5) {\Large{$C$}};

\draw[webs,-<-=.5] (0,-.5) to (3,-.5);

\draw[webs] (0,0) to[out=180,in=0] (-.5,-.5) node[left]{$k$};
\draw[webs,over] (0,-.5) to[out=180,in=0] (-.5,0);
\draw[webs] (0,1) to (-.5,1) node[left]{$i$};

\draw[webs] (3,.5) to[out=0,in=180] (3.5,0);
\draw[webs,over] (3,0) to[out=0,in=180] (3.5,.5);
\draw[webs] (3,-.5) to (3.5,-.5);
\draw[webs] (3.5,-.5) to[out=0,in=180] (4,0);
\draw[webs,over] (3.5,0) to[out=0,in=180] (4,-.5);

\draw[webs] (4,-.5) to (4.5,-.5);

\draw[webs] (4,0) to (4.5,0);
\draw[webs] (4.5,0) to[out=0,in=180] (5,-.5);
\draw[webs,over] (4.5,-.5) to[out=0,in=180] (5,0);
\draw[webs] (5,0) to[out=0,in=180] (5.5,-.5);
\draw[webs,over] (5,-.5) to[out=0,in=180] (5.5,0);
\draw[webs] (3.5,.5) to (5.5,.5);
\draw[webs] (5.5,0) to[out=0,in=180] (6,.5);
\draw[webs,over] (5.5,.5) to[out=0,in=180] (6,0);
\draw[webs] (5.5,-.5) to (6,-.5);
\draw[webs] (6,-.5) to[out=0,in=180] (6.5,0);
\draw[webs,over] (6,0) to[out=0,in=180] (6.5,-.5) node[right]{$j$};
\draw[webs] (3,1) to (6.5,1) node[right]{$l$};
\draw[webs] (6,.5) to (6.5,.5);
\end{tikzpicture}%
}
\\
\vspace{1cm}
{\Large{$\simeq$}}
\scalebox{.9}{%
\begin{tikzpicture}[anchorbase]
\draw[webs,biarr] (0,0) to (3,0);
    \draw[webs,->-=.1] (3,.5) to (1,.5);
    \draw[webs,biarr] (0,1) to (3,1);
    \filldraw[fill=white, draw=black]
        (.5,-.25) rectangle (2.5,1.25);
    \node at (1.5,.5) {\Large{$C$}};

\draw[webs,-<-=.5] (0,-.5) to (3,-.5);
\draw[webs] (0,0) to[out=180,in=0] (-.5,-.5) node[left]{$k$};
\draw[webs,over] (0,-.5) to[out=180,in=0] (-.5,0);
\draw[webs] (0,1) to (-.5,1) node[left]{$i$};

\draw[webs] (3,-.5) to[out=0,in=180] (3.5,0);
\draw[webs,over] (3,0) to[out=0,in=180] (3.5,-.5);
\draw[webs] (3,.5) to (3.5,.5);
\draw[webs] (3.5,.5) to[out=0,in=180] (4,0);
\draw[webs,over] (3.5,0) to[out=0,in=180] (4,.5);
\draw[webs] (3.5,-.5) to (4,-.5) node[right]{$j$};
\draw[webs] (3,1) to (4,1) node[right]{$l$};
\end{tikzpicture}%
{\Large{$=$}}
\begin{tikzpicture}[anchorbase]
\draw[webs,-<-=.5] (0,0) to (0,1.5) node[above]{$i$};
\draw[webs,->-=.5] (.5,0) to (.5,1.5) node[above]{$j$};

\draw[webs,-<-=.5] (1,0) to (1,.75);
\draw[webs] (1,.75) to[out=90,in=180] (1.25,1);
\draw[webs] (1.25,1) to[out=0,in=90] (1.5,.75);
\draw[webs,->-=.2] (1.5,0) to[out=90,in=270] (2,.75);
\draw[webs, over,->-=.2,->-=.8] (2,0) to[out=90,in=270] (1.5,.75);
\draw[webs,->-=.5] (2,.75) to (2,1.5);

\draw[webs,-<-=.5] (2.5,0) to (2.5,1.5) node[above]{$k$};
\draw[webs,->-=.5] (3,0) to (3,1.5) node[above]{$l$};
\end{tikzpicture}
{\huge{$\cdot$}}
\begin{tikzpicture}[anchorbase]
\draw[webs,biarr] (0,0) node[left]{$k$} to (3,0) node[right]{$j$};
    \draw[webs,->-=.1] (3,.5) to (1,.5);
    \draw[webs,biarr] (0,1) node[left]{$i$} to (3,1) node[right]{$l$};
    \filldraw[fill=white, draw=black]
        (.5,-.25) rectangle (2.5,1.25);
    \node at (1.5,.5) {\Large{$C$}};
\end{tikzpicture}%
}

\caption[blah blah blah]{A dial move in which the vertical strand passes under the cup and the only horizontal crossing acts nontrivially. The first equivalence is a sequence of braid relations, the second follows from Propositions \ref{prop: trace_local} and \ref{prop: markov}, and the third is a sequence of braid relations.}
\label{fig: nontrivial_dial_move_nontrivial_crossing}
\end{figure}

\FloatBarrier

In checking the remaining relations, we will often make use of the following consequence of the distant interchange, zigzag, and dial moves:

\begin{lemma} \label{lem: bending_around}
    Let $\beta$ be any oriented braid, and denote by $\beta'$ the oriented braid obtained from $\beta$ by rotation through an angle of $\pi$. Then we have the following natural equivalence of functors:

    
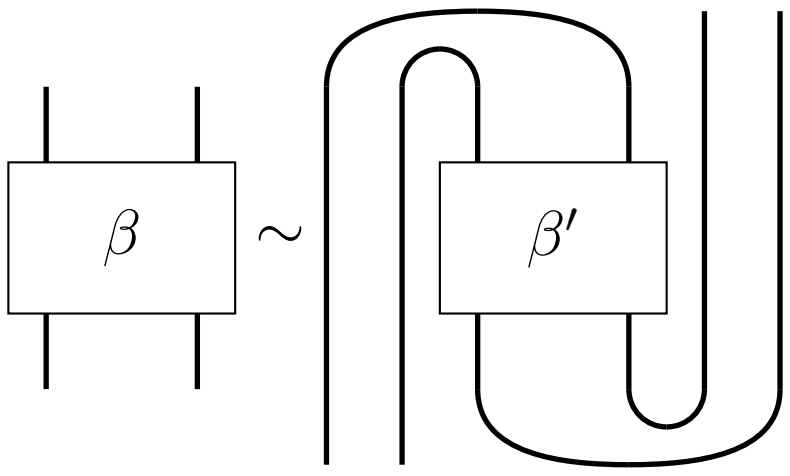
\begin{figure}[ht]
    \centering
        \begin{tikzpicture}[anchorbase]
            \draw[webs] (0,0) to (0,2);
            \draw[webs] (1,0) to (1,2);

            \filldraw[fill=white, draw=black]
            (-.25,.5) rectangle (1.25,1.5);
            \node at (.5,1) {\Large{$\beta$}};
        \end{tikzpicture}
        {\Large{$\sim$}}
        \begin{tikzpicture}[anchorbase]
            \draw[webs] (0,0) to (0,2);
            \draw[webs] (1,0) to (1,2);

            \filldraw[fill=white, draw=black]
            (-.25,.5) rectangle (1.25,1.5);
            \node at (.5,1) {\Large{$\beta'$}};

            \draw[webs] (-.5,2) to[out=90,in=180] (-.25,2.25);
            \draw[webs] (-.25,2.25) to[out=0,in=90] (0,2);
            \draw[webs] (-1,2) to[out=90,in=180] (0,2.5);
            \draw[webs] (0,2.5) to[out=0,in=90] (1,2);

            \draw[webs] (1,0) to[out=270,in=180] (1.25,-.25);
            \draw[webs] (1.25,-.25) to[out=0,in=270] (1.5,0);
            \draw[webs] (0,0) to[out=270,in=180] (1,-.5);
            \draw[webs] (1,-.5) to[out=0,in=270] (2,0);

            \draw[webs] (-1,-.5) to (-1,2);
            \draw[webs] (-.5,-.5) to (-.5,2);
            \draw[webs] (1.5,0) to (1.5,2.5);
            \draw[webs] (2,0) to (2,2.5);
        \end{tikzpicture}
    \end{figure}
\end{lemma}

\begin{proof} It suffices to consider a single crossing. This is a straightforward consequence of already established relations; see Figure \ref{fig: sliding_crossings} below. The same argument applies for any orientation of the two strands and either parity of the crossing.

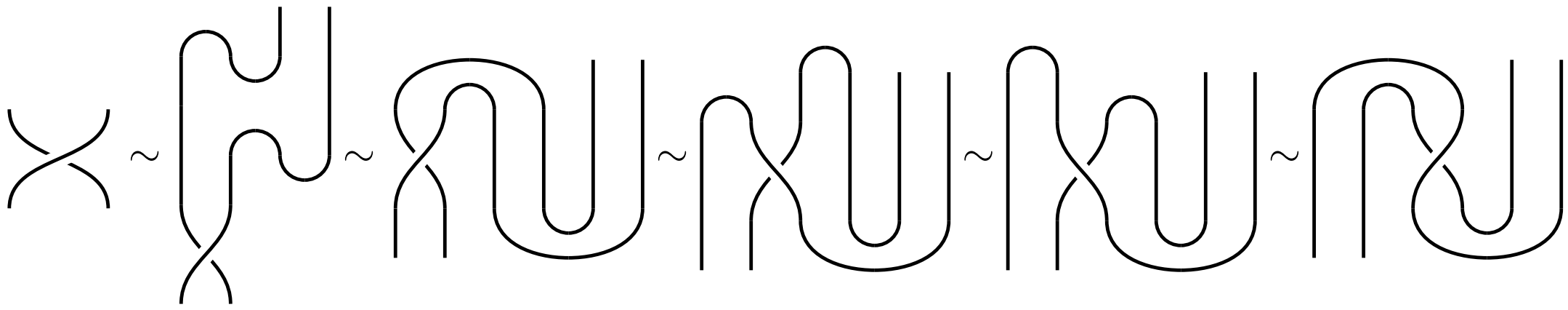
\begin{figure}[ht]
    \centering
    \scalebox{.9}{
    \begin{tikzpicture}[anchorbase]
        \draw[webs] (1,0) to[out=90,in=270] (0,1);
        \draw[webs,over] (0,0) to[out=90,in=270] (1,1);
    \end{tikzpicture}
    {\Large{$\sim$}}
    \begin{tikzpicture}[anchorbase]
        \draw[webs] (.5,0) to[out=90,in=270] (0,1);
        \draw[webs,over] (0,0) to[out=90,in=270] (.5,1);

        \begin{scope}[xshift=-.5cm]
        \draw[webs] (1,1) to (1,1.5);
        \draw[webs] (1,1.5) to[out=90,in=180] (1.25,1.75);
        \draw[webs] (1.25,1.75) to[out=0,in=90] (1.5,1.5);
        \draw[webs] (1.5,1.5) to[out=270,in=180] (1.75,1.25);
        \draw[webs] (1.75,1.25) to[out=0,in=270] (2,1.5);
        \draw[webs] (2,1.5) to (2,3);
        \end{scope}

        \draw[webs] (0,1) to (0,2);
        \begin{scope}[xshift=-1cm,yshift=1cm]
        \draw[webs] (1,1) to (1,1.5);
        \draw[webs] (1,1.5) to[out=90,in=180] (1.25,1.75);
        \draw[webs] (1.25,1.75) to[out=0,in=90] (1.5,1.5);
        \draw[webs] (1.5,1.5) to[out=270,in=180] (1.75,1.25);
        \draw[webs] (1.75,1.25) to[out=0,in=270] (2,1.5);
        \end{scope}
        \draw[webs] (1,2.5) to (1,3);
    \end{tikzpicture}
    {\Large{$\sim$}}
    \begin{tikzpicture}[anchorbase]
        \draw[webs] (.5,-.5) to (.5,0);
        \draw[webs] (1,-.5) to (1,0);
        \draw[webs] (1,0) to[out=90,in=270] (.5,1);
        \draw[webs,over] (.5,0) to[out=90,in=270] (1,1);

        \draw[webs] (1,1) to[out=90,in=180] (1.25,1.25);
        \draw[webs] (1.25,1.25) to[out=0,in=90] (1.5,1);
        \draw[webs] (.5,1) to[out=90,in=180] (1.25,1.5);
        \draw[webs] (1.25,1.5) to[out=0,in=90] (2,1);

        \draw[webs] (1.5,1) to (1.5,0);
        \draw[webs] (2,1) to (2,0);

        \begin{scope}[yscale=-1,xshift=1cm,yshift=-1cm]
        \draw[webs] (1,1) to[out=90,in=180] (1.25,1.25);
        \draw[webs] (1.25,1.25) to[out=0,in=90] (1.5,1);
        \draw[webs] (.5,1) to[out=90,in=180] (1.25,1.5);
        \draw[webs] (1.25,1.5) to[out=0,in=90] (2,1);
        \end{scope}

        \draw[webs] (2.5,0) to (2.5,1.5);
        \draw[webs] (3,0) to (3,1.5);
    \end{tikzpicture}
    {\Large{$\sim$}}
    \begin{tikzpicture}[anchorbase]
        \draw[webs] (0,-.5) to (0,1);
        \draw[webs] (.5,-.5) to (.5,0);

        \draw[webs] (.5,0) to[out=90,in=270] (1,1);
        \draw[webs,over] (1,0) to[out=90,in=270] (.5,1);

        \draw[webs] (1,1) to (1,1.5);

        \draw[webs] (1,1.5) to[out=90,in=180] (1.25,1.75);
        \draw[webs] (1.25,1.75) to[out=0,in=90] (1.5,1.5);
        \draw[webs] (0,1) to[out=90,in=180] (.25,1.25);
        \draw[webs] (.25,1.25) to[out=0,in=90] (.5,1);

        \draw[webs] (1.5,1.5) to (1.5,0);

        \draw[webs] (1.5,0) to[out=270,in=180] (1.75,-.25);
        \draw[webs] (1.75,-.25) to[out=0,in=270] (2,0);
        \draw[webs] (1,0) to[out=270,in=180] (1.75,-.5);
        \draw[webs] (1.75,-.5) to[out=0,in=270] (2.5,0);

        \draw[webs] (2,0) to (2,1.5);
        \draw[webs] (2.5,0) to (2.5,1.5);
    \end{tikzpicture}
    {\Large{$\sim$}}
    \begin{tikzpicture}[anchorbase]
        \draw[webs] (0,-.5) to (0,1.5);
        \draw[webs] (.5,-.5) to (.5,0);

        \draw[webs] (.5,0) to[out=90,in=270] (1,1);
        \draw[webs,over] (1,0) to[out=90,in=270] (.5,1);

        \draw[webs] (.5,1) to (.5,1.5);

        \draw[webs] (1,1) to[out=90,in=180] (1.25,1.25);
        \draw[webs] (1.25,1.25) to[out=0,in=90] (1.5,1);
        \draw[webs] (0,1.5) to[out=90,in=180] (.25,1.75);
        \draw[webs] (.25,1.75) to[out=0,in=90] (.5,1.5);

        \draw[webs] (1.5,1) to (1.5,0);

        \draw[webs] (1.5,0) to[out=270,in=180] (1.75,-.25);
        \draw[webs] (1.75,-.25) to[out=0,in=270] (2,0);
        \draw[webs] (1,0) to[out=270,in=180] (1.75,-.5);
        \draw[webs] (1.75,-.5) to[out=0,in=270] (2.5,0);

        \draw[webs] (2,0) to (2,1.5);
        \draw[webs] (2.5,0) to (2.5,1.5);
    \end{tikzpicture}
    {\Large{$\sim$}}
    \begin{tikzpicture}[anchorbase]
    \draw[webs] (0,-.5) to (0,1);
    \draw[webs] (.5,-.5) to (.5,1);

    \draw[webs] (0,1) to[out=90,in=180] (.75,1.5);
    \draw[webs] (.75,1.5) to[out=0,in=90] (1.5,1);
    \draw[webs] (.5,1) to[out=90,in=180] (.75,1.25);
    \draw[webs] (.75,1.25) to[out=0,in=90] (1,1);

    \draw[webs] (1,1) to[out=270,in=90] (1.5,0);
    \draw[webs,over] (1.5,1) to[out=270,in=90] (1,0);

    \draw[webs] (1.5,0) to[out=270,in=180] (1.75,-.25);
    \draw[webs] (1.75,-.25) to[out=0,in=270] (2,0);
    \draw[webs] (1,0) to[out=270,in=180] (1.75,-.5);
    \draw[webs] (1.75,-.5) to[out=0,in=270] (2.5,0);

    \draw[webs] (2,0) to (2,1.5);
    \draw[webs] (2.5,0) to (2.5,1.5);
    \end{tikzpicture}
    }

    \caption{Sliding a crossing around a zigzag. Here the first isotopy is a zigzag move, the second and fourth are distant interchanges, and the third and fifth are dial moves.}
    \label{fig: sliding_crossings}
\end{figure}

\FloatBarrier

\end{proof}

\subsubsection{Reidemeister I} \label{subsec: reidemeister_I}

There are again four distinct cases given by a choice of orientation and a loop on the left or right. We again proceed by direct computation. As before, we illustrate only one case; the rest are completely analogous.

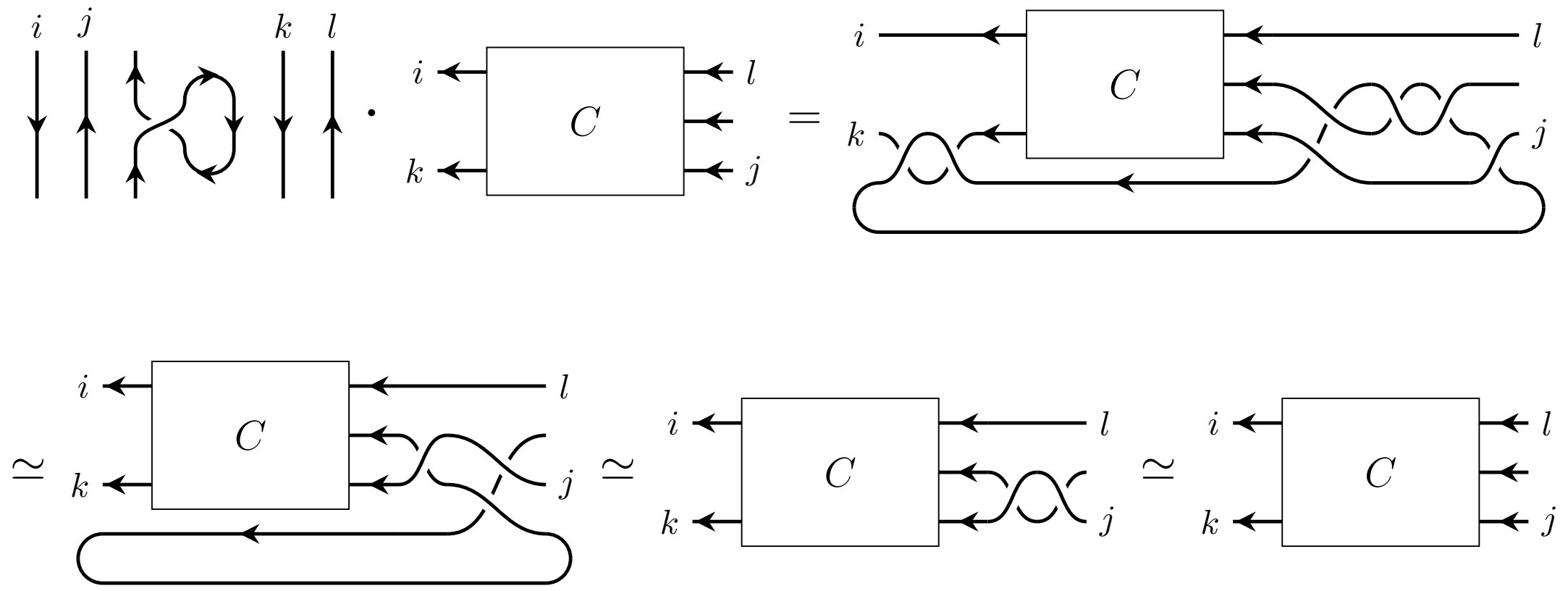
\begin{figure}[ht]
    \centering

    \scalebox{.9}{%
    \begin{tikzpicture}[anchorbase]
    \draw[webs,-<-=.5] (0,0) to (0,1.5) node[above]{$i$};
    \draw[webs,->-=.5] (0.5,0) to (0.5,1.5) node[above]{$j$};
    \draw[webs,-<-=.5] (2.5,0) node[below]{\textcolor{white}{$k$}} to (2.5,1.5) node[above]{$k$};
    \draw[webs,->-=.5] (3,0) to (3,1.5) node[above]{$l$};

    \draw[webs,->-=.5] (1,0) to (1,.5);
    \draw[webs,-<-=.95] (1.5,.5) to[out=270,in=180] (1.75,.25);
    \draw[webs] (1.75,.25) to[out=0,in=270] (2,.5);

    \draw[webs] (1.5,.5) to[out=90,in=270] (1,1);
    \draw[webs,over] (1,.5) to[out=90,in=270] (1.5,1);

    \draw[webs,->-=.5] (1,1) to (1,1.5);
    \draw[webs,->-=.95] (1.5,1) to[out=90,in=180] (1.75,1.25);
    \draw[webs] (1.75,1.25) to[out=0,in=90] (2,1);
    \draw[webs,->-=.5] (2,1) to (2,.5);

    \end{tikzpicture}%
    }
    {\huge{$\cdot$}}
    \scalebox{.9}{%
    \begin{tikzpicture}[anchorbase]
    \draw[webs,biarr] (0,0) node[left]{$k$} to (3,0) node[right]{$j$};
    \draw[webs,->-=.1] (3,.5) to (1,.5);
    \draw[webs,biarr] (0,1) node[left]{$i$} to (3,1) node[right]{$l$};
    \filldraw[fill=white, draw=black]
        (.5,-.25) rectangle (2.5,1.25);
    \node at (1.5,.5) {\Large{$C$}};
    \end{tikzpicture}%
    }
    {\Large{$=$}}
    \scalebox{.9}{%
    \begin{tikzpicture}[anchorbase]
    \draw[webs,biarr] (0,0) to (3,0);
    \draw[webs,->-=.1] (3,.5) to (1,.5);
    \draw[webs,biarr] (0,1) to (3,1);
    \filldraw[fill=white, draw=black]
        (.5,-.25) rectangle (2.5,1.25);
    \node at (1.5,.5) {\Large{$C$}};

    \draw[webs,-<-=.5] (0,-.5) to (3,-.5);
    \draw[webs] (3,-.5) to[out=0,in=180] (4,.5);
    \draw[webs,over] (3,0) to[out=0,in=180] (4,-.5);
    \draw[webs,over] (3,.5) to[out=0,in=180] (4,0);

    \draw[webs] (4,0) to[out=0,in=180] (4.5,.5);
    \draw[webs,over] (4,.5) to[out=0,in=180] (4.5,0);
    \draw[webs] (4.5,.5) to[out=0,in=180] (5,0);
    \draw[webs,over] (4.5,0) to[out=0,in=180] (5,.5);

    \draw[webs] (4,-.5) to (5,-.5);
    \draw[webs] (5,0) to[out=0,in=180] (5.5,-.5);
    \draw[webs,over] (5,-.5) to[out=0,in=180] (5.5,0) node[right]{$j$};
    \draw[webs] (5,.5) to (5.5,.5);
    \draw[webs] (3,1) to (5.5,1) node[right]{$l$};

    \draw[webs] (0,0) to[out=180,in=0] (-.5,-.5);
    \draw[webs,over] (0,-.5) to[out=180,in=0] (-.5,0);
    \draw[webs] (-.5,-.5) to[out=180,in=0] (-1,0) node[left]{$k$};
    \draw[webs,over] (-.5,0) to[out=180,in=0] (-1,-.5);
    \draw[webs] (0,1) to (-1,1) node[left]{$i$};

    \draw[webs] (-1,-.5) to[out=180,in=90] (-1.25,-.75);
    \draw[webs] (-1.25,-.75) to[out=270,in=180] (-1,-1);
    \draw[webs] (5.5,-.5) to[out=0,in=90] (5.75,-.75);
    \draw[webs] (5.75,-.75) to[out=270,in=0] (5.5,-1);
    \draw[webs] (-1,-1) to (5.5,-1);
    
    \end{tikzpicture}%
    }
    \\
    \vspace{1cm}
    {\Large{$\simeq$}}
    \scalebox{.9}{%
    \begin{tikzpicture}[anchorbase]
    \draw[webs,biarr] (0,0) node[left]{$k$} to (3,0);
    \draw[webs,->-=.1] (3,.5) to (1,.5);
    \draw[webs,biarr] (0,1) node[left]{$i$} to (3,1) ;
    \filldraw[fill=white, draw=black]
        (.5,-.25) rectangle (2.5,1.25);
    \node at (1.5,.5) {\Large{$C$}};

    \draw[webs] (3,.5) to[out=0,in=180] (3.5,0);
    \draw[webs,over] (3,0) to[out=0,in=180] (3.5,.5);
    \draw[webs] (3,-.5) to (3.5,-.5);
    \draw[webs] (3.5,-.5) to[out=0,in=180] (4.5,.5);
    \draw[webs,over] (3.5,0) to[out=0,in=180] (4.5,-.5);
    \draw[webs,over] (3.5,.5) to[out=0,in=180] (4.5,0) node[right]{$j$};
    \draw[webs] (3,1) to (4.5,1) node[right]{$l$};

    \draw[webs,->-=.5] (3,-.5) to (0,-.5);
    \draw[webs] (0,-.5) to[out=180,in=90] (-.25,-.75);
    \draw[webs] (-.25,-.75) to[out=270,in=180] (0,-1);
    \draw[webs] (4.5,-.5) to[out=0,in=90] (4.75,-.75);
    \draw[webs] (4.75,-.75) to[out=270,in=0] (4.5,-1);
    \draw[webs] (0,-1) to (4.5,-1);
    \end{tikzpicture}%
    }
    {\Large{$\simeq$}}
    \scalebox{.9}{%
    \begin{tikzpicture}[anchorbase]
    \draw[webs,biarr] (0,0) node[left]{$k$} to (3,0);
    \draw[webs,->-=.1] (3,.5) to (1,.5);
    \draw[webs,biarr] (0,1) node[left]{$i$} to (3,1) ;
    \filldraw[fill=white, draw=black]
        (.5,-.25) rectangle (2.5,1.25);
    \node at (1.5,.5) {\Large{$C$}};

    \draw[webs] (3,.5) to[out=0,in=180] (3.5,0);
    \draw[webs,over] (3,0) to[out=0,in=180] (3.5,.5);
    \draw[webs] (3.5,0) to[out=0,in=180] (4,.5);
    \draw[webs,over] (3.5,.5) to[out=0,in=180] (4,0) node[right]{$j$};
    \draw[webs] (3,1) to (4,1) node[right]{$l$};
    \end{tikzpicture}%
    }
    {\Large{$\simeq$}}
    \scalebox{.9}{%
    \begin{tikzpicture}[anchorbase]
    \draw[webs,biarr] (0,0) node[left]{$k$} to (3,0) node[right]{$j$};
    \draw[webs,->-=.1] (3,.5) to (1,.5);
    \draw[webs,biarr] (0,1) node[left]{$i$} to (3,1) node[right]{$l$} ;
    \filldraw[fill=white, draw=black]
        (.5,-.25) rectangle (2.5,1.25);
    \node at (1.5,.5) {\Large{$C$}};
    \end{tikzpicture}%
    }
    
    \caption{A Reidemeister I move. Here the first and third equivalences are sequences of braid relations, while the second is an application of Propositions \ref{prop: trace_local} and \ref{prop: markov}.}
    \label{fig: reidemeister_I}
\end{figure}

\FloatBarrier

\subsubsection{Reidemeister II} \label{subsec: reidemeister_II}

There are again 4 distinct cases given by a choice of orientation on each strand. When both strands are oriented the same direction, invariance follows immediately from the braid relations for Rouquier complexes. The case in which the overstrand is oriented upwards while the understrand is oriented downwards is handled by Lemma \ref{lem: trivial_horizontals} since both crossings act trivially. The remaining case follows in combination with Lemma \ref{lem: bending_around}, as illustrated below:

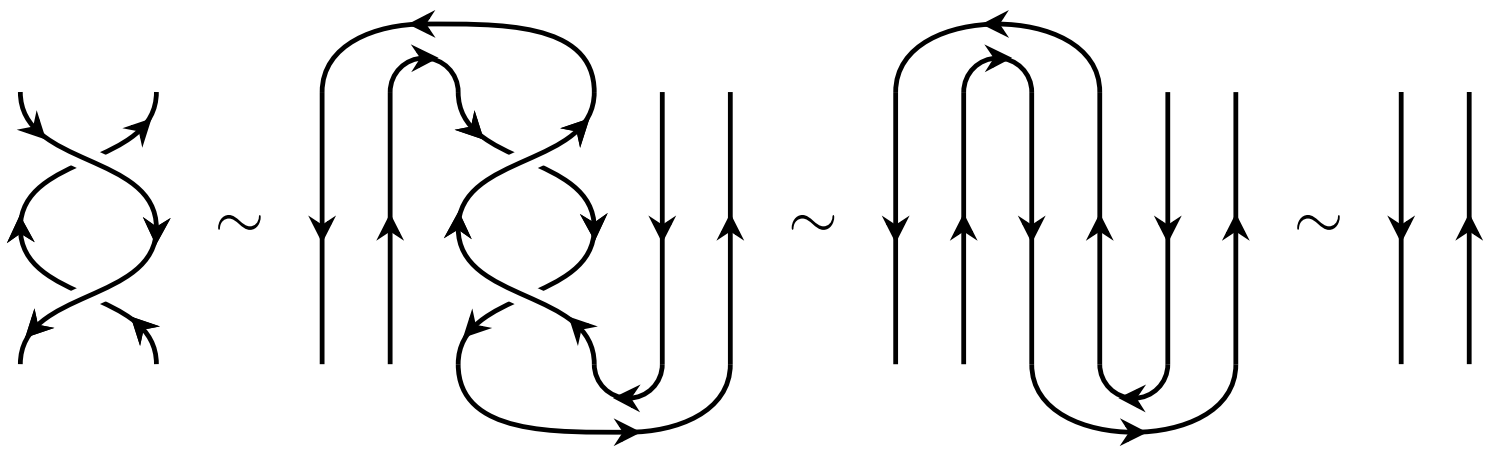
\begin{figure}[ht]
    \centering 

    \scalebox{.9}{%
    \begin{tikzpicture}[anchorbase]
    \draw[webs,->-=.8] (0,1) to[out=90,in=270] (1,2);
    \draw[webs,over,-<-=.8] (1,1) to[out=90,in=270] (0,2);
    \draw[webs,->-=.2,->-=.99] (1,0) to[out=90,in=270] (0,1);
    \draw[webs,over,-<-=.2,-<-=.99] (0,0) to[out=90,in=270] (1,1);
    \end{tikzpicture}%
    }
    \enspace {\Large{$\sim$}} \enspace 
    \scalebox{.9}{%
    \begin{tikzpicture}[anchorbase]
    \draw[webs,-<-=.5] (0,0) to (0,2);
    \draw[webs,->-=.5] (.5,0) to (.5,2);

    \draw[webs,->-=.99] (.5,2) to[out=90,in=180] (.75,2.25);
    \draw[webs] (.75,2.25) to[out=0,in=90] (1,2);
    \draw[webs,-<-=.99] (0,2) to[out=90,in=180] (.75,2.5);
    \draw[webs] (.75,2.5) to[out=0,in=90] (2,2);
    
    \draw[webs,->-=.8] (2,1) to[out=270,in=90] (1,0);
    \draw[webs,-<-=.8,over] (1,1) to[out=270,in=90] (2,0);
    \draw[webs,->-=.2,->-=.99] (1,2) to[out=270,in=90] (2,1);
    \draw[webs,over,-<-=.2,-<-=.99] (2,2) to[out=270,in=90] (1,1);

    \draw[webs,-<-=.99] (2,0) to[out=270,in=180] (2.25,-.25);
    \draw[webs] (2.25,-.25) to[out=0,in=270] (2.5,0);
    \draw[webs,->-=.99] (1,0) to[out=270,in=180] (2.25,-.5);
    \draw[webs] (2.25,-.5) to[out=0,in=270] (3,0);

    \draw[webs,-<-=.5] (2.5,0) to (2.5,2);
    \draw[webs,->-=.5] (3,0) to (3,2);
    \end{tikzpicture}%
    }
    \enspace {\Large{$\sim$}} \enspace 
    \scalebox{.9}{%
    \begin{tikzpicture}[anchorbase]
    \draw[webs,-<-=.5] (0,0) to (0,2);
    \draw[webs,->-=.5] (.5,0) to (.5,2);

    \draw[webs,->-=.99] (.5,2) to[out=90,in=180] (.75,2.25);
    \draw[webs] (.75,2.25) to[out=0,in=90] (1,2);
    \draw[webs,-<-=.99] (0,2) to[out=90,in=180] (.75,2.5);
    \draw[webs] (.75,2.5) to[out=0,in=90] (1.5,2);

    \draw[webs,->-=.5] (1,2) to (1,0);
    \draw[webs,-<-=.5] (1.5,2) to (1.5,0);

    \draw[webs,-<-=.99] (1.5,0) to[out=270,in=180] (1.75,-.25);
    \draw[webs] (1.75,-.25) to[out=0,in=270] (2,0);
    \draw[webs,->-=.99] (1,0) to[out=270,in=180] (1.75,-.5);
    \draw[webs] (1.75,-.5) to[out=0,in=270] (2.5,0);

    \draw[webs,-<-=.5] (2,0) to (2,2);
    \draw[webs,->-=.5] (2.5,0) to (2.5,2);
    \end{tikzpicture}%
    }
    \enspace {\Large{$\sim$}} \enspace 
    \scalebox{.9}{%
    \begin{tikzpicture}[anchorbase]
        \draw[webs,-<-=.5] (0,0) to (0,2);
        \draw[webs,->-=.5] (.5,0) to (.5,2);
    \end{tikzpicture}%
    }
    
    \caption{A wrong-way Reidemeister II move involving nontrivial horizontal crossings. The first and third isotopies follow from Lemma \ref{lem: bending_around}, while the second follows from Lemma \ref{lem: trivial_horizontals}.}
    \label{fig: reidemeister_II}
\end{figure}

\FloatBarrier

\subsubsection{Reidemeister III} \label{subsec: reidemeister_III}

Here there are 8 distinct cases, again given by a choice of orientation on each strand. We denote these cases by specifying the sequence of orientations on the input strands from left to right. The cases $+++$ and $---$ follow immediately from the braid relations for Rouquier complexes, while the cases $+--$ and $++-$ follow from Lemma \ref{lem: trivial_horizontals} by observing that only a single nontrivial crossing survives. Applying Lemma \ref{lem: bending_around}, we obtain the cases $-++$ and $--+$ from $+--$ and $++-$, respectively.

This leaves only the cases $+-+$ and $-+-$. Lemma \ref{lem: bending_around} again allows us to deduce one from the other, so we check only $+-+$. Here we compute directly using the decomposition of horizontal crossings from Figure \ref{fig: horizontal_crossing_decompositions}. We can apply Lemma \ref{lem: trivial_horizontals} once on each side to ignore the action of a trivial horizontal crossing. The two actions then differ only by braid relations on the right-hand side of the resulting register complex; see Figure \ref{fig: reidemeister_III_final} below.

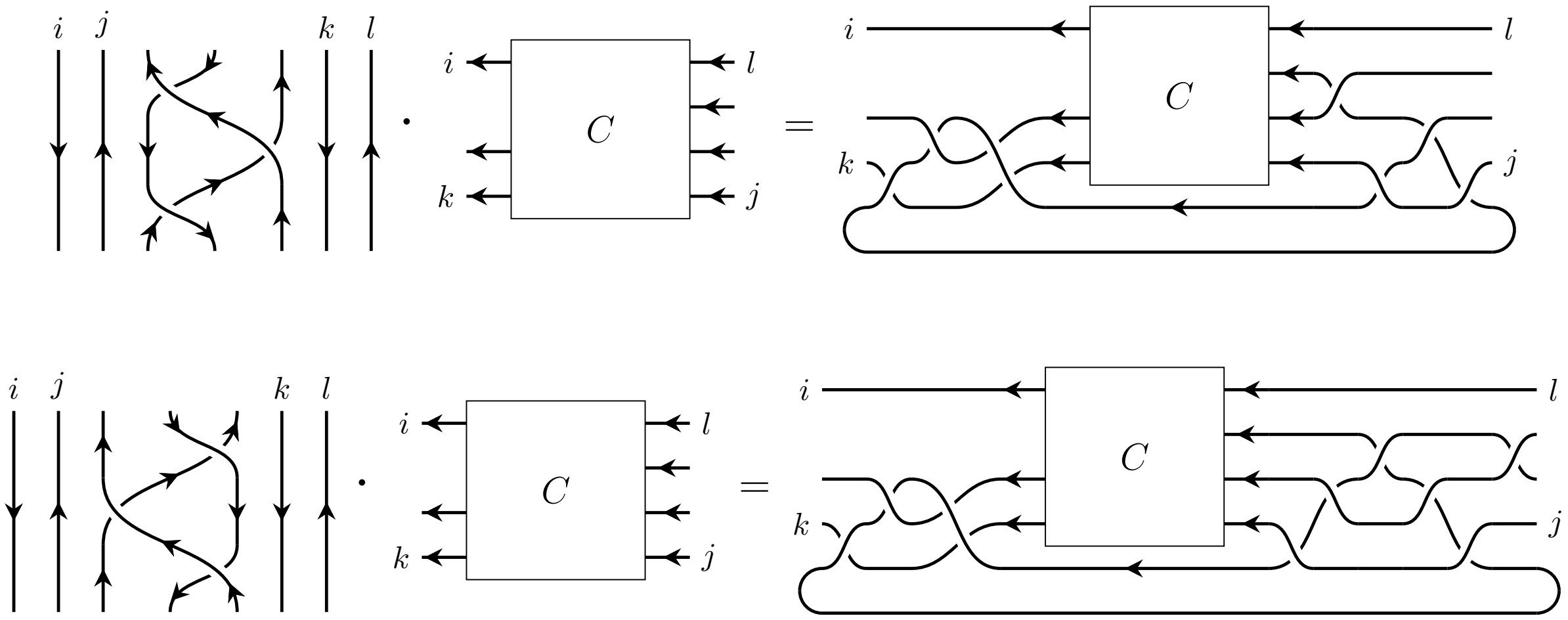
\begin{figure}[ht]
    \centering

    \scalebox{.9}{%
    \begin{tikzpicture}[anchorbase]
    \draw[webs,-<-=.5] (0,0) to (0,2.25) node[above]{$i$};
    \draw[webs,->-=.5] (.5,0) to (.5,2.25) node[above]{$j$};
    \begin{scope}[xshift=1cm,yshift=0]
    \draw[webs,->-=.1,->-=.5] (0,0) to[out=90,in=270] (1.5,1.5);
    \draw[webs,over,-<-=.2] (.75,0) to[out=90,in=270] (0,.75);
    \draw[webs,-<-=.5] (0,.75) to (0,1.5);
    \draw[webs,->-=.5] (1.5,0) to (1.5,.75);
    \draw[webs,-<-=.85] (0,1.5) to[out=90,in=270] (.75,2.25);
    \draw[webs,over,->-=.5,->-=.9] (1.5,.75) to[out=90,in=270] (0,2.25);
    \draw[webs,->-=.5] (1.5,1.5) to (1.5,2.25);
    \end{scope}
    \draw[webs,-<-=.5] (3,0) to (3,2.25) node[above]{$k$};
    \draw[webs,->-=.5] (3.5,0) to (3.5,2.25) node[above]{$l$};
    \end{tikzpicture}%
    }
    {\huge{$\cdot$}}
    \scalebox{.9}{%
    \begin{tikzpicture}[anchorbase]
    \draw[webs,biarr] (0,0) node[left]{$k$} to (3,0) node[right]{$j$};
    \draw[webs,biarr] (0,.5) to (3,.5);
    \draw[webs,-<-=.9] (.5,1) to (3,1);
    \draw[webs,biarr] (0,1.5) node[left]{$i$} to (3,1.5) node[right]{$l$};
    
    \filldraw[fill=white, draw=black]
        (.5,-.25) rectangle (2.5,1.75);
    \node at (1.5,.75) {\Large{$C$}};
    \end{tikzpicture}%
    }
    {\Large{$=$}}
    \scalebox{.9}{%
    \begin{tikzpicture}[anchorbase]
    \draw[webs,biarr] (0,0) to (3,0);
    \draw[webs,biarr] (0,.5) to (3,.5);
    \draw[webs,-<-=.9] (.5,1) to (3,1);
    \draw[webs,biarr] (0,1.5) to (3,1.5);
    
    \filldraw[fill=white, draw=black]
        (.5,-.25) rectangle (2.5,1.75);
    \node at (1.5,.75) {\Large{$C$}};

    \draw[webs,-<-=.5] (0,-.5) to (3,-.5);

    \draw[webs] (0,.5) to[out=180,in=0] (-1,0);
    \draw[webs] (0,0) to[out=180,in=0] (-1,-.5);
    \draw[webs,over] (0,-.5) to[out=180,in=0] (-1,.5);
    \draw[webs] (-1,.5) to[out=180,in=0] (-1.5,0);
    \draw[webs,over] (-1,0) to[out=180,in=0] (-1.5,.5);
    \draw[webs] (-1,-.5) to (-1.5,-.5);
    \draw[webs] (-1.5,-.5) to[out=180,in=0] (-2,0) node[left]{$k$};
    \draw[webs,over] (-1.5,0) to[out=180,in=0] (-2,-.5);
    \draw[webs] (-1.5,.5) to (-2,.5);
    \draw[webs] (0,1.5) to (-2,1.5) node[left]{$i$};

    \draw[webs] (-2,-.5) to[out=180,in=90] (-2.25,-.75);
    \draw[webs] (-2.25,-.75) to[out=270,in=180] (-2,-1);

    \draw[webs] (3,1) to[out=0,in=180] (3.5,.5);
    \draw[webs,over] (3,.5) to[out=0,in=180] (3.5,1);
    \draw[webs] (3.5,-.5) to[out=0,in=180] (4,0);
    \draw[webs,over] (3.5,0) to[out=0,in=180] (4,-.5);
    \draw[webs] (4,.5) to[out=0,in=180] (5,-.5);
    \draw[webs,over] (4,0) to[out=0,in=180] (4.5,.5);
    \draw[webs] (4.5,.5) to (5,.5);
    \draw[webs] (4,-.5) to (4.5,-.5);
    \draw[webs,over] (4.5,-.5) to[out=0,in=180] (5,0) node[right]{$j$};
    \draw[webs] (3,1.5) to (5,1.5) node[right]{$l$};
    \draw[webs] (3.5,1) to (5,1);
    \draw[webs] (3,0) to (3.5,0);
    \draw[webs] (3,-.5) to (3.5,-.5);
    \draw[webs] (3.5,.5) to (4,.5);

    \draw[webs] (5,-.5) to[out=0,in=90] (5.25,-.75);
    \draw[webs] (5.25,-.75) to[out=270,in=0] (5,-1);

    \draw[webs] (-2,-1) to (5,-1);
    
    \end{tikzpicture}%
    }
    \\
    \vspace{1cm}
    \scalebox{.9}{%
    \begin{tikzpicture}[anchorbase]
    \draw[webs,-<-=.5] (0,0) to (0,2.25) node[above]{$i$};
    \draw[webs,->-=.5] (.5,0) to (.5,2.25) node[above]{$j$};
    \begin{scope}[yscale=1,xscale=-1,xshift=-2.5cm]
    \draw[webs,-<-=.2] (.75,0) to[out=90,in=270] (0,.75);
    \draw[webs,-<-=.5] (0,.75) to (0,1.5);
    \draw[webs,->-=.5] (1.5,0) to (1.5,.75);
    \draw[webs,->-=.5,->-=.9] (1.5,.75) to[out=90,in=270] (0,2.25);
    \draw[webs,->-=.5,over] (1.5,1.5) to (1.5,2.25);
    \draw[webs,-<-=.85,over] (0,1.5) to[out=90,in=270] (.75,2.25);
    \draw[webs,->-=.1,->-=.5,over] (0,0) to[out=90,in=270] (1.5,1.5);
    \end{scope}
    \draw[webs,-<-=.5] (3,0) to (3,2.25) node[above]{$k$};
    \draw[webs,->-=.5] (3.5,0) to (3.5,2.25) node[above]{$l$};
    \end{tikzpicture}%
    }
    {\huge{$\cdot$}}
    \scalebox{.9}{%
    \begin{tikzpicture}[anchorbase]
    \draw[webs,biarr] (0,0) node[left]{$k$} to (3,0) node[right]{$j$};
    \draw[webs,biarr] (0,.5) to (3,.5);
    \draw[webs,-<-=.9] (.5,1) to (3,1);
    \draw[webs,biarr] (0,1.5) node[left]{$i$} to (3,1.5) node[right]{$l$};
    
    \filldraw[fill=white, draw=black]
        (.5,-.25) rectangle (2.5,1.75);
    \node at (1.5,.75) {\Large{$C$}};
    \end{tikzpicture}%
    }
    {\Large{$=$}}
    \scalebox{.9}{%
    \begin{tikzpicture}[anchorbase]
    \draw[webs,biarr] (0,0) to (3,0);
    \draw[webs,biarr] (0,.5) to (3,.5);
    \draw[webs,-<-=.9] (.5,1) to (3,1);
    \draw[webs,biarr] (0,1.5) to (3,1.5);
    
    \filldraw[fill=white, draw=black]
        (.5,-.25) rectangle (2.5,1.75);
    \node at (1.5,.75) {\Large{$C$}};

    \draw[webs,-<-=.5] (0,-.5) to (3,-.5);

    \draw[webs] (0,.5) to[out=180,in=0] (-1,0);
    \draw[webs] (0,0) to[out=180,in=0] (-1,-.5);
    \draw[webs,over] (0,-.5) to[out=180,in=0] (-1,.5);
    \draw[webs] (-1,.5) to[out=180,in=0] (-1.5,0);
    \draw[webs,over] (-1,0) to[out=180,in=0] (-1.5,.5);
    \draw[webs] (-1,-.5) to (-1.5,-.5);
    \draw[webs] (-1.5,-.5) to[out=180,in=0] (-2,0) node[left]{$k$};
    \draw[webs,over] (-1.5,0) to[out=180,in=0] (-2,-.5);
    \draw[webs] (-1.5,.5) to (-2,.5);
    \draw[webs] (0,1.5) to (-2,1.5) node[left]{$i$};

    \draw[webs] (-2,-.5) to[out=180,in=90] (-2.25,-.75);
    \draw[webs] (-2.25,-.75) to[out=270,in=180] (-2,-1);

    \draw[webs] (3,-.5) to[out=0,in=180] (4,.5);
    \draw[webs,over] (3,0) to[out=0,in=180] (3.5,-.5);
    \draw[webs] (3.5,-.5) to (4,-.5);
    \draw[webs] (3,.5) to (3.5,.5);
    \draw[webs, over] (3.5,.5) to[out=0,in=180] (4,0);
    \draw[webs] (4,1) to[out=0,in=180] (4.5,.5);
    \draw[webs,over] (4,.5) to[out=0,in=180] (4.5,1);
    \draw[webs] (4.5,.5) to[out=0,in=180] (5.5,-.5);
    \draw[webs,over] (4.5,0) to[out=0,in=180] (5,.5);
    \draw[webs] (5,.5) to (5.5,.5);
    \draw[webs] (4.5,-.5) to (5,-.5);
    \draw[webs, over] (5,-.5) to[out=0,in=180] (5.5,0);
    \draw[webs] (5.5,1) to[out=0,in=180] (6,.5);
    \draw[webs,over] (5.5,.5) to[out=0,in=180] (6,1);

    \draw[webs] (3,1.5) to (6,1.5) node[right]{$l$};
    \draw[webs] (3,1) to (4,1);
    \draw[webs] (4.5,1) to (5.5,1);
    \draw[webs] (4,0) to (4.5,0);
    \draw[webs] (4,-.5) to (4.5,-.5);
    \draw[webs] (5.5,0) to (6,0) node[right]{$j$};
    \draw[webs] (5.5,-.5) to (6,-.5);

    \draw[webs] (6,-.5) to[out=0,in=90] (6.25,-.75);
    \draw[webs] (6.25,-.75) to[out=270,in=0] (6,-1);

    \draw[webs] (-2,-1) to (6,-1);
    
    \end{tikzpicture}%
    }

    \caption{An oriented Reidemeister III move. Note that the resulting complexes are naturally homotopy equivalent via  a sequence of braid relations.}
    \label{fig: reidemeister_III_final}

\end{figure}

\FloatBarrier

\subsection{Useful Identities} \label{subsec: useful_identities}

Before moving on, we pause to compute the effect on $\Phi$ of a few especially useful operations on tangles.

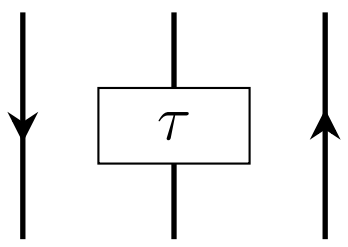
\begin{figure}[ht]
\centering

\begin{tikzpicture}
    \draw[webs,-<-=.5] (0,0) to (0,1.5);
    \draw[webs,->-=.5] (2,0) to (2,1.5);
    \draw[webs] (1,0) node[below]{$n$} to (1,1.5) node[above]{$m$};

    \filldraw[fill=white, draw=black]
        (.5,.5) rectangle (1.5,1);
        \node at (1,.75) {\Large{$\tau$}};
\end{tikzpicture}

\caption{$I(\tau)$ for an oriented $(m, n)$-tangle $\tau$.}
\label{fig: inclusion_tangles}
\end{figure}

\begin{definition} \label{def: tangle_inclusion}
    For any balanced oriented $(m,n)$-tangle $\tau$, we denote by $I(\tau)$ the balanced oriented $(m+2,n+2)$-tangle given by adding a downwards-oriented strand to the left of $\tau$ and an upwards-oriented strand to the right of $\tau$. See Figure \ref{fig: inclusion_tangles}.
\end{definition}

\begin{prop} \label{prop: monoidality}
    Let $C \in K^b(\extsbim_l)$ and $C' \in K^b(\extsbim_n)$ be given. Then for any balanced oriented $(2l,2m)$-tangle $\tau$, we have $ \Phi(I^n(\tau)) \big(C \boxtimes C' \big) = \big(\Phi(\tau)(C) \big) \boxtimes C'$.
\end{prop}

\begin{proof}
    We use the diagrammatic framework of Figure \ref{fig: general_framework}. The leftmost $n$ downwards-oriented vertical strands of the instruction tangle $I^n(\tau)$ connect directly to the left-hand boundary of the component $C'$ of the register complex; similarly, the rightmost $n$ upwards-oriented vertical strands connect to the right-hand boundary of $C'$. What remains is the component $\tau$ of the instruction tangle connected via some shuffle braid to the component $C$ of the register complex, but this exactly illustrates the action of $\tau$ on $C$.
\end{proof}

We can now interpret the algebraic inclusion functor $I$ on $K^b(\extsbim)$ in terms of concentric cups:

\begin{cor} \label{cor: concentric_cups}
    For each $n \geq 0$, let $\bigcap_n$ be the $(0,2n)$-tangle consisting of $n$ concentric counterclockwise-oriented cups. Then for each $C \in K^b(\extsbim_m)$, we have $I^n(C) = \Phi(I^m(\bigcap_n))(C)$.
\end{cor}

\begin{proof}
    Observe that $\bigcap_n = \bigcap_{n-1} \star I^{n-1}(\bigcap_1)$ for each $n$. The desired result then follows immediately from Proposition \ref{prop: monoidality} and the fact that $\Phi(\bigcap_1)(\mathbb{Z}) = \strand{1}$.
\end{proof}

We have a similar interpretation of the partial trace using concentric caps:

\begin{lemma} \label{lem: concentric_caps}
    For each $n \geq 0$, let $\bigcup_n$ be the $(2n,0)$-tangle consisting of $n$ concentric counterclockwise-oriented caps. Then for each $C \in K^b(\extsbim_m)$ and each $n \leq m$, we have $Tr^n(C) = \Phi(I^{m-n}(\bigcup_n))(C)$.
\end{lemma}

\begin{proof}
    Direct computation.
\end{proof}

\section{HOMFLY Homology} \label{sec: HHH}

\subsection{HOMFLY Homology} \label{subsec: HHH}

We now turn our attention to oriented $(0, 2n)$-tangles. Any such tangle $\tau$ is automatically balanced, so it defines a functor $\Phi(\tau) \colon K^b(\extsbim_0) \to K^b(\extsbim_n)$. Because $\Phi(\tau)$ is an exact additive functor, it is completely determined by the complex $\Phi(\tau)(\mathbb{Z}) \in K^b(\extsbim_n)$.

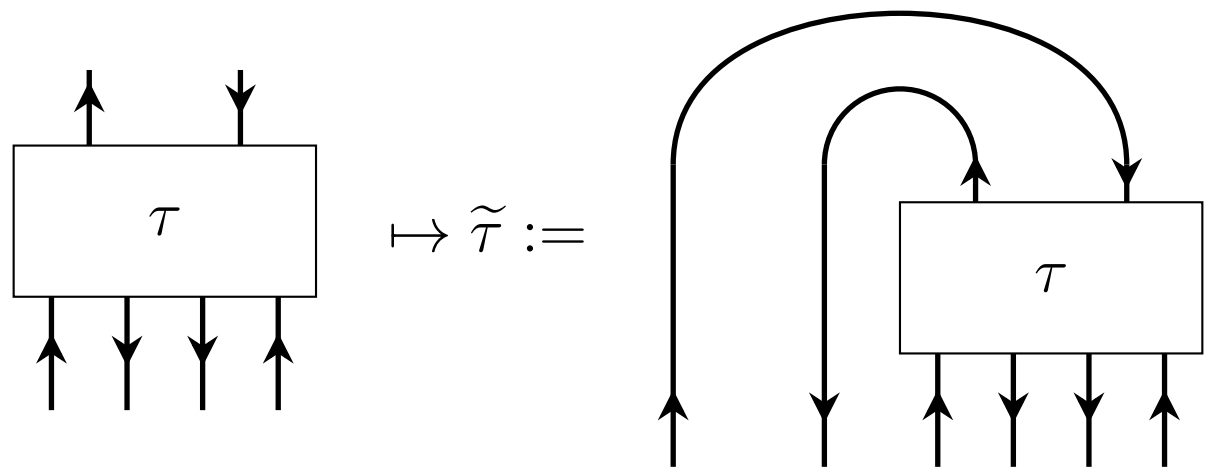
\begin{figure}[ht]
    \centering

    \begin{tikzpicture}[anchorbase]
        \draw[webs,->-=.4] (0,0) to (0,1);
        \draw[webs,-<-=.4] (.5,0) to (0.5,1);
        \draw[webs,-<-=.4] (1,0) to (1,1);
        \draw[webs,->-=.4] (1.5,0) to (1.5,1);

        \draw[webs,->-=.85] (.25,1) to (.25,2.25);
        \draw[webs,-<-=.85] (1.25,1) to (1.25,2.25);

        \filldraw[fill=white, draw=black]
        (-.25,.75) rectangle (1.75,1.75);
        \node at (.75,1.25) {\Large{$\tau$}};
    \end{tikzpicture}
    \quad {\Large{$\mapsto \tilde{\tau} :=$}} \quad
    \begin{tikzpicture}[anchorbase]
        \draw[webs,->-=.4] (0,0) to (0,1);
        \draw[webs,-<-=.4] (.5,0) to (0.5,1);
        \draw[webs,-<-=.4] (1,0) to (1,1);
        \draw[webs,->-=.4] (1.5,0) to (1.5,1);

        \draw[webs,->-=.95] (.25,1) to (.25,2);
        \draw[webs] (.25,2) to[out=90,in=0] (-.25,2.5);
        \draw[webs] (-.25,2.5) to[out=180,in=90] (-.75,2);
        \draw[webs,->-=.8] (-.75,2) to (-.75,0);
        \draw[webs,-<-=.95] (1.25,1) to (1.25,2);
        \draw[webs] (1.25,2) to[out=90,in=0] (-.25,3);
        \draw[webs] (-.25,3) to[out=180,in=90] (-1.75,2);
        \draw[webs,-<-=.8] (-1.75,2) to (-1.75,0);

        \filldraw[fill=white, draw=black]
        (-.25,.75) rectangle (1.75,1.75);
        \node at (.75,1.25) {\Large{$\tau$}};
    \end{tikzpicture}

    \caption{Bending an oriented $(2,4)$-tangle to produce an oriented $(0,6)$-tangle.}
    \label{fig: bending_tangles}
\end{figure}

More generally, given any oriented $(m, n)$-tangle $\tau$, there is a canonical\footnote{In principle there are many ways to bend a given oriented tangle \textit{diagram} to the left, but any isotopy of oriented $(m,n)$-tangle diagrams is sent to an isotopy of oriented $(0,m+n)$-tangle diagrams by this procedure.} way to produce an oriented $(0, m+n)$-tangle $\tilde{\tau}$ by bending the top boundary of $\tau$ to the left; see Figure \ref{fig: bending_tangles}. We can therefore assign a complex in $K^b(\extsbim)$ to any oriented tangle by the same argument. We collect these observations into the following definition:

\begin{definition} \label{def: tangle_hhh}
    Let $\tau$ be an oriented tangle with $2n$ endpoints. Then we denote by $HHH(\tau)$ the complex

    \[
    HHH(\tau) := \Phi(\tilde{\tau})(\mathbb{Z}) \in K^b(\extsbim_n).
    \]

    We call $HHH(\tau)$ the \textit{HOMFLY homology} of the oriented tangle $\tau$.
\end{definition}

Theorem \ref{thm: main_invariance} guarantees that $HHH(\tau)$ is a tangle invariant, but before stating a theorem regarding this invariance, we need to address an ambiguity in our framework. In Section \ref{subsec: otan}, all oriented tangle diagrams were assumed to be embedded in the thickened unit square $[0,1]^2 \times (-\epsilon, \epsilon)$, and all isotopies were assumed to fix the boundary. This introduces an auxiliary choice of a division of the endpoints of $\tau$ into top and bottom endpoint sequences. After bending around to obtain a $(0, 2n)$-tangle $\tilde{\tau}$, this choice is equivalent to a single choice of a leftmost (or rightmost) boundary point along the bottom boundary of $[0,1]^2$, rendering the division into a ``top" and ``bottom" collection of endpoints more or less arbitrary.

\begin{figure}[ht]
    \centering
    
    \begin{tikzpicture}[scale=1]
    \coordinate (top)   at (0,3.5);
    \coordinate (left)  at (-3.6,-3.4);
    \coordinate (right) at ( 3.6,-3.4);

    \draw[->, thick, draw=black!65]
        (-0.9,1.45) -- (-3.25,-2.25)
        node[midway, above, sloped, fill=white, inner sep=1pt] {Cut at $p$};

    \draw[->, thick, draw=black!65]
        (0.9,1.45) -- (3.6,-0.45)
        node[midway, above, sloped, fill=white, inner sep=1pt] {Cut at $q$};

    \begin{scope}[shift={(top)}]
        \draw[line width=1.4pt, draw=black!55] (0,0) circle (2);

        \draw[line width=1.2pt] (105:1.85) -- (105:2.15);
        \node at (105:2.35) {$p$};

        \draw[line width=1.2pt] (-160.7:1.85) -- (-160.7:2.15);
        \node at (-160.7:2.35) {$q$};

        \draw[thick, fill=white] (-0.9,-0.6) rectangle (0.9,0.6);
        \node at (0,0) {\Large{$\tau$}};

        \draw[webs,->-=.5] (-0.9,  0.45) .. controls (-1.15,  0.65) and (-1.35,  0.95) .. (-1.56,  1.25);
        \draw[webs,-<-=.5] (-0.9,  0.00) .. controls (-1.25,  0.00) and (-1.60,  0.00) .. (-2.00,  0.00);
        \draw[webs,->-=.5] (-0.9, -0.45) .. controls (-1.15, -0.65) and (-1.35, -0.95) .. (-1.56, -1.25);

        \draw[webs,->-=.5] ( 0.9,  0.45) .. controls ( 1.15,  0.65) and ( 1.35,  0.95) .. ( 1.56,  1.25);
        \draw[webs,-<-=.5] ( 0.9,  0.00) .. controls ( 1.25,  0.00) and ( 1.60,  0.00) .. ( 2.00,  0.00);
        \draw[webs,-<-=.5] ( 0.9, -0.45) .. controls ( 1.15, -0.65) and ( 1.35, -0.95) .. ( 1.56, -1.25);
    \end{scope}

    \begin{scope}[shift={(left)}]
        \draw[thick, fill=white] (-0.9,-0.6) rectangle (0.9,0.6);
        \node at (0,0) {\Large{$\tau$}};

        \draw[webs,->-=.5] (-0.9,  0.45) to[out=180,in=90] (-2.5,-1.5);
        \draw[webs,-<-=.5] (-0.9,  0.00) to[out=180,in=90] (-2,-1.5);
        \draw[webs,->-=.5] (-0.9, -0.45) to[out=180,in=90] (-1.5,-1.5);

        \draw[webs,->-=.5] ( 0.9,  0.45) to[out=0,in=90] (2.5,-1.5);
        \draw[webs,-<-=.5] ( 0.9,  0.00) to[out=0,in=90] (2,-1.5);
        \draw[webs,-<-=.5] ( 0.9, -0.45) to[out=0,in=90] (1.5,-1.5);
    \end{scope}

    \begin{scope}[shift={(right)}]
        \draw[thick, fill=white] (-0.9,-0.6) rectangle (0.9,0.6);
        \node at (0,0) {\Large{$\tau$}};

        \draw[webs,->-=.5]
            (-0.9, 0.45)
            .. controls (-1.45,0.75) and (-1.55,1.35) .. (-0.2,1.55)
            .. controls (1.35,1.75) and (3.0,0.55) .. (3.0,-0.85)
            -- (3.0,-1.5);

        \draw[webs,-<-=.5]
            (-0.9, 0.00)
            .. controls (-1.85,0.10) and (-2.55,1.35) .. (-1.05,2.15)
            .. controls (0.65,3.00) and (3.5,1.30) .. (3.5,-0.85)
            -- (3.5,-1.5);

        \draw[webs,->-=.5] (-0.9, -0.45) to[out=180,in=90] (-1.5,-1.5);

        \draw[webs,->-=.5] ( 0.9,  0.45) to[out=0,in=90] (2.5,-1.5);
        \draw[webs,-<-=.5] ( 0.9,  0.00) to[out=0,in=90] (2,-1.5);
        \draw[webs,-<-=.5] ( 0.9, -0.45) to[out=0,in=90] (1.5,-1.5);
    \end{scope}
    \end{tikzpicture}

    \caption{An oriented tangle $\tau$ embedded in a thickened disk. Cutting the boundary circle at the two marked points $p$ and $q$ produces two non-isotopic oriented $(0,2n)$-tangles.}
    \label{fig: disk_tangle}
\end{figure}
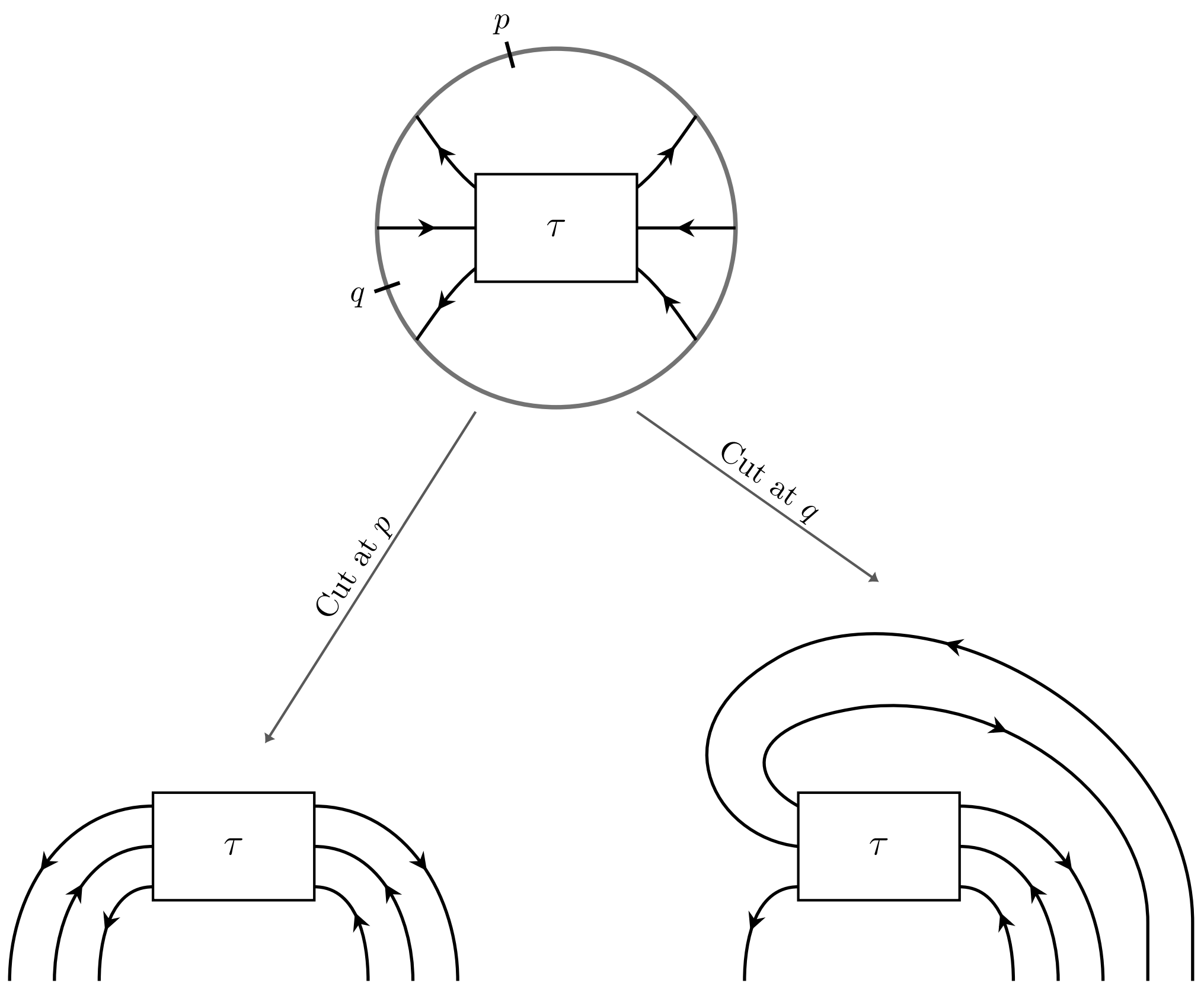

To avoid this artificial distinction, we can instead take the more typical approach of working with tangles embedded in the thickened disk $\mathbb{D}^2 \times (-\epsilon, \epsilon)$. The choice of a leftmost boundary point is then equivalent to a choice of marked point between two endpoints of $\tau$ along the boundary $S^1 = \partial \mathbb{D}^2$; see Figure \ref{fig: disk_tangle} below for an example. With this disclaimer out of the way, we can now state our primary invariance theorem:

\begin{theorem} \label{thm: main_invariance_HHH}
    Let $\tau$ and $\tau'$ be two oriented tangle diagrams embedded in $\mathbb{D}^2 \times (-\epsilon, \epsilon)$, and suppose $\tau$ and $\tau'$ depict isotopic oriented tangles. Then for each choice of marked point $p \in \partial \mathbb{D}^2$, we have $HHH(\tau) \cong HHH(\tau')$.
\end{theorem}

\FloatBarrier

One might hope that different choices of marked points are invisible at the level of $HHH$ so that $HHH$ becomes an invariant of oriented tangles in the thickened disk. This turns out to be too optimistic, but we can nevertheless analyze the effect of different choices of marked points on $HHH$ of the resulting tangles. We begin with a bit of setup:

\begin{prop} \label{prop: ccw_cup}
    Let $\tau$ be an oriented $(0, 2n)$-tangle, and consider the oriented $(0,2n+2)$-tangle $\tau'$ obtained by adding a counterclockwise cup over $\tau$. Then $HHH(\tau') = HHH(\tau) \boxtimes \strand{1}$.
\end{prop}

\begin{proof}
    We can express $\tau'$ as a vertical concatenation $\tau' \sim \bigcap \star I(\tau)$. Then by Proposition \ref{prop: monoidality} and functoriality of $\Phi$, we have

    \[
    HHH(\tau') = \Phi(I(\tau))(\strand{1}) = HHH(\tau) \boxtimes \strand{1}.
    \]
\end{proof}

We will also need an analogue of Proposition \ref{prop: ccw_cup} for clockwise cups. This requires a bit of additional notation: for each $n \geq 0$, we denote by $J_{n+1}$ the Rouquier complex assigned to the braid word $\sigma_n \sigma_{n-1} \dots \sigma_{2} \sigma_1^2 \sigma_{2} \dots \sigma_{n-1} \sigma_n$.

\begin{prop} \label{prop: clockwise_cup}
    Let $\tau$ be an oriented $(0,2n)$-tangle, and consider the oriented $(0,2n+2)$-tangle $\tau'$ obtained by adding a clockwise cup over $\tau$. Then $HHH(\tau') = J_{n+1} \star \left(HHH(\tau) \boxtimes \strand{1} \right)$.
\end{prop}

\begin{proof}
    We begin by applying the following isotopy:

    \FloatBarrier

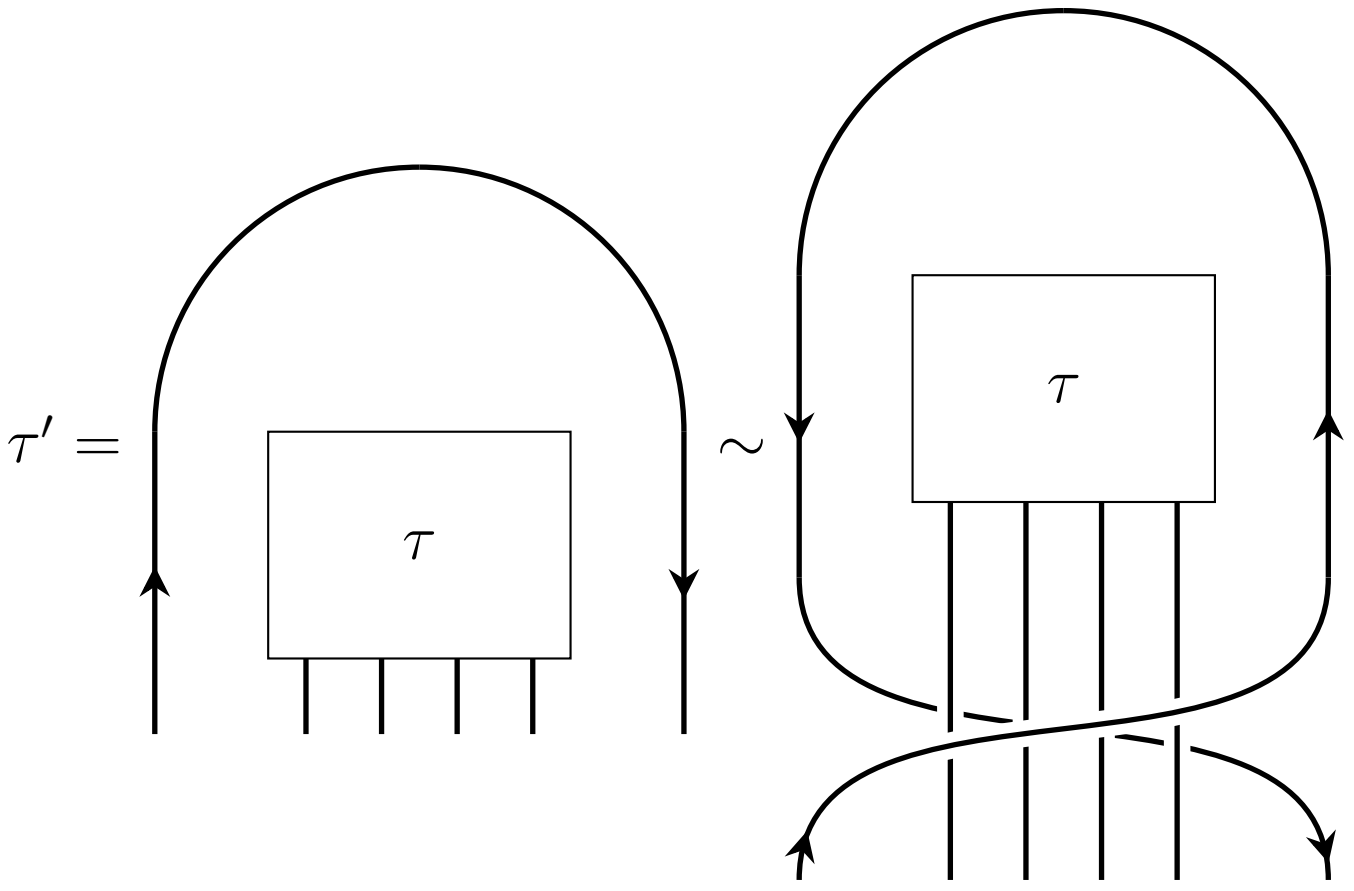
    \begin{figure}[ht]
    \centering

    {\Large{$\tau' = $}}
    \begin{tikzpicture}[anchorbase]
        \draw[webs] (0,0) to (0,2);
        \draw[webs] (.5,0) to (.5,2);
        \draw[webs] (1,0) to (1,2);
        \draw[webs] (1.5,0) to (1.5,2);

        \filldraw[fill=white, draw=black]
        (-.25,.5) rectangle (1.75,2);
        \node at (.75,1.25) {\Large{$\tau$}};

        \draw[webs,->-=.5] (-1,0) to (-1,2);
        \draw[webs] (-1,2) to[out=90,in=180] (.75,3.75);
        \draw[webs] (.75,3.75) to[out=0,in=90] (2.5,2);
        \draw[webs,->-=.5] (2.5,2) to (2.5,0);
    \end{tikzpicture}
    {\Large{$\sim$}}
    \begin{tikzpicture}[anchorbase]
        \draw[webs,-<-=.5] (-1,0) to (-1,2);
        \draw[webs] (-1,2) to[out=90,in=180] (.75,3.75);
        \draw[webs] (.75,3.75) to[out=0,in=90] (2.5,2);
        \draw[webs,-<-=.5] (2.5,2) to (2.5,0);

        \draw[webs,->-=.95] (-1,0) to[out=270,in=90] (2.5,-2);
        
        \draw[webs,over] (0,-2) to (0,2);
        \draw[webs,over] (.5,-2) to (.5,2);
        \draw[webs,over] (1,-2) to (1,2);
        \draw[webs,over] (1.5,-2) to (1.5,2);

        \filldraw[fill=white, draw=black]
        (-.25,.5) rectangle (1.75,2);
        \node at (.75,1.25) {\Large{$\tau$}};
        
        \draw[webs,over,-<-=.95] (2.5,0) to[out=270,in=90] (-1,-2);
    \end{tikzpicture}
    \end{figure}

    \FloatBarrier

    We compute $HHH(\tau')$ using the latter oriented tangle diagram. By Proposition \ref{prop: ccw_cup}, $\Phi$ applied to the portion of the diagram above the crossings takes $\mathbb{Z}$ to $HHH(\tau) \boxtimes \strand{1}$.
    
    Now, as the upwards-oriented boundary strand passes over the strands exiting $\tau$, any horizontal crossing in which it participates can be safely ignored by Lemma \ref{lem: trivial_horizontals}. The remainder of the crossings in which it participates are vertical crossings between upwards-oriented strands; since this strand passes over all of the upwards-oriented strands exiting $\tau$ from right to left, $\Phi$ of these crossings gives a total factor of $- \star F(\sigma_1 \sigma_{2} \dots \sigma_n)$.
    
    Similarly, as the downwards-oriented boundary strand passes under the strands exiting $\tau$, any horizontal crossing in which it participates must involve an upwards-oriented strand and can again be safely ignored by Lemma \ref{lem: trivial_horizontals}. The remainder of the crossings in which it participates are vertical crossings between downwards-oriented strands; since this strand passes under all of the downwards-oriented strands exiting $\tau$ from left to right, $\Phi$ of these crossings gives a total factor of $F(\sigma_n \sigma_{n-1} \dots \sigma_1) \star -$. After reassociating using Proposition \ref{prop: diagrammatics_isotopy} (proved below), we see that $\Phi$ of the bottom portion of this diagram contributes exactly $J_{n+1} \star -$.
\end{proof}

We can now describe how $HHH$ varies as we rotate the marked point past a boundary point of the tangle:

\begin{prop} \label{prop: tangles_in_disks}
    Let $\tau$ be an oriented tangle in $\mathbb{D}^2 \times (-\epsilon, \epsilon)$ with $2n$ endpoints, and consider the oriented $(0, 2n)$-tangles $\tilde{\tau}$ and $\tilde{\tau'}$ obtained by cutting the boundary circle in the two locations indicated in Figure \ref{fig: two_cuts} on the left and right, respectively. If $p$ is oriented inward, then $HHH(\tilde{\tau'}) \simeq HHH(\tilde{\tau}) \star F(\sigma_1^{-1} \sigma_2^{-1} \dots \sigma_{n-1}^{-1})$. Similarly, if $p$ is oriented outward, then $HHH(\tilde{\tau'}) \simeq F(\sigma_1 \sigma_2 \dots \sigma_{n-1}) \star HHH(\tilde{\tau})$. 

    \FloatBarrier
    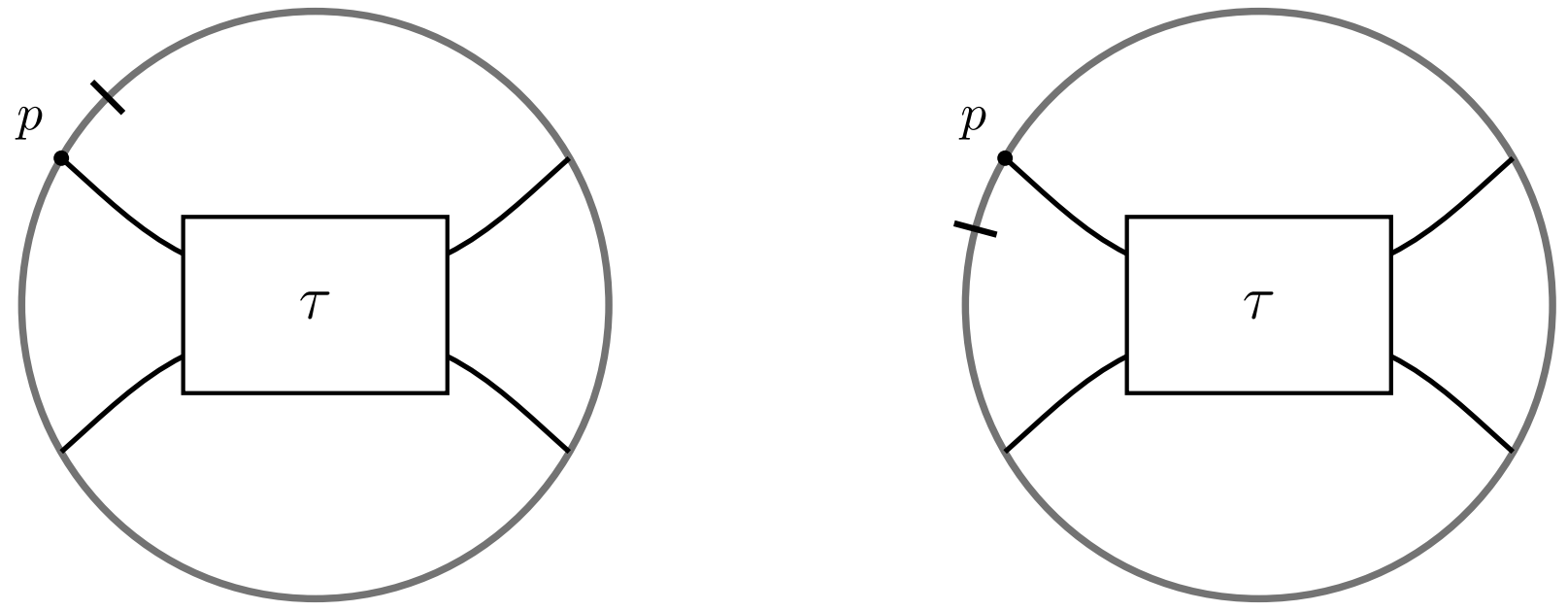
\begin{figure}[ht]
    \centering
    \begin{tikzpicture}[scale=1]
    \draw[line width=1.4pt, draw=black!55] (0,0) circle (2);

    \fill (-1.73, 1.00) circle (1.5pt);
    \node[above left] at (-1.73, 1.00) {$p$};

    \draw[thick, fill=white] (-0.9,-0.6) rectangle (0.9,0.6);
    \node at (0,0) {\Large{$\tau$}};

    \draw[webs] (-0.9,  0.35) .. controls (-1.2,  0.5) and (-1.45,  0.75) .. (-1.73,  1.00);
    \draw[webs] (-0.9, -0.35) .. controls (-1.2, -0.5) and (-1.45, -0.75) .. (-1.73, -1.00);
    \draw[webs] ( 0.9,  0.35) .. controls ( 1.2,  0.5) and ( 1.45,  0.75) .. ( 1.73,  1.00);
    \draw[webs] ( 0.9, -0.35) .. controls ( 1.2, -0.5) and ( 1.45, -0.75) .. ( 1.73, -1.00);

    \draw[line width=1.2pt] (135:1.85) -- (135:2.15);
    \end{tikzpicture}
    \hspace{2cm}
    \begin{tikzpicture}[scale=1]
    \draw[line width=1.4pt, draw=black!55] (0,0) circle (2);

    \fill (-1.73, 1.00) circle (1.5pt);
    \node[above left] at (-1.73, 1.00) {$p$};

    \draw[thick, fill=white] (-0.9,-0.6) rectangle (0.9,0.6);
    \node at (0,0) {\Large{$\tau$}};

    \draw[webs] (-0.9,  0.35) .. controls (-1.2,  0.5) and (-1.45,  0.75) .. (-1.73,  1.00);
    \draw[webs] (-0.9, -0.35) .. controls (-1.2, -0.5) and (-1.45, -0.75) .. (-1.73, -1.00);
    \draw[webs] ( 0.9,  0.35) .. controls ( 1.2,  0.5) and ( 1.45,  0.75) .. ( 1.73,  1.00);
    \draw[webs] ( 0.9, -0.35) .. controls ( 1.2, -0.5) and ( 1.45, -0.75) .. ( 1.73, -1.00);

    \draw[line width=1.2pt] (165:1.85) -- (165:2.15);
    \end{tikzpicture}

    \caption{An oriented tangle $\tau$ in the disk and boundary cuts on either side of a chosen boundary point $p$.}
    \label{fig: two_cuts}
    \end{figure}
    \FloatBarrier
\end{prop}

\begin{proof}
    Suppose first that $p$ is oriented inward. Then we can compute $HHH(\tilde{\tau'})$ using Proposition \ref{prop: ccw_cup} as follows:

    \FloatBarrier
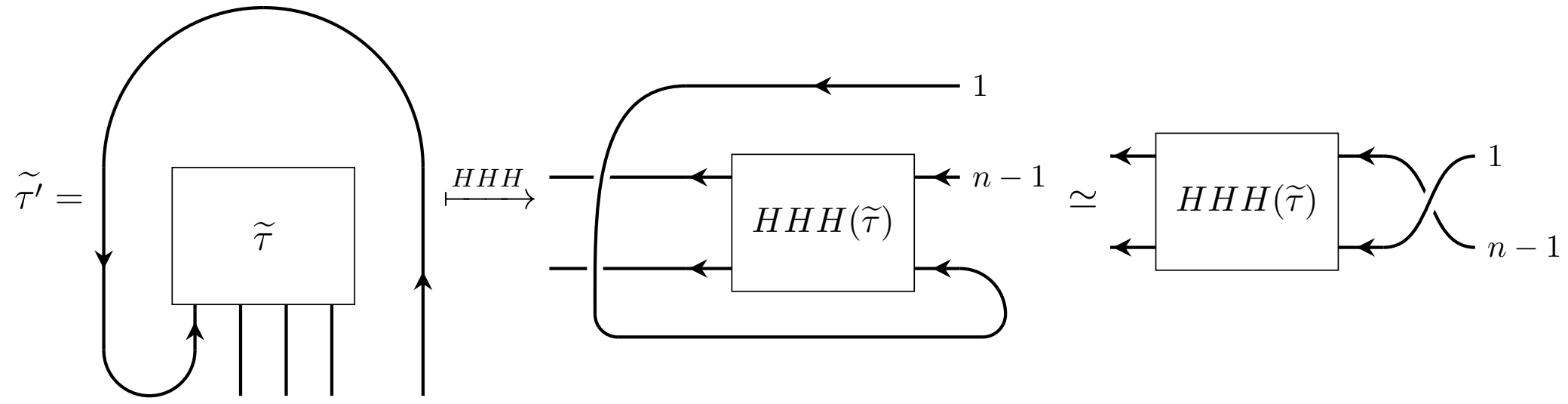
    \begin{figure}[ht]
    \centering
    \scalebox{.8}{
    {\Large{$\tilde{\tau'} = $}}
    \begin{tikzpicture}[anchorbase]
        \draw[webs,->-=.1] (0,0) to (0,2);
        \draw[webs] (.5,-.5) to (.5,2);
        \draw[webs] (1,-.5) to (1,2);
        \draw[webs] (1.5,-.5) to (1.5,2);

        \filldraw[fill=white, draw=black]
        (-.25,.5) rectangle (1.75,2);
        \node at (.75,1.25) {\Large{$\tilde{\tau}$}};

        \draw[webs] (0,0) to[out=270,in=0] (-.5,-.5);
        \draw[webs] (-.5,-.5) to[out=180,in=270] (-1,0);

        \draw[webs,-<-=.5] (-1,0) to (-1,2);
        \draw[webs] (-1,2) to[out=90,in=180] (.75,3.75);
        \draw[webs] (.75,3.75) to[out=0,in=90] (2.5,2);
        \draw[webs,-<-=.5] (2.5,2) to (2.5,-.5);
    \end{tikzpicture}
    {\Large{$\xmapsto{HHH}$}}
    \begin{tikzpicture}[anchorbase]
        \draw[webs,biarr] (0,0) to (3,0);
        \draw[webs,biarr] (0,1) to (3,1) node[right]{$n-1$};
        
        \filldraw[fill=white, draw=black]
        (.5,-.25) rectangle (2.5,1.25);
        \node at (1.5,.5) {\large{$HHH(\tilde{\tau})$}};

        \draw[webs,-<-=.5] (0,2) to (3,2) node[right]{$1$};

        \draw[webs] (0,0) to (-1.5,0);
        \draw[webs] (0,1) to (-1.5,1);
        \draw[webs,over] (0,2) to[out=180,in=90] (-1,-.5);

        \draw[webs] (-1,-.5) to[out=270,in=180] (-.75,-.75);
        \draw[webs] (3,0) to[out=0,in=90] (3.5,-.5);
        \draw[webs] (3.5,-.5) to[out=270,in=0] (3.25,-.75);
        \draw[webs] (-.75,-.75) to (3.25,-.75);
    \end{tikzpicture}
    {\Large{$\simeq$}}
    \begin{tikzpicture}[anchorbase]
        \draw[webs,biarr] (0,0) to (3,0);
        \draw[webs,biarr] (0,1) to (3,1);
        
        \filldraw[fill=white, draw=black]
        (.5,-.25) rectangle (2.5,1.25);
        \node at (1.5,.5) {\large{$HHH(\tilde{\tau})$}};

        \draw[webs] (3,1) to[out=0,in=180] (4,0) node[right]{$n-1$};
        \draw[webs,over] (3,0) to[out=0,in=180] (4,1) node[right]{$1$};
    \end{tikzpicture}
    }
    \end{figure}
    \FloatBarrier

    Here the homotopy equivalence follows from Proposition \ref{prop: diagrammatics_isotopy}. On the other hand, if $p$ is oriented outward, we can compute $HHH(\tilde{\tau'})$ using Proposition \ref{prop: clockwise_cup} as follows:

    \FloatBarrier
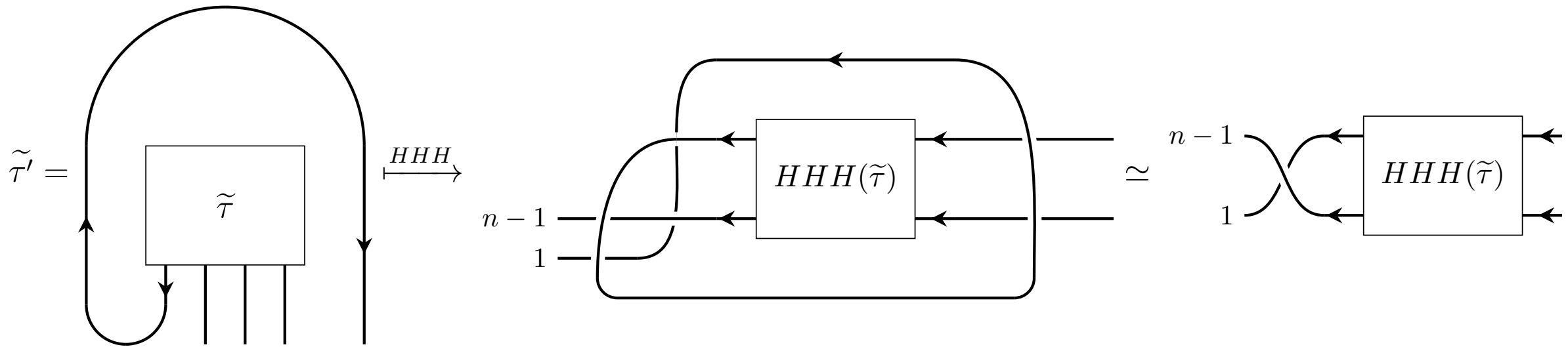
    \begin{figure}[ht]
    \centering
    \scalebox{.8}{
    {\Large{$\tilde{\tau'} = $}}
    \begin{tikzpicture}[anchorbase]
        \draw[webs,-<-=.1] (0,0) to (0,2);
        \draw[webs] (.5,-.5) to (.5,2);
        \draw[webs] (1,-.5) to (1,2);
        \draw[webs] (1.5,-.5) to (1.5,2);

        \filldraw[fill=white, draw=black]
        (-.25,.5) rectangle (1.75,2);
        \node at (.75,1.25) {\Large{$\tilde{\tau}$}};

        \draw[webs] (0,0) to[out=270,in=0] (-.5,-.5);
        \draw[webs] (-.5,-.5) to[out=180,in=270] (-1,0);

        \draw[webs,->-=.5] (-1,0) to (-1,2);
        \draw[webs] (-1,2) to[out=90,in=180] (.75,3.75);
        \draw[webs] (.75,3.75) to[out=0,in=90] (2.5,2);
        \draw[webs,->-=.5] (2.5,2) to (2.5,-.5);
    \end{tikzpicture}
    {\Large{$\xmapsto{HHH}$}}
    \begin{tikzpicture}[anchorbase]
        \draw[webs,biarr] (0,0) to (3,0);
        \draw[webs,biarr] (0,1) to (3,1);
        
        \filldraw[fill=white, draw=black]
        (.5,-.25) rectangle (2.5,1.25);
        \node at (1.5,.5) {\large{$HHH(\tilde{\tau})$}};

        \draw[webs,-<-=.5] (0,2) to (3,2);

        \draw[webs] (0,2) to[out=180,in=0] (-1,-.5);
        \draw[webs] (-1,-.5) to (-2,-.5) node[left]{$1$};
        \draw[webs,over] (0,1) to (-.5,1);
        \draw[webs,over] (0,0) to (-2,0) node[left]{$n-1$};
        \draw[webs,over] (-.5,1) to[out=180,in=90] (-1.5,-.75);

        \draw[webs] (3,0) to (5,0);
        \draw[webs] (3,1) to (5,1);
        \draw[webs,over] (3,2) to[out=0,in=90] (4,-.75);

        \draw[webs] (-1.5,-.75) to[out=270,in=180] (-1.25,-1);
        \draw[webs] (4,-.75) to[out=270,in=0] (3.75,-1);
        \draw[webs] (-1.25,-1) to (3.75,-1);
    \end{tikzpicture}
    {\Large{$\simeq$}}
    \begin{tikzpicture}[anchorbase]
        \draw[webs,biarr] (0,0) to (3,0);
        \draw[webs,biarr] (0,1) to (3,1);
        
        \filldraw[fill=white, draw=black]
        (.5,-.25) rectangle (2.5,1.25);
        \node at (1.5,.5) {\large{$HHH(\tilde{\tau})$}};

        \draw[webs] (-1,0) node[left]{$1$} to[out=0,in=180] (0,1);
        \draw[webs,over] (-1,1) node[left]{$n-1$} to[out=0,in=180] (0,0);
    \end{tikzpicture}
    }
    \end{figure}
    \FloatBarrier

    Here again the homotopy equivalence follows from Proposition \ref{prop: diagrammatics_isotopy}.
\end{proof}

\begin{example} \label{ex: trivial_tangle_in_disk}
    Consider the oriented tangle $\tau$ given by two parallel upwards-oriented vertical strands embedded in a disk. If we cut the boundary circle to the right of the strands to produce a $(0,4)$-tangle $\tilde{\tau}$ consisting of two concentric counterclockwise cups, then $HHH(\tilde{\tau}) = \strand{2}$ by Proposition \ref{prop: ccw_cup}. On the other hand, if we cut the boundary circle to the left of the strands to produce a $(0,4)$-tangle $\tilde{\tau'}$ consisting of two concentric \textit{clockwise} cups, then $HHH(\tilde{\tau'}) = FT_2$ by Proposition \ref{prop: clockwise_cup}. See Figure \ref{fig: parallel_strands_disk_tangle} below.

    \FloatBarrier
    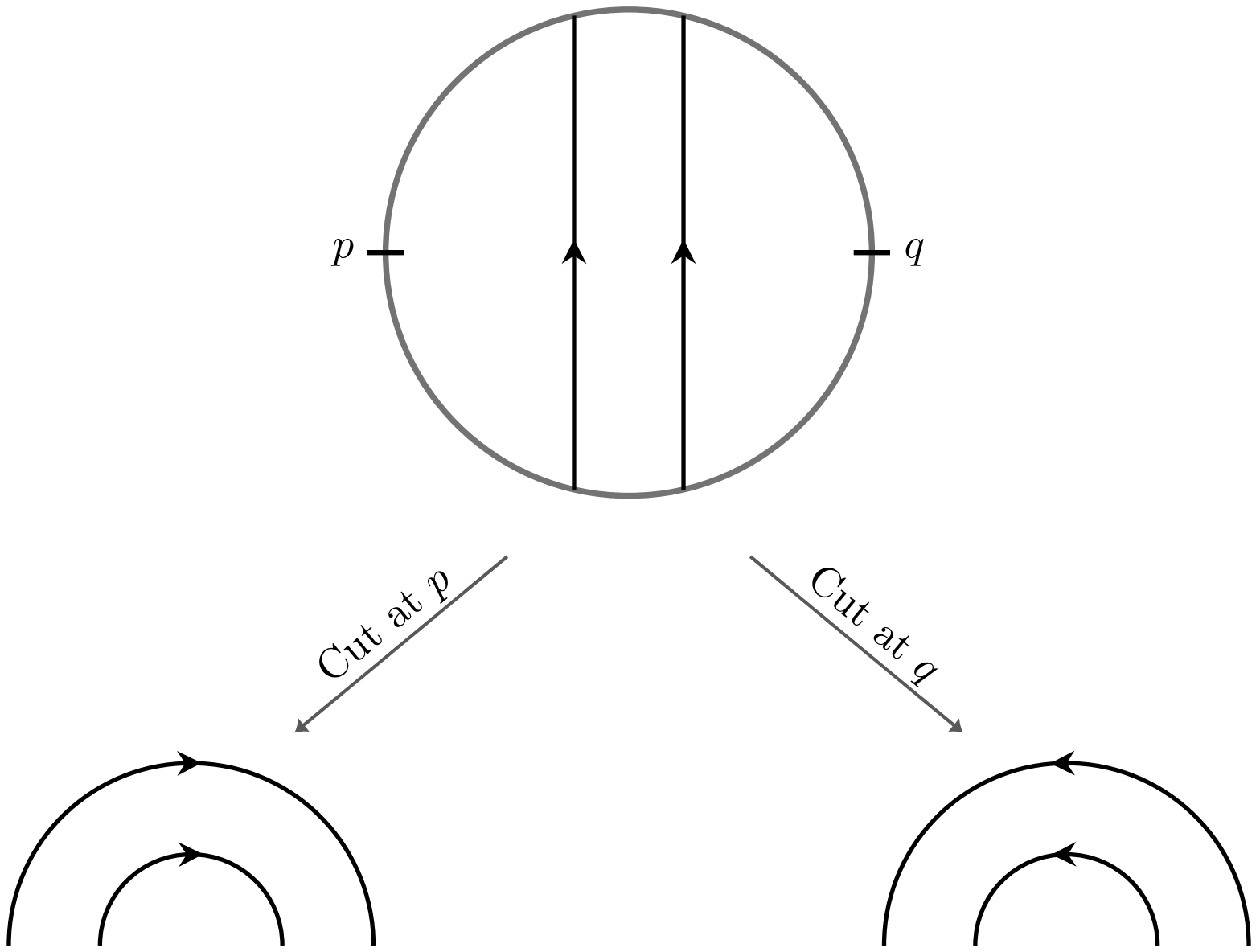
\begin{figure}[ht]
    \centering

    \begin{tikzpicture}[scale=1]
        \coordinate (top)   at (0,3.4);
        \coordinate (left)  at (-3.6,-2.3);
        \coordinate (right) at ( 3.6,-2.3);


        \draw[->, thick, draw=black!65]
            (-1.0,0.9)
            to
            node[midway, above, sloped, fill=white, inner sep=1pt] {Cut at $p$}
            (-2.75,-0.55);

        \draw[->, thick, draw=black!65]
            (1.0,0.9)
            to
            node[midway, above, sloped, fill=white, inner sep=1pt] {Cut at $q$}
            (2.75,-0.55);

        \begin{scope}[shift={(top)}]
            \draw[line width=1.4pt, draw=black!55] (0,0) circle (2);

            \draw[line width=1.2pt] (180:1.85) -- (180:2.15);
            \node at (180:2.35) {$p$};

            \draw[line width=1.2pt] (0:1.85) -- (0:2.15);
            \node at (0:2.35) {$q$};

            \draw[webs,->-=.5] (-0.45,-1.95) -- (-0.45,1.95);
            \draw[webs,->-=.5] ( 0.45,-1.95) -- ( 0.45,1.95);
        \end{scope}

        
        \begin{scope}[shift={(right)}]
            \draw[webs,-<-=.99] (-.75,0) to[out=90,in=180] (0,.75);
            \draw[webs]         (0,.75)  to[out=0,in=90]   (.75,0);

            \draw[webs,-<-=.99] (-1.5,0) to[out=90,in=180] (0,1.5);
            \draw[webs]         (0,1.5)  to[out=0,in=90]   (1.5,0);
        \end{scope}


        \begin{scope}[shift={(left)}]
            \draw[webs,->-=.99] (-.75,0) to[out=90,in=180] (0,.75);
            \draw[webs]         (0,.75)  to[out=0,in=90]   (.75,0);

            \draw[webs,->-=.99] (-1.5,0) to[out=90,in=180] (0,1.5);
            \draw[webs]         (0,1.5)  to[out=0,in=90]   (1.5,0);
        \end{scope}
    \end{tikzpicture}
    
    \caption{Two parallel upwards strands embedded in a disk. $HHH$ of the two oriented $(0,4)$-tangles resulting from different boundary cuts differ by a full twist.}
    \label{fig: parallel_strands_disk_tangle}
    \end{figure}
    \FloatBarrier
\end{example}

We can use Proposition \ref{prop: tangles_in_disks} to deduce various symmetries of $HHH$. For example, given any oriented tangle $\tau$ in the disk, we can consider the effect of rotating the marked point on the boundary circle by a complete rotation. Topologically, this operation preserves the underlying $(0,2n)$-tangle; algebraically, it corresponds to conjugating $HHH(\tau)$ by a full twist. This immediately gives the following:

\begin{cor} \label{cor: ft_central}
    For any oriented tangle $\tau$ with $2n$ endpoints, we have a natural homotopy equivalence

    \[
    HHH(\tau) \simeq FT_n \star HHH(\tau) \star FT_n^{-1} \in K^b(\extsbim_n).
    \]
\end{cor}

\FloatBarrier

\subsection{Coherence} \label{subsec: coherence}

Our next task is to justify the agreement of our topological construction of $HHH$ from Definition \ref{def: tangle_hhh} with the algebraic construction of $HHH_{alg}$ from Definition \ref{def: hhh_algebraic}. In fact, we prove a much stronger statement. Recall that $HHH_{alg}(\L)$ is defined using a presentation of $\L$ as a braid closure, which is a special type of coherent tangle diagram (see Definition \ref{def: coherent_tangle_diagram}). Any coherent tangle diagram $\tau$ with $2n$ endpoints encodes a complex $C \in K^b(\extsbim_n)$, and we can therefore make a systematic comparison between $C$ and the complex $HHH(\tau) \in K^b(\extsbim_n)$.

\begin{theorem} \label{thm: coherence}
    Let $\tau$ be a coherent tangle diagram encoding a complex $C \in K^b(\extsbim)$. Then $HHH(\tau) \simeq C$.
\end{theorem}

We proceed by showing that the topological operations which produce coherent tangle diagrams recover the corresponding algebraic operations on $K^b(\extsbim)$ after applying $HHH$. We begin with the inclusion of a strand:

\begin{prop} \label{prop: hhh_inclusion}
    Let $\tau$ be an oriented $(m, n)$-tangle with only positively oriented endpoints, and consider the oriented $(m+1, n+1)$-tangle $\tau'$ obtained by adding a single upwards-oriented vertical strand to the left of $\tau$. Then $HHH(\tau') = I(HHH(\tau))$.
\end{prop}

\begin{proof}
    Following the prescription of Figure \ref{fig: bending_tangles}, the oriented $(0, m+n+2)$-tangle $\tilde{\tau'}$ given by bending $\tau'$ around to the left is exactly the vertical concatenation of the oriented $(0, m+n)$-tangle $\tilde{\tau}$ with the oriented $(m+n, m+n+2)$-tangle $I^{\frac{1}{2}(m+n)}(\bigcap)$. By Corollary \ref{cor: concentric_cups} and functoriality of $\Phi$ with respect to vertical concatenation, we have 
    \[
    HHH(\tau') = \Phi\left(I^{\frac{1}{2}(m+n)}\left(\bigcap\right)\right)(HHH(\tau)) = I(HHH(\tau)).
    \]
\end{proof}

Next, we turn to partial trace:

\begin{prop} \label{prop: hhh_partial_closure}
    Let $\tau$ be an oriented $(m, n)$-tangle with only positively oriented endpoints, and consider the oriented $(m-1, n-1)$-tangle $\tau'$ obtained by closing the leftmost top and bottom endpoints of $\tau$. Then $HHH(\tau') = Tr(HHH(\tau))$.
\end{prop}

\begin{proof}
    The oriented $(0, m+n-2)$-tangle $\tilde{\tau'}$ given by bending $\tau'$ around to the left is exactly the vertical concatenation of the oriented $(0, m+n)$-tangle $\tilde{\tau}$ with the oriented $(m+n, m+n-2)$-tangle $I^{\frac{1}{2}(m+n-2)}(\bigcup)$. By Lemma \ref{lem: concentric_caps} and functoriality of $\Phi$ with respect to vertical concatenation, we have 
    \[
    HHH(\tau') = \Phi\left(I^{\frac{1}{2}(m+n-2)}\left(\bigcup\right)\right)(HHH(\tau)) = Tr(HHH(\tau)).
    \]
\end{proof}

Finally, we turn to braiding:

\begin{prop} \label{prop: hhh_braiding}
    Let $\tau$ be an oriented $(m,n)$-tangle with only positively oriented endpoints, and consider the oriented $(m,n)$-tangle $\tau'$ obtained by adding a single crossing $\sigma_i^{\pm}$ between the $i$ and $i+1^{st}$ strands along the top (resp. bottom) boundary of $\tau$. Then $HHH(\tau') \simeq F(\sigma_i^{\pm}) \star HHH(\tau)$ (resp. $HHH(\tau) \star F(\sigma_i^{\pm})$.
\end{prop}

\begin{proof}
    Consider the oriented $(0,m+n)$-tangle $\tilde{\tau'}$. If $\sigma_i^{\pm}$ was added below $\tau$, then $\tilde{\tau'}$ is exactly the vertical concatenation of $\tilde{\tau}$ with the vertical crossing $\sigma_i^{\pm}$. Otherwise, if $\sigma_i^{\pm}$ was added above $\tau$, we can pull this crossing around the new bends to become a vertical crossing between downwards-oriented strands, so $\tilde{\tau'}$ is isotopic to the concatenation of $\tilde{\tau}$ with a vertical crossing. In either case, the desired expression for $HHH(\tau')$ follows from functoriality of $\Phi$ and the action of vertical crossings from \S \ref{subsubsec: crossings}.
\end{proof}

\begin{proof}[Proof of Theorem \ref{thm: coherence}]
    Coherent tangle diagrams are by definition built from the empty tangle by iteratively adding strands on the left (using the inclusion functor $I$ of Definition \ref{def: adding_strand}), adding vertical crossings between adjacent upwards-oriented strands (using the Rouquier complex), and connecting the leftmost top and bottom endpoints of upwards-oriented strands (using the partial trace functor $Tr$ of Definition \ref{def: partial_trace}). Notice that any tangle diagram constructed in this way must have only positively oriented endpoints. By Propositions \ref{prop: hhh_inclusion}, \ref{prop: hhh_partial_closure}, and \ref{prop: hhh_braiding}, each of these topological operations on tangles is sent to the corresponding algebraic operation on $K^b(\extsbim)$.
\end{proof}

As special cases of Theorem \ref{thm: coherence}, we immediately recover the following:

\begin{cor} \label{cor: hhh_braids}
    Let $\beta$ be an upwards-oriented braid on $n$ strands, considered as an oriented $(n, n)$-tangle. Then $HHH(\beta) = F(\beta) \in K^b(\extsbim_n)$.
\end{cor}

\begin{cor} \label{cor: hhh_partial_braid_closure}
    Let $\tau$ be the oriented $(n-m,n-m)$-tangle obtained from closing the leftmost $m$ strands of an upwards-oriented braid $\beta$ on $n$ strands. Then $HHH(\tau) = Tr^m(F(\beta))$.
\end{cor}

\begin{cor} \label{cor: hhh_links}
    Let $\mathcal{L}$ be an oriented link, considered as a $(0,0)$-tangle. Then $HHH(\mathcal{L}) \cong HHH_{alg}(\L)$.
\end{cor}

\begin{proof}
    Pick a presentation of $\mathcal{L}$ as a braid closure, then apply Corollary \ref{cor: hhh_partial_braid_closure}.
\end{proof}

\subsection{Locality} \label{subsec: locality}

One of the primary benefits of promoting a link invariant to a tangle invariant is that this often enables \textit{local} computations. This provides (at least) a few huge advantages. First, actually computing the invariant assigned to any oriented link becomes vastly easier and faster when one can make local simplifications. The previous absence of locality has been a huge obstruction to any fast, systematic computation of HOMFLY homology for large knots to date. As one example, the local computations pioneered by Bar-Natan (and optimized by Jeremy Green) in \cite{BN07} could compute the Khovanov homology of the 48-crossing knot $T(8,7)$ in a few minutes (in 2005!). Compare these with the recent state of the art in HOMFLY homology computations by Nakagane and Sano \cite{NS25}, who in 2024 could only compute the HOMFLY homology of prime knots up to 11 crossings. Second, locality is a crucial tool in proving any sort of functoriality; see e.g. \cite{BN05} for the first application of this principle to Khovanov homology, \cite{ETW18} for an extension to colored $\mathfrak{gl}_N$ homology, and \cite{EK10a} for an application to HOMFLY homology restricted to certain cobordisms between braid diagrams.

At first glance, our construction of $HHH$ in Definition \ref{def: tangle_hhh} is inherently nonlocal; computing it for a given oriented tangle $\tau$ requires first bending over to produce an oriented $(0,2n)$-tangle $\tilde{\tau}$, then computing $\Phi$ of each local piece, and finally composing to obtain $HHH(\tau)$. On the other hand, our proof of Theorem \ref{thm: coherence} essentially states that $HHH$ of a very special class of oriented tangle diagrams can be computed directly from $\tau$ by first computing $HHH$ locally, then assembling via tensor products, inclusions, and partial traces. This suggests that our oriented tangle invariant might be amenable to local simplifications more generally.

Our main goal in this section is to show that our invariant can indeed be computed locally. We begin by observing a straightforward, but extremely useful, consequence of Theorem \ref{thm: coherence}:

\begin{prop} \label{prop: diagrammatics_isotopy}
    Let $\tau, \tau'$ be any coherent tangle diagrams encoding complexes $C, C' \in K^b(\extsbim)$. If $\tau$ and $\tau'$ are isotopic as oriented tangles, then there is a natural equivalence $C \simeq C'$.
\end{prop}

\begin{proof}
    By Theorem \ref{thm: coherence}, we have $HHH(\tau) = C$ and $HHH(\tau') = C'$. Any isotopy of oriented tangles realizing $\tau \sim \tau'$ in particular gives an isotopy of oriented tangles $\tilde{\tau} \sim \tilde{\tau'}$. Then by Theorem \ref{thm: main_invariance}, we obtain

    \[
    C = HHH(\tau) = \Phi(\tilde{\tau})(\mathbb{Z}) \simeq \Phi(\tilde{\tau'})(\mathbb{Z}) = HHH(\tau') = C'.
    \]
\end{proof}

Proposition \ref{prop: diagrammatics_isotopy} tells us that we can extend the usual diagrammatic calculus to treat \textit{any} oriented tangle diagram $\tau$ with positively oriented endpoints as encoding $HHH(\tau)$. Indeed, the prescription of Figure \ref{fig: general_framework} tells us that we can compute $HHH$ of any tangle by repeatedly applying isotopies to each elementary tangle composing $\tau$ until we arrive at a coherent tangle diagram $\tau^{coh}$ encoding $HHH(\tau)$ composed with a shuffle braid. When $\tau$ has only positively oriented endpoints, this shuffle braid will be trivial, so $\tau$ is isotopic to the coherent tangle diagram $\tau^{coh}$.

This perspective immediately reveals our main locality theorem:

\begin{theorem} \label{thm: upwards_locality}
    Let $\tau$ and $\tau'$ be oriented tangles with only positively oriented endpoints. Then we have a natural equivalence

    \[
    HHH(\tau \sqcup \tau') \simeq HHH(\tau) \boxtimes HHH(\tau').
    \]

    Further, if $\tau$ is an oriented $(\aa, \bb)$-tangle and $\tau'$ is an oriented $(\bb, \cc)$-tangle, then we have a natural equivalence

    \[
    HHH(\tau \star \tau') \simeq HHH(\tau) \star HHH(\tau').
    \]
\end{theorem}

\begin{proof}
    In either case, we can isotope $\tau$ and $\tau'$ to obtain coherent tangle diagrams $\tau^{coh}$ and $\tau'^{coh}$. Now both $\tau^{coh} \sqcup \tau'^{coh}$ and $\tau^{coh} \star \tau'^{coh}$ are coherent tangle diagrams isotopic to $\tau \sqcup \tau'$ and $\tau \star \tau'$, respectively\footnote{Recall that isotopies of oriented tangle diagrams fix the boundary, so isotopies $\tau \sim \tau^{coh}$ and $\tau' \sim \tau'^{coh}$ automatically glue to give an isotopy $\tau \star \tau' \sim \tau^{coh} \star \tau'^{coh}$.}, and we conclude by Theorem \ref{thm: coherence}.
\end{proof}

We also have locality for tangles with arbitrary boundary, though some care must be taken to account for orientations. We begin by treating locality under horizontal composition $\sqcup$.  For this, we will make use of a generalization of Proposition \ref{prop: monoidality}.

\begin{definition}
    Given $a, b \geq 0$, we denote by $Sh_{a, b}$ a cabled positive crossing of $a$ upwards-oriented strands over $b$ upwards-oriented strands; see Figure \ref{fig: shuffle_braid}.

    \begin{figure}[ht]
    \centering

    {\Large{$Sh_{a,b} :=$}}
    \begin{tikzpicture}[anchorbase]
        \draw[webs,->-=.2,->-=.8] (1.5,0) node[below]{$b$} to[out=90,in=270] (0,1.5) node[above]{$b$};
        \draw[webs,->-=.2,->-=.8,over] (0,0) node[below]{$a$} to[out=90,in=270] (1.5,1.5) node[above]{$a$};
    \end{tikzpicture}

    \caption{A positive shuffle braid $Sh_{a,b}$.}
    \label{fig: shuffle_braid}
    \end{figure}
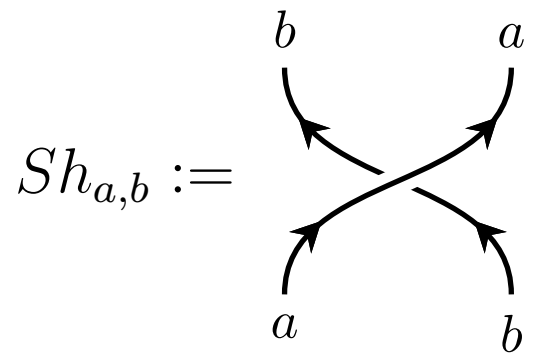
    \FloatBarrier
\end{definition}

\begin{prop} \label{prop: general_monoidality}
    Let $\tau$ be an oriented $(0,2n)$-tangle, and consider the oriented tangle $\tau'$ obtained from $\tau$ by adding $i$ downwards- and $j$ upwards-oriented vertical strands to the left of $\tau$ and $k$ downwards- and $l$ upwards-oriented vertical strands to the right of $\tau$. Assume $i + k = j + l$, so that $\tau'$ is balanced. Then for any complex $C \in K^b(\extsbim_{i + k})$, we have

    \[
    \Phi(\tau')(C) = (F(Sh_{n, k}) \boxtimes \strand{i}) \star (HHH(\tau) \boxtimes C) \star (F(Sh_{j, n}) \boxtimes \strand{l}).
    \]
\end{prop}

\begin{proof}
    Recalling the discussion of Warning \ref{warning: non_local}, we may safely assume that all of the downwards-oriented strands added to the left and right of $\tau$ live to the left of their upwards-oriented counterparts. With this convention in mind, we begin by applying the following isotopy:

    \begin{figure}[ht]
    \centering

    \begin{tikzpicture}[anchorbase]
        \draw[webs,-<-=.5] (0,0) to (0,1.5) node[above]{$i$};
        \draw[webs,->-=.5] (.5,0) to (.5,1.5) node[above]{$j$};
        \draw[webs,-<-=.5] (2.5,0) to (2.5,1.5) node[above]{$k$};
        \draw[webs,->-=.5] (3,0) to (3,1.5) node[above]{$l$};

        \draw[webs] (1.5,0) node[below]{$2n$} to (1.5,.75);

        \filldraw[fill=white, draw=black]
        (1,.5) rectangle (2,1);
        \node at (1.5,.75) {\Large{$\tau$}};
    \end{tikzpicture}
    \ {\large{$\sim$}} \ 
    \begin{tikzpicture}[anchorbase]
        \draw[webs,-<-=.5] (0,0) to (0,1.5);
        \draw[webs,-<-=.5] (.5,0) to (.5,1.5);
        \draw[webs,->-=.5] (2.5,0) to (2.5,1.5);
        \draw[webs,->-=.5] (3,0) to (3,1.5);

        \draw[webs] (1.5,0) to (1.5,.75);

        \filldraw[fill=white, draw=black]
        (1,.5) rectangle (2,1);
        \node at (1.5,.75) {\Large{$\tau$}};

        \draw[webs,-<-=.8] (0,1.5) to (0,3) node[above]{$i$};
        \draw[webs,-<-=.8] (.5,1.5) to[out=90,in=270] (2.5,3) node[above]{$k$};
        \draw[webs,->-=.8, over] (2.5,1.5) to[out=90,in=270] (.5,3) node[above]{$j$};
        \draw[webs,->-=.8] (3,1.5) to (3,3) node[above]{$l$};

        \draw[color=blue, dashed] (-.5,0) to (3.5,0);

        \draw[webs,->-=.8] (0,0) to (0,-1.5) node[below]{$i$};
        \draw[webs,->-=.8] (.5,0) to[out=270,in=90] (2.5,-1.5) node[below]{$k$};
        \draw[webs,over] (1.5,0) to (1.5,-1.5) node[below]{$2n$};
        \draw[webs,over] (2.5,0) to[out=270,in=90] (.5,-1);
        \draw[webs,-<-=.5] (.5,-1) to (.5,-1.5) node[below]{$j$};
        \draw[webs,-<-=.8] (3,0) to (3,-1.5) node[below]{$l$};
    \end{tikzpicture}
    \end{figure}

    Many of the additional crossings introduced in this isotopy, including in particular all of the crossings occurring above the dashed blue line, are horizontal crossings in which an upwards-oriented strand crosses over a downwards-oriented strand. These act trivially by Lemma \ref{lem: trivial_horizontals}. Observing that the bottom boundary of $\tau$ must consist of $n$ upwards- and $n$ downwards-oriented strands, the only nontrivial crossings are therefore a vertical cabled crossing of $j$ upwards-oriented strands over the $n$ upwards-oriented strands exiting $\tau$ and a vertical cabled crossing of the $n$ downwards-oriented strands exiting $\tau$ over $k$ downwards-oriented strands. The portion of the tangle above the dashed blue line then acts by $HHH(\tau) \boxtimes (-)$, while the remaining portion acts exactly by tensoring with the Rouquier complexes in the statement of the proposition.
\end{proof}

\begin{theorem} \label{thm: hhh_concatenation_general}

Let $\tau$ be an oriented $(\aa, \bb)$-tangle with $2m$ endpoints and $\tau'$ an oriented $(\cc, \dd)$-tangle. Suppose $\cc$ has $i$ positive and $j$ negative entries, while $\dd$ has $k$ positive and $l$ negative entries. Then

\[
HHH(\tau \sqcup \tau') \simeq (F(Sh_{m,l}) \boxtimes \mathbb{1}_{i}) \star (HHH(\tau) \boxtimes HHH(\tau')) \star (F(Sh_{j,m}) \boxtimes \mathbb{1}_k).
\]
\end{theorem}

\begin{proof}
    We begin by applying the following isotopy to the oriented $(\emptyset, r(\cc) \sqcup r(\aa) \sqcup \bb \sqcup \dd)$-tangle $\tilde{(\tau \sqcup \tau')}$:

    \begin{figure}[ht]
    \centering
    
    {\Large{$\tilde{(\tau \sqcup \tau')} =$}}
    \begin{tikzpicture}[anchorbase]
        \draw[webs] (.5,0) node[below]{$\bb$} to (.5,1.5);
        \draw[webs] (1.5,0) node[below]{$\dd$} to (1.5,1.5);

        \filldraw[fill=white, draw=black]
        (.1,.5) rectangle (.9,1.2);
        \node at (.5,.85) {\Large{$\tau$}};
        \filldraw[fill=white, draw=black]
        (1.1,.5) rectangle (1.9,1.2);
        \node at (1.5,.85) {\Large{$\tau'$}};

        \draw[webs] (.5,1.5) to[out=90,in=0] (0,2);
        \draw[webs] (0,2) to[out=180,in=90] (-.5,1.5);
        \draw[webs] (1.5,1.5) to[out=90,in=0] (0,2.5);
        \draw[webs] (0,2.5) to[out=180,in=90] (-1.5,1.5);

        \draw[webs] (-.5,1.5) to (-.5,0) node[below]{$r(\aa)$};
        \draw[webs] (-1.5,1.5) to (-1.5,0) node[below]{$r(\cc)$};
    \end{tikzpicture}
    \ {\Large{$\sim$}} \ 
    \begin{tikzpicture}[anchorbase]
        \draw[webs] (.5,0) node[below]{$\bb$} to (.5,1.5);

        \filldraw[fill=white, draw=black]
        (.1,.5) rectangle (.9,1.2);
        \node at (.5,.85) {\Large{$\tau$}};

        \draw[webs] (.5,1.5) to[out=90,in=0] (0,2);
        \draw[webs] (0,2) to[out=180,in=90] (-.5,1.5);

        \draw[webs] (-.5,1.5) to (-.5,0) node[below]{$r(\aa)$};

        \draw[webs] (1.5,0) node[below]{$\dd$} to (1.5,3.5);
        
        \filldraw[fill=white, draw=black]
        (1.1,2.5) rectangle (1.9,3.2);
        \node at (1.5,2.85) {\Large{$\tau'$}};
        
        \draw[webs] (1.5,3.5) to[out=90,in=0] (0,4.5);
        \draw[webs] (0,4.5) to[out=180,in=90] (-1.5,3.5);
        
        \draw[webs] (-1.5,3.5) to (-1.5,0) node[below]{$r(\cc)$};

        \draw[color=blue, dashed] (-2,2.25) to (2.5,2.25);
    \end{tikzpicture}
    \end{figure}

    Now, the portion of the right-hand oriented tangle above the dashed blue line is exactly $\tilde{\tau'}$. Meanwhile, the portion below is the oriented tangle obtained by adding $i$ downwards- and $j$ upwards-oriented strands to the left of $\tilde{\tau}$ and $l$ downwards- and $k$ upwards-oriented strands to the right of $\tilde{\tau}$. Applying $\Phi$, we first send $\mathbb{Z}$ to $HHH(\tau')$, then (using Proposition \ref{prop: general_monoidality}) to $(F(Sh_{m,l}) \boxtimes \mathbb{1}_{i}) \star (HHH(\tau) \boxtimes HHH(\tau')) \star (F(Sh_{j,m}) \boxtimes \mathbb{1}_k)$.
\end{proof}

As a special case of Theorem \ref{thm: hhh_concatenation_general}, we obtain a slight strengthening of Theorem \ref{thm: upwards_locality} under horizontal composition:

\begin{cor} \label{cor: hhh_concatentation_upwards}
    Let $\tau, \tau'$ be two oriented tangles, and suppose $\tau'$ has only positively oriented endpoints. Then

    \[
    HHH(\tau \sqcup \tau') \simeq HHH(\tau) \boxtimes HHH(\tau').
    \]
\end{cor}

\begin{proof}
    Observe that the shuffle braids in Theorem \ref{thm: hhh_concatenation_general} vanish.
\end{proof}

\begin{rem} \label{rem: canopolis}
By similar computations, one can give similar statements for the behavior of $HHH$ under arbitrary combinations of vertical and horizontal concatenations for tangles with arbitrarily oriented endpoints. One perspective on this locality is that $HHH$ is a \textit{planar algebra homomorphism} between the collection of oriented tangle diagrams and the category $K^b(\extsbim)$ equipped with a certain planar algebra structure in the sense of \cite{Jon21}. We postpone a detailed description of this structure for future work.
\end{rem}

\appendix
\section{Markov Theorem for Oriented Tangles} \label{app: markov}
In this Appendix we provide a topological perspective on our construction of $HHH$ by distilling from the main body of the paper an Alexander-style and Markov-style theorem for arbitrary oriented tangles\footnote{As opposed to the usual Alexander and Markov theorems for oriented links.}, giving a regular form for oriented tangles in terms of both coherent tangles and partial braid closures and a list of moves relating these. To the best of our knowledge, these are new results.

We begin with the version of these results involving coherent tangles:

\begin{theorem} \label{thm: markov_tangles_coherent}

    Any oriented tangle $\tau$ is isotopic to the product of a coherent tangle $\tau^{coh}$ and a shuffle braid (depending only on the orientations of endpoints of $\tau$ and choice of marked point in the boundary of the disk), and two coherent tangles $\tau^{coh}$ and $\tau'^{coh}$ are isotopic if and only if they are related by a series of the following moves:

    \begin{itemize}
        
        \item Inclusion locality\footnote{This is implicit in the monoidal structure of $\otan$.}: $\sigma_{i+1}^{\pm} \star I(-) \leftrightarrow I(\sigma_i^{\pm} \star -)$ or $I(-) \star \sigma_{i+1}^{\pm} \leftrightarrow I(- \star \sigma_i^{\pm})$;
        
        \item Trace locality\footnote{In the sense of Figure \ref{fig: trace_local}.}: $\sigma_i^{\pm} \star Tr(-) \leftrightarrow Tr(\sigma_{i+1}^{\pm} \star -)$ or $Tr(-) \star \sigma_i^{\pm} \leftrightarrow Tr(- \star \sigma_{i+1}^{\pm})$;
        
        \item Braid relations;

        \item Leftmost conjugation: $Tr^2(\sigma_1^{\pm} \star -) \leftrightarrow Tr^2(- \star \sigma_1^{\pm})$;

        \item Identity relation: $\sigma_1^{\pm} \star I^2(-) \leftrightarrow I^2(-) \star \sigma_1^{\pm}$;

        \item Stabilization: $\emptyset \leftrightarrow Tr(I(-) \star \sigma_1^{\pm})$;

        \item Exchange moves: $Tr(\sigma_1^{\pm} \star I(-) \star \sigma_1^{\mp}) \leftrightarrow I(Tr(-))$.
    \end{itemize}
\end{theorem}

\begin{proof}
    It is clear that each of these moves preserves the isotopy class of the resulting tangle, so it suffices to prove the converse. 
    
    Overall, the idea is this: the procedure described in Sections \ref{subsec: general_framework} and \ref{subsec: actions} scans over the tangle $\tau$, and at each cup, cap, or crossing updates a ``register complex" of Soergel bimodules. Each update is by a combination of the operations $\sigma_i \star -$, $- \star \sigma_i$, $I$, or $Tr$, however, so we can just as well scan over $\tau$ and instead update a ``register coherent diagram" in the same way -- in the end, we produce a coherent diagram $\tau^{coh}$, whose product with a shuffle braid is isotopic to $\tau$. This proves the existence portion of the result.

    A key aspect of this procedure is that, if $\tau$ is already coherent, then $\tau^{coh} = \tau$. In general, updating the register coherent diagram following a cup or a cap requires introducing extra crossings, but a coherent diagram is precisely one for which all of the cups and caps give just $I$ or $Tr$, with no extra crossings needed and the diagram remaining unchanged.

    So, to get the completeness of these moves, consider a series of tangles $\tau_0,\tau_1,\dots,\tau_n$ differing by a single Reidemeister move at each step, with $\tau_0$ and $\tau_n$ coherent but the intermediate ones arbitrary. Running the existence part of the proof at each step, we obtain a series of coherent tangles $\tau_i^{coh}$, with $\tau_i^{coh} = \tau_i$ for $i = 0,n$. The arguments of Section \ref{subsec: invariance} now show that, since $\tau_i$ and $\tau_{i+1}$ differ by a Reidemeister move, $\tau_i^{coh}$ and $\tau_{i+1}^{coh}$ differ by a series of the moves in the theorem statement -- what is directly shown is that the corresponding complexes of Soergel bimodules are homotopy equivalent, but the only properties of the action of $\sigma_i \star -$, $- \star \sigma_i$, $I$, and $Tr$ on $\extsbim$ used in the proof are this list of moves.
\end{proof}

\begin{remark}
   The shuffle braids are a bit vestigial in the statement of Theorem~\ref{thm: markov_tangles_coherent}, but are absolutely essential to its proof -- even when the tangles involved all have only positive boundary points, the intermediate tangles during the scan do not. 
\end{remark}

Seeing as coherent tangle diagrams are relatively non-standard constructions in the existing literature, we conclude by comparing to the more familiar class of partial braid closures. The latter can be distinguished within the collection of all coherent tangle diagrams using the encoding as a string of operations discussed above; they correspond exactly to strings of the form $Tr^m(\beta \star I^n(\emptyset))$ (abbreviated in the sequel as $Tr^m(\beta)$) for some $m \leq n$ and $n$-strand braid word $\beta$. Up to isotopy, no generality is lost in making this restriction, as we now show:

\begin{lemma} \label{lem: coherents_are_partial_closures}
    Every coherent tangle diagram is isotopic to a partial braid closure.
\end{lemma}

\begin{proof}
    We proceed by induction on the length $\ell(\tau)$ of a string constructing an arbitrary coherent diagram $\tau$. If $\ell(\tau) = 1$, then we must have $\tau = I(\emptyset)$, which is already a partial braid closure. Otherwise, $\tau$ must be built from some coherent tangle diagram $\tau'$ by applying the operations $I$, $Tr$, $\sigma_i^{\pm} \star -$, or $- \star \sigma_i^{\pm}$. By assumption, $\tau'$ is itself isotopic to a partial braid closure $Tr^m(\beta)$ for some non-negative integers $m \leq n$ and some $n$-strand braid $\beta$. We now treat each operation separately:

    \begin{itemize}
        \item If $\tau = Tr(\tau') \sim Tr^{m+1}(\beta)$, then $\tau$ is already isotopic to a partial braid closure.

        \item If $\tau = I(\tau')$, then $\tau \sim Tr^m(\sigma_m \sigma_{m-1} \dots \sigma_1 \star (\strand{1} \sqcup \beta) \star \sigma_1^{-1} \dots \sigma_{m-1}^{-1} \sigma_m^{-1})$ by repeated applications of exchange moves and inclusion locality (see Figure \ref{fig: inclusion_regularization} below).

        \item If $\tau = \sigma_i^{\pm} \star \tau'$, then $\tau \sim Tr^m(\sigma_{i+m}^{\pm} \star \beta)$ by repeated applications of trace locality.

        \item If $\tau = \tau' \star \sigma_i^{\pm}$, then $\tau \sim Tr^m(\beta \star \sigma_{i+m}^{\pm})$ by repeated applications of trace locality and the identity relation.
    \end{itemize} 
\end{proof}

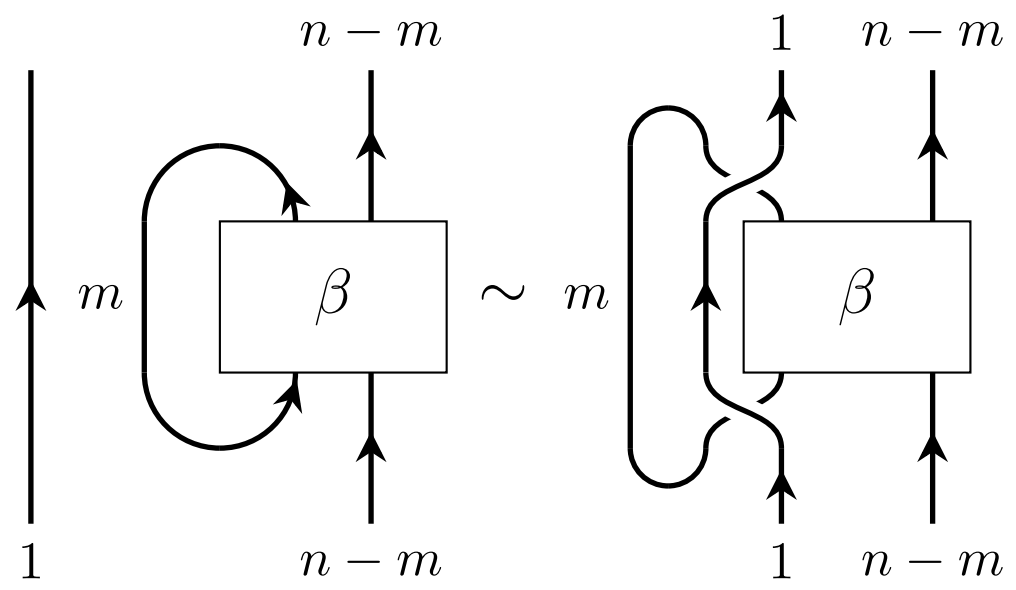
\begin{figure}[ht]
\centering

\begin{tikzpicture}[anchorbase]
\filldraw[fill=white, draw=black]
    (-.25,0) rectangle (1.25,1);
    \node at (.5,.5) {\Large{$\beta$}};

    \draw[webs,->-=.2] (.25,1) to[out=90,in=0] (-.25,1.5);
    \draw[webs] (-.25,1.5) to[out=180,in=90] (-.75,1);
    \draw[webs] (-.75,1) to node[midway,left]{$m$} (-.75,0);
    \draw[webs] (-.75,0) to[out=270,in=180] (-.25,-.5);
    \draw[webs,->-=.8] (-.25,-.5) to[out=0,in=270] (.25,0);

    \draw[webs,->-=.5] (.75,1) to (.75,2) node[above]{$n-m$};
    \draw[webs,->-=.5] (.75,-1) node[below]{$n-m$} to (.75,0);

    \draw[webs,->-=.5] (-1.5,-1) node[below]{$1$} to (-1.5,2);
\end{tikzpicture}
{\Large{$\sim$}}
\begin{tikzpicture}[anchorbase]
\filldraw[fill=white, draw=black]
    (-.25,0) rectangle (1.25,1);
    \node at (.5,.5) {\Large{$\beta$}};

    \draw[webs] (0,1) to[out=90,in=270] (-.5,1.5);
    \draw[webs,over] (-.5,1) to[out=90,in=270] (0,1.5);

    \draw[webs,->-=.5] (-.5,0) to (-.5,1);

    \draw[webs] (0,0) to[out=270,in=90] (-.5,-.5);
    \draw[webs,over] (-.5,0) to[out=270,in=90] (0,-.5);

    \draw[webs] (-.5,1.5) to[out=90,in=0] (-.75,1.75);
    \draw[webs] (-.75,1.75) to[out=180,in=90] (-1,1.5);
    \draw[webs] (-1,1.5) to node[midway,left]{$m$} (-1,-.5);
    \draw[webs] (-1,-.5) to[out=270,in=180] (-.75,-.75);
    \draw[webs] (-.75,-.75) to[out=0,in=270] (-.5,-.5);

    \draw[webs,->-=.5] (1,1) to (1,2) node[above]{$n-m$};
    \draw[webs,->-=.5] (1,-1) node[below]{$n-m$} to (1,0);

    \draw[webs,->-=.5] (0,1.5) to (0,2) node[above]{$1$};
    \draw[webs,->-=.5] (0,-1) node[below]{$1$} to (0,-.5);
\end{tikzpicture}

\caption{Applying $I(-)$ to a partial braid closure.}
\label{fig: inclusion_regularization}
\end{figure}

We pause to point out a few features of the proof of Lemma \ref{lem: coherents_are_partial_closures}. First, it is constructive: given a coherent tangle diagram $\tau$, it explicitly assigns non-negative integers $m \leq n$ and an $n$-strand braid $\beta$ satisfying $\tau \sim Tr^m(\beta)$. More precisely, each operation appearing in a string describing $\tau$ is assigned to a particular operation on partial braid closures, and composing these operations produces a partial braid closure which is isotopic to $\tau$.

We refer to this process as \textit{regularization} both at the level of diagrams and individual operations, so that the partial braid closure $Tr^m(\beta)$ is the regularization of the coherent diagram $\tau$. Note that the regularization of an individual operation appearing in a string describing $\tau$ depends only on the number of instances of the operation $Tr$ which have already appeared, and regularizing any string of operations preserves the number of instances of the operations $Tr$ and $I$. If $\tau$ happens to already be a partial braid closure, then the regularization of $\tau$ is exactly $\tau$.

Our final task in this appendix is to move one level higher, providing regularizations of the coherent Markov moves listed in Theorem \ref{thm: markov_tangles_coherent} and therefore a Markov-type theorem for partial braid closures.

\begin{theorem} \label{thm: markov_tangles_braid_closure}
    Let $\tau$ and $\tau'$ be two partial braid closures. Then $\tau$ and $\tau'$ are isotopic if and only if $\tau$ can be obtained from $\tau'$ by a finite sequence of moves of the following type:

    \begin{itemize}
        \item Braid relations;
        \item Conjugation through closed strands: $Tr^m(\sigma_i^{\pm} \star \beta) \leftrightarrow Tr^m(\beta \star \sigma_i^{\pm})$, where $i \leq m - 1$;
        \item Braid stabilization: $Tr^m(\beta \star \beta') \leftrightarrow Tr^{m+1}((\strand{1} \sqcup \beta) \star \sigma_1^{\pm} \star (\strand{1} \sqcup \beta'))$.
    \end{itemize}
\end{theorem}

\begin{proof}
    Let $\tau^{reg}$ and $\tau'^{reg}$ be two isotopic partial braid closures. Then in particular $\tau^{reg}$ and $\tau'^{reg}$ are coherent tangle diagrams, so by Theorem \ref{thm: markov_tangles_coherent}, they are related by a finite sequence of coherent Markov moves. Each coherent diagram appearing in this sequence has a regularization as a partial braid closure, and $\tau^{reg}$ and $\tau'^{reg}$ are preserved by regularization. It follows that the entire sequence of coherent Markov moves can be collapsed into a finite sequence of operations on partial braid closures, each of which is obtained by considering the regularization of the left- and right-hand sides of a single coherent Markov move.

    With this in mind, suppose $\tau$ and $\tau'$ are two coherent tangle diagrams related by a single coherent Markov move. By construction, we have $\tau = \varphi \psi \rho(\emptyset)$ and $\tau' = \varphi \psi' \rho(\emptyset)$ for some strings of operations $\varphi, \psi, \psi', \rho$, where $\psi$ and $\psi'$ appear on opposite sides of one of the coherent Markov moves listed in Theorem \ref{thm: markov_tangles_coherent}.

    Assume for the moment that $\varphi$ is empty, and let $Tr^m(\beta)$ be the regularization of $\rho(\emptyset)$. We analyze the effect of each coherent Markov move on the corresponding regularizations:

    \begin{itemize}
        \item \textit{Inclusion locality:} We analyze the first case, in which $\psi = \sigma_{i+1}^{\pm} \star I(-)$ and $\psi' =  I(\sigma_i^{\pm} \star -)$; the second follows similarly. The regularization of $\tau$ is then
        \begin{align*}
        \tau & = \sigma_{i+1}^{\pm} \star I(Tr^m(\beta)) \\
        & \sim \sigma_{i+1}^{\pm} \star Tr^m (\sigma_m \sigma_{m-1} \dots \sigma_1 \star (\strand{1} \sqcup \beta) \star \sigma_1^{-1} \dots \sigma_{m-1}^{-1} \sigma_m^{-1}) \\
        & \sim Tr^m (\sigma_{m + i + 1}^{\pm} \sigma_m \sigma_{m-1} \dots \sigma_1 \star (\strand{1} \sqcup \beta) \star \sigma_1^{-1} \dots \sigma_{m-1}^{-1} \sigma_m^{-1}),
        \end{align*}
        while the regularization of $\tau'$ is
        \begin{align*}
        \tau' & = I(\sigma_i^{\pm} \star Tr^m(\beta \star I^n(\emptyset))) \\
        & \sim I(Tr^m(\sigma_{i + m}^{\pm} \star \beta)) \\
        & \sim Tr^m (\sigma_m \sigma_{m-1} \dots \sigma_1 \sigma_{m + i + 1}^{\pm} \star (\strand{1} \sqcup \beta) \star \sigma_1^{-1} \dots \sigma_{m-1}^{-1} \sigma_m^{-1}).
        \end{align*}
        Up to the distant interchange braid relation $\sigma_m \sigma_{m-1} \dots \sigma_1 \sigma_{m + i + 1}^{\pm} \sim \sigma_{m + i + 1}^{\pm} \sigma_m \sigma_{m-1} \dots \sigma_1$, this move also has no effect after regularizing both sides.

        \item \textit{Trace locality:} We again analyze only the first case, in which $\psi = \sigma_i^{\pm} \star Tr(-)$ and $\psi' = Tr(\sigma_{i+1}^{\pm} \star -)$. In this case the regularization of both $\tau = \sigma_i^{\pm} \star Tr^{m+1} (\beta \star I^n(\emptyset))$ and $\tau' = Tr(\sigma_{i+1}^{\pm} \star Tr^m(\beta))$ is $Tr^{m+1}(\sigma_{i+m+1}^{\pm} \star \beta)$, so this move has no effect after regularizing both sides.

        \item \textit{Braid relations:} Say $\psi = \gamma \star -$ and $\psi' = \gamma' \star -$ for two isotopic $(n-m)$-strand braids $\gamma$ and $\gamma'$. Then the corresponding regularizations are
        \[
        \tau = \gamma \star Tr^m(\beta) \sim Tr^m((\strand{m} \sqcup \gamma) \star \beta)
        \]
        and
        \[
        \tau' = \gamma' \star Tr^m(\beta) \sim Tr^m((\strand{m} \sqcup \gamma') \star \beta).
        \]
        But $\strand{m} \sqcup \gamma$ and $\strand{m} \sqcup \gamma'$ are themselves isotopic braids, so up to braid relations, this move has no effect after regularizing both sides.

        \item \textit{Leftmost conjugation:} In this case the corresponding regularizations are
        \[
        \tau = Tr^2(\sigma_1^{\pm} \star Tr^m(\beta \star I^n(\emptyset))) \sim Tr^{m+2}(\sigma_{m+1}^{\pm} \star \beta)
        \]
        and
        \[
        \tau' = Tr^2(Tr^m(\beta \star I^n(\emptyset)) \star \sigma_1^{\pm}) \sim Tr^{m+2}(\beta \star \sigma_{m+1}^{\pm}).
        \]

        These two partial braid closures are related by conjugation through closed strands.

        \item \textit{Identity relation:} The two regularizations are
        \[
        Tr^m(\sigma_{m+1}^\pm \sigma_m\dots\sigma_1\sigma_{m+1} \dots \sigma_2 \star (\strand{2} \sqcup \beta) \star \sigma_2^{-1} \dots \sigma_{m+1}^{-1}\sigma_1^{-1} \dots \sigma_m^{-1})
        \]
        and
        \[
        Tr^m(\sigma_m\dots\sigma_1\sigma_{m+1} \dots \sigma_2 \star (\strand{2} \sqcup \beta) \star \sigma_2^{-1} \dots \sigma_{m+1}^{-1}\sigma_1^{-1} \dots \sigma_m^{-1} \sigma_{m+1}^\pm),
        \]
        which are braid-isotopic.

        \item \textit{Stabilization:} We analyze only the version $\emptyset \leftrightarrow Tr(I(-) \star \sigma_1)$; the other case is exactly analogous. Applying $Tr(I(-) \star \sigma_1)$ to $Tr^m(\beta)$ and regularizing, we obtain

        \begin{align*}
            \tau' & = Tr(I(Tr^m(\beta)) \star \sigma_1) \\
            & \sim Tr(Tr^m(\sigma_m \dots \sigma_1 \star (\strand{1} \sqcup \beta) \star \sigma_1^{-1} \dots \sigma_m^{-1}) \star \sigma_1) \\
            & \sim Tr^{m+1}(\sigma_m \dots \sigma_1 \star (\strand{1} \sqcup \beta) \star \sigma_1^{-1} \dots \sigma_m^{-1} \sigma_{m+1}).
        \end{align*}

        We now apply a sequence of braid relations and conjugations through closed strands, obtaining the following:

        
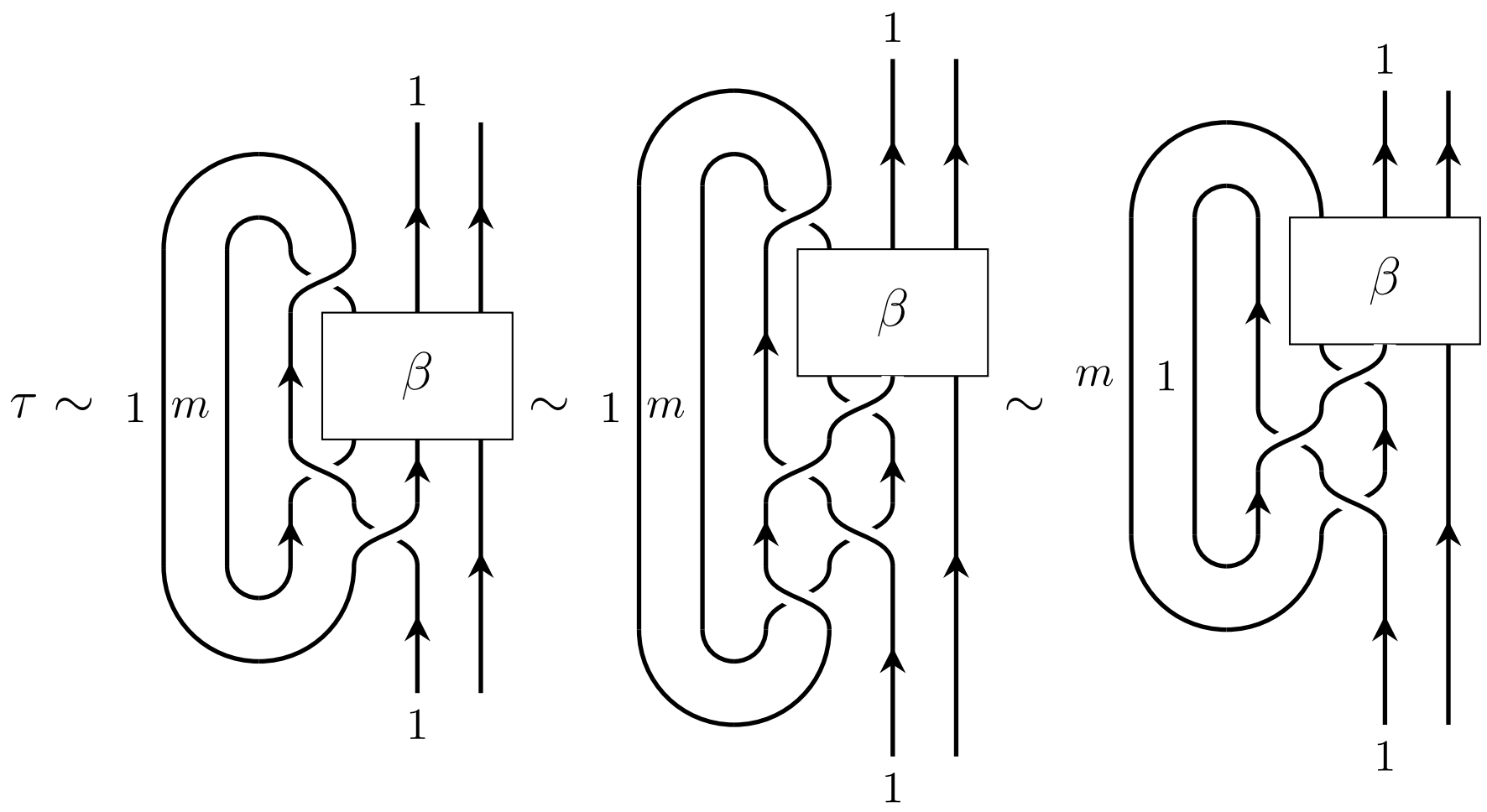
\begin{figure}[ht]
        \centering
    
        {\Large{$\tau \sim$}}
            \begin{tikzpicture}[anchorbase]
            \filldraw[fill=white, draw=black]
            (-.25,0) rectangle (1.25,1);
            \node at (.5,.5) {\Large{$\beta$}};
    
            \draw[webs] (0,1) to[out=90,in=270] (-.5,1.5);
            \draw[webs,over] (-.5,1) to[out=90,in=270] (0,1.5);
            \draw[webs,->-=.5] (.5,1) to (.5,2.5) node[above]{$1$};
            \draw[webs,->-=.5] (1,1) to (1,2.5);
    
            \draw[webs,->-=.5] (-.5,0) to (-.5,1);
            \draw[webs] (0,0) to[out=270,in=90] (-.5,-.5);
            \draw[webs,over] (-.5,0) to[out=270,in=90] (0,-.5);
            \draw[webs] (0,-.5) to[out=270,in=90] (.5,-1);
            \draw[webs,over] (.5,-.5) to[out=270,in=90] (0,-1);
    
            \draw[webs,->-=.5] (-.5,-1) to (-.5,-.5);
            \draw[webs,->-=.5] (.5,-.5) to (.5,0);
            \draw[webs,->-=.5] (.5,-2) node[below]{$1$} to (.5,-1);
            \draw[webs,->-=.5] (1,-2) to (1,0);
    
            \draw[webs] (-.5,1.5) to[out=90,in=0] (-.75,1.75);
            \draw[webs] (-.75,1.75) to[out=180,in=90] (-1,1.5);
            \draw[webs] (-1,1.5) to node[midway,left]{$m$} (-1,-1);
            \draw[webs] (-1,-1) to[out=270,in=180] (-.75,-1.25);
            \draw[webs] (-.75,-1.25) to[out=0,in=270] (-.5,-1);
    
            \draw[webs] (0,1.5) to[out=90,in=0] (-.75,2.25);
            \draw[webs] (-.75,2.25) to[out=180,in=90] (-1.5,1.5);
            \draw[webs] (-1.5,1.5) to node[midway,left]{$1$} (-1.5,-1);
            \draw[webs] (-1.5,-1) to[out=270,in=180] (-.75,-1.75);
            \draw[webs] (-.75,-1.75) to[out=0,in=270] (0,-1);
    
            \end{tikzpicture}
            {\Large{$\sim$}}
            \begin{tikzpicture}[anchorbase]
            \filldraw[fill=white, draw=black]
            (-.25,0) rectangle (1.25,1);
            \node at (.5,.5) {\Large{$\beta$}};
    
            \draw[webs] (0,1) to[out=90,in=270] (-.5,1.5);
            \draw[webs,over] (-.5,1) to[out=90,in=270] (0,1.5);
            \draw[webs,->-=.5] (.5,1) to (.5,2.5) node[above]{$1$};
            \draw[webs,->-=.5] (1,1) to (1,2.5);
    
            \draw[webs,->-=.5] (-.5,-.5) to (-.5,1);
            \draw[webs] (0,0) to[out=270,in=90] (.5,-.5);
            \draw[webs,over] (.5,0) to[out=270,in=90] (0,-.5);
            \draw[webs] (-.5,-.5) to[out=270,in=90] (0,-1);
            \draw[webs,over] (0,-.5) to[out=270,in=90] (-.5,-1);
            \draw[webs] (.5,-1) to[out=270,in=90] (0,-1.5);
            \draw[webs,over] (0,-1) to[out=270,in=90] (.5,-1.5);
            \draw[webs] (0,-1.5) to[out=270,in=90] (-.5,-2);
            \draw[webs,over] (-.5,-1.5) to[out=270,in=90] (0,-2);
    
            \draw[webs,->-=.5] (-.5,-1.5) to (-.5,-1);
            \draw[webs,->-=.5] (.5,-1) to (.5,-.5);
            \draw[webs,->-=.5] (.5,-3) node[below]{$1$} to (.5,-1.5);
            \draw[webs,->-=.5] (1,-3) to (1,0);
    
            \draw[webs] (-.5,1.5) to[out=90,in=0] (-.75,1.75);
            \draw[webs] (-.75,1.75) to[out=180,in=90] (-1,1.5);
            \draw[webs] (-1,1.5) to node[midway,left]{$m$} (-1,-2);
            \draw[webs] (-1,-2) to[out=270,in=180] (-.75,-2.25);
            \draw[webs] (-.75,-2.25) to[out=0,in=270] (-.5,-2);
    
            \draw[webs] (0,1.5) to[out=90,in=0] (-.75,2.25);
            \draw[webs] (-.75,2.25) to[out=180,in=90] (-1.5,1.5);
            \draw[webs] (-1.5,1.5) to node[midway,left]{$1$} (-1.5,-2);
            \draw[webs] (-1.5,-2) to[out=270,in=180] (-.75,-2.75);
            \draw[webs] (-.75,-2.75) to[out=0,in=270] (0,-2);
    
            \end{tikzpicture}
            {\Large{$\sim$}}
            \begin{tikzpicture}[anchorbase]
            \filldraw[fill=white, draw=black]
            (-.25,0) rectangle (1.25,1);
            \node at (.5,.5) {\Large{$\beta$}};
            
            \draw[webs,->-=.5] (.5,1) to (.5,2) node[above]{$1$};
            \draw[webs,->-=.5] (1,1) to (1,2);
    
            \draw[webs,->-=.5] (-.5,-.5) to (-.5,1);
            \draw[webs] (0,0) to[out=270,in=90] (.5,-.5);
            \draw[webs,over] (.5,0) to[out=270,in=90] (0,-.5);
            \draw[webs] (-.5,-.5) to[out=270,in=90] (0,-1);
            \draw[webs,over] (0,-.5) to[out=270,in=90] (-.5,-1);
            \draw[webs] (.5,-1) to[out=270,in=90] (0,-1.5);
            \draw[webs,over] (0,-1) to[out=270,in=90] (.5,-1.5);
    
            \draw[webs,->-=.5] (-.5,-1.5) to (-.5,-1);
            \draw[webs,->-=.5] (.5,-1) to (.5,-.5);
            \draw[webs,->-=.5] (.5,-3) node[below]{$1$} to (.5,-1.5);
            \draw[webs,->-=.5] (1,-3) to (1,0);
    
            \draw[webs] (-.5,1) to[out=90,in=0] (-.75,1.25);
            \draw[webs] (-.75,1.25) to[out=180,in=90] (-1,1);
            \draw[webs] (-1,1) to node[midway,left]{$1$} (-1,-1.5);
            \draw[webs] (-1,-1.5) to[out=270,in=180] (-.75,-1.75);
            \draw[webs] (-.75,-1.75) to[out=0,in=270] (-.5,-1.5);
    
            \draw[webs] (0,1) to[out=90,in=0] (-.75,1.75);
            \draw[webs] (-.75,1.75) to[out=180,in=90] (-1.5,1);
            \draw[webs] (-1.5,1) to node[midway,left]{$m$} (-1.5,-1.5);
            \draw[webs] (-1.5,-1.5) to[out=270,in=180] (-.75,-2.25);
            \draw[webs] (-.75,-2.25) to[out=0,in=270] (0,-1.5);
            \end{tikzpicture}
        \end{figure}
    
        \FloatBarrier

        Up to braid stabilization and a braid relation, this is exactly the partial braid closure $\tau$.

    \item \textit{Exchange moves:} We begin with the version of this move with $\psi = Tr(\sigma_1 \star I(-) \star \sigma_1^{-1})$ and $\psi' = I(Tr(-))$. Then the respective regularizations are
        \begin{align*}
        \tau & = Tr(\sigma_1 \star I(Tr^m(\beta)) \star \sigma_1^{-1}) \\
        & \sim Tr(\sigma_1 \star Tr^m(\sigma_m \sigma_{m-1} \dots \sigma_1 \star (\strand{1} \sqcup \beta) \star \sigma_1^{-1} \dots \sigma_{m-1}^{-1} \sigma_m^{-1}) \star \sigma_1^{-1}) \\
        & \sim Tr^{m+1}(\sigma_{m+1} \sigma_m \dots \sigma_1 \star (\strand{1} \sqcup \beta) \star \sigma_1^{-1} \dots \sigma_m^{-1} \sigma_{m-1}^{-1})
        \end{align*}
        and
        \begin{align*}
        \tau' & = I(Tr^{m+1}(\beta)) \\
        & \sim Tr^{m+1}(\sigma_{m+1} \sigma_m \dots \sigma_1 \star (\strand{1} \sqcup \beta) \star \sigma_1^{-1} \dots \sigma_m^{-1} \sigma_{m+1}^{-1}).
        \end{align*}
        Hence this version of the exchange move has no effect after regularizing both sides.
        
        On the other hand, taking $\psi = Tr(\sigma_1^{-1} \star I(-) \star \sigma_1)$ and $\psi'$ as before, the regularization of $\tau$ becomes

        \begin{align*}
        \tau & = Tr(\sigma_1^{-1} \star I(Tr^m(\beta)) \star \sigma_1) \\
        & \sim Tr(\sigma_1^{-1} \star Tr^m(\sigma_m \sigma_{m-1} \dots \sigma_1 \star (\strand{1} \sqcup \beta) \star \sigma_1^{-1} \dots \sigma_{m-1}^{-1} \sigma_m^{-1}) \star \sigma_1) \\
        & \sim Tr^{m+1} (\sigma_{m+1}^{-1} \sigma_m \dots \sigma_1 \star (\strand{1} \sqcup \beta) \star \sigma_1^{-1} \dots \sigma_m^{-1} \sigma_{m+1}).
        \end{align*}

        We can now apply a similar sequence of braid relations and conjugations through closed strands as in our analysis of stabilization to obtain the following:

        
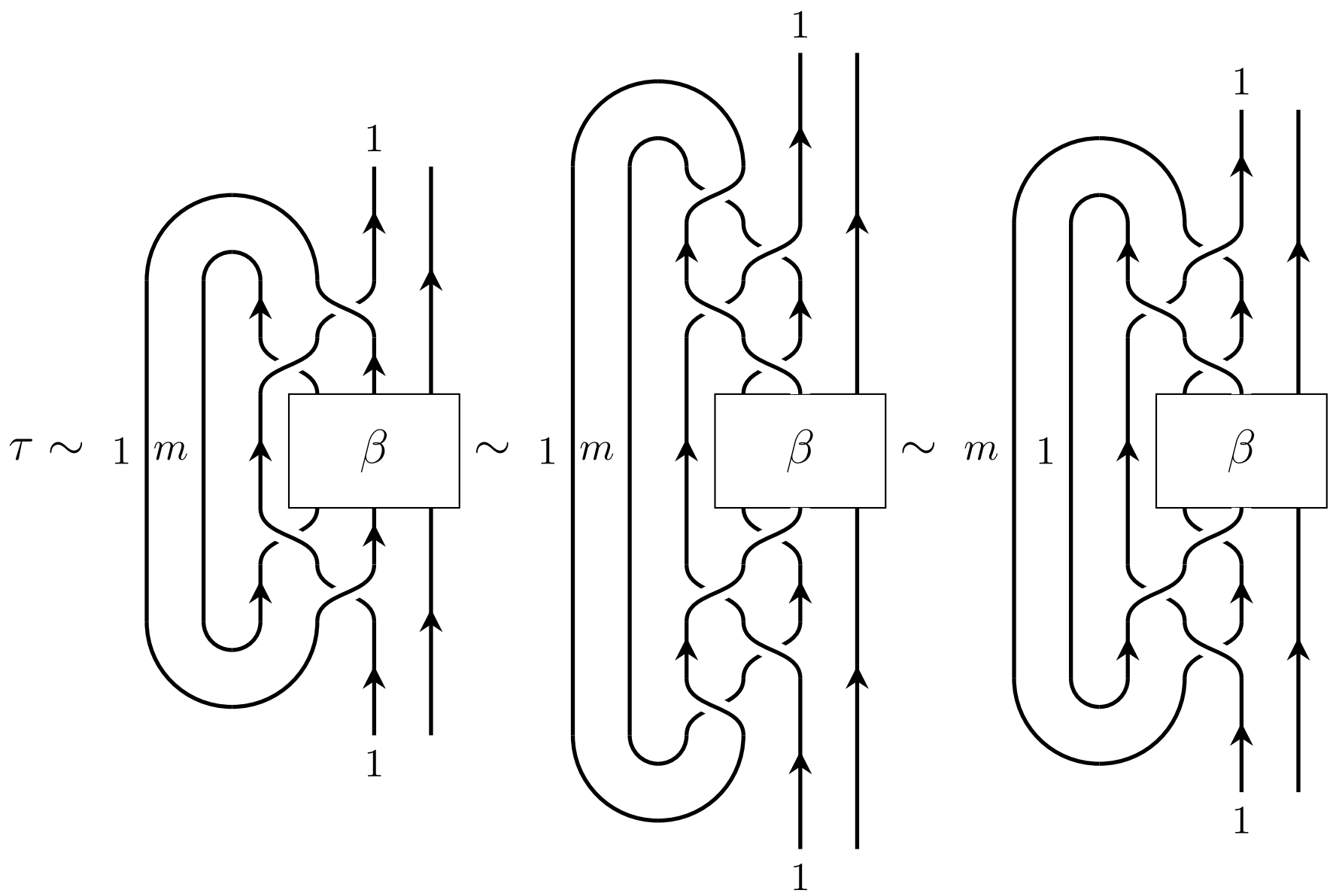
\begin{figure}[ht]
        \centering

        {\Large{$\tau \sim$}}
        \begin{tikzpicture}[anchorbase]
        \filldraw[fill=white, draw=black]
        (-.25,0) rectangle (1.25,1);
        \node at (.5,.5) {\Large{$\beta$}};
        
        \draw[webs] (0,1) to[out=90,in=270] (-.5,1.5);
        \draw[webs,over] (-.5,1) to[out=90,in=270] (0,1.5);
        \draw[webs,->-=.5] (.5,1) to (.5,1.5);
        \draw[webs,->-=.5] (1,1) to (1,3);

        \draw[webs] (0,1.5) to[out=90,in=270] (.5,2);
        \draw[webs,over] (.5,1.5) to[out=90,in=270] (0,2);
        \draw[webs,->-=.5] (-.5,1.5) to (-.5,2);
        \draw[webs,->-=.5] (.5,2) to (.5,3) node[above]{$1$};

        \draw[webs,->-=.5] (-.5,0) to (-.5,1);
        \draw[webs] (0,0) to[out=270,in=90] (-.5,-.5);
        \draw[webs,over] (-.5,0) to[out=270,in=90] (0,-.5);
        \draw[webs,-<-=.5] (.5,0) to (.5,-.5);

        \draw[webs] (0,-.5) to[out=270,in=90] (.5,-1);
        \draw[webs,over] (.5,-.5) to[out=270,in=90] (0,-1);
        \draw[webs,-<-=.5] (-.5,-.5) to (-.5,-1);
        \draw[webs,-<-=.5] (1,0) to (1,-2);
        \draw[webs,-<-=.5] (.5,-1) to (.5,-2) node[below]{$1$};

        \draw[webs] (-.5,2) to[out=90,in=0] (-.75,2.25);
        \draw[webs] (-.75,2.25) to[out=180,in=90] (-1,2);
        \draw[webs] (-1,2) to node[midway,left]{$m$} (-1,-1);
        \draw[webs] (-1,-1) to[out=270,in=180] (-.75,-1.25);
        \draw[webs] (-.75,-1.25) to[out=0,in=270] (-.5,-1);

        \draw[webs] (0,2) to[out=90,in=0] (-.75,2.75);
        \draw[webs] (-.75,2.75) to[out=180,in=90] (-1.5,2);
        \draw[webs] (-1.5,2) to node[midway,left]{$1$} (-1.5,-1);
        \draw[webs] (-1.5,-1) to[out=270,in=180] (-.75,-1.75);
        \draw[webs] (-.75,-1.75) to[out=0,in=270] (0,-1);
        \end{tikzpicture}
        {\Large{$\sim$}}
        \begin{tikzpicture}[anchorbase]
        \filldraw[fill=white, draw=black]
        (-.25,0) rectangle (1.25,1);
        \node at (.5,.5) {\Large{$\beta$}};
        
        \draw[webs,->-=.5] (-.5,-.5) to (-.5,1.5);

        \draw[webs] (0,1) to[out=90,in=270] (.5,1.5);
        \draw[webs,over] (.5,1) to[out=90,in=270] (0,1.5);
        \draw[webs] (-.5,1.5) to[out=90,in=270] (0,2);
        \draw[webs,over] (0,1.5) to[out=90,in=270] (-.5,2);
        \draw[webs] (.5,2) to[out=90,in=270] (0,2.5);
        \draw[webs,over] (0,2) to[out=90,in=270] (.5,2.5);
        \draw[webs] (0,2.5) to[out=90,in=270] (-.5,3);
        \draw[webs,over] (-.5,2.5) to[out=90,in=270] (0,3);

        \draw[webs,->-=.5] (-.5,2) to (-.5,2.5);
        \draw[webs,->-=.5] (.5,1.5) to (.5,2);
        \draw[webs,->-=.5] (.5,2.5) to (.5,4) node[above]{$1$};
        \draw[webs,->-=.5] (1,1) to (1,4);

        \begin{scope}[shift={(0,1)},yscale=-1]
        \draw[webs] (0,1) to[out=90,in=270] (.5,1.5);
        \draw[webs,over] (.5,1) to[out=90,in=270] (0,1.5);
        \draw[webs] (-.5,1.5) to[out=90,in=270] (0,2);
        \draw[webs,over] (0,1.5) to[out=90,in=270] (-.5,2);
        \draw[webs] (.5,2) to[out=90,in=270] (0,2.5);
        \draw[webs,over] (0,2) to[out=90,in=270] (.5,2.5);
        \draw[webs] (0,2.5) to[out=90,in=270] (-.5,3);
        \draw[webs,over] (-.5,2.5) to[out=90,in=270] (0,3);

        \draw[webs,-<-=.5] (-.5,2) to (-.5,2.5);
        \draw[webs,-<-=.5] (.5,1.5) to (.5,2);
        \draw[webs,-<-=.5] (.5,2.5) to (.5,4) node[below]{$1$};
        \draw[webs,-<-=.5] (1,1) to (1,4);
        \end{scope}

        \draw[webs] (-.5,3) to[out=90,in=0] (-.75,3.25);
        \draw[webs] (-.75,3.25) to[out=180,in=90] (-1,3);
        \draw[webs] (-1,3) to node[midway,left]{$m$} (-1,-2);
        \draw[webs] (-1,-2) to[out=270,in=180] (-.75,-2.25);
        \draw[webs] (-.75,-2.25) to[out=0,in=270] (-.5,-2);

        \draw[webs] (0,3) to[out=90,in=0] (-.75,3.75);
        \draw[webs] (-.75,3.75) to[out=180,in=90] (-1.5,3);
        \draw[webs] (-1.5,3) to node[midway,left]{$1$} (-1.5,-2);
        \draw[webs] (-1.5,-2) to[out=270,in=180] (-.75,-2.75);
        \draw[webs] (-.75,-2.75) to[out=0,in=270] (0,-2);
        \end{tikzpicture}
        {\Large{$\sim$}}
        \begin{tikzpicture}[anchorbase]
        \filldraw[fill=white, draw=black]
        (-.25,0) rectangle (1.25,1);
        \node at (.5,.5) {\Large{$\beta$}};
        
        \draw[webs,->-=.5] (-.5,-.5) to (-.5,1.5);

        \draw[webs] (0,1) to[out=90,in=270] (.5,1.5);
        \draw[webs,over] (.5,1) to[out=90,in=270] (0,1.5);
        \draw[webs] (-.5,1.5) to[out=90,in=270] (0,2);
        \draw[webs,over] (0,1.5) to[out=90,in=270] (-.5,2);
        \draw[webs] (.5,2) to[out=90,in=270] (0,2.5);
        \draw[webs,over] (0,2) to[out=90,in=270] (.5,2.5);

        \draw[webs,->-=.5] (-.5,2) to (-.5,2.5);
        \draw[webs,->-=.5] (.5,1.5) to (.5,2);
        \draw[webs,->-=.5] (.5,2.5) to (.5,3.5) node[above]{$1$};
        \draw[webs,->-=.5] (1,1) to (1,3.5);

        \begin{scope}[shift={(0,1)},yscale=-1]
        \draw[webs] (0,1) to[out=90,in=270] (.5,1.5);
        \draw[webs,over] (.5,1) to[out=90,in=270] (0,1.5);
        \draw[webs] (-.5,1.5) to[out=90,in=270] (0,2);
        \draw[webs,over] (0,1.5) to[out=90,in=270] (-.5,2);
        \draw[webs] (.5,2) to[out=90,in=270] (0,2.5);
        \draw[webs,over] (0,2) to[out=90,in=270] (.5,2.5);

        \draw[webs,-<-=.5] (-.5,2) to (-.5,2.5);
        \draw[webs,-<-=.5] (.5,1.5) to (.5,2);
        \draw[webs,-<-=.5] (.5,2.5) to (.5,3.5) node[below]{$1$};
        \draw[webs,-<-=.5] (1,1) to (1,3.5);
        \end{scope}

        \draw[webs] (-.5,2.5) to[out=90,in=0] (-.75,2.75);
        \draw[webs] (-.75,2.75) to[out=180,in=90] (-1,2.5);
        \draw[webs] (-1,2.5) to node[midway,left]{$1$} (-1,-1.5);
        \draw[webs] (-1,-1.5) to[out=270,in=180] (-.75,-1.75);
        \draw[webs] (-.75,-1.75) to[out=0,in=270] (-.5,-1.5);

        \draw[webs] (0,2.5) to[out=90,in=0] (-.75,3.25);
        \draw[webs] (-.75,3.25) to[out=180,in=90] (-1.5,2.5);
        \draw[webs] (-1.5,2.5) to node[midway,left]{$m$} (-1.5,-1.5);
        \draw[webs] (-1.5,-1.5) to[out=270,in=180] (-.75,-2.25);
        \draw[webs] (-.75,-2.25) to[out=0,in=270] (0,-1.5);
        \end{tikzpicture}
        \end{figure}

        \FloatBarrier

        Meanwhile, applying a similar sequence to $\tau'$, we obtain:

        
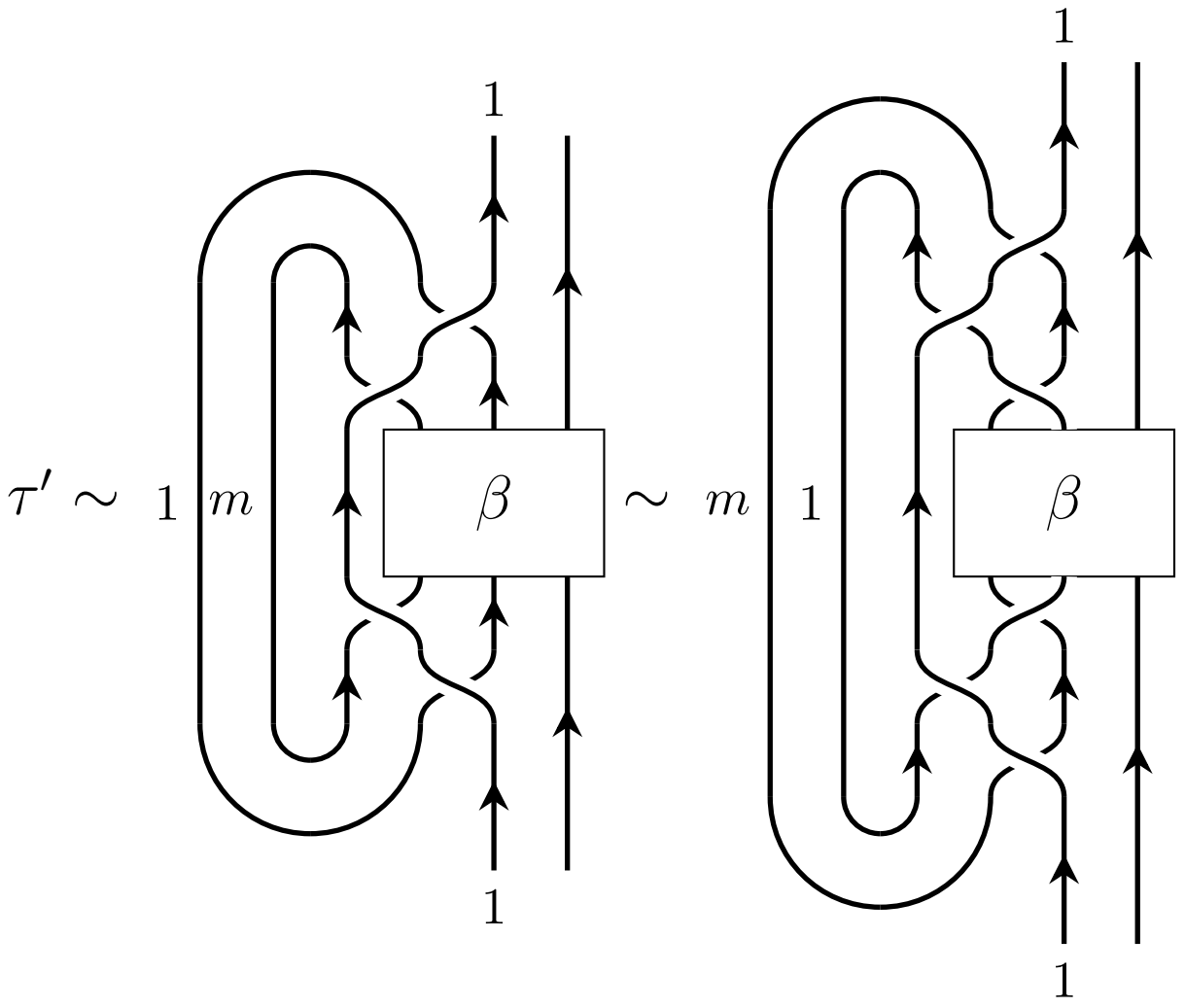
\begin{figure}[ht]
        \centering

        {\Large{$\tau' \sim$}}
        \begin{tikzpicture}[anchorbase]
        \filldraw[fill=white, draw=black]
        (-.25,0) rectangle (1.25,1);
        \node at (.5,.5) {\Large{$\beta$}};

        \draw[webs] (0,1) to[out=90,in=270] (-.5,1.5);
        \draw[webs,over] (-.5,1) to[out=90,in=270] (0,1.5);
        \draw[webs,->-=.5] (.5,1) to (.5,1.5);
        \draw[webs,->-=.5] (1,1) to (1,3);

        \draw[webs] (.5,1.5) to[out=90,in=270] (0,2);
        \draw[webs,over] (0,1.5) to[out=90,in=270] (.5,2);
        \draw[webs,->-=.5] (-.5,1.5) to (-.5,2);
        \draw[webs,->-=.5] (.5,2) to (.5,3) node[above]{$1$};

        \draw[webs,->-=.5] (-.5,0) to (-.5,1);
        \draw[webs] (0,0) to[out=270,in=90] (-.5,-.5);
        \draw[webs,over] (-.5,0) to[out=270,in=90] (0,-.5);
        \draw[webs,-<-=.5] (.5,0) to (.5,-.5);

        \draw[webs] (.5,-.5) to[out=270,in=90] (0,-1);
        \draw[webs,over] (0,-.5) to[out=270,in=90] (.5,-1);
        \draw[webs,-<-=.5] (-.5,-.5) to (-.5,-1);
        \draw[webs,-<-=.5] (1,0) to (1,-2);
        \draw[webs,-<-=.5] (.5,-1) to (.5,-2) node[below]{$1$};

        \draw[webs] (-.5,2) to[out=90,in=0] (-.75,2.25);
        \draw[webs] (-.75,2.25) to[out=180,in=90] (-1,2);
        \draw[webs] (-1,2) to node[midway,left]{$m$} (-1,-1);
        \draw[webs] (-1,-1) to[out=270,in=180] (-.75,-1.25);
        \draw[webs] (-.75,-1.25) to[out=0,in=270] (-.5,-1);

        \draw[webs] (0,2) to[out=90,in=0] (-.75,2.75);
        \draw[webs] (-.75,2.75) to[out=180,in=90] (-1.5,2);
        \draw[webs] (-1.5,2) to node[midway,left]{$1$} (-1.5,-1);
        \draw[webs] (-1.5,-1) to[out=270,in=180] (-.75,-1.75);
        \draw[webs] (-.75,-1.75) to[out=0,in=270] (0,-1);
        \end{tikzpicture}
        {\Large{$\sim$}}
        \begin{tikzpicture}[anchorbase]
        \filldraw[fill=white, draw=black]
        (-.25,0) rectangle (1.25,1);
        \node at (.5,.5) {\Large{$\beta$}};
        
        \draw[webs,->-=.5] (-.5,-.5) to (-.5,1.5);

        \draw[webs] (0,1) to[out=90,in=270] (.5,1.5);
        \draw[webs,over] (.5,1) to[out=90,in=270] (0,1.5);
        \draw[webs] (0,1.5) to[out=90,in=270] (-.5,2);
        \draw[webs,over] (-.5,1.5) to[out=90,in=270] (0,2);
        \draw[webs] (.5,2) to[out=90,in=270] (0,2.5);
        \draw[webs,over] (0,2) to[out=90,in=270] (.5,2.5);

        \draw[webs,->-=.5] (-.5,2) to (-.5,2.5);
        \draw[webs,->-=.5] (.5,1.5) to (.5,2);
        \draw[webs,->-=.5] (.5,2.5) to (.5,3.5) node[above]{$1$};
        \draw[webs,->-=.5] (1,1) to (1,3.5);

        \begin{scope}[shift={(0,1)},yscale=-1]
        \draw[webs] (0,1) to[out=90,in=270] (.5,1.5);
        \draw[webs,over] (.5,1) to[out=90,in=270] (0,1.5);
        \draw[webs] (0,1.5) to[out=90,in=270] (-.5,2);
        \draw[webs,over] (-.5,1.5) to[out=90,in=270] (0,2);
        \draw[webs] (.5,2) to[out=90,in=270] (0,2.5);
        \draw[webs,over] (0,2) to[out=90,in=270] (.5,2.5);

        \draw[webs,-<-=.5] (-.5,2) to (-.5,2.5);
        \draw[webs,-<-=.5] (.5,1.5) to (.5,2);
        \draw[webs,-<-=.5] (.5,2.5) to (.5,3.5) node[below]{$1$};
        \draw[webs,-<-=.5] (1,1) to (1,3.5);
        \end{scope}

        \draw[webs] (-.5,2.5) to[out=90,in=0] (-.75,2.75);
        \draw[webs] (-.75,2.75) to[out=180,in=90] (-1,2.5);
        \draw[webs] (-1,2.5) to node[midway,left]{$1$} (-1,-1.5);
        \draw[webs] (-1,-1.5) to[out=270,in=180] (-.75,-1.75);
        \draw[webs] (-.75,-1.75) to[out=0,in=270] (-.5,-1.5);

        \draw[webs] (0,2.5) to[out=90,in=0] (-.75,3.25);
        \draw[webs] (-.75,3.25) to[out=180,in=90] (-1.5,2.5);
        \draw[webs] (-1.5,2.5) to node[midway,left]{$m$} (-1.5,-1.5);
        \draw[webs] (-1.5,-1.5) to[out=270,in=180] (-.75,-2.25);
        \draw[webs] (-.75,-2.25) to[out=0,in=270] (0,-1.5);
        \end{tikzpicture}
        \end{figure}

        \FloatBarrier

        These two partial braid closures differ only by a local application of the exchange relation $Tr^{m+1}(\sigma_1 \star (\strand{1} \sqcup \beta) \star \sigma_1^{-1}) \leftrightarrow Tr^{m+1}(\sigma_1^{-1} \star (\strand{1} \sqcup \beta) \star \sigma_1)$. In fact, this version of the exchange relation can be realized by a sequence of braid relations, conjugations through closed strands, and braid stabilizations; see Figure~\ref{fig: ExchangeViaStabilization}.
    \end{itemize}

    \begin{figure}[ht]
        \centering
        \resizebox{\textwidth}{!}{
\begin{tikzpicture}[anchorbase]
    \draw[webs] (1,-1.65) to (1,3.9);
    \draw[webs] (.25,-1.65) to (.25,-.75);
    \draw[webs] (.25,3) to (.25,3.9);
    \draw[webs] (-.5,0) to (-.5,2.25);
    \draw[webs] (-1.25,-.75) to (-1.25,3);
    \draw[webs] (-1.25,3) to[out=90,in=180] (-.875,3.375) to[out=0,in=90] (-.5,3);
    \draw[webs] (-1.25,-.75) to[out=270,in=180] (-.875,-1.125) to[out=0,in=270] (-.5,-.75);
    \draw[webs] (-.5,2.25) to[out=90,in=270] (.25,3);
    \draw[webs,over] (.25,2.25) to[out=90,in=270] (-.5,3);
    \draw[webs] (.25,-.75) to[out=90,in=270] (-.5,0);
    \draw[webs,over] (-.5,-.75) to[out=90,in=270] (.25,0);
    \filldraw[fill=white, draw=black] (0,0) rectangle (1.25,2.25);
    \node at (.625,1.125) {\Large{$\beta$}};
\end{tikzpicture}
{\Large{$\sim$}}
\begin{tikzpicture}[anchorbase]
    \draw[webs] (1,-4.5) to (1,4.5);
    \draw[webs] (.25,-4.5) to (.25,-3);
    \draw[webs] (.25,-2.25) to (.25,-1.5);
    \draw[webs] (.25,3) to (.25,4.5);
    \draw[webs] (-.5,0) to (-.5,2.25);
    \draw[webs] (-1.25,-3) to (-1.25,-2.25);
    \draw[webs] (-1.25,-1.5) to (-1.25,3);
    \draw[webs] (-2,-3) to (-2,3);
    \draw[webs] (-2.75,-3) to (-2.75,3);
    \draw[webs] (-2,3) to[out=90,in=180] (-1.625,3.375) to[out=0,in=90] (-1.25,3);
    \draw[webs] (-2,-3) to[out=270,in=180] (-1.625,-3.375) to[out=0,in=270] (-1.25,-3);
    \draw[webs] (-2.75,-3) to[out=270,in=180] (-1.625,-4.125) to[out=0,in=270] (-.5,-3);
    \draw[webs] (-2.75,3) to[out=90,in=180] (-1.625,4.125) to[out=0,in=90] (-.5,3);
    \draw[webs] (-1.25,-2.25) to[out=90,in=270] (-.5,-1.5);
    \draw[webs,over] (-.5,-2.25) to[out=90,in=270] (-1.25,-1.5);
    \draw[webs] (.25,-1.5) to[out=90,in=270] (-.5,-.75);
    \draw[webs,over] (-.5,-1.5) to[out=90,in=270] (.25,-.75);
    \draw[webs] (-.5,2.25) to[out=90,in=270] (.25,3);
    \draw[webs,over] (.25,2.25) to[out=90,in=270] (-.5,3);
    \draw[webs] (.25,-.75) to[out=90,in=270] (-.5,0);
    \draw[webs,over] (-.5,-.75) to[out=90,in=270] (.25,0);
    \draw[webs] (-.5,-3) to[out=90,in=270] (.25,-2.25);
    \draw[webs,over] (.25,-3) to[out=90,in=270] (-.5,-2.25);
    \filldraw[fill=white, draw=black] (0,0) rectangle (1.25,2.25);
    \node at (.625,1.125) {\Large{$\beta$}};
\end{tikzpicture}
{\Large{$\sim$}}
\begin{tikzpicture}[anchorbase]
    \draw[webs] (1,-3.75) to (1,4.5);
    \draw[webs] (.25,-3.75) to (.25,-.75);
    \draw[webs] (.25,3) to (.25,4.5);
    \draw[webs] (-.5,.75) to (-.5,2.25);
    \draw[webs] (-1.25,-.75) to (-1.25,0);
    \draw[webs] (-1.25,.75) to (-1.25,3);
    \draw[webs] (-2,-2.25) to (-2,3);
    \draw[webs] (-2.75,-2.25) to (-2.75,3);
    \draw[webs] (-2,3) to[out=90,in=180] (-1.625,3.375) to[out=0,in=90] (-1.25,3);
    \draw[webs] (-2,-2.25) to[out=270,in=180] (-1.625,-2.625) to[out=0,in=270] (-1.25,-2.25);
    \draw[webs] (-2.75,-2.25) to[out=270,in=180] (-1.625,-3.375) to[out=0,in=270] (-.5,-2.25);
    \draw[webs] (-2.75,3) to[out=90,in=180] (-1.625,4.125) to[out=0,in=90] (-.5,3);
    \draw[webs] (-.5,-.75) to[out=90,in=270] (.25,0);
    \draw[webs,over] (.25,-.75) to[out=90,in=270] (-.5,0);
    \draw[webs] (-.5,-2.25) to[out=90,in=270] (-1.25,-1.5);
    \draw[webs,over] (-1.25,-2.25) to[out=90,in=270] (-.5,-1.5);
    \draw[webs] (-.5,2.25) to[out=90,in=270] (.25,3);
    \draw[webs,over] (.25,2.25) to[out=90,in=270] (-.5,3);
    \draw[webs] (-.5,-1.5) to[out=90,in=270] (-1.25,-.75);
    \draw[webs,over] (-1.25,-1.5) to[out=90,in=270] (-.5,-.75);
    \draw[webs] (-1.25,0) to[out=90,in=270] (-.5,.75);
    \draw[webs,over] (-.5,0) to[out=90,in=270] (-1.25,.75);
    \filldraw[fill=white, draw=black] (0,0) rectangle (1.25,2.25);
    \node at (.625,1.125) {\Large{$\beta$}};
\end{tikzpicture}
{\Large{$\sim$}}
\begin{tikzpicture}[anchorbase]
    \draw[webs] (1,-3) to (1,6);
    \draw[webs] (.25,-3) to (.25,-1.5);
    \draw[webs] (.25,-.75) to (.25,0);
    \draw[webs] (.25,3) to (.25,6);
    \draw[webs] (-.5,-.75) to (-.5,1.5);
    \draw[webs] (-1.25,-1.5) to (-1.25,1.5);
    \draw[webs] (-1.25,2.25) to (-1.25,3);
    \draw[webs] (-2,-1.5) to (-2,4.5);
    \draw[webs] (-2.75,-1.5) to (-2.75,4.5);
    \draw[webs] (-2,4.5) to[out=90,in=180] (-1.625,4.875) to[out=0,in=90] (-1.25,4.5);
    \draw[webs] (-2,-1.5) to[out=270,in=180] (-1.625,-1.875) to[out=0,in=270] (-1.25,-1.5);
    \draw[webs] (-2.75,-1.5) to[out=270,in=180] (-1.625,-2.625) to[out=0,in=270] (-.5,-1.5);
    \draw[webs] (-2.75,4.5) to[out=90,in=180] (-1.625,5.625) to[out=0,in=90] (-.5,4.5);
    \draw[webs] (-.5,-1.5) to[out=90,in=270] (.25,-.75);
    \draw[webs,over] (.25,-1.5) to[out=90,in=270] (-.5,-.75);
    \draw[webs] (-.5,3) to[out=90,in=270] (-1.25,3.75);
    \draw[webs,over] (-1.25,3) to[out=90,in=270] (-.5,3.75);
    \draw[webs] (-.5,2.25) to[out=90,in=270] (.25,3);
    \draw[webs,over] (.25,2.25) to[out=90,in=270] (-.5,3);
    \draw[webs] (-.5,3.75) to[out=90,in=270] (-1.25,4.5);
    \draw[webs,over] (-1.25,3.75) to[out=90,in=270] (-.5,4.5);
    \draw[webs] (-1.25,1.5) to[out=90,in=270] (-.5,2.25);
    \draw[webs,over] (-.5,1.5) to[out=90,in=270] (-1.25,2.25);
    \filldraw[fill=white, draw=black] (0,0) rectangle (1.25,2.25);
    \node at (.625,1.125) {\Large{$\beta$}};
\end{tikzpicture}
{\Large{$\sim$}}
\begin{tikzpicture}[anchorbase]
    \draw[webs] (1,-3) to (1,6.75);
    \draw[webs] (.25,-3) to (.25,-1.5);
    \draw[webs] (.25,-.75) to (.25,0);
    \draw[webs] (.25,3.75) to (.25,4.5);
    \draw[webs] (.25,5.25) to (.25,6.75);
    \draw[webs] (-.5,-.75) to (-.5,2.25);
    \draw[webs] (-1.25,-1.5) to (-1.25,3.75);
    \draw[webs] (-1.25,4.5) to (-1.25,5.25);
    \draw[webs] (-2,-1.5) to (-2,5.25);
    \draw[webs] (-2.75,-1.5) to (-2.75,5.25);
    \draw[webs] (-2,5.25) to[out=90,in=180] (-1.625,5.625) to[out=0,in=90] (-1.25,5.25);
    \draw[webs] (-2,-1.5) to[out=270,in=180] (-1.625,-1.875) to[out=0,in=270] (-1.25,-1.5);
    \draw[webs] (-2.75,-1.5) to[out=270,in=180] (-1.625,-2.625) to[out=0,in=270] (-.5,-1.5);
    \draw[webs] (-2.75,5.25) to[out=90,in=180] (-1.625,6.375) to[out=0,in=90] (-.5,5.25);
    \draw[webs] (-.5,-1.5) to[out=90,in=270] (.25,-.75);
    \draw[webs,over] (.25,-1.5) to[out=90,in=270] (-.5,-.75);
    \draw[webs] (.25,2.25) to[out=90,in=270] (-.5,3);
    \draw[webs,over] (-.5,2.25) to[out=90,in=270] (.25,3);
    \draw[webs] (-.5,4.5) to[out=90,in=270] (.25,5.25);
    \draw[webs,over] (.25,4.5) to[out=90,in=270] (-.5,5.25);
    \draw[webs] (.25,3) to[out=90,in=270] (-.5,3.75);
    \draw[webs,over] (-.5,3) to[out=90,in=270] (.25,3.75);
    \draw[webs] (-1.25,3.75) to[out=90,in=270] (-.5,4.5);
    \draw[webs,over] (-.5,3.75) to[out=90,in=270] (-1.25,4.5);
    \filldraw[fill=white, draw=black] (0,0) rectangle (1.25,2.25);
    \node at (.625,1.125) {\Large{$\beta$}};
\end{tikzpicture}
{\Large{$\sim$}}
\begin{tikzpicture}[anchorbase]
    \draw[webs] (1,-1.65) to (1,3.9);
    \draw[webs] (.25,-1.65) to (.25,-.75);
    \draw[webs] (.25,3) to (.25,3.9);
    \draw[webs] (-.5,0) to (-.5,2.25);
    \draw[webs] (-1.25,-.75) to (-1.25,3);
    \draw[webs] (-1.25,3) to[out=90,in=180] (-.875,3.375) to[out=0,in=90] (-.5,3);
    \draw[webs] (-1.25,-.75) to[out=270,in=180] (-.875,-1.125) to[out=0,in=270] (-.5,-.75);
    \draw[webs] (.25,2.25) to[out=90,in=270] (-.5,3);
    \draw[webs,over] (-.5,2.25) to[out=90,in=270] (.25,3);
    \draw[webs] (-.5,-.75) to[out=90,in=270] (.25,0);
    \draw[webs,over] (.25,-.75) to[out=90,in=270] (-.5,0);
    \filldraw[fill=white, draw=black] (0,0) rectangle (1.25,2.25);
    \node at (.625,1.125) {\Large{$\beta$}};
\end{tikzpicture}}
    \caption{Accomplishing a single exchange move via (de)stabilizations, braid isotopies, and conjugation through closed strands}\label{fig: ExchangeViaStabilization}
    \end{figure}
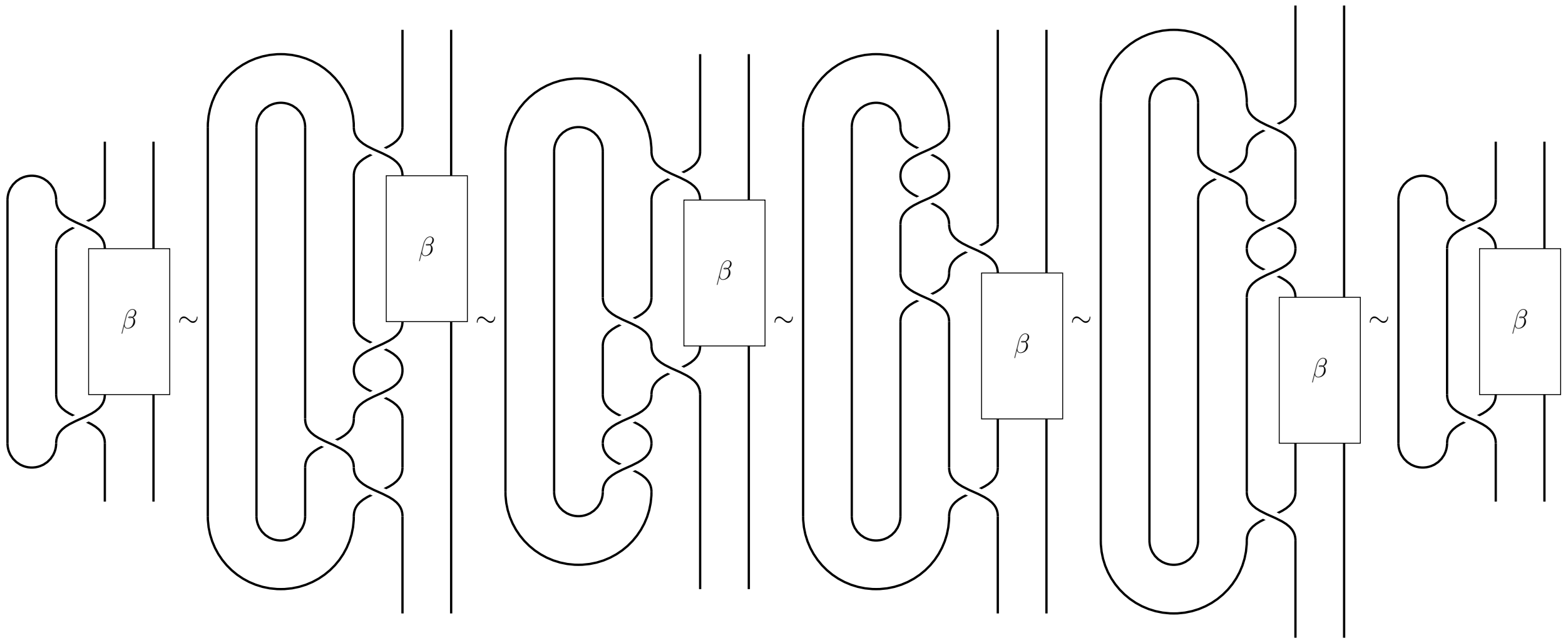

    It remains to examine the case when $\varphi \neq 0$. To do this, we consider how the partial braid moves in the theorem statement behave under application of the operations $Tr$, $I$, $\sigma_i^{\pm} \star -$, and $- \star \sigma_i^{\pm}$. We analyze each such move in turn:

    \begin{itemize}
        \item \textit{Braid relations}: Say $\beta$ and $\beta'$ are two isotopic $n$-strand braids, and consider the two partial braid closures $Tr^m(\beta)$ and $Tr^m(\beta')$. After applying any of $Tr(-)$, $I(-)$, $\sigma_i^{\pm} \star -$, or $- \star \sigma_i^{\pm}$ to both and regularizing, it is easy to see that the resulting partial braid closures are still related by braid relations.

        \item \textit{Conjugation through closed strands:} Consider the two partial braid closures $Tr^m(\sigma_i^{\pm} \star \beta)$ and $Tr^m(\beta \star \sigma_i^{\pm})$ for some $i \leq m - 1$. If we apply $Tr(-)$ to both, the results are clearly still related by conjugation through closed strands. If we apply $\sigma_j^{\pm} \star -$ or $- \star \sigma_j^{\pm}$ to both, then after regularizing and applying distant interchange of crossings, the two resulting partial braid closures are again related by conjugation through closed strands. Finally, if we apply $I(-)$ to both, then after regularizing, the crossing $\sigma_i^{\pm}$ can slide under the new crossings, and the resulting partial braid closures are again related by conjugation through closed strands.

        \item \textit{Braid stabilization:} Consider the two partial braid closures $Tr^m(\beta)$ and $Tr^{m+1}((\strand{1} \sqcup \beta) \star \sigma_1^{\pm})$. If we apply $Tr(-)$ to both, the resulting partial braid closures are still related by braid stabilization. The same is true of the operations $\sigma_j^{\pm} \star -$ and $- \star \sigma_j^{\pm}$, as the braids in the resulting partial braid closures are modified identically and the crossing $\sigma_1^{\pm}$ can be slid to the bottom by distant interchange. Finally, if we apply $I(-)$ to both, then after regularizing, the resulting partial braid closures are also still related by braid relations and a braid stabilization.
    \end{itemize}
\end{proof}

\begin{remark}
    As an immediate consequence of Theorem \ref{thm: markov_tangles_braid_closure}, we obtain a new proof of the classical Markov theorem for braid closures: two braids $\beta$ and $\beta'$ have isotopic closures if and only if they are related by a finite sequence of braid relations, conjugations, and stabilizations.
\end{remark}

\FloatBarrier

\bibliographystyle{alpha}
\bibliography{Tangles_bib.bib}

\end{document}